\documentclass{amsart}
\usepackage{amsmath} 
\usepackage{amssymb}
\usepackage{mathrsfs}

\newbox\noforkbox \newdimen\forklinewidth
\setbox0\hbox{$\textstyle\bigcup$}
\setbox1\hbox to \wd0{\hfil\vrule width \forklinewidth depth \dp0
                        height \ht0 \hfil}
\setbox\noforkbox\hbox{\box1\box0\relax}
\def\unionstick{\mathop{\copy\noforkbox}\limits}
\def\nonfork#1#2_#3{#1\unionstick_{\textstyle #3}#2}
\def\nonforkin#1#2_#3^#4{#1\unionstick_{\textstyle #3}^{\textstyle #4}#2}     
\setbox0\hbox{$\textstyle\bigcup$}
\setbox1\hbox to \wd0{\hfil{\sl /\/}\hfil}
\setbox2\hbox to \wd0{\hfil\vrule height \ht0 depth \dp0 width
                                \forklinewidth\hfil}
\newbox\doesforkbox
\setbox\doesforkbox\hbox{\box1\box0\relax}
\def\nunionstick{\mathop{\copy\doesforkbox}\limits}

\def\fork#1#2_#3{#1\nunionstick_{\textstyle #3}#2}
\def\forkin#1#2_#3^#4{#1\nunionstick_{\textstyle #3}^{\textstyle #4}#2}

\newtheorem{theorem}{Theorem}[section] 
\newtheorem{claim}[theorem]{Claim}
\newtheorem{lemma}[theorem]{Lemma} 
\newtheorem{ml}[theorem]{Main Lemma} 
 
\newtheorem{conclusion}[theorem]{Conclusion}
\newtheorem{observation}[theorem]{Observation}

\theoremstyle{definition}
\newtheorem{definition}[theorem]{Definition}
\newtheorem{preliminaries}[theorem]{Preliminaries}
\newtheorem{tcd}[theorem]{The Construction Definition}
\newtheorem{example}[theorem]{Example}
\newtheorem{fact}[theorem]{Fact}
\newtheorem{context}[theorem]{Context}
\newtheorem{construction}[theorem]{Construction}
\newtheorem{convention}[theorem]{Convention}
\newtheorem{tcc}[theorem]{The Context Continued}

\newtheorem{discussion}[theorem]{Discussion}

\theoremstyle{remark}
\newtheorem{remark}[theorem]{Remark}
\newtheorem{question}[theorem]{Question}
\newtheorem{notation}[theorem]{Notation}

\newcommand{\st}{{\rm such that}}

\newcommand{\Th}{{\rm Th}}
\newcommand{\tp}{{\rm tp}}

\newcommand{\otp}{{\rm otp}}
\newcommand{\sat}{{\rm sat}}
\newcommand{\subadd}{{\rm subadd}}

\newcommand{\ps}{{\rm ps}}
\newcommand{\ns}{{\rm ns}}
\newcommand{\pk}{{\rm pk}}
\newcommand{\ak}{{\rm ak}}
\newcommand{\cd}{{\rm cd}}
\newcommand{\ex}{{\rm ex}}
\newcommand{\br}{{\rm br}}
\newcommand{\ms}{{\rm ms}}
\newcommand{\tms}{{\rm tms}}
\newcommand{\isb}{{\rm isb}}
\newcommand{\hor}{{\rm hor}}
\newcommand{\ver}{{\rm ver}}
\newcommand{\dis}{{\rm dis}}
\newcommand{\dcl}{{\rm dcl}}

\newcommand{\mt}{{\rm mt}}
\newcommand{\mg}{{\rm mg}}
\newcommand{\uf}{{\rm uf}}
\newcommand{\Av}{{\rm Av}}
\newcommand{\CH}{{\rm CH}}
\newcommand{\dms}{{\rm dms}}
\newcommand{\dws}{{\rm dws}}
\newcommand{\prob}{{\rm prob}}
\newcommand{\wm}{{\rm wm}}
\newcommand{\wmg}{{\rm wmg}}
\newcommand{\po}{{\rm po}}
\newcommand{\pot}{{\rm pot}}
\newcommand{\poe}{{\rm poe}}
\newcommand{\ind}{{\rm ind}}
\newcommand{\even}{{\rm even}}
\newcommand{\av}{{\rm av}}
\newcommand{\fs}{{\rm fs}}

\newcommand{\na}{{\rm na}}
\newcommand{\true}{{\rm true}}

\newcommand{\false}{{\rm false}}
\newcommand{\acl}{{\rm acl}}
\newcommand{\ids}{{\rm ids}}
\newcommand{\qf}{{\rm qf}}
\newcommand{\Dl}{{\rm Dl}}

\newcommand{\truth}{{\rm truth}}

\newcommand{\nsp}{{\rm nsp}}
\newcommand{\lev}{{\rm lev}}

\newcommand{\at}{{\rm at}}

\newcommand{\tr}{{\rm tr}}
\newcommand{\pr}{{\rm pr}}

\newcommand{\incr}{{\rm increasing}}

\newcommand{\id}{{\rm id}}

\newcommand{\Mod}{{\rm Mod}}

\newcommand{\dom}{{\rm dom}}

\newcommand{\sep}{{\rm sp}}

\newcommand{\Dom}{{\rm Dom}}

\newcommand{\Rang}{{\rm Rang}}
\newcommand{\rang}{{\rm rang}}

\newcommand{\rest}{{\restriction}}

\newcommand{\iif}{{\rm if}}
\newcommand{\wilog}{{\rm without loss of generality}}
\newcommand{\Wilog}{{\rm Without loss of generality}}

\newcommand{\then}{{\underline{then}}}
\newcommand{\when}{{\underline{when}}}

\newcommand{\Then}{{\underline{Then}}}

\newcommand{\If}{{\underline{if}}}
\newcommand{\Iff}{{\underline{iff}}}
\newcommand{\mn}{{\medskip\noindent}}
\newcommand{\sn}{{\smallskip\noindent}}

\newcommand{\bbF}{{\mathbb F}}
\newcommand{\gb}{{\mathfrak b}}
\newcommand{\gA}{{\mathfrak A}}
\newcommand{\gB}{{\mathfrak B}}
\newcommand{\gK}{{\mathfrak K}}
\newcommand{\gk}{{\mathfrak k}}
\newcommand{\gC}{{\mathfrak C}}
\newcommand{\cC}{{\mathscr C}}
\newcommand{\cL}{{\mathscr L}}

\newcommand{\cD}{{\mathscr D}}
\newcommand{\cE}{{\mathscr E}}

\newcommand{\cH}{{\mathscr H}}
\newcommand{\cJ}{{\mathscr J}}

\newcommand{\cI}{{\mathscr I}}

\newcommand{\bbL}{{\mathbb L}}
\newcommand{\bbI}{{\mathbb I}}
\newcommand{\bbN}{{\mathbb N}}

\newcommand{\bbP}{{\mathbb P}}
\newcommand{\bbR}{{\mathbb R}}

\newcommand{\cP}{{\mathscr P}}
\newcommand{\gp}{{\mathfrak p}}
\newcommand{\varp}{{\varepsilon}}
\newcommand{\gs}{{\mathfrak s}}

\newcommand{\cS}{{\mathscr S}}

\newcommand{\gt}{{\mathfrak t}} 
\newcommand{\cU}{{\mathscr U}}

\newcommand{\cY}{{\mathscr Y}}

\newcommand{\cf}{{\rm cf}}

\newcount\skewfactor
\def\mathunderaccent#1#2 {\let\theaccent#1\skewfactor#2
\mathpalette\putaccentunder}
\def\putaccentunder#1#2{\oalign{$#1#2$\crcr\hidewidth
\vbox to.2ex{\hbox{$#1\skew\skewfactor\theaccent{}$}\vss}\hidewidth}}

\newenvironment{PROOF}[2][\proofname.]
   {\begin{proof}[#1]\renewcommand{\qedsymbol}{\textsquare$_{\rm #2}$}}
   {\end{proof}}

\begin{document}

\title {Compactness of the Quantifier on ``Complete embedding of BA's"}
\author {Saharon Shelah}
\address{Einstein Institute of Mathematics\\
Edmond J. Safra Campus, Givat Ram\\
The Hebrew University of Jerusalem\\
Jerusalem, 91904, Israel\\
 and \\
 Department of Mathematics\\
 Hill Center - Busch Campus \\ 
 Rutgers, The State University of New Jersey \\
 110 Frelinghuysen Road \\
 Piscataway, NJ 08854-8019 USA}
\email{shelah@math.huji.ac.il}
\urladdr{http://shelah.logic.at}
\thanks{The author thanks Alice Leonhardt for the beautiful typing.
First version: January 1993; May, 1993; corrected May 7, 1993;
next correction Nov 1993, Jan 2002; 
next correction introduced Nov 2002; next Dec 2002, LAST May 2003.
Publication 482.}




\date {January 13, 2016}

\begin{abstract}
We try to build, provably in ZFC, for a first order $T$ a model
in which any isomorphism between two Boolean algebras is definable. 
The problem, compared to [SH:384], is with pseudo-finite Boolean
algebras. A side benefit is that we do not use Skolem function 
(which do not matter for proving compactness of logics but still 
are of interest).
Let $\lambda$ be $2^\mu$ if regular and its successor otherwise. 
Model theoretically we investigate notions of bigness of types, usually 
those are ideals of the set of formulas in a model, definable in 
appropriate sense.  We build a model of cardinality $\lambda^+$ by a sequence
of models $M_\alpha$ of cardinally $\lambda$ for $\alpha < \lambda^+$, 
each $M_\alpha$ equips with a sequence $\langle (M_{\alpha,i},a_{\alpha,i},
\Omega d_{\alpha,i}):i \in S_I \subseteq \lambda \rangle$,
with $M_{\alpha,i}$ is of cardinality $< \lambda,\prec$-increasing
continuous with $i,\Omega_{\alpha,i}$ a bigness notion defined using 
parameters from  $M_{\alpha,i}$ and $a_{\alpha,i}$
realized in $M_{\alpha,i+1}$ over $M_{\alpha,i}$ a 
$\Omega_{\alpha,i}$-big type. As $ \alpha $ increase,
not only $M_\alpha$ increase, but this extra structure increasing
modulo a club of $\lambda$, this is why we have insisted on
$\lambda$ being regular.

This can be considered as a way to omit types of cardinality
$\lambda$, which in general is hard.  The fact that $\lambda$ is not 
too much larger than $\mu$ help us to guarantee that any
 possible automorphism of structures be
defined in $M= \cup \{M_\alpha:\alpha < \lambda^+\}$
by approximations of cardinal $\mu$ and so we can enumerate them all.

The bigness notions involved has to relate to the kind of structures we are
interested in interpreting in $M$, e.g. dense for linear orders and
subsets of $\cP(n)$ of large cardinality for $n$ pseudo finite.
During the construction for each $\alpha$, some $b_\alpha \in
M_{\alpha + 1}$ realizes a big type
over $M_\alpha$ for an appropriate bigness notion.
We have to guarantee that the bigness notions
used in the horizontal direction (that is $\alpha < \lambda^+$)
and the bigness notions used in the vertical direction
(that is for $i < \lambda$) do not interact.
So we have to prove that enough pairs of bigness
notions are so called orthogonal:
if $p_\ell(x_\ell) \in {\bold S}(M)$ is $\Omega_{\ell}$-big
then we can find $p(x_0,x_1) \in {\bold S}^2(M)$
extending both such that it says that ``$x_\ell$ is $\Omega_\ell$-big
over $M+ x_{1-\ell}"$.
\end{abstract}

\maketitle
\numberwithin{equation}{section}
\setcounter{section}{-1}
\newpage

\section {Introduction}

We continue here the attempt of extracting and strengthening the purely
combinatorial content of \cite[Ch.VIII,\S2]{Sh:a} = \cite[Ch.VIII,\S2]{Sh:c}
(i.e., many models for unsuperstable theories), as done in
\cite[\S2]{Sh:E59}, \cite[\S1,\S2]{Sh:331} (and \cite[\S3]{Sh:331}); 
(so also \cite{Sh:136}) or \cite{Sh:E60}.

In \cite{Sh:384} we succeed to get ``complicated models'' by
omitting small (e.g. countable) types, so that building a model 
of size $\lambda^+$
by a sequence of $\lambda^+$ approximation each of size $\lambda$,
$\lambda$ regular we suffice to guess, e.g. when $\lambda =
\lambda^{\aleph_0}$. But we have been stuck on the
problem of automorphism of pseudo finite Boolean algebras. 
Here we use a different approach, building a model of size
$\lambda^+,\lambda$ is, e.g. $(2^\kappa)^+$; so we can enumerate all
subsets of size $\kappa$, and instead of guessing automorphism on
$\gB_\alpha$ we try to make the model code them by a subset of size
$\kappa$, so we can enumerate them.  See general
construction in \S4, our specific construction in \S5.
 
The model $\gB_{\alpha+1}$ is build over ${\gB}_\alpha$ as
an \incr\ sequence of length $\lambda$ of approximations, each a type 
$p^\alpha_i(\langle \bar x_j: j<i\rangle)$ over $\gB_\alpha$
of cardinality $< \lambda$ for $i<\lambda$, restricted by being ``big" in
appropriate sense.  Bigness notions are defined in \S1, bigness 
notion of general type are investigated in \S2 and more specific ones in \S3.

But how do we omit types? Generally we do not know, types of some special
forms we know: we represent $\gB_\alpha$ as $\langle \bar a^\alpha_i:i<
\lambda\rangle$, and demand that for a club of $\alpha < \lambda$, in 
$\gB_{\alpha+1}$ the type $\tp(\bar a^\alpha_i,\bigcup\limits_{j< i}(\bar a
^\alpha_j\cup \bar x^\alpha_j),\gB_{\alpha+1})$ is big in appropriate sense.

To a large extent here we continue \cite{Sh:72}, \cite{Sh:73} rather than
\cite{Sh:107}. In the later we use a general omitting types theorem for
$\lambda^+$, quite powerful but it depends on $(\Dl)_\lambda$ 
(hence necessarily $\lambda=\lambda^{<\lambda}$). 
In \cite{Sh:72}, \cite{Sh:73} we use a special way to
omit types: we build a model of cardinality $\lambda^+$, by an
increasing chain of $\gB_\alpha$ for $\alpha< \lambda$, and the
omitting of types in stage $\alpha$ has the form: the type is
represented by a stationary $S \subseteq \lambda$ and
$\langle \langle a^\beta_i: i< \omega\rangle:
\beta\in S\rangle$ with $a^\beta_i$ or $a^\beta_i \in \gB_\alpha$ and
we ``promise" that for every $\beta \ge \alpha$ and finite 
$A \subseteq {\gB}$, $\{\beta\in S: \langle a^\beta_i: i<
\omega\rangle$ is not discernible over $A\}$ is not
stationary. Such properties are preserved in any limit stage, even of small
cofinality, the problematic case. For
wider framework we use ``bigness of types'' as in \cite{Sh:107}, but here the
restriction of bigness act in two ways: ``horizontally'', building
${\gB}_\alpha$ by a sequence of $\Omega_\alpha$-big types over 
$\gB_\alpha$, and ``vertically'', preserving: for $\alpha< \beta<
\lambda^+$ and any finite $A \subseteq {\gB}_\alpha$, for a club of
$\delta< \lambda$, the element $a^\alpha_i$ realizes a 
$\Gamma^*_i$-big type over $A$.
To be able to do it we need the so called ``orthogonality".
See more in \cite{Sh:800} and see history in \cite[\S0]{Sh:384}.

This paper was supposed to be Ch.XI to the book
``Non-structure" and probably will be if it materializes, it has been
circulated and lectured on since 1993.

Our main results are on models of $T$ which is first order complete
coding enough set theory, see clause (B) of \ref{a2}.  This is enough
for proving compactness of first order logic extended by suitable
second order quantifiers.  Probably we can get all those results for
any (first order complete) $T$, i.e. as in \ref{a2}(A), using reducts
of global bigness notions but this is delayed.

The intentions were \cite{Sh:E58} (revising \cite{Sh:229}) for
Ch.I, and \cite{Sh:421} for Ch.II and
\cite{Sh:E59} for Ch.III and
\cite{Sh:309} for Ch.IV and
\cite{Sh:363} for Ch.V  and
\cite{Sh:331} for Ch.VI  and
\cite{Sh:511} for Ch.VII and
\cite{Sh:E60}, a revision of \cite{Sh:128}, for Ch.VIII,
and \cite{Sh:E62}, for the appendix and \cite{Sh:757}, \cite{Sh:384},
\cite{Sh:482} and \cite{Sh:800}, for Ch. IX, X, XI, XII respectively.
References like \cite[3.7=Lc2]{Sh:E62} means that c2 is the label of
3.7 in \cite{Sh:E62}, will only help the author if changes in the
paper \cite{Sh:E62} will change the number. 

\begin{notation}
\label{x2}
1) Let $A+{\bar a}$ be $A\cup\Rang (\bar a)$, similarly $A+a$, 
$A+{\bar a}+b +C$.

\noindent
2) Let $\cL$ denote a logic, $\tau$ a vocabulary, $\cL(\tau)$
the language i.e. the set of $\cL$-formulas in the vocabulary 
$\tau,\tau_M$ is the vocabulary of the model $M$; $\cL(\tau,A)$ means
we add all members of $A$ as individual constants.

\noindent
3)  Let ${\bbL}$ be first order logic.

\noindent
4) $\cL(\dot{\bold Q})$ means we add to to the logic $\cL$ the quantifier 
$\dot{\bold Q}$.

\noindent
5) Let $T$ denote a theory, first order if not said otherwise, 
usually complete.

\noindent
6) For $T$, $\tau(T)=\tau_T$ is its vocabulary, $\cL(T)= \cL(\tau(T))$ 
the corresponding language (first order for $\bbL(T) = \bbL(\tau_T)$)

\noindent
7) If $A \subseteq M$, $M\models T$ then
\mn
\begin{enumerate}
\item[$(a)$]  $T[A] = \{\varphi(\bar a):
\bar a \in {}^{\omega >}A,\ \ M \models \varphi [\bar a]\}$
\sn
\item[$(b)$]  so $\tau(T[A])=\tau(T)\cup A$
\sn
\item[$(c)$]  $\acl(A,M) = \{b \in M:\tp(b,A,M)$ is algebraic,
  i.e. some formula in it is reazlied by finitely many elements$\}$.
\end{enumerate}
\mn
8) We say $\gp$ is a type definition over $N$ \when \,:
\mn
\begin{enumerate}
\item[$(a)$]  $\gp$  is an ultrafilter on ${}^\alpha N$, and if 
$N \subseteq  A \subseteq  M,N \prec M$ then ${\gp}^A = 
\{\varphi(\bar x,\bar a):\bar a \in {}^{\omega >}A$ and 
$\{\bar b \in {}^\alpha N: M\models \varphi[\bar b,\bar a]\} \in
\gp\}$ or
\sn
\item[$(b)$]  $\gp$  is a function from
$\{\langle \varphi(\bar x,\bar y),q(\bar y)\rangle:\varphi$ a formula, 
$q(\bar y)$  a complete type over $N\}$ to $\{$truth,false$\}$ and if 
$N \subseteq  A \subseteq M,N \prec M$ then 
${\gp}^A = \{\varphi(\bar x,\bar a):\bar a \in {}^{\omega >}A$
and ${\gp}(\langle \varphi(\bar x,\bar y),\tp(\bar a,N,M)\rangle) =
\truth\}$.
\end{enumerate}
\mn
8A) Above we say ${\gp}$ is of kind $\nonfork{}{}_{\fs}$
or of kind $\nonfork{}{}_{\nsp}$ respectively.
In the first section we use an arbitrary compact logic ${\cL}$
but the reader may concentrate on ${\bbL}$, first order logic.
\end{notation}
\newpage

\section {Bigness notions: basic definitions and properties.}

For a complete first order theory $T$, a $\ell$-bigness notion
$\Gamma$ ($\ell$ for local) 
is a scheme defining for every model $M$ of $T$, an ideal of the set
of formulas $\varphi(\bar x_\Gamma, \bar a)$, $\bar a \subseteq M$ .

We are interested in such ideals preserved by elementary
embedding. Such notions play crucial role in our construction of
models. This section is soft - just giving definitions and easy
applications. If this section is too abstract, the reader can read
parallelly \S2, \S3 which deal with examples.

The reader may concentrate on the case $\gk = ({\bold K},\le)$
is the class of models of a complete  first order theory $T$ with Skolem
functions, $\le = \prec,\gC=\gC_T$ a monster model 
and on the case of simple $\ell$-bigness notion.

\begin{context}
\label{a2}
Let $\tau$ be a vocabulary, $\gk = (\bold K_{\gk},\le_{\gk}) 
= (\bold K,\le)$ a 
class of $\tau$-models so $M$ means $M \in \bold K_{\gk}$ 
(usually the class of models of a fixed first order theory, $M \le_{\gk} N$ iff
$M \prec N$), $\cL$ a logic. Always satisfaction of ${\cL}$-formulas is
preserved by extensions and the pair $({\gk},{\cL})$ is compact; see below. 
Usually $({\gk},{\cL})$ is one of the following:
\mn
\begin{enumerate}
\item[$(A)$]  $T$ is complete first order, $\bold K$ 
the class of models of $T$ i.e.
${\bold K} = \Mod(T)$ (of course $\tau = \tau(T)$-models) and $M \le_{\gk} N
\Leftrightarrow  M \prec N,\cL = \bbL$ first order logic
and $\gC = \gC_T$ is a monster model $\bold S^\alpha(A,M,\gk) = \bold
S^\alpha_{\gk}(A,M) = \{p:p$ a complete type over $A$ in $M\}$
\sn
\item[$(B)$]  like (A), but we may denote $T$ by $T^*$ and it has a model
$\gC^*$, an expansion of $({\cH}(\chi^*),\in,<^*)$ where
$\chi^*$ is strong limit cardinal, $<^*$ a well ordering of
${\cH}(\chi)$ see \cite[2.1]{Sh:384}. Then $\gC$ 
denote a ``monster'' model of $T^*$ 
and $\dot e$ denote the membership inside it, 
i.e. $\dot e^{\gC^*} = \in \rest{\cH}(\chi)$. 
We call such $T$ ``of set theory character". Saying ``$n$" we mean
a true natural number but also its interpretation in $\gC$, we
use $\dot n,\dot m,\dot k$ for members of $\bbN^{\gC}$.  So many times
it is better to deal with ``sets" not classes.
\sn
\item[$(C)$]   $T$ is a universal first order theory with amalgamation,
$\le_{\gk}$ is being a submodel (so $\gk = (\Mod_T,\subseteq))$ and
$\cL = \bbL$ and $\bold S^\alpha (A,M,\bold K) = 
\{\tp_{\qf}(\bar{a},A,N):M \le_{\gk} N,N\models T,\bar a \in {}^\alpha N\}$
\sn
\item[$(D)$]  $T$ is a universal first order theory, ${\bold K}$ 
is the family of
existentially closed models of $T,\le_{\gk}$ is being a submodel, 
let $\cL(\tau) = \Sigma(\tau) = 
\{\varphi:\varphi$ an existential first order formula
in the vocabulary $\tau\}$ and so $\bold S^\alpha_{\gk} (A,M) =
\{\tp_{\Sigma(\tau)}(\bar{a},N,N):M \le_{\gk} N,N \models T, 
\bar a \in {}^\alpha N\}$. 
\end{enumerate} 
\end{context}

\begin{definition}
\label{a5}
1)  For $A \subseteq M, \bar{a}$ a sequence from
$M$ let
\mn
\begin{enumerate}
\item[$(a)$]  $\tp_{\cL}(\bar a, A, M)=\{\varphi(\bar x, \bar b): \bar b\in
{}^{\omega>} A,\ell g(\bar x) = \ell g(\bar a), 
M \models \varphi[\bar a,\bar b]$ and
$\varphi\in \cL(\tau_M)\}$
\sn
\item[$(b)$]  ${\bold S}^\alpha_{\cL} (A,M)=\{p:p$ is a maximal set of 
${\cL}(\tau_M)$-formulas with free variables among 
${\bar x}=\langle x_i:i<\alpha\rangle$ and parameters from A and every
  finite subset is realized in some $N$, satisfies $M \le_{\gk} N\}$;
if $\cL = \bbL$ we may omit it; the length of 
$\bar{a}$, i.e. of $\bar{x}$ is not necessarily finite.  Writing $p
\in \bold S_{\cL}(A,M)$ means for some $\alpha$ clear from the context
\sn
\item[$(c)$]  $\bold S^\alpha_{\cL}(A,M) = \{\tp(\bar a,A,N):M
  \le_{\gk} N$ and $\bar a \in {}^\alpha N\}$.
\end{enumerate}
\mn
2) $({\gK}, {\cL})$ is $\cL$-compact means:
\mn
\begin{enumerate}
\item[$(a)$]   if a set $\Delta$ of $\cL$-formulas
with parameters in $M \in \bold K$ and is finitely satisfiable in $M$
\then \, it is realized in some $N$, $M \le_{\gk} N \in \bold K$
\sn
\item[$(b)$]  if $p \in \bold S^\alpha_{\cL} (A,M)$ then\footnote{is this
  not covered by clause(a)? not always! see, e.g. \ref{a7}(B)}
 $p$ is realized in some $N$ extending $M$ (i.e. $M \le_{\gk} N \in
 {\bold K}$).
\end{enumerate}
\mn
We let $\gC= \gC_T$ be a monster.
\end{definition}

\begin{discussion}
\label{a7}
We may consider several further general contexts:
\mn
\begin{enumerate}
\item[$(A)$]   The class of models of $T$, a complete countable theory in
$\cL = \bbL(\bold Q)$ (with $\bold Q x$ 
interpreted as $(\exists^{\ge\lambda}x)$, usually
$\lambda=\aleph_1$), with elimination of quantifiers for simplicity,
$M \le_{\gk} N$ if $M \prec N$ and $M \models \neg \bold Q x \varphi
(x,\bar a) \Rightarrow \varphi (N,\bar a)\subseteq M$
\sn
\item[$(B)$]   $T$ is first order complete, $D \subseteq D(T) =: 
\{\tp (\bar a,\emptyset,M):M \models T,\bar a \in {}^{\omega>}M\}$, 
and $D$ is good (see \cite{Sh:3}), $\gC$ a monster model,
i.e. $(D,\bar\kappa)$-sequence-homogeneous model, $K_{\gk} = \{M:M$ a
model of $T$ such that every finite $\bar a \subseteq M$ realizes a
type from $D\}$ and $\le_{\gk} = \prec$; in this case as well as in
(C), it is natural to use global bigness notions
\sn
\item[$(C)$]   $\bold K$ is a universal class (i.e. $M \in \gk$ iff for 
every $\bar a \in {}^{\omega>} M,M \restriction c \ell_M(\bar a)
\in \bold K$), the relation $\le_{\gk}$ is being submodel (i.e. locally finite 
models of a universal theory), 
\[
\bold S_\cL(A,M,\gk) = \{\tp_{\cL}(\bar a,A,N):M \le_{\gk} N \in
{\bold K},
\bar a \subseteq N\}
\]

\[
\bold S_{\cL}^{\alpha}(A,M,\gk) = \{\tp_{\cL}(\bar a,A,N):M \le_{\gk} 
N \in \bold K,
\bar a \in {}^{\alpha} N\}.   
\]
\end{enumerate}
\mn
E.g. the class of locally finite gruopos (or existentially closed
ones), see \cite{Sh:312}
\mn
\begin{enumerate}
\item[$(D)$]    For some class $\bold K' 
\subseteq \bold K$ and partial order $\le'$ on
$\bold K'$ the union of a $\leq'$-directed system of models from
$\bold K'$
(were $\bold K'$, $\le'$ closed under isomorphism satisfying natural
  conditions)
\sn
\item[$(E)$]   Abstract elementary class (amalgamation is not demanded).
\end{enumerate}
\mn
Those contexts will not be used, but we may remark on them or give examples.
\end{discussion}

\begin{definition}
\label{a11}
1) We call $\Gamma$ a $l$-bigness (=local bigness) notion for
$(\gk,\cL)$ (with set of parameters $A_\Gamma \subseteq M_\Gamma \in
\bold K$; for
simplicity, we usually restrict ourselves to $\{M \in \bold K:M_\Gamma 
\le_{\gk} M\}$ or $\gC$ is a monster model for $T, A_\Gamma \subseteq
\gC$ is ``small") \If \, (it gives a sequence $\bar x=\bar x_\Gamma$ of
variables of length $\alpha(\Gamma)$, in the usual case
singleton $x$ or at least finite and):
\mn
\begin{enumerate}
\item[$(a)$]   for every $M \in \bold K$ (such that $M_\Gamma \le_{\gk} M$),
$\Gamma_M^-= \Gamma^{-} (M)$ is a subset of the family of
formulas $\varphi(\bar x,\bar a)$, $\varphi \in \cL(\tau)$, $\bar
a \subseteq M$, and $\Gamma_M=\Gamma_M^+ =\Gamma^+(M)$  is the complement
of $\Gamma^-(M)$ inside this family
\sn
\item[$(b)$]  $\Gamma_M^-$ is preserved by automorphisms of $M$ over
$A_\Gamma$
\sn
\item[$(c)$]  $\Gamma_M^-$ is a proper ideal, i.e.
\sn
\begin{enumerate}
\item[($\alpha$)] if $M \models (\forall\bar x)[\varphi(\bar x,\bar a)
  \rightarrow \psi(\bar x,\bar b)]$ and $\psi(\bar x,\bar b)\in\Gamma_M^-$
\then \, $\varphi(\bar x,\bar b) \in \Gamma_M^-$
\sn
\item[($\beta$)]  if $\varphi_1(\bar x,\bar a_1)$, $\varphi_2(\bar x,
\bar a_2) \in \Gamma_M^-$ \then \, $\varphi_1(\bar x,\bar
a_1)\lor\varphi_2(\bar x,\bar a_2) \in \Gamma_M^-$
\sn
\item[$(\gamma)$]  $\Gamma^+_M \ne \emptyset$.
\end{enumerate}
\end{enumerate}
\mn
2)  Assume $\bar x_\Gamma$ is finite \then \, $\Gamma$ is called 
non-trivial if, when $M_{\Gamma} \le_{\gk} M \in \bold K$.
\mn
\begin{enumerate}
\item[$(*)$]   $\bar x=\bar a \in\Gamma_M^-$.
\end{enumerate}
\mn
3)  We call members of $\Gamma_M^-$ ``$\Gamma$-small in $M$",
members of $\Gamma_M^+,``\Gamma$-big in $M$". We may write
$M \models (\bold Q^\Gamma \bar x) \varphi(\bar x,\bar a)$ for ``$\varphi(\bar
x,\bar a)$ is $\Gamma$-big in $M$" and $M \models (\bold Q_\Gamma\bar
x)\varphi(\bar x,\bar a)$ for ``$\varphi(\bar x,\bar a)$ is
$\Gamma$-small in $M$", (in this notation, $\Gamma_1 \perp \Gamma_2$ 
defined below essentially means that $\bold Q^{\Gamma_1},
\bold Q^{\Gamma_2}$ commute).

\noindent
4)  A $\Gamma$-big type $p(\bar x)$ is a set of formulas $\psi(\bar
x,\bar a)$, any finite conjunction of which is $\Gamma$-big.
\end{definition}

\begin{example}
\label{a14}
Let $\chi$ be an infinite cardinal, and let $\exists^{\ge \chi}\bar
x\varphi(\bar x,\bar y)$ be the formula which says ``at least $\chi$
pairwise disjoint sequences $\bar x$ satisfy $\varphi(\bar x,\bar y)$" and
$\cL = \bbL(\exists^{\ge\chi})$. Consider a theory $T$
in $\cL(\tau)$, \wilog \, every formula is equivalent to a 
predicate, and $T' = T \cap \bbL$ (so not exactly in the
context \ref{a2}(1) for $T$).  This naturally defines a local notion
$\Gamma_T$ of bigness for $T'$, for $M$ a model of $T,\varphi(x,\bar
a)$ is $\Gamma$-big iff $M \models R_\varphi(\bar a)$ where $(\forall
\bar y)(\exists^{\ge \chi} \bar x,\varphi(x,\bar y) \equiv R(y)) \in T$.
\end{example}

\begin{convention}
\label{a17}
1) We will, abusing notation, first define bigness notions
 and only then prove they are bigness notions.

\noindent
2) As we shall deal here only with invariant bigness notions 
[see definition below] we may ``forget" this adjective.
\end{convention}

\begin{definition}
\label{a20}
Let $\Gamma$ be a local bigness notion for $(\gk,\cL)$.

\noindent
1) We say that $\Gamma$ is weakly invariant \If \,  for every
$\varphi(\bar x,\bar y) \in \cL(\tau),\bar a \in M,M \le_{\gk} N$ (models in
$\bold K$) we have: $\varphi(\bar x,\bar a)$ is $\Gamma$-big in $M$ iff
$\varphi(\bar x,\bar a)$ is $\Gamma$-big in $N$.

\noindent
2)  We say $\Gamma$ is invariant \If \, $\tp_\cL (\bar a,A_\Gamma,M)$ 
and $\varphi(\bar x,\bar y)$ determine whether $\varphi(\bar x,\bar a)$ is
$\Gamma$-big (in $M$).

\noindent
3) We say that $\Gamma$ is $\lambda$-strong [or $\lambda$-co strong]
\If \, for every $M \in \bold K$ and $\varphi(\bar x,\bar a)$ which is 
$\Gamma$-big [or $\Gamma$-small] in $M$
there is $\tau_\varphi \subseteq \tau_{\gk},|\tau_\varphi|<\lambda$ such that:
$\varphi (\bar x,{\bar a}')$ is $\Gamma$-big [or $\Gamma$-small] in
$M'$ whenever ${\bar a}'\subseteq M' \in \bold K$, ${\bar a}'$ realized in
$M'$ the $\cL (\tau_\varphi)$-type which $\bar a$ realize over $A_\Gamma$ in
$M\restriction\tau_\varphi$.

\noindent
4) We say $\Gamma$ is very $\lambda$-strong  \If \,  
for every $\varphi = \varphi(\bar x,\bar y) \in \cL(\tau_{\gk})$ 
there is $\tau_\varphi \subseteq \tau_{\bold K},|\tau_\varphi|<\lambda$ 
such that: for every $\varphi(\bar x, \bar a), 
\bar a \subseteq M \in \bold K$ the type 
$\tp_{\cL}(\tau_\varphi)(\bar a,A_\Gamma,M \restriction \tau_\varphi)$
determine if $\varphi(\bar x,\bar a)$ is $\Gamma$-small or $\Gamma$-big in $M$.
\end{definition}

\begin{definition}
\label{a23}
Let $\Gamma$ be a local invariant bigness notion.

\noindent
1)  We say that $\Gamma$ is $\lambda$-simple [or $\lambda$-co-simple]
\If \, for every $M \in \bold K,\bar a \subseteq M$,
and $\varphi$ such that $\varphi(\bar x,\bar a)$ is $\Gamma$-big [or
$\Gamma$-small] in $M$ there is $q \subseteq \tp_{\cL} 
(\bar a,A_\Gamma,M),|q|<\lambda$ such that:

\[
[\bar a' \in M' \in \bold K \, \& \, \bar a' \text{ realizes } q
\text{ in } M' \Rightarrow \varphi(\bar x,\bar a') \text{ is }
\Gamma \text{-big} \text{ [or } \Gamma \text{-small] in } M'].
\]

\mn
If $\lambda=\aleph_0$ we may omit it.

\noindent
2  $\Gamma$ is very $\lambda$-simple \If \, for every
$\varphi(\bar x,\bar y)$ there is a set $\Delta$ of $<\lambda$ formulas
of the form $\psi(\bar y; \bar a) \in \cL (\tau, A_\Gamma)$, with
$\bar a\subseteq A_\Gamma$ such that: if $M \in \bold K$ and

\[
\Delta \cap \tp_{{\cL}} (\bar b^1, A_\Gamma, M)=\Delta \cap
\tp_{\cL}(\bar b^2, A_\Gamma, M)
\]

\mn
\then \,:
$\varphi(\bar x, \bar b^1)$ is $\Gamma$-big iff $\varphi (\bar x,
\bar b^2)$ is $\Gamma$-big.
If $\lambda=\aleph_0$ we may omit it.

\noindent
3) We say $\Gamma$ is uniformly $\lambda$-simple [or $\lambda$-co-simple]
\when \, for any $\varphi(\bar x, \bar y)$ with $\lg(\bar x)=
\alpha(\Gamma)$ there is a type $q_\varphi(\bar y)$ over $A_\Gamma$
such that: for any relevant $M,\bar a$ we have
$\varphi(\bar x,\bar a)$ is $\Gamma$-big iff $\bar a$ realizes
$q_\varphi$ [or: $\varphi(\bar x, \bar a)$ is $\Gamma$-small iff $\bar
a$ realizes $q_\varphi$].
\end{definition}


\begin{claim}
\label{a27}
Let $\Gamma$ be a local bigness notion.

\noindent
1) If $\lambda > |A_\Gamma|+ |{\cL}|$, $\cL$ has finite
occurrence\protect\footnote{this means that every $\varphi \in {\cL}
(\tau)$ for some finite $\tau' \subseteq \tau$; note $|{\cL}|$ 
is the number of sentences up to
renaming the predicates and function symbols, so $|{\cL} (\tau)| \le
|\tau|+ |{\cL}| + \aleph_0.$} \then
\mn
\begin{enumerate}
\item[$(a)$]   $\Gamma$ is $\lambda$-strong \underline{ if and only
  if} $\Gamma$ is $\lambda$-simple 
\sn
\item[$(b)$]   $\Gamma$ is $\lambda$-co-strong \underline{if and only
  if} $\Gamma$ is $\lambda$-co-simple
\sn
\item[$(c)$]   $\Gamma$ is very $\lambda$-strong \underline{if and
  only if} $\Gamma$ is very $\lambda$-simple.
\end{enumerate}
\mn
2) 
\mn
\begin{enumerate}
\item[$(a)$]   If $\Gamma$ is very $\lambda$-simple \then \, $\Gamma$ is
$\lambda$-simple and $\lambda$-co-simple
\sn
\item[$(b)$]   If $\Gamma$ is very $\lambda$-strong \then \,
$\Gamma$ is $\lambda$-strong and $\lambda$-co-strong
\sn
\item[$c)$]  If $\Gamma$ is $\lambda$-simple 
\then \, $\Gamma$ is $\lambda$-strong
\sn
\item[$(d)$]  If $\Gamma$ is $\lambda$-co-simple \then \, $\Gamma$
is $\lambda$-co-strong
\sn
\item[$(e)$]  If $\Gamma$ is very $\lambda$-simple \then \,
$\Gamma$ is very $\lambda$-strong
\sn
\item[$(f)$]  If $\lambda_1 <\lambda_2$ and $\Gamma$ is
$\lambda_1$-strong \then \, $\Gamma$ is $\lambda_2$-strong; similarly for
``$\lambda_\ell$-co-strong", ``very $\lambda_\ell$-strong",
``$\lambda_\ell$-simpleq", ``$\lambda_\ell$-co-simple", ``very
$\lambda_\ell$-simple"
\sn
\item[$(g)$]  if $\lambda > |T|+\aleph_0$ \then \, the
corresponding strong and simple properties are equal.
\end{enumerate}
\mn
3) If $\Gamma$ is $\lambda$-simple and co-$\lambda$-simple
\then \, $\Gamma$ is very simple (using the logic ${\cL}$ being
compact).

\noindent
4)  If $A \subseteq M$, $p$ a $\Gamma$-big type over $A$ in $M$ \then
\, we can find a $\Gamma$-big $q \in \bold S_\cL(A,M)$ extending $p$ 
using the logic ${\cL}$ being compail.

\noindent
5) Parallel of parts (1), (2), (3) hold for global bigness
notion defined below.

\noindent
6) If $\lambda > |{\cL} (\tau)|+ |A_\Gamma|$  \then \,: 
$\Gamma$ is very $\lambda$-strong, very $\lambda$-simple.

\noindent
7) If $\Gamma$ is uniformly $\lambda$-simple \then \, $\Gamma$
is co-simple. If $\lambda$ is uniformly $\lambda$-co-simple then $\Gamma$ is
simple.
\end{claim}

\begin{PROOF}{\ref{a27}}
By the definitions (and compactness when demanded).
\end{PROOF}

\begin{definition} 
\label{a31}
1) We say $\Gamma$ is a  $g$.($=$global) bigness notion for $(\gk,\cL)$
(with set parameters $A_\Gamma\subseteq M_\Gamma \in \bold K$) \If \, 
(it gives a sequence $\bar x=\bar x_\Gamma$ of variables and):
\mn
\begin{enumerate}
\item[$(a)$]  for every $M \in \bold K$ (satisfying $M_\Gamma \le_{\gk} M$
  as usual), 
$\Gamma_M$ is a family of types $p(\bar x)$ such that for some 
$A$, $A_\Gamma \subseteq A \subseteq M$ and
 $p(\bar x)\in \bold S_\cL^{\alpha(\Gamma)}(A,M)$
\sn
\item[$(b)$]  \underline{local character}: if $M \in \bold K,A_\Gamma\subseteq
A\subseteq M$, $p(\bar x)\in \bold S_\cL(A,M)$ \then \, 
$p(\bar x) \in \Gamma_M \Leftrightarrow (\forall \text{ finite }
B\subseteq A) [p(\bar x)\restriction (A_\Gamma\cup B)\in\Gamma_M]$
\sn
\item[$(c)$]  \underline{the extension property}: if $A_\Gamma\subseteq
A\subseteq B\subseteq M$, $p \in \bold S^{\alpha(\Gamma)}_\cL(A,M)$ 
is in $\Gamma_M$ \then \, some extension $q \in \bold S_\cL(B,M)$ of $p$ is in
$\Gamma_M$
\sn
\item[$(d)$]  \underline{existence}: if $A_\Gamma\subseteq A\subseteq M$
then there is a $\Gamma$-big $p \in \bold S^{\alpha(\Gamma)}(A, M)$. 
[We can close the family under restriction thus allowing 
``$p \in \bold S^{\alpha(\Gamma)}(A, M)$ is
$\Gamma$-big" though $A_\Gamma\not\subseteq A$.]
\end{enumerate}
\mn
2) We define ``$\Gamma$, a $g$.bigness notion is weakly
invariant/invariant" as in definition \ref{a20}(1),(2) above.

\noindent
3)  A g.bigness notion $\Gamma$ is $\lambda$-strong [or $\lambda$-co-strong]
\If \, for every $\Gamma$-big [or $\Gamma$-small] 
$p(\bar x, \bar a)=\tp_{\cL} (\bar b,\bar a,M)
\in \bold S^{\alpha(\Gamma)}_\cL (\bar a\cup A_\Gamma M)$, for 
some $\tau_p\subseteq\tau$ of cardinality 
$<\lambda$, we have: $M \le N,\tp_{\cL \restriction \tau_p}
(\bar b \char 94 \bar a, A_\Gamma,M) =
\tp_{\cL \restriction\tau_p}(\bar a',A_\Gamma,N)$
implies $\tp_{{\cL}}(\bar b'\bar a',N)$ is $\Gamma$-big [or
  $\Gamma$-small].

\noindent
4) We say $\Gamma$ is very $\lambda$-strong \If \, for any $n<\omega$ for some
$\tau_n$ for any $n$ and $p(\bar{x}) \in \bold S^{\alpha(\Gamma)} 
(A_\Gamma, M), q(\bar{y}) \in \bold S^n (A_\Gamma,M)$ with 
$\ell g(\bar{y})=n,$ there is $\tau_{p,q}
\subseteq \tau({\bold K})$ of cardinality $<\lambda$ such that if for
$\ell=1,2$ we have $M \le N_\ell,\bar{a}_\ell,\bar{b}_\ell
\in N_\ell$ realizing $q,p$ respectively in $N_\ell$ and 
$\tp_{\cL}(\tau_{p,q}) (\bar a_1 \char 94 \bar b_1,A_\Gamma,N_1) 
= \tp_{\cL}(\tau_{p,q}) ({\bar a}_2 \char 94 \bar b_2,A_\Gamma,N_2)$
  \then \, $\tp_{\cL} (\bar{b}_1,\bar{a}_1 \cup A_\Gamma,N)$ is $\Gamma$-
big \Iff \, $\tp_{\cL}(\bar{b}_2,\bar{a}_2 \cup A_\Gamma, N_2)$ is
$\Gamma$-big.

\noindent
5) A $g$.bigness notion $\Gamma$ is very $\lambda$-simple
\If \, for every $m$ there is a set $\Delta$ of $<\lambda$ 
formulas $\psi(\bar y)$ (with parameters from $A_\Gamma$) 
such that: if $\bar a \in {}^m M$, 
``$\tp_{\cL} (\bar b,\bar a\cup A_\Gamma,M)$ 
is $\Gamma$-big" depend just on $\tp_\Delta(\bar b \char 94 
\bar a,A_\Gamma,M)$.

\noindent
6) A g.bigness notion $\Gamma$ is $\lambda$-simple [or
$\lambda$-co-simple] \If \, $\tp_{\cL}(\bar b,\bar a\cup A_\Gamma,M)\in 
\bold S^{\alpha}_\cL (\bar a\cup A_\Gamma, M)$
is $\Gamma$-big [is $\Gamma$-small] in $M$ implies that
for some $q \subseteq \tp_{\cL}(\bar b \char 94 \bar a,
A_\Gamma, M)$ of cardinality $<\lambda$ we have
$\bar b \char 94 \bar a' \in M' \in \bold K$ realizes $q \Rightarrow \tp_{\cL}
(\bar b',\bar a' \cup A_\Gamma,M')$ is $\Gamma$-big [or $\Gamma$-small] 
in $M'$ (for example inconsistent).

\noindent
7)  We say $\Gamma$ is a semi-$g$. bigness notion \If \, above  we omit the
local character.
\end{definition}

\begin{claim}
\label{a34}
1) Every $\ell$-bigness notion is a $g$-bigness notion
(when we restrict ourselves to complete $\cL$-types over sets including
$A_\Gamma$; we do not always bother to make the distinction).

\noindent
2) If an $\ell$-bigness notion $\Gamma$ is 
$\lambda$-strong/co-$\lambda$-strong/$\lambda$-simple/
co-$\lambda$-simple/weakly invariant/ invariant \then \, as a
$g$-bigness notion it satisfies the corresponding property.
\end{claim}

\begin{PROOF}{\ref{a34}}
Easy.
\end{PROOF}

\begin{definition}
\label{a37}
1) If $\bar\Gamma = \langle\Gamma^i = \Gamma_i:i<\alpha\rangle$ is a  
sequence of $g$-bigness
notion (with $\bar x_{\Gamma_i}$ pairwise disjoint for notational
simplicity), we consider $\bar \Gamma$ also as a $g$-bigness notion by:
$\bar x_{\bar\Gamma}=\langle\bar x_{\Gamma_i}:i<\alpha\rangle$
(formally - their concatenation), $A_\Gamma = \cup\{A_{\Gamma_i}:i <
\alpha\}$ and:
for $p \in \bold S^{\alpha(\bar\Gamma)}_\cL(A,M),A_\Gamma\subseteq A$ we have:
 $p$ is $\bar\Gamma$-big \underline{if and only if}  $p = 
p(\ldots,\bar x_{\Gamma_i},\ldots)_{i<\alpha}$, and
whenever $M \le_{\gk} N \in \bold K,
\langle\bar a_i:i<\alpha\rangle$ realizes $p$
we have $\tp_\cL(\bar a_i,A \cup \bigcup\limits_{j<i}
\bar a_j,M)\in\Gamma^i_N$ for each $i<\alpha$.

\noindent
2) Similarly when $A_{\Gamma^i} \subseteq A_\Gamma \cup \{\bar
x_{\Gamma_j}:j <i\}$ defined naturally.
\end{definition}

\begin{claim}
\label{a41}
1) If $\bar\Gamma$ is a sequence of invariant $g$-bigness notions
\then \, $\bar\Gamma$ is itself an invariant $g$-bigness notion.

\noindent
2) If $\bar\Gamma$ is a sequence of [very] $\lambda$-[co-]strong $g$-bigness
notions and $\ell g(\bar{\Gamma}) < \cf(\lambda)$ \then \, 
$\bar\Gamma$ is a [very] $\lambda$-[co-]strong $g$-bigness notion.

\noindent
3) If $\bar\Gamma$ is a sequence of [very] $\lambda$-[co-]simple
$g$-bigness notions and $\ell g(\bar \Gamma) < \cf(\lambda)$ \then \,
$\bar\Gamma$ is a [very] $\lambda$-[co-]simple $g$.bigness notion.
\end{claim}

\begin{PROOF}{\ref{a41}}
Easy.
\end{PROOF}

\begin{remark}
\label{a42}
1) What about $\bar\Gamma = \langle \Gamma_t:t \in I\rangle,I$ a linear
order which is not a well ordering?  If the $\Gamma_t$'s are very
simple, this is O.K.

\noindent
2) If $A_{\Gamma_t} \backslash A_\Gamma$ is finite for every $t \in
T$, there are no restrictions on $I$.
\end{remark}

\noindent
Now we turn to the central relation here between bigness notions
here-orthogonality.
\begin{definition}
\label{a44}
1)  Let $\Gamma_1,\Gamma_2$ be two $g$.bigness notions for
$(\gk,\cL)$, for the sequences of variables $\bar x^1,\bar x^2$
respectively (maybe infinite). We say that $\Gamma_1,\Gamma_2$ 
are orthogonal (or say $\Gamma_1$ is orthogonal to $\Gamma_2$, or say
$\Gamma_1 \perp \Gamma_2$) \If \, for any model $M \in \bold K$,
$A\subseteq M$, and sequences $\bar a^1,\bar a^2\in M$
of length $\ell g(\bar x^1),\ell g(\bar x^2)$ respectively 
such that $\tp_{\cL}(\bar a^l,A,M)$ is $\Gamma_l$-big for 
$l=1,2,$ \underline{there are} an $\le_{\gk}$-extension $N$ of $M$, and 
sequences $\bar b^1,\bar b^2\in N$ of length
$\ell g(\bar x^1),\ell g(\bar x^2)$ respectively such that 
for $l=1,2$ the sequence $\bar b^l$ realizes $\tp_{\cL}(\bar a^l,A,N)$ and
$\tp_{\cL}(\bar b^l,A\cup\bar b^{3-l}, N)$ is $\Gamma_l$-big. Similarly
``for $T$".

\noindent
2) In part (1) we say $\Gamma_1,\Gamma_2$ are nicely orthogonal
or we say $\Gamma_1$ is nicely orthogonal to $\Gamma_2$, or
we write $\Gamma_1 \perp_n \Gamma_2$, \If \,: adding to the
assumption $A_{\Gamma_1} \cup A_{\Gamma_2} \subseteq A = \acl_M A$ we can
add to the conclusion $\acl_M(A \cup \bar b^1) \cap \acl_M
(A\cup\bar b^2)=A$ ($\acl$ stands for algebraic closure, i.e.
$\acl_M(A) = \{b \in M$: for some $\bar a \subseteq A\subseteq M$ and $\varphi
(y, \bar x)$ we have  $M \models \varphi [b,\bar a]$
and $M \models(\exists^{<n} y) \varphi (y, \bar a)$
for some finite $n\}$).
\end{definition}

\begin{remark}
\label{a47}
1) If $\Gamma_l$ has parameters in $M_{\Gamma_l}$ and some
$M^* \in \bold K$ such that $M_{\Gamma_1} \le_{\gk} M^*$ and
$M_{\Gamma_2} \le_{\gk} M^*$ is given, in \ref{a44} 
by $M(\in \bold K)$ we mean for $\le_{\gk}$-extension of it, otherwise 
we look at any $\le_{\gk}$-extensions of
$M_{\Gamma_1},M_{\Gamma_2}$.

\noindent
2) Under context \ref{a2}, orthogonal is equivalent to nicely
orthogonal and every $\Gamma$ is nice, see below Definitions \ref{a50} and
Claim \ref{a53}(4).

\noindent
3) In fact also in \ref{a44}(3) we are demanding $A_{\Gamma_1} \cup
A_{\Gamma_2} \subseteq A$. However in \ref{a31}(1) we can 
weaken $A_\Gamma \subseteq A \subseteq M$ to 
$A_\Gamma \cup A \subseteq M$, and in \ref{a31}(1)(b) replace $p(\bar{x})
\restriction B$, etc. with minor changes.
\end{remark}

\begin{definition}
\label{a50}
$\Gamma$ is {\bf nice} when: if $p \in \bold S_\cL(A,M)$ is in $\Gamma_M$, and
$\alpha$ an ordinal, $A=\acl_MA$ \then \,  in some $N,M \le_{\gk}
N \in \bold K$ we can find $\bar a^i\subseteq N$ for $i<\alpha$ such that:
\mn
\begin{enumerate}
\item[$(a)$]   $\bar a^i$ realizes $p$ (in $N$)
\sn
\item[$(b)$]  $\tp(\bar a^i,A \cup \bigcup\limits_{j<i} \bar a^j)\in
\Gamma_N$ 
\sn
\item[$(c)$] $\acl_N(A\cup \bar a^i) \setminus \acl (A)$ are 
pairwise disjoint (for $i<\alpha$).
\end{enumerate}
\end{definition}

\noindent
We now give some basic properties of those notions.
\begin{claim}
\label{a53}
1)  If $\bar\Gamma^\ell=\langle\Gamma_i^\ell:i<\alpha_\ell\rangle$
$(\ell=1,2)$ are two sequences of $g$-bigness notions and $(\forall
i<\alpha_1)(\forall j<\alpha_2)[\Gamma_i^1 \perp\Gamma_j^2]$ 
\then \, $\bar\Gamma^1 \perp \bar\Gamma^2$ (on such $g$-bigness
notions, see definition \ref{a37}).

\noindent
2)  If $\Gamma_1,\Gamma_2$ are nicely orthogonal $g$-bigness notion
\then \,  they are orthogonal.

\noindent
3) If $\Gamma_1,\Gamma_2$ are orthogonal $g$-bigness
notion and each $\Gamma_\ell$ is invariant and at least one is nice
\then \,  $\Gamma_1,\Gamma_2$ are nicely orthogonal.

\noindent
4) If $\Gamma$ is an invariant $g$-bigness notion \then \, 
$\Gamma$ is nice.

\noindent
5) If during a proof of the orthogonality of $\Gamma_1,\Gamma_2$
we are given $\Gamma_\ell$-big $p_{\ell} \in
\bold S^{\alpha(\Gamma_\alpha)}_{\cL}(A, M^*)$ for 
$\ell=1,2$ we can replace $M^*$ by any $N^*$ such that 
$M^* \le_{\gk} N^*$ and $A$, $p_1$, $p_2$ by $A'$,
$A \subseteq A' \subseteq N^*$ and $p'_1$, $p'_2$ respectively such that
${p_{\ell}}' \in \bold S^{\alpha(\Gamma_\alpha)}_{\cL}(A,M')$ 
extend $p_{\ell}$ and is $\Gamma_{\ell}$-big.

\noindent
6) If $\bar\Gamma=\langle\Gamma_i:i<\alpha\rangle$, each
$\Gamma_i$ is an invariant $g$-bigness notion and nice \then \, 
so is $\bar\Gamma$.
\end{claim}

\begin{PROOF}{\ref{a53}}
E.g.

\noindent
3) Say $\Gamma_1$ is nice. Let $p^\ell$ for $\ell=1,2$ be a complete
$\Gamma_\ell$-big  type over $A$ in $M$, $A_{\Gamma_1}\cup
A_{\Gamma_2} \subseteq A \subseteq M$.
Find $N  \in \bold K,M \le_{\gk} N$ and $\langle\bar a^1_{i}: 
i<\lambda^+\rangle$ (where $\lambda$ is the supimum on the number 
of $\cL$-formulas over a set of cardinality $\le |A|+\aleph_0$), such that:
$\bar a^1_{i}\subseteq N$ realizes $p^1$ and $\langle  \acl_N(A\cup\bar
a^1_{i}) \setminus \acl_N (A):i<\lambda^+ \rangle$
are pairwise disjoint and $\tp(\bar a^1_i, A\cup \bigcup\limits_{j<i}
\bar a^1_j, N)$ is $\Gamma_1$-big (possible as $\Gamma_1$ is nice).
Choose by induction on $i \le \lambda^+$ a type $p^2_i \in 
\bold S^{\alpha(\Gamma_2)}_\cL(A \cup 
\bigcup\limits_{j<i}\bar a^1_j, N)$ such that:
\mn
\begin{enumerate}
\item[$(\alpha)$]   $p^2_i$ is $\Gamma_2$-big
\sn
\item[$(\beta)$]   if $\bar a^2$ realizes $p^2_i$ in $N'$, $N \le_{\gk} N'$
then $\tp_{\cL}(\bar a^1_i, (A \cup \bigcup\limits_{j<i}\bar a^1_j)\cup
\bar a^2, N')$ is $\Gamma_1$-big
\sn
\item[$(\gamma)$]   $p^2_i$ is increasing continuously
\sn
\item[$(\delta)$]  $p^2_0=p^2$.
\end{enumerate}
\mn
[Why possible?  For $i=0$, trivial: for $i=j+1$ we can take care of clauses
$(\alpha)+(\beta)$ as $\Gamma_1 \perp \Gamma_2$; for $i$ limit use Definition
\ref{a31}(1)(b)).

\Wilog \, some $\bar b\subseteq N$ realizes $p^2_{\lambda^+}$. Choose $i$
such that $\acl_N(A\cup\bar a^{1}_i) \setminus \acl_M (A)$ is disjoint to
$\acl_N(A\cup\bar b)$. Now $N,\bar{b},\bar{a}_i^1$ are as required.

\noindent
4) Clearly to prove that $\Gamma$ is nice, it suffice to prove:
\mn
\begin{enumerate}
\item[$(*)$]  for $\Gamma$-big, $p \in
  \bold S^{\alpha(\Gamma)}_\cL(A,M)$, 
and $A \subseteq B \subseteq M$ we can find $\bar a$, $N$ such that $M
\le_{\gk} N$, $\bar a\subseteq N$, $\tp_\cL(\bar a, B, M)$ is a $\Gamma$-big 
type extending $p$ and $\acl_N(A+\bar a)\cap \acl_N(B)=\acl (A)$.
\end{enumerate}
\mn
(just use it repeatedly).

To prove $(*)$, let $\lambda$ be as in the proof of part (3) above,
\wilog \, $A = \acl_M(A)$; by the compactness (and basic properties 
of algebraicity) we can find $M'$, $M \le_{\gk} M'$ and for
$i<\lambda^+$, elementary mapping $f_i$ (for $M'$) such that $\Dom(f_i)=B$,
$f_i\restriction A=\id$

\[
[i<j<\lambda^+\Rightarrow \acl_M (f_i(B))\cap \acl_M(f_j(B))=\acl (A)]
\]

\mn
[Why?  We should consider

\[
\{x_b \ne c: b \subseteq B \smallsetminus A, c \in M\} \cup
 \{\varphi(x_{b_1}\ldots,x_{b_n},\bar{a}):\bar{a} \subseteq A,M
\vDash \varphi [b_1,\ldots,b_n, \bar{a}]
\]

\mn
it is finitely satisfiable in $M$ by the definition of algebraic and the finite
$\triangle$- system lemma).
Then we can find $N$, $\bar a$ such that $M' \le N$, $\tp_L( \bar a,
\bigcup\limits_{i<\lambda^+} f_i(B), N)$ extend $p$ and is $\Gamma$-big.
Clearly for some $i<\lambda^+$, $\acl_N(f_i(B))\cap \acl_N(A+\bar a)=\acl
(A)$, and by invariance we are done.]
\end{PROOF}

\noindent
Now we consider a quite general scheme for defining bigness notion; as an
example see \ref{c2}.
\begin{definition}
\label{a56}
1) Let $\gs$ be a first order theory. A co-pre $\gs$-bigness
notion scheme\footnote{In principle we should
denote schemes by a different letter, so in the definition we 
use $\bold \Gamma$ but usually we do not} $\Gamma$ is a sentence
(in possibly infinitary logic) called $\psi_{\bold \Gamma}$ in the
vocabulary $\tau(\gs)\cup\{P\},P$ has arity
$\ell g(\bar x_{\bf \Gamma})$ but is not in $\tau({\gs})$, 
($\bar x_{\bf \Gamma}$ finite for simplicity-otherwise we
should have $P_w$ for every finite $w \subseteq \ell g(\bar x_{\bf \Gamma})$).
We may write $\psi_{{\bf \Gamma}}(P)$, (treating $P$ as a variable;
we shall use $P^*$, if $P$ is already occupied).

\noindent
2)  We call $\bar\varphi$ an interpretation with parameters of $\gs$
in a model $M^*\in \bold K$ \If \,:
$\bar\varphi=\langle\varphi_R(\bar y_R,\bar a_R) :R\in\tau(\gs)\rangle$
where $\varphi_R \in \cL(\tau_{\bold K})$ and $\tau(\gs)$ is the
vocabulary of $\gs$ including equality\footnote{We may add:
  $\varphi_=(x,y,\bar a_=) \rightarrow x=x$.}
 (for each sort of $\tau(\gs)$),
treating also function symbols as predicates, so $R$ is interpreted as
$\{\bar b:M^*\models\varphi_R(\bar b,\bar a_R)$,
$\ell g(\bar b) = \ell g(\bar y_R)$ $(=$ arity of $R)\}$.  The
interpreted model has universe $\{b:\varphi(M \models
\varphi_=[b,b,\bar a_=]\}$, if $\gs$ is multi-sort we have equality
for each sort. Of course, we assume
that (it holds in the cases we are considering) if
$M^* \le_{\gk} M$ and  $\bigcup\limits_{R\in\tau}\bar a_R\subseteq N
\le_{\gk} M$ then $N$ inherits the interpretation.
The interpreted model is called $M[\bar\varphi]$ or $M^{[\bar \varphi]}$ 
and we demand that it is 
a model\protect\footnote{We use equivalence classes as elements,
equality is interpreted as equivalence relation and we will not take the
trouble of dividing by it;
alternatively we can have $\gs$ not first order} of $\gs$; we demand further 
(for $\bar \varphi$ to be an interpretation of ${\gs}$ in $M$) that if 
$M^{[\bar \varphi]}$ is a model of ${\gs}$ and $M \le_{\gk} N$ then 
$N^{[\bar \varphi]}$ is a model of ${\gs}$ and $M^{[\bar \varphi]} \le_{\gk}
N^{[\bar \varphi]}$, in the context \ref{a27}(A) we get $M^{[\bar \varphi]}
\prec N^{[\bar \varphi]}$.

\noindent
3) For a co-pre $\gs$-bigness notion scheme ${\bf \Gamma} = 
\psi_{\bf \Gamma}$ and interpretation $\bar\varphi$
of $\gs$ in $M^*\in \bold K$, we define ${\bf \Gamma}[\bar\varphi]
= {\bf \Gamma}[\bar{\varphi},M^*]$, the
$\bar\varphi$-derived local bigness notion, as follows: given
$M \in \bold K$ such that $M^* \le_{\gk} M$, $\vartheta(\bar x,\bar b)$ 
is ${\bf \Gamma}[\bar\varphi]$-small in $M$ if for any quite saturated 
(see below) $N^*,M\prec N^*$ letting $P = \{\bar a:
N^*\models\vartheta[\bar a,\bar b]$ and $\bar a\subseteq
N^*[\bar\varphi]\}$ (in the relevant sorts, of course) we have
$(N^*[\bar\varphi],P)\models\psi_{\bf \Gamma}$.

\noindent
3A)  The ``quite saturated" means:
\mn
\begin{enumerate}
\item[$(a)$]  if $(\gk,\cL)= (\Mod(T),\bbL),T$ first order, means
$\lambda^+$-saturated where $\lambda=|T|+ |\tau(\gs)|+\aleph_0$
\sn
\item[$(b)$] if ${\bold K}$ is an a.e.c. with amalgamation (the last follows by
compactness if ${\cL}$ is non-trivial), we use
$\lambda^+$-model-homogeneous universal where $\lambda = {\rm
  L.S.}(\bold K) + |\tau(\gs)|$ This is needed for invariance to  hold.
\end{enumerate}
\mn
3B) For a model $N$ of $\gs$ and the identity interpretations we
define $\Gamma_N \subseteq {\cP}(N)$ as above.

\noindent
4)  We omit the ``pre" \If \, every ${\bf \Gamma}[\bar\varphi,M^*]$ is a
$\ell$-bigness notion (usually but not always for our fixed $\bold K$).
[If this holds for some $\gs$, we write $*$.]

$\Gamma$ is derived from ${\bf \Gamma}$ if it is
of the form ${\bf \Gamma} [ \bar{\varphi}, M^*]$ for some
$\bar{\varphi},M^*$.

\noindent
5) We say a property holds for ${\bf \Gamma}$ \If \, it holds for every
derived $\ell$. bigness notion (so we demand that $\bar\varphi$ is 
an interpretation of ${\gs}$ in $M$; see part (3)) 
in which not every formula is small.
We may put the ``co-'' before big.

\noindent
6)  We omit the ``co" if in part (3) we replace
$\Gamma[\bar{\varphi}]$-small by $\Gamma[\bar{\varphi}]$-big (so there 
is no real need for both notions).

\noindent
7)  We say ``$\Gamma$ is a co-${\gs}$-bigness notion, it is of the form 
${\bf \Gamma}[\bar \varphi,M]$, for ${\bf \Gamma}$ a 
co-${\gs}$-bigness notion scheme,
$\bar \varphi$ and interpretation of ${\gs}$ in $M$.

\noindent
8) We can define global ${\bf \Gamma}$ parallely.
\end{definition}

\begin{observation}
\label{a59}
Assume in Definition \ref{a56} that $\Gamma^*={\bf \Gamma}[\bar\varphi]$
for $\bold K$, and ${\bf \Gamma}$ is a co-pre $\gs$-bigness notion scheme.

\noindent
1)  It has the parameters from the interpretation $\bar \varphi$ 
i.e. $A_\Gamma$ is the set of parameters appearing in $\bar \varphi$.

\noindent
2)  ${\bf \Gamma}[\bar\varphi]$ is invariant.
\end{observation}
\begin{definition}
\label{a60}
1)  A local bigness notion $\Gamma$ (for $(\gk,\cL)$) is $\lambda$-presentable
\If \, $\bar x=\bar x_\Gamma$, and for each $\varphi(\bar
x_\Gamma,\bar y)
\in \cL$, for some set of $\cL$-formulas $\Delta_\Gamma$ of the form
$\vartheta(\ldots,\bar x^i,\ldots)_{i \in I_{\Gamma}} \in \cL(\tau(\gk))$
with $\lambda>|I_\Gamma|,|\triangle|$ (we may have parameters) where 
for $i \in I_\Gamma,\ell g(\bar x^i) = \ell g(\bar x_\Gamma)$ and
for $M \in \bold K,\bar b\subseteq M$ we have $|\Delta_\Gamma|<\lambda$ and:
$\varphi(\bar x,\bar b)$ is $\Gamma$-big in $M$ \underline{if and only if}
 the set $\{\varphi(\bar x^i,\bar b):i\in I\}\cup\Delta_\Gamma$ is finitely
satisfiable in $M$.

\noindent
2) If we omit $\lambda$ we mean $\lambda=\aleph_0$. \Wilog \,
$|I_\Gamma| \le |\Delta_\Gamma|$.

\noindent
3)  We define similarly a $\lambda$-co-presentable bigness notion
scheme, i.e. we replace above ``$\Gamma$-big" by ``$\Gamma$-small".
\end{definition}

\begin{claim}
\label{a62}
1)  If $\Gamma$ is $\gs$-[co-]bigness notion \then \,  it is
invariant.

\noindent
2)  If $\Gamma$ is a $\lambda$-presentable local bigness scheme
\then \, $\Gamma$ is invariant very $\lambda$-strong uniformly $\lambda$-simple
(hence co-simple) local bigness notion.

\noindent
3)  If $\Gamma$ is $\lambda$-co-presentable local bigness scheme
\then \, $\Gamma$ is invariant presentable, very $\lambda$-strong,
uniformly $\lambda$-co-simple (hence simple) local bigness notion.

\noindent
4) For ``very $\lambda$-strong" we can replace $\lambda$ by
$\aleph_0 + |\tau(\Delta_\Gamma)|$.

\noindent
5)  In parts 2), 3), if $\lambda =\aleph_0$ we get ``very simple" =
``very co-simple".
\end{claim}

\begin{PROOF}{\ref{a62}}
Easy (for the co-simple/simple in parts (2),(3), use \ref{a27}(7).
\end{PROOF}

\begin{claim}
\label{a65}
1) If $\langle \bar x^i : i\in I_\Gamma\rangle$ and a set $\Delta$ of
$< \lambda$ formulas in $\cL(\tau(\gk))$ and the variable $\bar x^i$
($i\in I_\Gamma$) are given, $\ell g(\bar x^i) = \alpha$ for $i \in
I_\Gamma$, and we try to define a $\lambda$-presentable 
$\Gamma$, i.e. ``$\Gamma$-big formula
in $M \in \bold K$", i.e. as in \ref{a60}(1) (so $\alpha(\Gamma) = \alpha$),
\then \, for some $\gs$ (a set of $< \lambda$ formulas of
$\cL(\tau(\gk))$ possibly with parameter) $\Gamma$ is a derived case of the
scheme pre $\gs$-bigness notion scheme.

\noindent
2) Similarly for co-representable and co-pre $\gs$-bigness notion scheme.
\end{claim}

\begin{convention}
\label{a68}
If we define $\Gamma^x_{\bar a}$ for every appropriate sequence
$\bar a$ of parameters (from some $M \in \bold K$) \then \,  we call
$\Gamma^x$ a scheme, and $\Gamma^x_{\bar a}$ an instance of this scheme.
Abusing notation we may call ${\bf \Gamma}[\bar{\varphi}]$ from
\ref{a56} a case of ${\bf \Gamma}$; in definition
\ref{a56} we may write $\Gamma$ instead of ${\bf \Gamma}$.

See more in \ref{e8} (and \ref{e27}-\ref{e41}).
\end{convention}
\newpage

\section {General examples of bigness notions}

We deal here with example of bigness notions which are general, i.e. no
formulas play special roles. Consider a model $M$, very saturated and
$A\subseteq M$, and type $p(\bar x)$ over $A$ in $M$ and we describe
the examples defined below.  We say $p(\bar x)$ is 
$\Gamma^{\tr}$-big always, $\Gamma^{\na}$-big if some $\bar a$
realizing $p$ is disjoint to $\acl(M)$. This will be helpful in
guaranteeing no undesirable algebraicity.

We say $p= p(\langle \bar x_i: i< \alpha\rangle)$ is
$\Gamma^{\rm ids}_\alpha$-big if some/any sequence 
$\langle \bar a_i: i< \alpha\rangle$
realizing $p$ is indiscernible over $A$; this is helpful in omitting
types of a sort: if we know that $\bold f$ 
(say an outside automorphism we would
like to ``kill'') satisfies $f(a) = b \notin \acl(A+a)$ and $\langle b_n:
n< \omega\rangle$ is indiscernible over $A+a$, if we guarantee it is
still indiscernible over $B+a$, $A\subseteq B$, this helps to ``kill"
$\bold f$.

Then we consider $\Gamma^{\av}_{\dot D,\bar a}$ where $p(x)$ is
$\Gamma^{\av}_{\dot D,\bar a}$-big iff every formula $\varphi \in p$
(or finite conjunction of such formulas) is satisfied by a 
$\dot D$-positive set of $a_i$. This is used, e.g. when
adding a very small non-standard natural number.

For each we are interested in its simplicity etc. and in
orthogonality. Now for theories $T$ 
 with enough set theory coded in then we can
more easily define bigness notion, so we may expand $T$ to such $T^+$,
define there $\Gamma^+$ and see what it induce on $T$, this is
promising, but as not presently used, we say little (see \ref{a41}, 
\ref{a47}).  We consider also
weakening the local property of $g$-bigness notion.

\begin{context}
\label{b2}
$T$ is first order complete, $\tau=\tau(T)$, $L =  
\bbL(\tau)$, and all $\Gamma$ here are (by \ref{a53}(4)) nice; so $\gk
= (\Mod_T,\prec),\le_{\gk} = \prec$ and $\gC$ a monster for $T$.
So $\tp$ will mean $\tp_{\bbL}$ and $\le_{\gk}$ is $\prec$ and $M,N$ 
are models of $T$.
\end{context}

\begin{definition}
\label{b5}
1) $\Gamma^0 = \Gamma^{\tr}$, the trivial bigness notion is defined by:
$\varphi(x,\bar a)$ is $\Gamma^{\tr}$-big in $M$ \underline{if and
  only if} ($\bar a\in M,\varphi\in \bbL_\tau$ and)
$M\models\exists x\varphi(x,\bar a)$.

\noindent
2) $\Gamma^1 = \Gamma^{\na}$, the  non-algebraicity bigness notion is
defined by:
$\varphi(x,\bar a)$ is $\Gamma^{\na}$-big in $M$ \underline{if and
  only if} ($\bar a\in M$, $\varphi\in \cL$ and): $M \models
\exists^{\ge n}x \varphi(x,\bar a)$ for every natural number $n$.

\noindent
3) $\Gamma^1_\alpha = \Gamma^{\na}_\alpha$, the
$\alpha$-non-algebraicly bigness notion, is defined by (where 
$\bar x = \bar x_{[\alpha]}\langle x_i:i< \alpha\rangle):
\varphi(\bar x, \bar a)$ is $\Gamma^1_\alpha$-big in
$M$, a model of $T$ \underline{if and only if} 
$\{\varphi(\langle x^k_i: i< \alpha\rangle, \bar a): k< \omega\}\cup
\{x^k_i \ne x^m_i:i<\alpha, k< m< \omega\}$ is finitely satisfiable in $M$.
\end{definition}

\begin{claim}
\label{b8}
1)  $\Gamma^{\tr}$ is a nice, presentable (invariant) $\ell$-bigness notion,
orthogonal to every (invariant) $g$-bigness notions
$\Gamma$ (trivially sometimes) and $A_{\Gamma^{\tr}} = \emptyset$.

\noindent
2) $\Gamma^{\na}$ is $\aleph_1$-presentable invariant $\ell$.bigness notion
orthogonal to every invariant $\Gamma$ and
$A_{\Gamma^{na}_\alpha}=\emptyset$.

\noindent
3) $\Gamma^{\rm na}_\alpha$ is $(\aleph_1 + |\alpha|^+)$-presentable
invariant $\ell$-bigness notion orthogonal to every invariant $\Gamma$ and
$A_{\Gamma^\tr_\alpha}=\emptyset$. 
\end{claim}

\begin{PROOF}{\ref{b8}}
By the proof of \ref{a53}(4).
\end{PROOF}

\begin{definition}
\label{b11}
$\Gamma^2_\alpha=\Gamma_\alpha^{\rm ids}$ (for $\alpha\geq\omega$), the
indiscernibility bigness notion, in the variables 
$\bar x = \langle x_i:i<\alpha\rangle$ (or each $x_i$ replaced 
by a sequence of length $n$ (possibly infinite) - does not matter) 
is defined by: $\varphi(\bar x,\bar a)$ is $\Gamma_\alpha^{\ids}$-big 
in a model $M$ of $T$ \underline{if and only if}
$\varphi(\bar x,\bar a)\cup\{\psi(x_{i_0},\ldots,x_{i_{n-1}},\bar a) 
\equiv \psi(x_{j_0},\ldots,x_{j_{n-1}},\bar a):\psi \in \bbL(\tau)$ 
and $i_0 < i_1 < \cdots < i_n < \alpha$ and 
$j_0 < j_1 < \cdots < j_n < \alpha \}$ is finitely satisfiable in $M$.
\end{definition}

\begin{claim}
\label{b14}
1) $\Gamma_\alpha^{\ids}$ is a $(|\tau|+\aleph_0)^+$-presentable
$\ell$-bigness notion, $A_{\Gamma^{\ids}_\alpha}=\emptyset$.

\noindent
2) $\Gamma_\alpha^{\ids}$ is orthogonal to any $g$-bigness notion.
\end{claim}

\begin{PROOF}{\ref{b14}}
For part (1), the proof of the (relatively) non-trivial part is contained
in the proof of part (2) which is rephrased and proved in \ref{b17} below.
\end{PROOF}

\begin{lemma}
\label{b17}
Suppose $M$ is $\lambda$-saturated (or just $\lambda$-compact),
$p$ is a $\Gamma$-big $\alpha(\Gamma)$-type over
$M$, $\Gamma$ any invariant g.bigness notion for $T$,
$|p|<\lambda$, and the set of parameters
appearing in $p$ is $\subseteq A$ and $A_\Gamma \subseteq A$. 
Let $\dot{\bold I} = \{\bar a_i:i<\alpha\}\subseteq M$ be an infinite
indiscernible sequence over $A$, \then \,  we can find a
$\Gamma$-big type $q\in \bold S^{\alpha(\Gamma)}(A \cup \bold I)$, such
that $p\subseteq q$, $q$ is $\Gamma$-big, and: if $\bar a$ realizes
$q$, then $\dot{\bold I}$ is an indiscernible sequence over $A\cup\bar a$.
\end{lemma}

\begin{remark}
\label{b18}
We can do this to $\dot{\bold I}$ whose index-set is any 
infinite (linearly) ordered set.
\end{remark}

\begin{PROOF}{\ref{b17}}
For notational simplicity assume $\alpha=\alpha(\Gamma)<\omega$ (see below).
We can replace $M$ by any $\mu^+$-saturated elementary extension
for any $\mu$ and similarly $\alpha$ can be increased.  So \wilog \, 
$\mu=\beth_{(2^\chi)^+}$, $\chi=\lambda+|T| +|A| + |\alpha(\Gamma)|$
and $\alpha = \mu$.  We can extend $p$ to
some $\Gamma$-big $p_1\in \bold S^\alpha(A \cup \dot{\bold I})$
and assume $\bar a\in M$ realizes $p_1$. Expand $M$ to $N$ by making
all elements of $A \cup \bar a$ into individual constants, and 
making the set $R^N = \dot{\bold I}$ and the order
$<^N = \{\langle \bar a_i,\bar a_j\rangle:i<j<\mu\}$ into relations of $M$.
The fact that $\bar a$ realizes over $A \cup \dot{\bold I} = 
A \cup R^N$ a $\Gamma$-big complete $\lg (\bar a)$-type,
can be expressed by omitting some types (remember the ``local
character", i.e. Definition \ref{a31}(b)).

By Morley theorem on the Hanf number of omitting types,
(see, e.g. \cite[Ch.VII,\S5]{Sh:a} = \cite[Ch.VII,\S5]{Sh:c}), there 
is a model $N'$, elementarily equivalent to $N$ and omitting all
the types over $\emptyset$ that $N$ omits, such that in $R^{N'}$ there
is an infinite indiscernible sequence $\dot{\bold J}$ (even in the vocabulary
 of $N'$).  As $\Gamma$ has local character (see \ref{a31} clause (b))
necessarily $\bar a$ realizes a $\Gamma$-big complete 
$\ell g({\bar a})$-type over $A \cup \dot{\bold J}$ in $N' \rest
\tau_{\gk}$.  Now we can compute $q$ from $\tp(\bar a,A \cup
\dot{\bold J})$ in the $\tau_{\bold K}$-reduct of $N'$. 
If $\alpha(\Gamma) \ge \omega$ then in $N$ we have also the
partial function defined by $F_\zeta^N,F_\epsilon^N(a_{i,0}) =
a_{i,\zeta}$ for $\zeta < \alpha(\Gamma)$ and $<^N= 
\{(a_{i,\zeta},a_{j,0}):i<j<\mu\}$.
\end{PROOF}

\begin{remark}
\label{b20}
We can generalize this to other cases where we have  a generalization of
Ramsey theorem (for \ref{b17} if $\Gamma$ is very simple) or Erd\"os-Rado
theorem with enough colours (for \ref{b17} in the general case)
\cite[\S2]{Sh:E59}.
\end{remark}

\begin{definition}
\label{b23}
1) We define $\Gamma=\Gamma_{\dot D,\bar a}^{\av}$
(the averaging $\ell$-bigness notion) where $\bar a=\langle\bar
a_\beta:\beta<\alpha\rangle$ is a sequence
of sequences from some $M_0=M_\Gamma,\dot D$ a filter on $\alpha$, as follows
($\ell g(\bar x) = \ell g(\bar a_\beta)$ is constant): 
$\varphi(\bar x,\bar a)$ is $\Gamma$-big in the model $M$
\underline{if and only if} ($\bar a\subseteq M,M_\Gamma \le_{\gk} 
M \in \bold K$ and)
$\{\beta<\alpha:M\models\varphi[\bar a_\beta,\bar a]\}\not=\emptyset 
\mod \dot D$.
(So we call $\Gamma_{\dot D,\bar a}^{\av}$ an instance of the scheme
$\Gamma_{\dot D^{\av}}$.)

\noindent
2)  We say $\dot D,\bar a$ is non-trivial if for some $\zeta < 
\ell g(\bar x)$ and for every finite $u \subseteq \ell g(\bar a)$, 
the set $\{\beta< \alpha:a_{\beta,\zeta} \notin \{a_{\gamma,\zeta}:
\gamma\in u\}\}$ belongs to $\dot D$. If $\bar a_i= \langle a_i
\rangle$ we may write $a_i$ instead $\bar a_i\}$.
\end{definition}

\begin{claim}
\label{b27}
If $\bar a_i \in M_{\Gamma} \in \bold K$ for $i<\alpha$ and $\dot D$ a
$\lambda$-complete filter on $\alpha$, and $\Gamma =
\Gamma_{\dot D,\bar a}^{\av}$ \then \, 
\mn
\begin{enumerate}
\item[$(1)$]   $A_\Gamma = \bigcup\limits_{i<\alpha} \bar a_i $ (but also
$\bigcup\limits_{i\in Y} \bar a_i$ for $Y \in \dot D$ is O.K.)
\sn
\item[$(2)$]  $\Gamma$ is an invariant, very $\aleph_0$-strong,
very $|\alpha|^+$-simple $\ell$-bigness notion
\sn
\item[$(3)$]  If $\dot D$ is an ultrafilter, \then \, $\Gamma$ is
orthogonal to any uniformly $\lambda$-co-simple $\ell$-bigness notion
\sn
\item[$(4)$]  $\Gamma^{\av}_{\dot D,\bar a}$ is non-trivial \If \,
  $\dot D$ is non-trivial.
\end{enumerate}
\end{claim}

\begin{PROOF}{\ref{b27}}
For part (2) the proof of the (relatively) non-trivial part is contained in
the proof of part (3) which is rephrased and proved in \ref{a31} below.
\end{PROOF}

\begin{lemma}
\label{b31}
Suppose
\mn
\begin{enumerate}
\item[$(a)$]   $p$ is a $\Gamma$-big type over $A$ in $M$, and $\bar a
  \in {}^\zeta M$, $A_\Gamma\subseteq A\subseteq M$
\sn
\item[$(b)$]  $\bar a_i \in {}^\zeta A$ for each $i<\alpha$, and $\dot D$ is a
$\lambda$-complete filter over $I$
\sn
\item[$(c)$]  for any formula $\varphi(\bar x,\bar b)$ with parameters
from $A$, if $\models \varphi[\bar a,\bar b]$ \then \, 
$\{i\in I:M\models\varphi[{\bar a}_i,\bar b]\}\in \dot D$,
(hence $\tp(\bar a,A,M)$ is $\Gamma^{\av}_{\dot D,\langle{\bar a}_i:
i<\alpha\rangle}$-big; if $\dot D$ is an ultrafilter this is an 
equivalent formulation)
\sn
\item[$(d)$]  $\Gamma$ is a uniformly $\lambda$-co-simple notion of
$\ell$-bigness.
\end{enumerate}
\mn
We can conclude that we can extend $p$ to a $\Gamma$-big type $q\in
\bold S^m(A\cup\bar a)$ such that for any formula 
$\varphi(\bar x,\bar y,\bar b)$: \If \, $\bar b\in A$,
and $\{i<\alpha:\varphi(\bar x,\bar a_i,\bar b) \in q\} \in \dot D$, 
\then \, $\varphi(\bar x,\bar a,\bar b) \in q$ (i.e. if 
$\bar c$ realizes $q$ then $\tp(\bar a,A \cup \bar c, M)$ is 
$\Gamma^\av_{\dot D,\bar a}$-big).
\end{lemma}

\begin{PROOF}{\ref{b31}}
Let $r_{\varphi(\bar y,\bar z)}(\bar z)$ be a type of cardinality
$< \lambda$ (over $A_\Gamma$) such that: 
$\varphi(\bar y,\bar b)$ is $\Gamma$-small \underline{if and only if}
$\bar b$ realizes $r_{\varphi(\bar y,\bar z)}(\bar z)$ (in $N$ whenever 
$M \le_{\gk} N,\bar b \in N$ of course; exist as $\Gamma$ is $\lambda$-co-simple).
\Wilog \, $p \in \bold S^{\alpha(\Gamma)}(A,M)$, and now we define $q$, by

\[
q=q(\bar x,\bar y) = p(\bar x) \cup q_0 (\bar y) \cup q_1(\bar x, \bar y)\cup
q_2(\bar x,\bar y)
\]

where

\[
q_0(\bar y) = \{\varphi (\bar{y}, \bar{b}): \bar{b} \in A,\text{ and }
 M \vDash \varphi[\bar{a},\bar{b}]\} 
\]

\[
q_1(\bar x,\bar y) = \left\{\varphi(\bar x,\bar y,\bar b):\bar b\in A,\
\{i<\alpha:\varphi(\bar x,\bar a_i,\bar b)\in p(\bar{x})\} \in \dot D \right\}
\]

and

\[
q_2(\bar x,\bar y) = \{\neg \varphi(\bar x,\bar y,\bar b):\bar b \in A
\text{ and } \varphi(\bar x,\bar a,\bar b) \text{ is }
\Gamma\text{-small}\}.
\]

\mn
By the hypothesis on $\bar a$
\mn
\begin{enumerate}
\item[$(*)_1$]   $q$ extend $q_0(\bar y)= \{\varphi(\bar y,\bar b):\bar b
\in A,M \models \varphi[\bar a,\bar b]\}$ and
\sn
\item[$(*)_2$]  $p(\bar x) \subseteq q_1(\bar x,\bar y) \subseteq q$ 
\end{enumerate}
\mn
and lastly
\mn
\begin{enumerate}
\item[$(*)_3$]  every finite $q' \subseteq q$ is realized in $M$.
\end{enumerate}
\mn
Why $(*)_3$ holds?  As $p,q_0$ are complete types over $A,\dot D$ 
is a filter and the set of $\Gamma$-small formulas form an ideal 
clearly $p({\bar x}),q_0(\bar y),q_1({\bar x},\bar y),q_2
({\bar x},\bar y)$ are closed under
conjunctions hence \wilog \, $q' = \{\varphi({\bar x},\bar b)\} \cup
\{\varphi_0 ({\bar y},b_0)\} \cup \{\varphi_1(\bar x,\bar y,\bar b_1)\} \cup
\{\varphi_2(\bar x,\bar y,\bar b_2)\}$ where $\varphi(\bar x,\bar b) 
\in p(\bar x),\varphi_0({\bar y},b_0) \in q_0(\bar y)$ and 
$\varphi_\ell(\bar x,\bar y,\bar b_\ell) \in q_\ell$ for $\ell=1,2$. 
As $\varphi_2(\bar x,\bar y,\bar b_2)\in q_2$ necessarily 
$\bar a \char 94 \bar b_2$ realizes a type $r_{\psi(\bar x;\bar y
\char 94 \bar z)}(\bar y,\bar z)$ when we let $\varphi_2 = \neg \psi$ 
which satisfies: if $\bar a' \char 94 \bar b'$ realizes the type 
$r_{\psi(\bar x,\bar y \char 94 \bar z)} (\bar y,\bar z)$ then 
$\psi(\bar x,\bar a',\bar b)$ is $\Gamma$-small. 

Now

\[
\vartheta(\bar y,\bar z) \in r_{(\bar x;\bar y \char 94 \bar z)} 
(\bar y,\bar z) \Rightarrow M \models \vartheta
[\bar a,\bar b_2] \Rightarrow I_{\vartheta(\bar y,\bar z)} =: \{i<
\alpha:\models \vartheta[\bar a_i,\bar b_2]\} \in \dot D.
\]

\mn
But by assumption (c) as $\dot D$ is $\lambda$-complete  
${\cU}^* = \bigcap \{I_{\vartheta(\bar y,\bar z)}:\vartheta(\bar
y,\bar z) \in r_{\psi(x,\bar y \char 94 \bar z)}\}$
belongs to $\dot D$.  Now $\varphi_1(\bar x,\bar y,\bar b_1) \in q_1(\bar x,
\bar y)$ hence ${\cU}' =: \{i<\alpha:\varphi_1(\bar x,\bar a_i,\bar
b_1) \in p(\bar x)\} \in \dot D$.
As $\varphi_0 (\bar y,\bar b) \in q_0(\bar y)$ clearly ${\cU} =
\{i<\alpha:\models \varphi_0[{\bar a}_i,\bar{b}]\} \ne \emptyset
\mod \dot D$.

We conclude that $({\cU}\cap {\cU}^*)\cap {\cU}' \ne \emptyset \mod
\dot D$ hence there is $i \in {\cU} \cap {\cU}' \cap {\cU}^*$.

Now we can find $N,{\bar c}$ such that $M \le N$ and ${\bar c} \in N$
such that $\tp({\bar c},A \cup \bar a_i,N)$ is $\Gamma$-big and extend
$p(\bar x)$ exists by \ref{a27}(4).  We claim that
${\bar c} \char 94 \bar a_i$ realizes $q'$ in $N$.  Why?  First $N \models
\varphi[{\bar c},b]$ as $\varphi({\bar x},b) \in p(\bar x)
\subseteq \tp(\bar c,A \cup \bar a, N)$.   Second, $N \models \varphi_0
[{\bar a}_i,\bar b]$ as $i \in {\cU}$.  Third, $N \models
\varphi_1 [ \bar c,{\bar a}_i,b]$ as $i \in {\cU}'$ which implies
$\varphi_1 (\bar x,{\bar a}_i,{\bar b}_1) \in p(\bar x)
\subseteq \tp({\bar c},A,N) \subseteq \tp({\bar c},A \cup \bar a_i,N)$
and ${\bar c}$ realizes this type.  Fourth, $N \models \varphi_2[\bar
c,{\bar a}_i,{\bar b}_2]$ which holds as $\bar{a}_i \char 94 \bar{b}_2$
realizes the type $r_{\psi(\bar{x,y \char 94 z})} (\bar y,\bar z)$ hence
$\psi (\bar x,\bar a_i,\bar b_2)$ is $\Gamma$-small but $\tp(\bar{c},A \cup
\bar{a}_i,N)$ is $\Gamma$-big in $N$ so necessarily $\neg \psi (\bar x,a_i,b_2)
\notin \tp(\bar{c},A \cup \bar{a}_i,N)$ hence $\neg \psi (x,\bar{a}_i,
\bar{b}_2) \in \tp(\bar{c},A \cup \bar{a}_i,\bar{b}_2)$ so $N \models
\neg \psi [ \bar{c},{\bar a}_i,{\bar b}_2]$ but 
$\varphi_2 (\bar x,\bar y,\bar z)= \neg \psi(\bar x,\bar y,\bar z)$ 
so we are done.  So $\bar{c} \char 94 \bar{a}_i$ really 
realizes $q'$, so $q'$ is realized as promised in $(*)_3$.

So $q$ is indeed finitely satisfiable (in $M$) but $q \ge p_0(\bar y)
\in \bold S^\zeta(A,M)$  hence even $q (\bar x,\bar a)$ is 
finitely satisfiable so let $\bar c^*$ realizes $q(\bar x,a)$ in $N^*$
where $M \le N^*$ and let $q^* = \tp(\bar c,A \cup \bar a, N^*)$. By
our choice $q_2(\bar x,\bar a)\cup q_1(\bar x,\bar a) 
\subseteq q^*$ so clearly $q (\bar x,\bar a)$ is also 
$\Gamma$-big.  Obviously it extends $p$, and satisfies the
conclusion of the lemma. 
\end{PROOF}

\begin{claim}
\label{b34}
In \ref{b31}, if $\dot D$ is an ultrafilter \then \,  the
conclusion of \ref{b31} is valid even if we (seemingly) 
weaken the demand (c) to
\mn
\begin{enumerate}
\item[$(c)^-$]  $\tp(\bar a,A,M)$ is $\Gamma^{\av}_{\dot D,
\langle\bar a_i:i<\alpha\rangle}$-big.
\end{enumerate}
\end{claim}

\begin{PROOF}{\ref{b34}}
Should be clear.
\end{PROOF}

\begin{observation}
\label{b37}
1)  If $\Gamma_1$, $\Gamma_2$ are local bigness notions, \then \, the following
are equivalent:
\mn
\begin{enumerate}
\item[$(a)$]   $\Gamma_1 \perp \Gamma_2$
\sn
\item[$(b)$]  for any $\Gamma_1$-big, $p_1=\tp({\bar a}^1,A,M)$ 
and $\Gamma_2$-big $p_2 = \tp({\bar a}^2,A,M)$
the following set of formulas is finitely satisfiable (in $M$)
\begin{equation}
\begin{array}{clcr}
p_1 ({\bar x}^1) \, \& \, &\cup p_2({\bar x}^2) \cup 
\{\neg\varphi({\bar x}^1,{\bar x}^2,\bar b):{\bar b} \subseteq
A,\varphi({\bar x}^1,{\bar a}^2,\bar b \text{ is }
\Gamma_1\text{-small}\} \\
  &\cup \{\neg\varphi({\bar x}^1,{\bar x}^2,\bar b):b \subseteq
A,\varphi({\bar a}^1,{\bar x}^2,\bar b) \text{ is }
\Gamma_2\text{-small}\}
\end{array}
\end{equation}
\sn
\item[$(c)$]   Assume $M_{\Gamma^\ell}\le N,A_{\Gamma^1}
\cup A_{\Gamma^2}\subseteq A \subseteq N,N$ is $|A|^+$-saturated, 
$\bar b_\ell\in N,p_{\ell} = \tp(\bar b_\ell,A,N) \in \Gamma^\ell_N$. 
We can find $\bar b'_1,\bar b'_2 \in N$ such that $\bar b'_\ell$ 
realizes $p_\ell$ ($=\tp(\bar b_\ell,A,N)$) and
$\tp(\bar b'_\ell, A \cup\{\bar b'_{3-\ell}\},N)$ is
$\Gamma_\ell$-big for $\ell =1,2$
\sn
\item[$(c)'$]  there are $N \supseteq A_* \supseteq A_{\Gamma_1} \cup 
A_{\Gamma_2}$ such that $N$ is $|A_*|^+$-saturated such that if $B
\subseteq N$ is finite and $A = A_* \cup B$ then the second sentence
in cluase (c) holds.
\end{enumerate}
\mn
2)  If $\Gamma_1$, $\Gamma_2$ are co-simple local bigness notions
\then \, we can add:
\mn
\begin{enumerate}
\item[$(d)$]   if $p_1$ is $\Gamma^1$-big, $p_2$ is $\Gamma^2$-big,
$\vartheta_1({\bar x}^1)\in p_1,\vartheta_2({\bar x}^2)\in p_2$,
$\vartheta_1({\bar x}^1)$ witness $\varphi_1({\bar x}^2;{\bar x}^1)$ is
$\Gamma^2$-small, $\vartheta_2({\bar x}^2)$ witness $\varphi_2({\bar x}^1;
{\bar x}^2)$ is $\Gamma^1$-small \then \, 
\newline
$\{\vartheta_1({\bar x}^1),\vartheta_2({\bar x}^2),\neg \varphi_1({\bar x}^2;
{\bar x}^1),\neg \varphi_2({\bar x}^1;{\bar x}^2)\}$ is consistent; of course
$\vartheta_1,\vartheta_2,\varphi_1,\varphi_2$ may have parameters. We
can guarantee $\varphi_{\ell}({\bar x}^{3-\ell};{\bar a}'_\ell)$
is $\Gamma^\ell$-small for any ${\bar a}'_{\ell}$.
\end{enumerate}
\end{observation}

\begin{PROOF}{\ref{b37}}
Easy.
\end{PROOF}

\begin{definition}
\label{b41}
Let $\tau_1 \subseteq \tau_2$ be vocabularies, $T_\ell$ a complete theory
in $\bbL(\tau_\ell)$, $T_1\subseteq T_2$.

\noindent
1)  If $\Gamma^2$ is a local bigness notion for $T_2$, we define
$\Gamma^1=\Gamma^2\restriction \tau_1$ by: for $N_1$ a model of $T_1,
\Gamma^1_{N_1} = \{\varphi(\bar x_{\Gamma^2},\bar a): \bar a\subseteq
N_1$ and for some $N_2 \models T_2,N_1 \prec N_2\restriction \tau_1$
 we have $\varphi(\bar x_{\Gamma^2},\bar a) \in \Gamma^2_{N_2}\}$.

\noindent
2) If $\Gamma^2$ is a global bigness notion for $T_2$,
$\Gamma^1=\Gamma^2\restriction \tau_1$ is defined by: $p$, a 
complete type over $A$ in $N_1$ is $\Gamma_1$-big (in $N_1$) \If \,
for some $N_2 \models T_2$, $N_1\prec N_2\restriction \tau_1$ and
$p$ can be extended to a complete type over $A$ in $N_2$ which is
$\Gamma^2$-big.
\end{definition}

\begin{remark} 
\label{b44}
Note that above $\Gamma^2 \restriction \tau_1$ is not a priori 
a bigness notion.
\end{remark}

\begin{claim}
\label{b47}
Let $\tau_1,\tau_2,T_1,T_2$ be as in Definition \ref{b41}.

\noindent
1)  The two parts of Definition \ref{b41} are compatible.

\noindent
2)  In Definition \ref{b41}(1) if $\Gamma^2$ is a local bigness notion
for $T_2$ \then \, $\Gamma_1$ is a local bigness notion for
$T_1$ (and is invariant). 
Similarly (see Definition \ref{a31}(7)) for global bigness notion.

\noindent
3)  $\Gamma^{\tr},\Gamma^{\na},\Gamma^{\ids}_\alpha,\Gamma^{\av}_{\dot
  D,\bar a}$, commute with the restriction operation.

\noindent
4) In Definition \ref{b41}(1), if $\Gamma^2$ is a $\lambda$-strong/
co-$\lambda$-strong/very $\lambda$-strong/co-simple local bigness
notion for $T_2$ \then \,  $\Gamma^1$ is a $\lambda$-strong/very
$\lambda$-strong/co-simple local bigness notion for $T_1$.

\noindent
5) Assume $\Gamma'_2,\Gamma''_2$ are global (or local) bigness notions
for $T_2$ and $\Gamma'_1= \Gamma'_2 \restriction \tau_1,\Gamma''_1=
\Gamma''_2 \restriction \tau_2.$   If $\Gamma'_1,\Gamma''_2$ are 
orthogonal \then \, $\Gamma'_1,\Gamma''_1$ are orthogonal.

\noindent
6) ${\bf \Gamma}$ is a co-$\gs$-bigness notion scheme, 
$\bar\varphi$ an interpretation of $\gs$ in a model $M_2$ of 
$T_2$, $\Gamma^2= {\bf \Gamma} [\bar{\varphi},M_2]$ (see Definition
\ref{a56}) and $\Gamma^1= \Gamma^2 \restriction \tau_1$
the relevant formulas $\varphi_i$ belong to $\bbL(\tau(T_1))$,
\then \, $\Gamma^1= {\bf \Gamma} [\bar \varphi, M_2 \restriction \tau_1]$.
Similarly for $\lambda-[co-]$ representable.

\noindent
7) If in Definition \ref{b41}, $\Gamma^2$ is a $\kappa$-weakly global bigness
notion \then \, so is $\Gamma^1$ (see Definition \ref{b54} below).
\end{claim}

\begin{remark}
\label{b50}
Why in \ref{b47}(4) we have ``co-simple" and not simple?  The point is
that in Definition \ref{b41} we have:
\mn
\begin{enumerate}
\item[$\bullet$]  $\varphi(\bar x_{\Gamma_2},\bar a)$ is
  $\Gamma_1$-big in $M_1 \models T_1$ \Iff \, there is $M_2 \models
  T_2$ such that $M_1 \prec M_2 \rest \tau(T_1)$ and $\varphi(\bar
  x_{\Gamma_2},\bar a)$ is $\Gamma_2$-big in $M_2$.
\end{enumerate}
\mn
So for $\Gamma_1$-small we have to say ``for every $M_2$ ...".
\end{remark}

\begin{PROOF}{\ref{b47}}
Straightforward, e.g.

\noindent
4) For notational simplicity let $A_\Gamma = \emptyset$.
\medskip

\noindent
\underline{Case 1}:  $\Gamma_1$ is $\lambda$-strong.

Assume, $\varphi(\bar x,\bar a)$ is a $\Gamma_1$-big formula in 
the model $M_1$ of $T_1$, so for some model $M_2$ of $T_2$ we 
have $M_1 \prec M_2 \restriction \tau_1$ and
$\varphi(\bar x,\bar a)$ is
$\Gamma_2$-big in $M_2$, hence for some $\tau^*_2 \subseteq \tau_2$ of
cardinality $< \lambda$ we have: if $M^*_2$ is a model of 
$T_2,\bar{a}^*_2 \in M^*_2$ realizes
$\tp(\bar{a},\emptyset,M_2 \restriction \tau^*_2)$
then $\varphi(\bar{x},\bar{a}^*)$ is $\Gamma_2$ big in $M^*_2$. We shall
show that $\tau^*_1= \tau_1 \cap \tau^*_2$
is as required, so assume that $M^*_1$ is a model of $T_1$ and $\bar{a}^*_1
\in M^*_1$ realizes $\tp(\bar{a},\emptyset,M_1 \restriction
\tau^*_1)$. 
By Robinson lemma there is a model $M^{**}_2$ of $T_2$ such that
$M^*_1 \prec M^{**}_2 \restriction \tau_1$ and $\bar{a}^*_1$ realizes
$\tp(\bar{a},\emptyset,M^*_2 \restriction \tau^*_2)$ 
hence $\varphi(\bar x,\bar a^*_1)$ is $\Gamma_1$-big.
\medskip

\noindent
\underline{Case 2}:  co $\lambda$-strong. The same proof.
\medskip

\noindent
\underline{Case 3}:  very $\lambda$-strong.  Just easier.
\medskip

\noindent
\underline{Case 4}:  simple.  Like case 5.
\medskip

\noindent
\underline{Case 5}:  co-simple (not co-$\lambda$-simple!).

Assume $\varphi(\bar{x}, \bar{a})$ is $\Gamma_1$-small in $M_1$, a model of
$T_1$. Let $T'=T_2 \cup \{\vartheta(\bar y):\vartheta(\bar y) \in
\bbL(\tau(T_1))$ and $M_1 \models 
\vartheta [\bar a]\} \cup \{\neg \psi (\bar y):\psi \in {\Bbb L}
(\tau_2),$ and if $M_2 \models \psi[\bar{b}], M_2$ a model of $T_2$ then
$\varphi(\bar{x},\bar{b})$ is $\Gamma_2$-small\}, clearly this 
set of (first order) formulas has no model hence is
inconsistent, but the second set in the
union is closed under conjunctions and also the third (as $\neg \psi_1
(\bar{y}) \wedge \neg \psi_2 (\bar{y})$ is $\neg(\psi_1(\bar{y})
\vee \psi_2 (\bar{y}))$.  So for some $\vartheta(\bar{y})
\in \tp(\bar{a},\emptyset, M_1)$ and $\psi(\bar{y})$ we have $T_2 \cup \{
\vartheta (\bar{y}),\neg \psi(\bar{y})\}$ is inconsistent and 
$[M_2 \models \psi [\bar{b}] \Rightarrow \varphi(\bar{x,b})$ 
is $\Gamma_2$-small] for every $\bar b \in M_2, M_2$ a model
of $T_2$. So $\vartheta (\bar{y})$ is as required. 
\end{PROOF}

\begin{definition}
\label{b54}
$\Gamma$ is a $(<\kappa)$-weakly global bigness notion for $T$ \If \,:
in Definition \ref{a31}(1) we weaken clause (b) to:
\mn
\begin{enumerate}
\item[$(b)^-$]  for $\lambda<\kappa$, the odd player has a winning
strategy in the following game: the game lasts $\lambda+1$ moves,
in the $\alpha$-th move a $\Gamma$-big $p_\alpha \in
\bold S^{\alpha(\Gamma)}(A_\alpha,M_\alpha)$ such that $\alpha< \beta
\Rightarrow M_\alpha \le M_\beta \, \& \, A_\alpha \subseteq A_\beta \, \& \,
p_\alpha = p_\beta \restriction A_\alpha$,
the even/odd player choosing for $\alpha$ even/odd. The even player
wins if he has no legal move for some $\alpha \le \lambda$. Otherwise
the odd player wins.
\end{enumerate}
\mn
Let $\alpha^*$-weak mean $(<\alpha^*+1)$-weak.
\end{definition}
\newpage

\section {Specific examples of bigness notion schemes}

We deal with bigness notions for which some formulas have special
roles. One family of natural ones are variants of a subsets of a partially
ordered sets which are somewhere dense (i.e. for some $x$ for every
$y> x$ there is $z> y$ in the set). By looking at intervals of a
linear order we can get
as a special instance the case of dense linear orders;  note this
density has a different meaning, those are important for automorphism
of ordered field not treated here.
Another family of natural ones (considered even earlier) are connected to
independence: outside a small set every combination
is possible (for this we need the strong independence property); the main
example here is a member of an atomic Boolean algebra such that except
for a small set of atoms we have total freedom which ones to put
inside and which ones to put outside. The main case here is having a
pseudo finite set $a$ as a parameter, and inside 
$\gC$ (see \ref{a2}(B)) we say $x \subseteq a$ is big 
if $|x|/|a| \geq c_i$ for each $i<\delta$ where $c_i \in 
(0,1)^{\gB}_{\bbR}$ is increasing with $i$. So ``finitary" 
theorems enter like the law of large numbers.  
We are in particular interested in the case $a={\cP}(a_1)$, and in 
particular if for some $\bold q \in (0,1)^{\gB}_{\bbR}$ we
give to $b \subseteq a$ the weight $\bold q^{|b|} \times
(1-\bold q)^{|a\setminus b|}$.

Here we usually do not mention ``in $\gC$" as it is obvious.

\begin{definition}
\label{c2}
1)  Let ${\gt}^{\po}$ be the first order theory such a structure $(A,<,R)$
satisfies ${\gt}^{\po}$ \underline{if and only if} $A \ne \emptyset,<$ 
a partial order, $R$ a symmetric two-place relation satisfying
$xRy \rightarrow \neg(\exists z)(x \le z \wedge y \le z)$, to which
$({}^{\omega>} 2,\triangleleft,\ntriangleleft)$ can be embedded (where if we
omit $R$ it means $xRy := \neg(\exists z)[x \le z \wedge y \le z)$.

\noindent
2) Let ${\gt}^\poe$ be as ${\gt}^\po$ adding to the theory $x < y
\Rightarrow (\exists z)[x<z \wedge z Ry]$ and $\forall x\exists
y(x<y)$.

\noindent
2A) Let ${\gt}^{\pot}$ be as ${\gt}^\po$ adding to the theory 
$(\forall x)(\exists y_1,y_2)[x<y_1 \, \& \, x<y_2 \, \& \, y_1 R y_2]$.

\noindent
3) Let $\Gamma^{\po}$ be the following pre-${\gt}^{\po}$-bigness notion 
scheme $\psi(P)$ (see Definition \ref{a56}):
for $M$ a model of ${\gt}^{\po}$, and $P \subseteq M$: $M \models \psi[P]$
says that the following is finitely satisfiable in $(M,P)$: 

\begin{equation}
\begin{array}{clcr}
\{P(\bar x_\eta):\eta\in {}^{\omega>}2\} &\cup\{x_\eta<x_\nu:
\eta \triangleleft \nu \in {}^{\omega>}2\} \\
  &\cup\{x_{\eta \char 94 \langle 0 \rangle} R x_{\eta \char 94 
\langle 1 \rangle}:\eta \in {}^{\omega>}2\}.
\end{array}
\end{equation}
\mn
4)  Let $\Gamma^{\poe}$ be the following pre-${\gt}^{\poe}$-bigness scheme:
for $M$, a model of ${\gt}^{\poe}$ and $P \subseteq M$: $\psi[P]$ says
$(\exists x)(\forall y)(\exists z\in P)[x<y\to y< z]$ (this means
somewhere dense).

\noindent
5) $\Gamma^\pot$ is defined like $\Gamma^\poe$. 
\end{definition}

\begin{remark}
\label{c5}
A natural example of \ref{c2}(4) (more exactly, a model 
of ${\gt}^\poe$) is the set of open intervals of a dense linear 
order ordered by inverse (strict) inclusion with $(a,b) R (a',b')$ 
iff $(a,b) \cap (a',b')=\emptyset$.
\end{remark}

\begin{claim}
\label{c8}
1) $\Gamma^{\po}$ is a ${\gt}^{\po}$-bigness notion scheme and is
$\aleph_1$-presentable so by \ref{a62}(2), (3)
is invariant, very $\aleph_0$-strong co-simple and uniformly
$\aleph_1$-simple local bigness notion.

\noindent
2)  $\Gamma^{\poe}$ is a ${\gt}^{\po}$-bigness notion scheme and
is presentable (so by \ref{a62} is invariant, very simple). 
Similarly for $\Gamma^\pot$.

\noindent
3) For $M \models {\gt}^{\poe}$ we have $\Gamma_M^{\po}= \Gamma_M^{\poe}$. If
$\bar\varphi$ is an interpretation of ${\gt}^{\poe}$ in a model $M$ then
$\Gamma_M^{\po}[\bar\varphi] = \Gamma_M^{\poe}[\bar\varphi]$. Similarly
for $\Gamma^\pot$.

\noindent
4) ${\gt}^{\poe} \bot {\gt}^\tr$ and ${\gt}^\tr \bot {\gt}^\po$.
\end{claim}

\begin{definition}
\label{c11}
1)  ${\gt}^{\ind}$ is the theory saying on a model $M,M=(P,Q,R)$,
$P,Q$ are disjoint (are two sorts (or if you prefer -- two unary predicates)), 
such that $R \subseteq P \times Q,R$ has the strong 
independence property (see Definition \ref{c14} below).

\noindent
2)  We define a pre-${\gt}^{\ind}$-bigness notion scheme
$\Gamma=\Gamma^{\ind}$ as follows: $\psi_{\Gamma}(P^*)$ says 
$P^*\subseteq Q$ and $(P,Q,R,P^*)$ satisfies:
for every $n$ it is $\Gamma_R$-$n$-big which means that for some finite
$A\subseteq P,P^*$ has $n$-independence outside
$A$ ($A$ is called a $\Gamma$-$n$-witness), which means:
for every pairwise distinct $a_0,\ldots,a_{2n-1}\in P\setminus A$, for some
$c \in P^*$ we have $\forall \ell < 2n \Rightarrow 
[a_\ell Rc]^{\iif(\ell < n)}$; so $\psi_\Gamma$ is not first 
order because we have said ``some finite $A$" but $\psi \in
\bbL_{\aleph_1,\aleph_0}$.
\end{definition}

\begin{definition}
\label{c14}
1) Let $P$, $Q$ be one place predicates, $R$ a two-place predicate
and suppose the theory $T$ contains the formula $(\forall xy)[xRy \rightarrow
P(x)\, \& \, Q(y)]$.

We say $R$ (more exactly $P,Q,R$) has the strong independence property (for
or in the theory $T$) \If \,:
\mn
\begin{enumerate}
\item[$(a)$]  $P$, $Q$, $R$ as above
\sn
\item[$(b)$]  for every $n<\omega$, $M\models T$, and pairwise distinct
$a_1,\ldots,a_{2n}\in P^M$ there is $c\in Q^M$ such that: $a_iRc$ iff
$i \le n$.
\end{enumerate}
\mn
2) We say $R$ i.e. $(P,Q,R)$ has comprehension, i.e.

\[
\forall x, y\exists z\ [Q(x) \wedge Q(y) \wedge x \ne y \rightarrow
P(z) \wedge (zRx \equiv \neg zRy)].
\]
\end{definition}

\begin{example}
\label{c17}
The following are examples of theories $T$ implying $(P,Q,R)$,
i.e. $(P^{\gC(H)},Q^{\gC(T)},R^{\gC(T)})$ has the strong
independence property and comprehension.

\noindent
1)  $T=$ true arithmetic, that is the theory of $\bbN = (\omega,+,\times,0,1)$.
Let

$P(x):x$ is prime-

$Q(x):x>0$ not divisible by any square of a prime,

$xRy:x$ divides $y$, and $x$ is prime and $Q(y)$.

\noindent
2)  $T$ as above

$P(x):x=x$.

$xRy:y$ codes a sequence in which $x$ appears (using a fix coding).

\noindent
3) $T=$ the first order theory of infinite atomic Boolean algebras

$P(x):x$ an atom,

$Q(x):x=x$,

$xRy:x \le y$
\end{example}

\begin{claim}
\label{c20}
$\Gamma^{\ind}$ is a ${\gt}^{\ind}$-bigness notion scheme (hence
invariant) and is very $\aleph_1$-strong and $\aleph_1$-co-strong 
(but not uniformly).
\end{claim}

\begin{PROOF}{\ref{c20}}
Assume $\bar\varphi$ is a ${\gt}^{\ind}$-interpretation, and let
$M[\bar\varphi]=(P,Q,R)$.
Concerning the $\aleph_1$-co-simple, note that ``$\varphi(x, \bar a)$
is $\Gamma_M$-small" iff for some $m$, we have:
$\bar a$ realizes the type

\begin{equation}
\begin{array}{clcr}
q_m(\bar y) =: \{&\neg(\exists y_0) \ldots (\exists y_{n-1})
(\forall z_0) \ldots (\forall z_{m-1}) \\
  &[\bigwedge\limits_{\ell<n}\, \bigwedge\limits_{k<m} (y_\ell \ne z_k
  \, \& \,\bigwedge\limits_{\ell<k<m} z_\ell \ne z_k \, \& \, 
\bigwedge\limits_{k<m} P(z_k) \rightarrow \\
  &\bigwedge\limits_{w\subseteq m} (\exists x)(Q(x) \, \& \,  
\varphi(x,\bar y) \, \& \, \bigwedge\limits_{k<m}
\varphi(x,y_k)^{\iif(k\in w)})]: n < \omega\}.
\end{array}
\end{equation}

\mn
We should check $\Gamma_M = (\Gamma_{M[\bar\varphi]}^{\ind})_M$ 
satisfies ``$\Gamma_M$ is a proper ideal" (the other conditions 
are obvious). So we should check $(\alpha),(\beta),(\gamma)$ of 
\ref{a11}(1) (c) in order to show
that ``$\Gamma_M$ is a proper ideal".

So $(\alpha)+(\gamma)$ obviously hold. How about $(\beta)$ i.e.
$\varphi_1\lor\varphi_2$? Suppose $\varphi_1$ is not $\Gamma_R-n$-big and
$\varphi_2$ is not $\Gamma_R-m$-big. 
Let $A\subseteq M$ be finite and we shall show that $A$ cannot be a
$\Gamma_R-(n+m)$-witness for $\varphi_1 \vee \varphi_2$ 
(see Definition \ref{c11}(2)).

As $\varphi_1$ is $\Gamma_M-n$-small, $A$ cannot be a witness hence 
there are $a_0,\ldots,a_{2n-1} \in P^M \setminus A$ with no repetition 
such that

\[
\neg(\exists x)\bigl(\varphi(x,\bar a) \wedge 
\bigwedge\limits_{i<2n} [a_iRx]^{\iif(i< n)}\bigr).
\]

\mn
Now let $A' = A\cup\{a_0,\ldots,a_{2n-1}\}$ it cannot be a 
$\Gamma_R-m$-witness for $\varphi_2$. So
there are $b_0,\ldots,b_{2m-1}\in P^M\setminus A'$ with no repetition such that

\[
\neg(\exists x)\Bigl(\varphi(x,\bar a) \wedge
\bigwedge\limits_{i<2n}[b_iRx]^{\iif(i<m)}\Bigr).
\]

\mn
Clearly $a_0,\ldots,a_{2n-1},b_0,\ldots,b_{2m-1}\in P^M\setminus A$
are pairwise distinct and

\[
\neg(\exists x)\Bigl((\varphi_1\lor\varphi_2)\, \& \, 
\bigwedge\limits_{i<2n}[a_iRx]^{\iif(i<n)} \, \& \,
\bigwedge\limits_{i<2n}[b_iRx]^{\iif(i<m)}\Bigr).
\]

\mn
So $A$ is not a $\Gamma_R-(n+m)$-witness for $\varphi_1 \vee \varphi_2$.
\end{PROOF}

\begin{definition}
\label{c23}
Let $T^*$ be as in \ref{a2}(2)\footnote{ we need just some
schemes}, and $M^* \prec \gC$ be a model of $T^*$.
Let $a$, $c_i$ ($i<\delta$, $\delta$ a limit ordinal) be
members of $M^*$ such that in $M^*$:
\mn
\begin{enumerate}
\item[$(\alpha)$]   $a$ is a ``finite set"
\sn
\item[$(\beta)$]  $|a|\ge n$ for every true natural number $n$
\sn
\item[$(\gamma)$]  $c_i$ is a rational (or even real), $0<c_i<{1\over n}$
for every true natural number $n$
\sn
\item[$(\delta)$]  $2c_i<c_{i+1}$ and $c_i< c_j$ for $i<j < \delta$
\end{enumerate}
\mn
1) Let $\bar c = \langle c_i:i<\delta\rangle$.
We define the local bigness notion $\Gamma=\Gamma^{\ms}_{a,\bar c}$
as follows:  $\varphi(\bar x,\bar b)$ is $\Gamma$-big \underline{if
  and only if} $M^* \models ``|\{x:x$ is a member of 
$a,\varphi(x,\bar b)\}|/|a|$ is $>c_i$" for every $i<\delta$.

\noindent
2)  Let $\Gamma_\delta^{\ms}$ be the scheme whose instances are
$\Gamma_{a,\bar c}^{\ms}$, where $a$ and ${\bar c}=\langle c_i:
i<\delta\rangle$ are as above for $T^*$. 
Let $\Gamma^{\ms}$ mean $\bigcup\limits_{\delta} \Gamma^{\ms}_\delta$.

\noindent
3)  We say ``the smallness of $\varphi$ is witnessed by $c_i$" 
\If \, the quotient in part (1) is $<c_i$.

\noindent
4)  If $a,\bar c$ satisfy $(\alpha)-(\delta)$ we say that
$\bar c$ is an increasing sequence for $a$. It is called O.K. for $a$
if also
\mn
\begin{enumerate}
\item[$(\epsilon)$]   $|a|\times c_i > 1$ for every $i < \delta$.
\end{enumerate}
\mn
We may say that $\bar c$ is O.K. for $a$ or $a$, $\bar c$ are O.K.

\noindent
5)  We say $\bar c$ is wide \If \, for every $i<\delta$ and $n< \omega$
we have $c_{i+1}/ c_i > n$.
We say that $a$, $\bar c$ are wide for $M^*$ or $\bar c$ wide for $a$
in $M^*$) \If \, clauses $(\alpha) - (\delta)$ above and $\bar c$ is wide.
\end{definition}

\begin{remark} 
\label{c26}
On $\Gamma^{\ms}_{a,\bar c_1},\Gamma^{\ms}_{a,\bar c_2}$ being 
equivalent see \ref{c29}(4).
\end{remark}

\begin{claim}
\label{c29}
1) If for $T^*,M^*,a,\bar c$ are as in definition \ref{c23}
(so clauses $(\alpha)--(\delta)$ holds), \then \, 
$\Gamma=\Gamma^{\ms}_{a,\bar c}$ is a uniformly
$|\delta|^+$-simple $\ell$-bigness notion (with set of 
parameters $\{a\}\cup\bar c$) hence $\Gamma$ is co-simple.  If in addition
 $\bar c$ is O.K. for $a$ then $\Gamma$ is not trivial. 
If $\bar{c}$ is not O.K. then $\Gamma$ is trivial.

\noindent
2)  Suppose $T^*,M^*$ are as in Definition \ref{c23}, and for
$\ell=1,2$ $a^\ell,\bar c^{\ell}=\langle c^\ell_\alpha:\alpha<\delta
\rangle$ are as in Definition \ref{c23}(A) (for $T^*,M^*$),
and $\Gamma^\ell = \Gamma^{\ms}_{a^\ell,\bar c^\ell}$. \Then \, 
$\Gamma^1$, $\Gamma^2$ are orthogonal.

\noindent
3)  Let $\dot D$ be a filter say on $\kappa$; the bigness notion 
$\Gamma_\delta^{\ms},\Gamma^{\av}_{\dot D}$ are orthogonal
if $\cf(\delta)>\aleph_0$ (not $\kappa$!); also $\Gamma^{\ms}_{a,\bar c},
\Gamma^{\av}_{\dot D}$ are orthogonal if $\bar c$ is wide.

\noindent
4) If in $M^*,a,{\bar c}^{\ell} = \langle c_i^{\ell}:
i<\delta_\ell\rangle$ is as in Definition \ref{c23} 
for $\ell=1,2$ and $(\forall i<\delta_1)(\exists
j<\delta_2)[c_i^1 <c_j^2]$ and $(\forall j<\delta_2)(\exists
i<\delta_1)[c^2_j<c^1_i]$ \then \, $\Gamma^{\ms}_{a,{\bar c}^1}
= \Gamma^{\ms}_{a,{\bar c}^2}$.  Hence for $a,\bar c^1$ as above if
$\cf(\delta) > \aleph_0$ then for some $\bar c^2$ wide for $\bar a$,
we have $\Gamma^{\ms}_{a,\bar c^1_*} = \Gamma^{\ms}_{a,c^2}$.

\noindent
5) $\Gamma^{\ms}_{a,\bar c}$ is non-trivial (i.e. no algebraic
type is $\Gamma^{\ms}_{a,\bar c}$-big) \If \,
$\bar c$ is wide or just O.K. for $a$.

\noindent
6) If $a,\bar c=\langle c_i:i < \delta\rangle$ is as in Definition \ref{c23}
and $\cf(\delta) >\aleph_0$ or just $\omega \omega$ divides $\delta$
\then \, $\langle c_{\omega i}:\omega i < \delta\rangle$ is wide.
\end{claim}

\begin{PROOF}{\ref{c29}}
1) Note that $M^*\models$ ``$2c_i < c_{i+1}$". So
assume $\varphi(x,\bar b)= \varphi_1(x,\bar b_1) \vee \varphi_2(x,
\bar b_2)$ and $\varphi_\ell(x, \bar b_\ell)$ is $\Gamma$-small for
$\ell=1,2$. So for $\ell=1, 2$ for some $i_\ell<\delta$ the formula
$\varphi_\ell(x,\bar b_\ell)$ being $\Gamma$-small is witnessed
by $i_\ell$ so $M^* \models |\{x \dot e E a:
\varphi_\ell(x,\bar b_\ell)\}| \le c_{i_\ell} \times |a|$.

Hence

\begin{equation}
\begin{array}{clcr}
M^* \models |\{x \dot e a:\varphi(x,\bar b)\}| &\le
\sum\limits^{2}_{\ell=1} |\{x \dot e a:
\varphi_\ell(x,\bar b_\ell)\}| \\
   & \le \sum\limits^{2}_{\ell+1} c_{i_\ell}\times |a|\le 2 
c_{\max\{i_1,i_2\}}\times |a| \\
   & < c_{\max\{i_1,i_2\}+1} \times |a|
\end{array}
\end{equation}

\mn
so $\varphi(x,\bar b)$ is $\Gamma$-small as witnessed by $\max\{i_1,
i_2\} +1 < \delta$. The other facts are even easier.

\noindent
2) This is really a discrete version of Fubini theorem but we
shall elaborate. \Wilog \, $\bar{c}^\ell$ is O.K. for $\ell=1,2$.

Let $A_{\Gamma_1} \cup A_{\Gamma_2} \subseteq A \subseteq M^*$ and for
$\ell=1,2$ let $p_\ell(\bar x_{\Gamma_\ell})$ be a $\Gamma_\ell$-big type
over $A,p_\ell(\bar x_{\Gamma_\ell}) \in \bold S^{\alpha(\Gamma_\ell)}(A)$.

As each $\Gamma_\ell$ is co-simple let
$q := p_1(x)\cup p_2(y)\cup\{\psi_1(y,\bar d) \rightarrow 
\neg\varphi_1(x,y,\bar d):\bar d\subseteq A,\psi_1(y,\bar z)$ 
witness $\varphi_1(x;y,\bar z)$ is $\Gamma_1$-small$\} \cup \{\psi_2(x,\bar
d) \rightarrow \neg\varphi_2(y,x,\bar d):\bar d \subseteq A$ and 
$\psi_2(x,\bar z)$ witness $\varphi_2(y;x,\bar z)$ is $\Gamma_2$-small$\}$.

By \ref{a37}(1)(b) it suffice to prove that this set of formulas
is finitely satisfiable in $M^*$, assume not. So we have
(increasing the sequences of parameters from $M^*$ noting $p_\ell(\bar
x_{\Gamma_\ell})$ is closed under conjunctions) $\vartheta_1(x,\bar
d^*) \in p_1(x),\vartheta_2(y,\bar d^*) \in p_2(y),\psi_{1,k}(y,\bar
d^*) \longrightarrow \neg \varphi_{1,k}(x,y,\bar d^*)$
for $k<k_1$ from the third term in the union with smallness witnessed
by $c^1_{\alpha_k}$ (see \ref{c23}(3)) and $\psi_{2,k}(x,\bar d^*)
\longrightarrow \neg\varphi_{2,k}(y,x,\bar d^*)$
for $k<k_2$ from the fourth term in the union with smallness
witnessed by $c^2_{\beta_k}$ (see \ref{c23}(3)). Note
$k_1,k_2$ are true natural numbers. Choose $\alpha(*)<\delta_1$
such that $(\forall k<k_1)[\alpha_k+k_1<\alpha(*)]$, and choose
$\beta(*)<\delta_2$ such that $(\forall k<k_2)[\beta_k+k_2<\beta(*)]$.
Let (recalling $\dot e$ is membership in $M^*$'s sense)
$Z=\{(x,y):x \dot e a^1 \, \& \, y \dot e a^2 \, \& \, \vartheta_1(x,\bar
d^*)\, \& \, \vartheta_2(y,\bar d^*)\}$.   
So $Z$ is (representable) in $M^*$ (we do not distinguish). 

Let (all is $M^*$'s sense):

\[
Z_1 = \{(x,y) \dot e Z:(\exists k<k_1)[\psi_{1,k}(y,\bar d^*)\, \& \,
\varphi_{1,k}(x;y,{\bar d}^*)]\}
\]

\[
Z_2 = \{(x,y) \dot e Z:(\exists k<k_2)[\psi_{2,k}(x,d^*)\, \& \,
\varphi_{2,k}(y;x,{\bar d}^*)]\}.
\]

\mn
So
\mn
\begin{enumerate}
\item[$(a)$]  $Z=Z_1\cup Z_2$ (by the ``assume not" above)
\sn
\item[$(b)$]  for every $y,|Z_1^{[y]}|/ |a^1| \le c^1_{\alpha(*)}$ where
$Z_1^{[y]}=\{x:(x,y) \dot e Z_1\}$.
\end{enumerate}
\mn
[Why?  As a union of $k_1$ sets, each with
$\le c^1_{\alpha(*)-{k_1}} \times |a^1| \le \frac{1}{2^{k_1}}\times
c^1_{\alpha(*)}\times |a^1|$ members has $\leq c^1_{\alpha
(*)}\times |a^1|$ members.]

And similarly
\mn
\begin{enumerate}
\item[$(c)$]  for every $x,|Z_2^{[x]}|/ |a^2|\le c^2_{\beta(*)}$ where
$Z_2^{[x]}=\{y:(x,y) \dot e Z_2\}$.
\end{enumerate}
\mn
There is in $M^*$ a set $X$, such that:

\[
x \dot e X \Rightarrow x \dot e a^1 \, \& \,
\vartheta_1(x,\bar d^*)
\]

\[
(c^1_{\alpha(*)+2})|a^1|+1 \ge |X| \ge (c^1_{\alpha(*)+2})\times |a^1|.
\]

\mn
There is in $M^*$'s sense a set $Y$ such that

\[
y \dot e Y \Rightarrow y \dot e a^2 \, \& \, \vartheta_2(y,\bar d^*)
\]

\[
(c^2_{\beta(*)+2})|a^2|+1 \ge |Y|> (c^2_{\beta(*)+2})|a^2|.
\]

\mn
By the choice of $Z$ clearly $X \times Y\subseteq Z$. So in $M^*$ (for the
fifth line recall $\bar{c}$ is O.K., for the seventh line 
recall $2c_i < c_{i+1}$)

\begin{equation}
\begin{array}{clcr}
(c^1_{\alpha(*)+2}|a^1|)(c^2_{\beta(*)+2}|a^2|) &\le |X|\times|Y|=
|X\times Y| \le |Z_1 \cap (X\times Y)|+|Z_2\cap (X\times Y)| \\
  &\le |X|(c^1_{\alpha(*)}|a^1|) + |Y|(c^2_{\beta(*)}|a^2|) \\
  &\le (c^1_{\alpha(*)} |a^1|)(c^2_{\beta(*)+2} |a^2|+1)
+(c^1_{\alpha(*)+2} |a^1|+1)(c^2_{\beta(*)} |a^2|) \\
  &\le (c^1_{\alpha(*)} |a^1|)(2c^2_{\beta(*)+2}|a^2|)
+(2c^1_{\alpha(*)+2}|a^1|)(c^2_{\beta(*)} |a^2|) \\
  &= (2c^1_{\alpha(*)}|a^1|)(c^2_{\beta(*)+2}|a^2|)
+(c^1_{\alpha(*)+2}|a^1|)(2c^2_{\beta (\ast)}|a^2|) \\
  &< (\frac{1}{2} c^1_{\alpha(*)+2}|a^1|)(c^2_{\beta(*)+2}|a^2|)
+(c^1_{\alpha(*)+2}|a^1|)(\frac{1}{2} c^2_{\beta(*)+2}|a^2|) \\
  = c^1_{\alpha(*)+2} |a^1|c^2_{\beta(*)+2} |a^2|
\end{array}
\end{equation}

\mn
contradiction.

\noindent
3) Let $\Gamma_1 = \Gamma^{\av}_{\dot D,\langle b^*_\epsilon:
\epsilon < \kappa\rangle},\Gamma_2 = \Gamma^{\ms}_{a,\bar c},
{\bar c}=\langle c_i:i<\delta\rangle$.

By (4)+(6) below (which does not depend on \ref{c29}(3) ) it
suffices to prove the second case i.e. prove orthogonality assuming
$(\forall n<\omega)(\forall i < \delta) ``n< c_{i+1}/c_i$" 
(as if $\cf (\delta)>\omega$, letting
${\bar c}'=\langle c_{\omega\times j}:\omega j <\delta\rangle$
we have $\Gamma^{\ms}_{a,\bar c}$, $\Gamma^{\ms}_{a,{\bar c}'}$ are
equal and ${\bar c}'$ satisfies the requirement above).

Given $A \subseteq M^*$ such that $b^*_\epsilon,a,c_j \in A$ (for 
$\epsilon < \kappa,j<\delta$), and $\Gamma_{\ell}$-big 
$p_{\ell} = \tp(b_{\ell},A,M^*)$ (for $\ell=1,2$), 
possibly increasing $M^*$ we can find $c\in M^*$ such that
$M^*\models$ ``$c$ a natural number, $n<c,2^{(c^n)} < \log_2
(c_{i+1}/c_i)$" for every $n<\omega,i<\delta$, and we can find 
a pseudo-finite set $d\in M^*$, such that: 
$M^* \models ``|d|=c \, \& \, b^*_{\epsilon} \dot e d$" 
for $\epsilon < \kappa$.

Let $A_0=A,A_1=A_0+d$ and let $p'_2 \in \bold S(A_1, M^*)$ be a 
$\Gamma_2$-big extension of $p_2$.
We shall show now that (\wilog \,) there are $e^*_\varphi\in M^*$
(for the $\varphi=\varphi (x,{\bar z}_\varphi)\in \bbL(T^*[A_1])$) such
that $M^*\models ``e^*_\varphi\subseteq {}^{n(\varphi)} d$" 
where $n(\varphi)= \ell g({\bar z}_\varphi)$ and
\mn
\begin{enumerate}
\item[$(*)$]   $p^*_2 = p'_2 \cup \{(\forall {\bar z} \dot e 
{}^n d)[{\bar z} \dot e e^*_\varphi \Leftrightarrow \varphi(x,\bar z)]: 
\varphi = \varphi(x,\bar z) \in \bbL_\tau(T^*[A_1])\}$ is $\Gamma_2$-big.
\end{enumerate}
\mn
Why $(*)$ holds?  As $\Gamma^{\ms}_{a,\bar c}$ is a local co-simple 
bigness notion, we can also replace $p'_2$ by one formula say 
$\vartheta_2(x,{\bar b}^*)$, and consider only $\varphi_1,\ldots,
\varphi_n \in \bbL(T^*[A_1])$ for some $n<\omega$. 
We can find $\varphi=\varphi(x,\bar z)$ such that for parameters from
${}^{n(\varphi)}d$ we get all the instances of $\varphi_1,\ldots,\varphi_n$ by 
increasing $\bar z$,  hence \wilog \, we can consider just
one $\varphi=\varphi(x,{\bar z}_\varphi)$. Let
$n = n(\varphi) = \ell g({\bar z}_\varphi)$, then $\vartheta(M^*,{\bar b}^*)$
has size $\ge c_\alpha|a|$ (for every $\alpha<\delta$, internal sense), it is
divided to $\le |\cP({}^{n(\varphi)}d)|$ parts according to the $\varphi$-type
over $d$, so the largest one (internal sense) is as required (or find
a right $e^*_\varphi$ in a $|\delta|^+$-saturated extension of
$M^*$) so really $(*)$ holds.

So we can find $p''_2 \in \bold S(A_2,M^*)$ extending $p^*_2$ which is
$\Gamma_2$-big where $A_2 = A_1\cup\{e^*_\varphi:\varphi\in \bbL(T^*[A_1])\}$.

Let $b'_2 \in M^*$ realize $p''_2,b'_1$ realize a
$\Gamma_1$-big $p''_1 \in \bold S(A_2+b_2',M^*)$ extending $p_1$ 
so clearly $[x \dot e d] \in \bar p''_1$.
They are as required (think or see \ref{c86}(3)).

\noindent
4), 5), 6) Trivial.                    
\end{PROOF}

\begin{remark}
You may wonder whether we can weaken the demand on $T$, still
demanding that $a$ behave like a finite set. Certainly we can,
e.g. by using restriction, see Definition \ref{a41}. 
We can do it in a more finely tuned way, we hope to
deal with it elsewhere.
\end{remark}

\begin{definition}
\label{c32}
1) Suppose that $T^*,\gC$ are as in \ref{a2}(B),
$a$ in $\gC$, $a$ is a pseudo finite set, $\dot w$ 
a function from $a$ to $[0,1]_{\bbR}$ such that 
$\gC \models ``\sum\limits_{x \dot e a}
\dot w(x)=1"$, and $\bar c = \langle c_i :i < \delta\rangle$ 
($\delta$ a limit ordinal) is an increasing sequence for $a$.
For every $a'\subseteq a$ let $\dot w(a')$ be
$\sum\limits_{x \dot e a'} \dot w(x)$, in $\gC$'s sense; 
if confusion may arise we shall write $\dot w(\{x\})$ for 
$x \dot e a$. 

We define $\Gamma=\Gamma^{\wm}_{a,\dot w,\bar c}$, a local 
bigness notion, by: $\psi(x,\bar b)$ is $\Gamma$--small \underline{if
  and only if} for some $i<\delta,\gC
\models ``\dot w(\{x \dot e a: \psi(x,\bar b)\}) \le c_i"$.

\noindent
2)  Assume above that in $\gC$ we have:
$a^+ = \cP(a),\bold q \dot e (0,1)_{\bbR}$.  Let
$\dot w_{\bold q} = \dot w_{a,\bold q}:a^+ \rightarrow
[0,1]_{\bbR}$ be defined by

\[
\dot w(b) = \bold q^{|b|}(1-\bold q)^{|a\setminus b|}.
\]

\mn
We let $\Gamma^{\wmg}_{a,\bold q,\bar c} = \Gamma^{\wm}_{\cP(a)},
\dot w_{\bold q,\bar c}$. If $\bold q = c_0$ we write 
$\Gamma^{\wmg}_{a,\bar c}$, we always assume $\bold q \le 1/2$.

\noindent
3) We say $\bar c$ is O.K. for $(a,\dot w)$ or $(a,\dot w,\bar c)$ is
O.K. \If \, for some $i< \ell g(\bar c)$ we have $\gC \Vdash
 ``(\forall x \dot e a) \dot w(x)< c_i$", we 
normally assume that this holds for $i=0$.

\noindent
4) We say $\bar c$ is wide for $(a,\dot w)$ or ($a,\dot w,\bar c$) is
wide \If \, $(\forall i < \ell g(\bar c))(\forall n < \omega)[c_{i+1}/c_i>n]$.
\end{definition}

\begin{remark}
\label{c35}
Remember $\ln (1-x)\sim\ -x$ for $x\in (0,1/2)$, more exactly 
for the natural logarithm, $|\ln (1-x) + x|< x^2$
for every $x\in (0,\frac{1}{2})$.  Why?  Because by the Taylor series,
$\ln(1-x) = (x + \frac{x^2}{2} + \frac{x^3}{3} + \ldots)$ and
$\frac{x^3}{3} + \frac{x^4}{4} + \ldots < \frac{1}{3} x^3 
+ \frac{1}{3} x^3(1+x + x^2 + \ldots) = \frac{1}{3} x^3 + 
\frac{1}{3} x^3 \frac{1}{1-x} \le \frac{1}{3} x^3
\frac{1}{1-1/2} = x^3 < \frac{1}{2} x^2)$.
\end{remark}

\begin{claim}
\label{c38}
1) For $T^*,M^*,a,\dot w,\bar c$ as in Definition \ref{c32}(A),
\then \, $\Gamma=\Gamma^{\wm}_{a,\dot w,\bar c}$ 
is a uniformly $|\lg(\bar c)|^+$-simple $\ell$-bigness
notion (with set of parameters $\{a\} \cup \{\dot w\} \cup 
\bar c$ so $\Gamma$ is co-simple).

\noindent
2) Suppose $T^*$, $M^*$ are as in Definition \ref{c32}(A), and for
$\ell=1,2$ $a^\ell,\dot w^\ell,\bar c^{\ell}=\langle
c^\ell_\alpha:\alpha<\delta\rangle$ are as in Definition \ref{c32}(A)
(for $T^*$, $M^*$), and $\Gamma^\ell =
\Gamma^{\wm}_{a^\ell,w^\ell,\bar c^\ell}$. \Then \,  $\Gamma^1$,
$\Gamma^2$ are orthogonal.

\noindent
3) $\Gamma_\delta^{\wm},\Gamma^{\av}_{\dot D}$ are orthogonal
\If \, $\cf(\delta)>\aleph_0;\Gamma^{\wm}_{a,\dot w,\bar c}$
$\Gamma^{\av}_{\dot D}$ are orthogonal if 
$\bar c$ is wide or $\omega \omega$ divides $\delta$.

\noindent
4)  If in $M^*,{\bar c}^{\ell} = \langle c_i^{\ell}:
i<\delta_\ell\rangle$ for $\ell=1,2$ and $(\forall i<\delta_1)(\exists
j<\delta_2)[c_i^1 <c_j^2]$ and $(\forall j<\delta_2)(\exists
i<\delta_1)[c^2_j<c^1_i]$ \then \,  $\Gamma^{\wm}_{a,\dot t,{\bar c}^1}
=\Gamma^{\wm}_{a,\dot w,{\bar c}^2}$.

\noindent
5)  $\Gamma^{\wm}_{a,\dot w,\bar c}$ is not trivial (i.e. if $p(x)$ is
$\Gamma^{\wm}_{a,\dot w,\bar c}$-big \then \, $p(x)$ is not an
algebraic type) \underline{if and only if} $(a,\dot w,\bar c)$ is
O.K. 

\noindent
6) $\Gamma^{\ms}_{a,\bar c}$ is a special case of 
$\Gamma^{\wm}_{a,\dot w,\bar c}$.
\end{claim}

\begin{PROOF}{\ref{c38}}
Like \ref{c29} except that in part (2) in the case $(a,\dot w,
\bar{c})$ is not O.K we have to take more care (and this case is not used).
\end{PROOF}

\noindent
A ``dual" notion is
\begin{definition}
\label{c41}
Let $T^*$ be as in \ref{a2}(B)\footnote{ we need just some
schemes}, and $M^*$ be a model of $T^*$.
Let $a$, $c_i$ ($i<\delta$, $\delta$ a limit ordinal) be
members of $M^*$ such that in $M^*$:
\mn
\begin{enumerate}
\item[$(\alpha)$]   $a$ is a ``finite" set
\sn
\item[$(\beta)$]   $|a|\ge n$ for every true natural number $n$
\sn
\item[$(\gamma)$]  $c_i$ is a ``rational", $0 < c_i < \frac{1}{n}$ for every
true natural number $n$
\sn
\item[$(\delta)$]  $c_i> 2c_{i+1}$ and $i<j \Rightarrow c_i> c_j$.
\end{enumerate}
\mn
1)  Let $\bar c=\langle c_i:i<\delta\rangle$.
We define the $\ell$-bigness notion $\Gamma=\Gamma^{\dms}_{a,\bar c}$
as follows:  $\varphi(\bar x,\bar b)$ is $\Gamma$-big \underline{if
  and only if}  $M^* \models ``|\{x:x$ is a member of
$a,\varphi(x,\bar b)\}|/|a|$ is $>c_i$" for some $i<\delta$.

\noindent
2) Let $\Gamma_\delta^{\dms}$ be the scheme whose instances are
$\Gamma_{a,\bar c}^{\dms}$, where $a$ and ${\bar c}=\langle c_i :
i<\delta\rangle$ as above for $T^*$.  Let $\Gamma^{\dms}$ mean 
$\bigcup\limits_{\delta} \Gamma^{\dms}_\delta$.

\noindent
3)  We shall say ``the bigness of $\varphi$ is witnessed by $c_i$" 
\If \, the quotient above is $>c_i$.

\noindent
4) $\bar c$ is decreasing sequence if: clauses $(\alpha)$ -- $(\delta)$ above
hold.  It is d-O.K. for $a$ if $\frac{n}{|a|} < c_i< \frac{1}{n}$
for $n \le \omega$. It is wide if $c_i/c_{i+1} > n$ for
$i< \delta$, $n<\omega$.
\end{definition}

\begin{claim}
\label{c44}
1)  For $T^*,M^*,a,\bar c$ as in Definition \ref{c41},
\then \,  $\Gamma=\Gamma^{\dms}_{a,\bar c}$ is a uniformly
$|\delta|^+$-co-simple $\ell$-bigness notion (with set of parameters 
$\{a\}\cup\bar c$) hence is simple;
$\Gamma$ is non-trivial if $\bar{c}$ is d-O.K..

\noindent
2) Suppose $T^*,M^*$ are as in Definition \ref{c41}, and for
$\ell=1,2$ $a^\ell,\bar c^{\ell} = \langle c^\ell_i:i < 
\delta\rangle$ are as in Definition \ref{c41}
(for $T^*$, $M^*$), and $\Gamma^\ell =
\Gamma^{\dms}_{a^\ell,\bar c^\ell}$. \Then \, $\Gamma^1,\Gamma^2$ 
are orthogonal.

\noindent
3)  $\Gamma_\delta^{\dms},\Gamma^{\av}_{\dot D}$ are orthogonal
if $\cf(\delta)> |\Dom(\dot D)]$ for any filter $\dot D$.

\noindent
4) If in $M^*,{\bar c}^{\ell} = \langle c_i^{\ell}:
i<\delta_\ell\rangle$ for $\ell=1,2$ and $(\forall i<\delta_1)(\exists
j<\delta_2)[c_j^2 <c_i^1]$ and $(\forall j<\delta_2)(\exists
i<\delta_1)[c^1_i<c^2_j]$ \then \,  $\Gamma^{\dms}_{a,{\bar c}^1}
=\Gamma^{\dms}_{a,{\bar c}^2}$.

\noindent
5)  In Definition \ref{c41}, if $\cf(\delta)>\aleph_0$ 
or just $\omega \omega$ divides $\delta$ then $\bar c'=
\langle c_{\omega i}:\omega i < \delta\rangle$ is wide and 
$\bar c,\bar c'$ are like $\bar c^1,\bar c^2$ in part (4).

\noindent
6) Assume that $T^*,M^*,a^1,\bar c^1$ are as in Definition 
\ref{c23}, $\Gamma_1 = \Gamma^{ms}_{\bar a^1,\bar c^1}$, and 
$(T^*,M^*),a^2,\bar c^2$ are as in Definition \ref{c41} and 
$\Gamma_2=\Gamma^{\dms}_{\bar a^2,\bar c^2}$. \Then \, 
$\Gamma_1,\Gamma_2$ are orthogonal. 
\end{claim}

\begin{PROOF}{\ref{c44}}
1) Note that $M^* \models ``c_i > 2c_{i+1}"$.
So assume $\varphi(x,\bar b) = \varphi_1(x,\bar b_1) \vee 
\varphi_2(x,\bar b_2),\varphi_\ell(x,\bar b_\ell)$ is 
$\Gamma$-small for $\ell=1,2$.  
So for every $i< \delta,\ell=1,2$ we have $M^* \models
``|\{x \dot e a:\varphi_\ell(x,\bar b_\ell)\}| \le c_{i+1}
\times |a|"$ hence $M^* \models ``|\{x \dot e a:
\varphi(x,\bar b)\}|/|a| \le 2 c_{i+1}< c_i$; as this holds 
for each $i$, clearly $\varphi(x,\bar b)$ is $\Gamma$-small. 
The other parts in the Definition are even easier.

\noindent
2) This is really a discrete version of Fubini theorem.

Consider, assuming $b_\ell$ realizes $p_\ell \in \bold S(A,M^*)$ which is
$\Gamma^{\dms}_{a_\ell,\bar c^\ell}$-big, 
$q := p_1(x) \cup p_2(y) \cup \{\neg\varphi_1(x,y,\bar d):\bar d 
\subseteq A,\varphi_1(x;b_\ell,\bar d)$ is $\Gamma_1$-small$\} \cup 
\{\neg\varphi_2(y,x,\bar d):\bar d \subseteq A$ and 
$\varphi_2(y;b_1,\bar d)$ is $\Gamma_2$-small$\}$.

By \ref{a37}(1) it suffice to prove that this set of formulas
is finitely satisfiable in $M^*$, toward contradiction assume not. So we have
(increasing the sequences of parameters from $M^*$ and recalling
$p_\ell$ is closed under conjunctions) $\vartheta_1(x,\bar d^*) \in 
p_1(x),\vartheta_2(y,\bar d^*) \in p_2(y),\neg\varphi_{1,k}(x,y,\bar d^*)$
for $k<k_1$ from the third term in the union and
$\neg\varphi_{2,k}(y,x,\bar d^*)$ for $k<k_2$ from the fourth 
term in the union. Note $k_1,k_2$ that are true natural numbers. 
Choose $\alpha(*)<\delta_1$
such that $\vartheta_1(x, \bar d^*)$ being $\Gamma_1$-big is witnessed by
$c^1_{\alpha(*)}$, and choose
$\beta(*)<\delta_2$ such that $\vartheta_2(y,\bar d^*)$ being
$\Gamma_2$-big is witnessed by $c^2_{\beta(*)}$. \Wilog \,  
for every $a \in \vartheta_1(M^*,\bar d^*),M^* \vDash 
``|\{y:\varphi_{2,k}(y;a,\bar d^*)\}|/|a^2| < 
c^2_{\beta(*)+k_2+2}"$ and for every $b \in \vartheta_2(M^*,\bar d),
M^* \vDash ``|\{x:\varphi_{1,k}(y;x,\bar d^*)\}|/|a^1| <
c^1_{\alpha(*)+k_1+2}"$.

Let (recall $\dot e$ is membership in $M^*$'s sense)
$Z = \{(x,y):x \dot e a^1 \, \& \, y \dot e a^2\ \&\ \vartheta_1(x,\bar
d^*) \, \& \, \vartheta_2(y,\bar d^*)\}$. 
So $Z$ is (representable) in $M^*$ (we
do not distinguish). 

Let

\[
Z_1 = \{(x,y) \dot e Z:(\exists k<k_1)\varphi_{1,k}(x;y,{\bar
  d}^*)\}
\]

\[
Z_2 = \{(x,y) \dot e Z:(\exists k<k_2)\varphi_{2,k}(y;x,{\bar
  d}^*)\}.
\]

\mn
So
\mn
\begin{enumerate}
\item[$(a)$]  $Z=Z_1 \cup Z_2$ (by the ``assume not" above)
\sn
\item[$(b)$]  for every $y,|Z_1^{[y]}|/|a_1| \le c^1_{\alpha(*)+2}$ where
$Z_1^{[y]}=\{x:(x,y)\in Z_1\}$.
\end{enumerate}
\mn
[Why?  As a union of $k_1$ sets each with
$<c^1_{\alpha(*) + {k_1}+2} \times |a^1| \le \frac{1}{2^{k_1}}\times
c^1_{\alpha(*)+2}\times |a^1|$ members has $\le c^1_{\alpha
(*)+2}\times |a^1|$ members.]

And similarly
\mn
\begin{enumerate}
\item[$(c)$]   for every $x,|Z_2^{[x]}|/ |a_2| \le c^2_{\beta(*)+2}$ where
$Z_2^{[x]}=\{y:(x,y)\in Z_2\}$.
\end{enumerate}
\mn
There is in $M^*$'s sense a set $X$, such that:

\[
x \dot e X \Rightarrow x \dot e a^1 \, \& \,
\vartheta_1(x,\bar d^*)
\]

\[
(c^1_{\alpha(*)})|a^1|+1 \ge |X|> (c^1_{\alpha(*)})|a^1|.
\]

\mn
There is in $M^*$'s sense a set $Y$ such that

\[
y \dot e Y \Rightarrow y \dot e a^2 \, \& \,
\vartheta_2(y,\bar d^*)
\]

\[
(c^2_{\beta(*)})|a^2|+1 \ge |Y| \ge (c^2_{\beta(*)})|a^2|.
\]
By the choice of $Z$ clearly $X\times Y\subseteq Z$. 

So in $M^*$

\begin{equation*}
\begin{array}{clcr}
(c^1_{\alpha(*)}|a^1|)(c^2_{\beta(*)}|a^2|) &< |X| 
\times |Y| = |X \times Y| \\
  &\le |Z_1\cap (X\times Y)|+|Z_2\cap (X\times Y)| \le
(c^1_{\alpha(*)+2}|a^1|)|Y| + |X|(c^2_{\beta(*)+2}|a^2|) \\
  &\le c^1_{\alpha(*)+2} |a^1|(c^2_{\beta(*)} |a^2|+1) +
(c^1_{\alpha(*)} |a^1|+ (c^2_{\beta(*)+2} |a^2|) \\
  &\le (c^1_{\alpha(*)+2} |a^1|)(2c^2_{\beta(*)}|a^2|)
+(2c^1_{\alpha(*)}|a^1|)(c^2_{\beta(*)+2} |a^2|) \\
  &= (2c^1_{\alpha(*)+2}|a^1|)(c^2_{\beta(*)}|a^2|)+(c^1_{\alpha(*)}
|a^1|)(2c^2_{\beta(*)+2}|a^2|) \\
  &< (\frac{1}{2} c^1_{\alpha(*)}|a^1|)(c^2_{\beta(*)}|a^2|)
+(c^1_{\alpha(*)}|a^1|)(\frac{1}{2} c^2_{\beta(*)}|a^2|) \\
  &= c^1_{\alpha(*)} |a^1|c^2_{\beta(*)} |a^2|
\end{array}
\end{equation*}

\mn
contradiction.
\noindent
3) Let $\Gamma_1 = \Gamma^{\av}_{\dot D,\langle b^*_\epsilon:
\epsilon<\kappa\rangle}$, $\Gamma_2=\Gamma^{\dms}_{a,\bar c},\bar c
=\langle c_i:i<\delta\rangle$.

Let $p_\ell = \tp(b_\ell,A,M^*)$ be $\Gamma_\ell$-big for $\ell=1,
2$; \wilog \, $b^*_\epsilon,a,c_i \in A$ for $\epsilon < \kappa,i<
\delta$;  let $\dot D_1$ be an ultrafilter 
on $\kappa$ extending $\dot D$. Now, possibly increasing $M^*$
\wilog \, $b_1$ realizes $\Av_{\dot D_1}(\langle b^*_\epsilon: \epsilon<
\kappa\rangle, A+b_2, M^*)$ so clearly $\tp(b_1, A+b_2, M^*)$ is
$\Gamma_1$-big. Now assume $\varphi(x,b_1,\bar d) \in 
\tp(b_2,A+b_1,M^*)$ hence ${\cU} =: \{\epsilon <\kappa:
\varphi(x,b^*_\epsilon,\bar d) \in \tp(b_2,A,M^*)\} \in \dot D_1$, clearly for 
$\epsilon \in {\cU}$ the formula $\varphi(x, b^*_{\epsilon},\bar d)$
is $\Gamma_2$-big (belonging to $\tp(b_2, A,M^*)$ which is 
$\Gamma_2$-big), so for $\epsilon \in {\cU}$, let
$i_\epsilon = \min\{i: \varphi(x, b^*_\epsilon, \bar d)$ is
$\Gamma^*_2$-big as witnessed by $c_i\}$; it is well defined, so let
$i(*)=\sup\{i_\epsilon:\epsilon\in {\cU}\}$; now $i(*) < \delta$ as 
$\cf(\delta) > \kappa \ge |{\cU}|$, and so easily
$\varphi(x, b_1, \bar d)$ is $\Gamma_2$-big being  witnessed by $c_{i(*)}$.

\noindent
4),5)  Easy.

\noindent
6) We let $A_{\Gamma_\ell}\subseteq A\subseteq M^*$ and 
$p_\ell = \tp(b_\ell,A,M)$ is $\Gamma_\ell$-big for $\ell=1,2$, 
\wilog \, $M$ is $|A|^+$-saturated. We can find a $||M^*||^+$-saturated 
$N^*$ \st \, $M^*<N^*$. We can find $\langle c_n^3:n<\omega\rangle$ in 
$N^*$ \st \,: $N^* \models ``c^3_n < c^3_{n+1} < c^2_\alpha"$ for $n<\omega,
\alpha < \ell g(\bar c^2)$ \st \, if $c \in (0,1)^{M^*}_{\bbR}, 
M^* \models ``c < c^2_\alpha"$ for $\alpha < \ell g(\bar c^2)$ then 
$N^* \models ``c < c^3_n"$.  Now $p_2$ is $\Gamma^{\ms}_{a^2,\bar
  c^3}$-big and $\Gamma^{\ms}_{a^2,\bar c^3} \perp
\Gamma^{\ms}_{a^1,\bar c^1}$, and so we can find in $N^3$ elements 
$b'_1,b'_2$ realizing $p_1,p_2$ respectively \st \, 
$\tp(b'_1,A+b'_2,N^*)$ is $\Gamma^{\ms}_{a^1,\bar c^1}$-big and 
$\tp(b_1,A+b_1,N^*)$ is $\Gamma_{a^2,\bar c^3}$-big: so 
$b'_1,b'_2$ exemplify the desirable result.    
\end{PROOF}

\begin{definition}
\label{c47}
Suppose $T^*,\gC,a,\bar c$ are as in Definition
\ref{c41}, so is a ``finite" and $\dot w$ is a function 
from $a$ to $[0,1]_{\bbR}$ (in $\gC$-s sense such that
$\gC \models \sum \{\dot w(x):x \in a\}=1$).

\noindent
1) Let $\Gamma=\Gamma^{\dws} = \Gamma^{\dws}_{a,\dot w,\bar c}$ 
be the following local bigness notion: $\varphi(x,\bar b)$ is 
$\Gamma$-small \underline{if and only if} for every $i<\delta$
we have $\gC \models ``c_i > \Sigma\{\dot
w(x):\varphi(x,\bar b) \, \& \, x \in a\}$.

\noindent
2) We say $\bar c$ is d-O.K. for $a,\dot w$ \underline{if and only if}
for $i < \ell g(\bar c)$ and every true natural number $n$ and 
$b_1,\ldots,b_n \in a$ we have $M^* \models \dot w(b_1) + \ldots +
\dot w(b_n) <c_i < 1/n$ (this is retained in any elementary extension 
of $M^*$).

\noindent
3)  We define $\Gamma^{\dms}_{a,\bold q,\bar c}$ parallely to 
\ref{c32}(2).
\end{definition}

\begin{claim}
\label{c50}
The parallel of \ref{c44} holds for $\Gamma^{\dws}$.
\end{claim}

\begin{PROOF}{\ref{c50}}
Similar to the proof of \ref{c44}.
\end{PROOF}

\begin{claim}
\label{c53}
1) Suppose $M^* \models T^*$ and (in $M^*$) we have: $(a,\bar c)$ 
is increasing and O.K., $A \subseteq M^*,\{a\}\cup \bar c \subseteq A,
c_\beta/c_{\beta+2} \le c_{\beta+n}$ (for every $\beta$ for some 
$n = n_\beta$, e.g. $n=3$) and 
$p = \tp(b,A,M^*)$ is $\Gamma^{\ms}_{a,\bar c}$-big, 
and $a=\cP(a_1)$.  Assume further $A \subseteq
A'\subseteq M^*$, and $B^* := \{d:d \in \acl(A'),d \in {}^{M^*} a_1$,
but for no $a_2$ do we have:
$a_2 \in \acl(A),M^* \models ``d \dot e a_2 \, \& \, a_2\subseteq a_1
\,\& \, (\forall \alpha < \ell g(\bar c))[2^{-|a_2|} \ge c_\alpha]"\}$.
If $h:B^* \rightarrow \{\true,\false\}$ \then \,:
\mn 
\begin{enumerate}
\item[$(\alpha)$]  $p' = p \cup \{[d \in x]^{\iif(h(d))}:d \in B^*\}$ is
$\Gamma^{\ms}_{a,\bar c}$-big,
\end{enumerate}
\mn
moreover for some $b'$ realizing $p'$ in some $N,M^* \le_{\gk} N$  we have
\mn
\begin{enumerate}
\item[$(\beta)$]  $\tp(b',A',N)$ is $\Gamma^{\ms}_{a,\bar c}$-big
\sn
\item[$(\gamma)$]  $\acl_N(A+b') \cap \acl_N(A') = \acl_N(A)$.
\end{enumerate}
\mn
2) If in addition $\log_2(1/c_\beta)/\log_2(1/c_{\beta+1})<c_{\beta+2}$ for 
$\beta<\delta$ \then \, we can add:
\mn
\begin{enumerate}
\item[$(\delta)$]  if $d \in B^*$, \then \, for no $a_2$ do we have: 
$a_2 \in \acl(A+b')$ and $M^* \models``d \dot \epsilon a_2 \subseteq 
a_1 \, \& \, (\forall \alpha < \ell g(\bar c))[2^{-|a_2|}\ge c_\alpha]"$.
\end{enumerate}
\end{claim}

\begin{PROOF}{\ref{c53}}
1) The proof of part 2) is similar only slightly harder, so read it.

\noindent
2) \Wilog \, $A = \acl(A)$ and $A'=\acl  (A')$. 
It suffices to prove $(\alpha)+(\delta)$, meaning:
\mn
\begin{enumerate}
\item[$\boxtimes$]  the following is a $\Gamma^{\ms}_{a,\bar c}$-big type

\begin{equation*}
\begin{array}{clcr}
q^*(x) :=  &p(x) \cup \{[d \in x]^{\iif(h(d))}:d\in B^*\} 
\cup \{\neg d \in \sigma(x):d\in B^* \text{ and}\\
  & \sigma(-) \text{ a term with parameters from } A \text{ such that} \\
  & p(x) \vdash ``\sigma(x) \text{ a finite set } \& \, 2^{-|\sigma(x)|}\ge
   c_\beta" \text{ for every } \beta < \ell g(\bar c)\}.
\end{array}
\end{equation*}
\end{enumerate}
\mn
Why is $\boxtimes$ enough?  As then we can get
$(\beta)$ by the extension property of 
``$\Gamma^{\ms}_{a,\bar c}$-big" types, and to get clause $(\gamma)$
we replace $A',h$ by $\bigcup\limits_{i<\lambda} A'_i,
\bigcup\limits_{i<\lambda} h_i$ where $\lambda> |T| + |A'| +
\aleph_0,{\bold f}_i$ is an elementary mapping with domain $\acl(A')$
in $N$ where $M^* \le N,{\bold f}_i \restriction \acl(A) = 
\id_{\acl(A)}$ for $i< \lambda$ and $\langle \bold f_i(\acl(A')) 
\setminus \acl(A): i< \lambda\rangle$ are pairwise disjoint, 
and $A'_i = \bold f_i(A')$ and $h_i = h \circ \bold f^{-1}_i$.

Now for simplicity assume $\Rang(\bold f_i) \subseteq M^*$. Now we apply
$\boxtimes$ with $A'' = \cup\{A_i:i<\lambda\},h' = 
\cup\{h_i:i<\lambda\}$ instead $A,h$, easily the desired 
conclusion follows; so indeed it is enough to prove $\boxtimes$.
(We suppress below parameters from $A$.)

Let ${\bf \Sigma} = \{\sigma(x): \sigma$ a term with parameters from $A$ such
that $M^* \vDash ``(\forall x)[|2^{-|\sigma(x)|}| \ge c_\beta \, \& \,
\sigma(x)$ a finite subset of $a_1$]" for every $\beta< \delta\}$. 
Note that if $\sigma_1(x),\sigma_2(x) \in {\bf \Sigma}$ and 
$\sigma(x) := \sigma_1(x) \cup \sigma_2(x)$, in $M^*$'s sense, 
then $\sigma(x) \in {\bf \Sigma}$ as
$\beta < \delta \Rightarrow 2^{-|\sigma(x)|} = 2^{-|\sigma_1(x)|} \cdot
2^{-\sigma_2(x)} \ge c_{\beta +n} \times c_{\beta +2} \ge c_\beta$ 
when $n=n_\beta$ is as in the assumption (see \ref{c53}(1)).
So if the type $q^*$ is not $\Gamma^{\ms}_{a, \bar c}$-big then for some
$\vartheta (x)\in p$, $k<\omega$, $\beta(*)< \delta$, $\sigma(x)\in
{\bf \Sigma}$ and distinct $d_0,\ldots,d_{k-1}\in B^*$, the (finite) type

\[
q = \{\vartheta(x)\} \cup \{[d_{\ell} \dot e x]^{h(d_{\ell})}: 
\ell<k\} \cup \{\neg [d_\ell \dot e \sigma(x)]:\ell<k\}
\]

\mn
is $\Gamma^{\ms}_{a,\bar c}$-small, say as witnessed by $c_{\beta(*)}$
where $\beta (*) < \ell g(\bar c)$, that is $c_{\beta(*)}$ is 
a witness for the conjunction  of $q$ so for some
$\psi(x_0,\ldots, x_{k-1})$ with parameters from $A$ we have 
$M^* \models \psi[d_0,\ldots, d_{k-1}]$ and $\psi$ state the 
$\Gamma^{\ms}_{a,\bar c}$-smallness of $q$ by $c_{\beta(*)}$ i.e.
\mn
\begin{enumerate}
\item[$(*)$]  $\gC \models (\forall y_0,\ldots, y_{k-1})[ \psi
(y_0,\ldots, y_{k-1}) \rightarrow
|\{x \in a:\vartheta (x) \, \& \, (\forall \ell<k)
[y_{\ell} \in x]^{\iif(h(d_{\ell}))}\ \&\ (\forall \ell<k)\neg d_\ell \in 
\sigma(x) \}| / |a|\leq c_{\beta(*)}]$.
\end{enumerate}
\mn
\Wilog

\[ 
\gC \models (\forall y_0,\ldots, y_{k-1})[\psi(y_0,
\ldots,y_{k-1}) \rightarrow \bigwedge_{\ell<m<k} y_\ell \ne
y_m \, \& \, \bigwedge_{\ell< k} y_\ell \dot e a_1].
\]

\mn
Let $\breve{f}^* \in {\gC}$ be a maximal family of pairwise disjoint
$k$-tuples satisfying $\psi$, it exists - see \ref{a2}(B) and
\wilog \, $\breve{f}^* \in A'$; let $\gC \models 
``\dot n := |\breve{f}^*|"$.

Clearly

\begin{equation*}
\begin{array}{clcr}
M^* \models ``&(\forall y_0 \ldots y_{k-1})[\psi(y_0,\ldots,
y_{k-1})\rightarrow (\exists z_0,\ldots,z_{k-1}) \\
   &(\langle z_0,\ldots,z_{k-1}\rangle \in \breve{f}^* \, \& \, \{y_0,\ldots,
y_{k-1}\}\cap\{z_0,\ldots,z_{k-1}\} \ne \emptyset)].
\end{array}
\end{equation*}
\mn

Hence for some $m<k$ and $k(*) < k$ we have $d_{k(*)} \in d^*_m$
where for each $\ell< k$ we let $d^*_\ell := \{d'_\ell:\bar d'
\dot e \breve{f}^*\}$ (in $M^*$'s sense). So $d^*_m \in \acl(A)$ and
$M^* \models ``|d^*_m| = \dot n"$, but as $d_{k(*)} \in d^*_m \, \& \,
d_{k(*)} \in B^*$, by the definition of $B^*$ we have 
$2^{-\dot n} = 2^{-|d^*_m|} < c_{\gamma(*)}$ for some
$\gamma(*) < \lg(\bar c)$. \Wilog \, $\gamma(*) = \beta(*) +1$.

Choose $\alpha(*) < \delta$ such that
$c_{\beta(*)}/c_{\beta(*)+2} < c_{\alpha(*)}$.

Let $a^* := \{x \in a:\vartheta(x)\}$, so as $\vartheta(x)$ is 
$\Gamma^{\ms}_{a,\bar c}$-big clearly $|a^*|/|a| \ge c_{\alpha}$ 
for each $\alpha$.  Also by $(*)$ above and the choice of $\breve{f}$ 
we have $|\{(\bar y,x):\bar y \dot e \breve{f}^*,x \dot e a^*$ and
$(\forall \ell<k)[y_{\ell} \dot e x]^{\iif(h(d_{\ell}))}$ 
and $(\forall \ell<k)[\neg d_\ell \dot e \sigma(x)]\}|$
$\le \dot n \times |a| \times c_{\beta (*)}$;
let us define $a'=:\{x \dot e a:|\{{\bar y} \dot e
\breve{f}^*:(\forall \ell <k)[y_{\ell} \dot e x]^{h(d_{\ell})}$ 
and $(\forall \ell<k)[\neg d_\ell \dot e \sigma(x)]\}|/\dot n 
\ge c_{\beta(*)+2}\}$, by the previous inequality $|a^*\cap a'| \le
(\dot n \times |a| \times c_{\beta(*)})/(\dot n \times 
c_{\beta (*)+2}) = |a|\times(c_{\beta (*)}/c_{\beta(*)+2}) \le 
|a|\times c_{\alpha(*)} < \frac{1}{2}|a^*|$.
Hence $|a\setminus a'|\ge |a^*\setminus a'|\ge (1/2)|a^*|$.

Now we shall show that $a\setminus a'\subseteq b$ where

\[
b = \{x \dot e a: |\{\bar y \dot e \breve{f}^*:
\neg(\forall \ell<k)[y_\ell \dot e x]^{(\iif(h(d_\ell))}\}|/\dot
n< c_{\beta(*)+3}\}.
\]

\mn
This holds as if $x \dot e a\setminus a'$ then 
\mn
\begin{enumerate}
\item[$(i)$]   $|\{\bar y \dot e \breve{f}^*:
(\forall \ell<k)[y_\ell \dot e x]^{\iif(h(d_\ell))}$ and $(\forall
\ell<k)[\neg d_\ell \dot e \sigma(x)\}| /\dot n < c_{\beta(*)+2}$ \\
(as $x \dot e a\setminus a'$)

and
\sn
\item[$(ii)$]   $|\{\bar y \dot e \breve{f}^*:(\exists \ell<k)
[y_\ell \dot e \sigma(x)$]$\}|/\dot n \le |\sigma(x)|/\dot n 
\le (\log_2(1/c_{\beta(*)}))/$
\newline
$\dot n < \log_2
(1/c_{\beta(*)})/(\log_2(1/c_{\gamma(*)}))=
\log_2((1/c_{\beta(*)})/
\log_2 (1/c_{\beta(*)+1}) \le c_{\beta(*)+2}$.
\end{enumerate}
\mn
[Why?  The first inequality as $\breve{f}$ is a set of pairwise 
disjoint $k$ types for the second inequality recall that $2^{-\sigma(x)}\ge
c_\beta$ holds by the choice of $\Sigma$ and for the third, 
$2^{-\dot n} < c_{\gamma(*)}$ was noted above the fourth 
(equally holds as $\gamma(*)= \beta(*)+1$ and for the last 
an assumption of part (2).]
\mn
\begin{enumerate}
\item[$(iii)$]   $2 c_{\beta(*)+2} \le c_{\beta(*)+3}$.
\end{enumerate}
\mn
By the previous paragraph $|a \smallsetminus a'| \ge 
(1/2) |a^*|$ hence $M^* \models ``|b| \ge (1/2)|a^*|$.

Now for the ``random variable" $x \dot e a$ the
events $(\forall \ell<k)[d'_{\ell} \dot e x]^{\iif(h(d'_{\ell}))}$
for $\langle d'_\ell :\ell<k\rangle\in \breve{f}^*$
each has probability $2^{-k}$ and they are independent and their number
is $\dot n$, so the probability that only $\le c_{\beta(*)+3}\times
\dot n$ of them occur for $x$ is sufficiently small by the law of 
large numbers, which mean (see e.g \cite[pg.29]{Sp87}
recall $\bold e$ is the basis of natural
logarithm) that the probability is, for some $n < \omega$ 

\begin{equation*}
\begin{array}{clcr}
\le \bold e^{-2(2^{-k}\dot n-c_{\beta(*)+3}\dot )^2/\dot n} &\le 
\bold e^{[-2^{-(2k+2)}\dot n]} \\
    & < 2^{-2^{-(2k+3)}\dot n} \\
    &\le 2^{-2^{-(2k+3)}|d^*_m|} \\
    & =(2^{-|d^*_n|})^{2^{-(2k+3)}} \\
    & <(c_{\gamma(*)})^{2^{-(2k+3)}}\le c_{\gamma(*)+n}
\end{array}
\end{equation*}

\mn
(last inequality as for every $\gamma <\ell g(\bar c)$ for some 
$n$ we have $c_\gamma \le (c_{\gamma(*)})^{2^{-(2k+3)}}$ hence $c_\gamma \le
c_{\gamma+n} \times c_{\gamma+n}$ so for some $n$ depending on $k$
and $\gamma(*)$, we are O.K.).  So $|b| \le |a|\times \prob (x \notin b)
<|a| \times c_{\gamma(*)+n} < \frac{1}{4} \times
c_{\gamma(*)+n+2}\times |a| \le \frac{1}{2}\times |a^*|$.

Together we get a contradiction.
\end{PROOF}

\noindent
We could have said something on $\breve{f}^*$
\begin{observation}
\label{c56}
Assume $M^* \models T^*,A \subseteq A' \subseteq M \prec {\gC}$,
$\dot n \in A = \acl(A),k<\omega$, for $\ell <k$ we have $d_{\ell}\in A'$,
$a_1\in A$ and $M^*\models$ ``$a_1$ is finite, $d_{\ell} \dot e a_1"$ 
and $M^* \models ``\dot n$ a natural number $>0"$ and for every 
$a_2 \in A$ and $\ell<k$ we have $\gC \models ``d_{\ell} \dot e
a_2\subseteq a_1 \rightarrow |a_2| \ge \dot n"$.  \Then \, 
we can find $\breve{f} \in \gC$ such that for every $\varphi
(\bar x) \in \tp(\langle d_0,\ldots,d_{k-1}\rangle,A,\gC)$ we
have:
\mn
\begin{enumerate}
\item[$(*)$]  $\gC \models ``\breve{f}$ is a set of $\dot n$ 
pairwise disjoint $k$-tuples from $a_1$, each satisfying $\varphi(\bar
x)"$.
\end{enumerate}
\end{observation}

\begin{PROOF}{\ref{c56}}
The properties of $\breve{f}$ can be represented as realizing a $k$-type,
so as $\tp(\langle d_0,\ldots,d_{k-1}\rangle,A,{\gC})$
is closed under finite conjunction, it is enough to find $f \in
{\gC}$ satisfying $(*)$, for one given $\varphi (\bar x) \in 
\tp(\langle d_0,\ldots,d_{k-1}\rangle, A, {\gC})$.

In ${\gC}$ there is a  set $\breve{f}^*$ which is (in
${\gC}$'s sense) a maximal set of pairwise disjoint $k$-types
$\subseteq a_1$ satisfying $\varphi (\bar x)$.
As $T^*$ has Skolem functions \wilog \, $\breve{f}^* \in \acl(A)=A$. If
${\gC} \models |\breve{f}^*| \ge \dot n$ we are done, so assume not; by
$\breve{f}^*$'s maximality ${\gC} \models$ ``some $d_{\ell}$ 
appear in one of the $k$-types from $\breve{f}^*$, say as the
$\ell(*)$-th member of this sequence'', so $d_{\ell}$ satisfies
$d_{\ell} \in \breve{f}^*_{\ell(*)} = \{y:y$ is the $\ell (*)$-th member of
some $k$-type from $\breve{f}^*\}$, but $\breve{f}^*_{\ell(*)}\in 
\acl(A)=A,\models ``|\breve{f}^*_{\ell(*)}| = |\breve{f}^*|< \dot n"$, 
contradicting an assumption.
\end{PROOF}

\begin{claim}
\label{c59}
Suppose $A \subseteq {\gC}$, ${\bar c}^{\ell}\cup\{a^\ell,
\dot w^\ell\}\subseteq A=\acl (A)$,
$p_{\ell}=\tp(b_{\ell}, A, \gC)\in {\bold S}(A,\gC)$ is
$\Gamma^{\wm}_{a^{\ell},\dot w^\ell,{\bar c}^{\ell}}$-big, for 
$\ell=1,2,\dot w^1$ constant, $A^2=\acl(A+b_2),B^* = 
\{d\in A^2:d \in M a^1$ but for every $a_2in A$ such that $d \in a_2$ 
we have $(\exists \alpha)(2^{-|d_2|}\le c_{\alpha}^1)\}$ and
$h:B^* \rightarrow \{\true,\false\}$.  

\Then \, we can find ${b'}_1$ such that:
\mn
\begin{enumerate}
\item[$(\alpha)$]  $b'_1$ realizes $p_1$
\sn
\item[$(\beta)$]  $\acl(A+b'_1) \cap \acl(A+b_2)=\acl (A) (=A)$
\sn
\item[$(\gamma)$]  $\tp(b'_1,A+b_2)$ is 
$\Gamma^{\wm}_{a^1,\dot w^1,\bar c^1}$-big
\sn
\item[$(\delta)$]  $\tp(b_2,A+b'_1)$ is 
$\Gamma^{\wm}_{a^2,\dot w^2{\bar c}^2}$-big
\sn
\item[$(\epsilon)$]  for $d \in B^*,\gC \models [d
  \dot e b'_1]^{h(d)}$
\sn
\item[$(\zeta)$]  every $d \in B^*$ is still ``large" over $\acl(A+b)$ 
(as in the definition of $B^*$).
\end{enumerate}
\end{claim}

\begin{remark}
\label{c62}
Note that the extra constraints in clause $(\zeta)$ are on
$b'_1$ only. If $|\Dom(h)|=1$ simpler better bound suffice,
 otherwise we use $\Delta$-system. This claim is used in fourth 
case Stage D proof of \ref{e8} below. Even if we would have failed, in
\ref{e8}, stage D, we can use the weak diamond argument.
\end{remark}

\begin{PROOF}{\ref{c59}}
We can repeat proof of \ref{c29}(2) + \ref{c53}. By ``the
local character of the demand'' (really a variant of \ref{a37})
we can replace $(\epsilon)$ by
\mn
\begin{enumerate}
\item[$(*)$] $(\forall \ell<k) [d_{\ell} \in  b'_1]^{h(d_{\ell})}$ for
some $k<\omega$ and $d_0,\ldots,d_{k-1}\in B^*$.
\end{enumerate}
\end{PROOF}

\begin{definition}
\label{c65}
For $M^*,a,\bar c$ as in Definition \ref{c23} and $\kappa$
we define $\Gamma=\Gamma^{\ms,\kappa}_{a,\bar c}$ as follows. 
We let ${\bar x}_{\Gamma}=\langle x_i:i < \kappa\rangle$, 
and $p(\bar x_\Gamma)$ is $\Gamma$-big \underline{if and only if} 
for any finite $q\subseteq p$ and  finite
$w \subseteq \kappa$ such that [$x_i$ appears in $q \Rightarrow i\in w$] 
we have

\[
M^* \models ``|\{\langle d_i: i\in w\rangle:M^* \models
\bigwedge q[\langle d_i: i\in w\rangle]\, \& \, \bigwedge\limits_{i\in w} 
d_i \in a\}|/|a|^{|w|} \ge c_\alpha"
\]

\mn
for each $\alpha < \ell g(\bar{c})$.
\end{definition}

\begin{remark}
We can develop this (as in \cite{Sh:107}) and used
it in the stage D of the proof of \ref{e8}.
Also we can define a parallel of $\Gamma_{\ind}$
from Definition \ref{c11} for $\kappa$.
\end{remark}

\begin{claim}
\label{c68} 
Assume $a \in \gC$ is pseudo finite such that $|a|\ge n$ for 
$n<\omega$ and $\gC$ satisfies:
$\bold p \in (0,1)_{\bbR},\bar c = \langle c_i:i < \delta\rangle,
c_i \in (0,1)_{\bbR},\bigwedge\limits_{i<j} 2c_i < c_j$, 
and let $\dot w_{\bold p}:\cP(a) \rightarrow [0,1)$ be as in 
Def. \ref{c32}(2).

1) \Then\ $\cP(a),\dot w_{\bold p},\bar c$ satisfy the requirements in
Definition \ref{c32}(1) on $a$, $\dot w$, $\bar c$.
The bigness notion $\Gamma^{\wmg}_{a,\bold q,\bar c}$ is non-trivial 
\If \, for each $n<\omega$ we have $\bold q^{|a|}<1/n$ and 
$(1-\bold q)^{|a|}<c_\alpha$ for some $\alpha<\lg(\bar c)$, 
equivalently for every $n<\omega,\ln(1/\bold q) > n/|a|$ and 
$\ln (1/(1-\bold q))> \ln(1) c_\alpha)/|a|$.

\noindent
2) If the type $p$ is $\Gamma^\wmg_{a,\bold p,\bar c}$-big, \then \, for no
$i<\lg(\bar c)$ and $b$ do we have, 
$\gC \models ``b \subseteq a,|b| \le \log(1-c_i)/\log (1-\bold
p)"$ and ``$[x \cap b \ne \emptyset]"' \in p$.

\noindent
3) If the type $p$ is $\Gamma^\wmg_{a,\bold p,\bar c}$-big and
$\bigvee\limits_{i<\lg(\bar c)} \bold p < c_i$ \then\, for no 
$d \in \gC$ do we have, ``$[d \in x]" \in p$.
\end{claim}

\begin{PROOF}{\ref{c68}}
1) Check.

\noindent
2) As $|b| \le \frac{\ln(1-c_i)}{\ln(1-\bold p)}$ and as $\ln(1-\bold p)<0$
clearly $|b|\ln(1-\bold p) \ge \ln(1-c_i)$ hence $\ln((1-\bold p)^{|b|}) \ge
\ln(1-c_i)$ hence $(1-\bold p)^{|b|} \ge (1-c_i)$. So
clearly $\dot w_{\bold p}(\{x \in \cP(a):x \cap b = \emptyset\})
= (1-\bold p)^{|b|} \ge 1-c_i$ hence $\dot w_{\bold p}(\{x \in
{\cP}(a):x \cap b \ne \emptyset\}) = 1-\dot w_{\bold p}(\{x \in
{\cP}(a):x \cap b =\emptyset\}) \le c_i$, hence
``$[x \cap b \ne \emptyset]" \in p$ gives easy contradiction.

\noindent
3) Because if $\bold p < c_i$ then $1-c_{i+3}> 1-\bold q$ hence 
${\ln}(1-c_{i+3}) > {\ln}(1-\bold q)$ hence as $\ln(1-\bold q)$ is negative
$|\{d\}| = 1 < \log (1-c_{i+3})/\log (1-\bold p)$ and apply part (2).
\end{PROOF}

\begin{claim}
\label{c71}
1)  Assume $p$, a type over $A \subseteq {\gC}$, 
is $\Gamma^\wmg_{a,\bold p,\bar c}$-big and $a,\bold p$ (hence 
$\dot w_{\bold p}),\bar c$ as in \ref{c32}(2),
$\{a,\bold p\} \cup {\bar c}\subseteq A=\acl(A)$. Suppose $\bar c'$ is
wide for $a$ and $e^* \in a$ and $\tp(e^*,A)$ is 
$\Gamma^\ms_{a,\langle c'_i:i<\lg(\bar c')\rangle}$-big
where $c'_0 \ge \ln (c_0)/(|a|\times \ln (1-\bold p))$ and 
assume $i<\lg(\bar c) \Rightarrow \ln (1/c_i)< c_{i+1}/c_i$. 
\Then \, $p(x \cup\{e^*\}) \cup \{\neg(e^* \in x)\}$
is $\Gamma^\wmg_{a,\bold p,\bar c}$-big.

\noindent
2) If $A=\acl(A),\{a,\bold p,d^*,e^*\}\in A,{\bar c}^1,
{\bar c}^2 \subseteq A,\tp(d^*,A)$ is $\Gamma^\wmg_{a,\bold p,
{\bar c}^1}$-big and $\tp(e^*,A)$ is $\Gamma^\ms_{a,\bar c^2}$-big and
$(*)$ below holds \then \,  we can find $d'$, $e'$ such that:
\mn
\begin{enumerate}
\item[$(\alpha)$]  $\tp(e',A+d')$ is $\Gamma^{\ms}_{a, \bar c^2}$-big extending
$\tp(e^*,A)$
\sn
\item[$(\beta)$]  $e' \in d'$
\sn
\item[$(\gamma)$]  $\tp(d'\setminus\{e'\}, A+e')$ is
$\Gamma^\wmg_{a,\bold p,\bar c^1}$-big
\sn
\item[$(\delta)$]  $\tp(d',A+e')$ nicely extend $\tp(d^*,A)$
\end{enumerate}
\mn
provided that
\mn
\begin{enumerate}
\item[$(*)$]   $(i) \quad \bar c^\ell=\langle c^\ell_i: i<\delta^\ell\rangle$ 
are as in Definition \ref{c32} (in particularly, O.K. and 

\hskip25pt  wide) for ${\cP}(a)$ (if $\ell=1$) or for $a$ 
(if $\ell=2$)\footnote{This is an over-kill but suffice}
\sn
\item[${{}}$]  $(ii) \quad \bold p \in (0,1)_{\bbR}$,
and ``$\bold p < 1/n$'' for each real natural number $n$
\sn
\item[${{}}$] $(iii) \quad \bold e !^{-\frac{1}{8}\bold p^2\times
c^2_{j_2}\times |a|} < \frac{1}{2}c^1_{j_1}$ for every $j_1<\delta_1$,
$j_2<\delta_2$\footnote{if $\delta_1=\delta_2$ then the situation is
simpler; $\bold e !$ is the bases of natural logarithm}
(or at least for 

\hskip25pt  some $(j_1, j_2)$ (by the monotonicity) hence 
we can omit the $\frac{1}{2}$)
\sn
\item[${{}}$]  $(iv) \quad c^2_i \le \bold p \times c^2_{i+1}$.
\end{enumerate}
\end{claim}

\begin{PROOF}{\ref{c71}}
1)  \Wilog \, $p \in {\bold S}(A,\gC)$.  As 
$\Gamma^\wmg_{a,\bold p,\bar c}$ is co-simple, $p$
closed under finite conjunctions, clearly if the conclusion fails
then for some $i < \lg (\bar c)$ and $\varphi (\bar x) \in p$ 
(suppressing parameters from $A$) we have:
\mn
\begin{enumerate}
\item[$(*)$]  $\varphi(x\cup \{e^*\})\cup \{\neg(e^* \in x)\}$ is
$\Gamma^\wmg_{a,\bold p,\bar c}$-small as witnessed by $c_i$.
\end{enumerate}
\mn
Similarly there is $\psi(y)\in \tp(e^*,A)$ (suppressing parameters
from A) such that for every $e \in \psi (\gC)$, we have $e
\in a$ and $\varphi(x \cup \{e\})\cup\{\neg(e \in x)\}$ is 
$\Gamma^\wmg_{a,\bold p,\bar c}$-small as witnessed by $c_i$.

Let $b^1=:\{x \in \cP(a): \varphi (x)\}$ (in $\gC$), and
so as $\varphi(x)\in p$ clearly 
\mn
\begin{enumerate}
\item[$(*)_0$]  $\dot w_{\bold p} (b^1) \ge c_j$ for every $j<\lg (\bar c)$.
\end{enumerate}
\mn
Let $b^2=:\{y \in a: \psi(y)\}$ (in $\gC$), now there is  $b^3$
 such that $\gC \models "b^3 \subseteq b^2\, \& \, 
|b^3|\sim\frac{\ln c_i}{\ln (1-\bold p)}$"
(possible by the assumption on $e^*$ as $c'_j \ge c'_0 \ge \ln(c_0)/(|a| \times
\ln(1-\bold p)) \ge \ln (c_i)/(|a|\times \ln (1-\bold p))$); 
pedantically we should
say $\frac{\ln c_i}{\ln (1-\bold q)}\le |b^3|<\frac{\ln c_1}{\ln 
(1-\bold q)}+1$ and complicate the computations a little (in ($(*)_1$). 
Note (as $\ln(1-\bold p)\sim -\bold p$ and $0 < \bold q < 1/2$), that
$\bigl(\frac{\bold q}{1-\bold q}\bigr)\times \frac{-1}{\ln(1-\bold q)}$ is in
the interval $(\frac{1}{2}, 2)$, and as $\ln(1/c_i)<c_{i+1}/c_i$ 
holds by an assumption we have
\mn
\begin{enumerate}
\item[$(*)_1$]  $\bigl(\frac{\bold q}{1-\bold q}\bigr)|b^3| \le 
\bigl(\frac{\bold q}{1-\bold q}\bigr)\times \frac{\ln c_i}{\ln
  (1-\bold q)}= \bigl(\frac{\bold q}{(1-\bold q)}\bigl)
\times \frac{-1}{\ln(1-\bold q)}\times \ln (1/c_i) \le 
2 \times \ln(1/c_i) \le 2\times (c_{i+1}/c_i)=(2 c_{i+1}/c_i)
\le \frac{c_{i+2}}{c_i}$.
\end{enumerate}
\mn
Now in $\gC$, by the choice of $\dot w_\bold p$ and of $b_3$
we have:
\mn
\begin{enumerate}
\item[$(*)_2$]  $\dot w_{\bold p}(\{x \subseteq a:x\cap b^3=\emptyset\})=
(1-\bold p)^{|b^3|} = \bold e !^{|b_3|\times\ln (1-\bold q)} \le c_i$
\end{enumerate}
\mn
and for every $e \in b^3$ ( as $e \in b^2$ i.e. $\models \psi[e]$ and
the choice of $\psi$)
\mn
\begin{enumerate}
\item[$(*)_3$]  $\dot w_{\bold p}(\{x\setminus\{e\}:x \in b^1,\ e \in
  x\}) \le c_i$
\end{enumerate}
\mn
and by the choice of $\dot w_{\bold p}$ we have
\mn
\begin{enumerate}
\item[$(*)_4$]  $\dot w_{\bold p}(\{x:x \in  b^1$ and $e \in x\})=
\bigl(\frac{\bold p}{1-\bold p}\bigr) 
\dot w_{\bold p}(\{x \setminus\{e\}:x \in 
b^1$ and $e \in x\})$
\end{enumerate}
\mn
hence (use logic, logic, $(*)_4,(*)_2 + (*)_3,(*)_1$,
requirement on $\bar c$ and $(*)_0$ respectively):
\mn
\begin{enumerate}
\item[$(*)_5$]  $\dot w_{\bold p}(b^1) = w_{\bold p}\bigl(\{x:x \in 
b^1$ and $x \cap b^3=\emptyset\}\cup \bigcup\limits_{e \in b^3}
\{x: x \in b^1$ and $e \in  x\}\bigr) \le \dot w_{\bold p}
(\{x:x \in b^1$ and $x \cap b^3 = \emptyset\}) +
\sum\limits_{e \in b^3} \dot w_{\bold p}(\{x: x \in b^1$ and
$e \in x\}) = \dot w_{\bold p}(\{x:x \in b^1$ and  $x \cap b^3=\emptyset\})+
\sum\limits_{e \in b^3} \bigl(\frac{\bold p}{1-\bold p}\bigr)
\dot w_{\bold p}(\{x \setminus \{e\}:x \in b^1$ and $e \in x\}) \le c_i+
|b^3| \times \frac{\bold p}{1-\bold p} \times c_i \le c_i 
+ c_{i+2}<c_{i+3} < \dot w_{\bold p}(b^1)$
\end{enumerate}
\mn
Contradiction.

\noindent
2) Let $\dot w$ be constantly $\frac{1}{|a|}$ on $a$.

Note that we can ignore the ``nicely", i.e. clause $(\delta)$ (Why? As
then we let $\lambda = (|T|+|A|+\aleph_0)^+$, and by
induction on $\zeta< \lambda$ we choose $d_\zeta$ such that:
\mn
\begin{enumerate}
\item[$(a)$]  $p_\zeta = \tp(d_\zeta, A \cup \{d_\xi:\xi <
\zeta\})$ is $\Gamma^\wmg_{a,\bold p,\bar{c}^1}$-big, increasing with
$\zeta$
\sn
\item[$(b)$]  $p_\zeta$ nicely extend $\tp(d^*,A)$
\sn
\item[$(c)$]  $e^* \in d_\zeta$
\sn
\item[$(d)$]  $\tp(d_\zeta \setminus \{e^*\}, A \cup \{d_\xi:
\xi<\zeta\}+ e^*)$ is $\Gamma^{\wmg}_{a,\bold p,\bar c^1}$-big
\sn
\item[$(e)$]   $\tp(e^*,A+\{d_\xi:\xi< \zeta\})$ is
$\Gamma^\ms_{a,\bar{c}^2}$-big.
\end{enumerate}
\mn
If we succeed then for some $\zeta,\tp(d_\zeta, A+e^*)$ nicely
extend $\tp(d^*,A)$ and we are done. For each $\zeta$ we can choose
$p_\zeta$ satisfying (a) + (b) as $\Gamma^\wmg_{a,\bold p,\bar c^1}$
is nice, and then $d_\zeta$ by the claim (without the ``nicely").

Let $p_1 = \tp(d^*,A),p_2 = \tp(e^*,A)$; as the bigness
notions are uniformly $\infty$-simple, if the conclusion fail then by
\ref{a37}(2) we can find $i_\ell < \lg(\bar c_\ell),
\varphi_\ell (x) \in p_\ell$ for $\ell=1,2$ (suppressing
parameters from $A$) such that in $\gC,
\varphi_1(x_1) \rightarrow x_1 \in {\cP}(a),
\varphi_2(x_2) \rightarrow x_2 \in a$ and: 
\mn
\begin{enumerate}
\item[$(*)_0$]  $\varphi_1(x_1)\, \& \, \varphi_2(x_2)\, \& \, x_2 \in x_1
\rightarrow \psi_1(x_1\setminus\{x_2\},x_2) \vee \psi_2(x_2,x_1)$
\end{enumerate}
\mn
where for every $e \in a$ the formula  $\psi_1(x_1 \setminus \{e\},e)$ 
is small for $\Gamma^\wmg_{a,\bold p,\bar c^1}$ as witnessed by $c^1_{i_1}$,
and for every $d \in \bold q(a)$ the formula
$\psi_2(x_2,d)$ is small for $\Gamma^\ms_{a,\bar c^2}$ as witnessed by
$c^2_{i_2}$ (and we suppress parameters from $A$).

Let $j_\ell \ge i_\ell+5$ (and $j_\ell < \delta_\ell$), and the pair
$(j_1,j_2)$ is as in $(*)(iii)$ from \ref{c71}(2) 
and trivially $c^\ell_{i_\ell} \le 
\frac{1}{3 2} \bold q \times c^\ell_{j_1}$. 

Choose non-empty $b_1 \subseteq \{x \in \cP(a):\varphi_1(x)\}$ 
(in ${\gC}$) such that:
\mn
\begin{enumerate}
\item[$(*)_{1,1}$]  $c_{j_1}^1 \le \dot w_{\bold p}(b_1)$, and
$d \in b_1 \Rightarrow \dot w_{\bold p} (b_1 \setminus\{d\}) < c^1_{j_1}$.
\end{enumerate}
\mn
[Why $b_1$ exists?  Choose $b_1 \subseteq \{x \in {\cP}(a):
\varphi_1(x)\}$ such that $c^1_{j_1} \le \dot w_{\bold p}(b_1)$, 
now there is at least one: $\{x \in {\cP}(a):\varphi_1(x)\}$, and 
$\emptyset$ fails this, and $d \in b_1 \Rightarrow |b_1/\{d\}|<|b_1| \,
\& \, \dot w_{\bold q}(b_1/\{d\}) < \dot w_\bold q(b_1)$ hence 
there is such $b_1$ of minimal cardinality and it is as required].

Similarly choose $b_2 \subseteq \{y \in a:\varphi_2(y)\}$ in $\gC$,
such that:
\mn
\begin{enumerate}
\item[$(*)_{1,2}$]  $c^2_{j_2} \le \dot w(b_2)$ and $e \in b_2
\Rightarrow \dot w(b'_2 \setminus\{e\})< c^2_{j_2}$
\end{enumerate}
\mn
So (as $d \in b_1 \Rightarrow \cdot \bold q(\{d\}) < c^1_0,e \in 
b_2 \Rightarrow \dot w(\{e\}) < c^2_0$) easily 
\mn
\begin{enumerate}
\item[$(*)_{1,*}$]  $c^1_{j_1} \le \dot w_{\bold q}(b_1)< 2 c^1_{j_1}$
and $c^2_{j_2} \le \dot w(b_2)< 2 c^2_{j_2}$ and recall 
$|b_2|/|a| = \dot w(b_2)$.
\end{enumerate}
\mn
Now in $\gC$ by the definition of $\dot w_{\bold p}$ we have
\mn
\begin{enumerate}
\item[$(*)_2$]  $\dot w_{\bold p}(\{d \in \cP(a):\dot w(d \cap
b_2)<(\bold p/2)|b_2|/|a|\}) \le \bold e^{-\frac{1}{2}(\bold
  p/2)^2|b_2|} \le {\gt}^{-\frac{1}{8} \bold p^2\times c^2_{j_2}\times |a|}$
\end{enumerate}
\mn
[Why? the second equality (i.e. the left side) is just noting
$c^2_{j_2} \le \dot w(b_2) = \sum\limits_{e \in b_2} \dot w(e)=|b_2|\times
(1/|a|)$.
For the first inequality, $\dot w_{\bold p}$ of the set is 
just the probability of $d$ satisfying $|d \cap b_2|/|a| = 
\dot w(d \cap b_2) \le (\bold p/2)|b_2|/|a|$, where
$d$ is gotten by throwing a coin for each $e \in a$ to decide 
whether $e \in d$ with the probability of yes being $\bold p$; so
the $e \in a \setminus b_2$ are irrelevant and we apply the law of
large numbers see e.g. \cite[p.29]{Sp87}.]

For every $ e \in b_2$ (as $x_\ell \in b_\ell \Rightarrow 
\varphi_\ell(x_\ell)$ and by the choice of $i_1$ and $\psi_1$)
\mn
\begin{enumerate}
\item[$(*)_3$ ] $\dot w_{\bold p}(\{x \setminus\{e\}:x \in b_1,e \in x,\psi_1
(x \setminus\{e\},e)\}) \le c^1_{i_1}$
\end{enumerate}
\mn
and by the choice of $\dot w_{\bold p}$ we have
\mn
\begin{enumerate}
\item[$(*)_4$]  $\dot w_{\bold p}(\{x: x \in b_1$ and $e \in x$ and
$\psi_1 (x \setminus\{e\},e)\}) = 
\frac{\bold p}{1-\bold p} \dot w_{\bold p}d (\{x \setminus\{e\}:x \in 
b_1,\ e \in x,\ \psi_1 (x \setminus\{e\},e)\})$
\end{enumerate}
\mn
By the choice of $\psi_2$, for every $d \in b_1$
\mn
\begin{enumerate}
\item[$(*)_5$]  $\dot w(\{y \in a:\psi_2(y;d)\}) \le c^2_{i_2}$. 
\end{enumerate}
\mn
Let $b^*_1 = \{d \in \cP (a):\dot w(d \cap b_2)=|d \cap b_2|/|a|<
(\bold p/2)\times |b_2|/ |a|\}$.
Note that by $(*)_2$ and $(iii)$ of $(*)$
\mn
\begin{enumerate}
\item[$(*)_6$] $\dot w_{\bold p}(b_1^*) \le \frac{1}{2} c^1_{j_1}$
\end{enumerate}
\mn
Now (in ${\gC}$) on the one hand:
\mn
\begin{enumerate}
\item[$\oplus_1$]  $\Sigma\{\dot w_{\bold p}(d) \times \dot w(e):d
  \in b_1,\ e \in b_2,\ e \in d\} \ge$\\

$\sum\limits_{d \in b_1 \setminus b_1^*}(\Sigma\{\dot w_{\bold p}(d)
\times \dot w(e):e \in b_2, e \in d\})=$\\

$\sum\limits_{d \in b_1 \setminus b^*_1} \dot w_{\bold p}(d)
\times \dot w(d \cap b_2) \ge \sum\limits_{d \in b_1 \setminus b_1^*} 
\dot w_{\bold p}(d) \times ((\bold p/2)|b_2|/|a|) \ge $\\

$(\bold p/2) \times c^2_{j_2} \times \sum\limits_{d \in b_1 \setminus b_1^*}
\dot w_{\bold p}(d) = (\bold p/2) \times c^2_{j_2} \times [\dot
  w_{\bold p}(b_1) - \dot w_{\bold p}(b_1^*)] \ge$\\

$(\bold p/2) \times c^2_{j_2} \times [\dot w_{\bold p}(b_1) -
\frac{1}{2}c^1_{j_1}] \ge$\\

$(\bold p/2) \times c^2_{j_2} \times [c^1_{j_1}-\frac{1}{2}c^1_{j_1}]
\ge$\\

$(\bold p/2) \times c^2_{j_2} \times \frac{1}{2}c^1_{j_1} =
\frac{1}{4} \bold p \times c^1_{j_1}\times c^2_{j_2}$.
\end{enumerate}
\mn
But on the other hand
\mn
\begin{enumerate}
\item[$\oplus_2$]  $\Sigma\{\dot w_{\bold p}(d) \times \dot w(e):d
  \in b_1,\ e \in b_2,\ e \in d\} \le$\\

$\Sigma\{\dot w_{\bold p}(d) \times \dot w(e):d \in b_1,\ e \in b_2,\
e \in d$ and $\psi_1(d \setminus \{e\},e)\}+$\\

$\Sigma\{\dot w_{\bold p}(d) \times \dot w(e):d \in b_1,\ e \in b_2,\
e \in d$ and $\psi_2(e,d)\}$\\

$=\sum\limits_{e \in b_2} \dot w(e) \times 
\dot w_{\bold p}(\{d:d \in b_1,\ e \in d$ and $\psi_1(d
\setminus\{e\},e)\})$\\

$\qquad + \sum\limits_{d \in b_1} \dot w_{\bold p}(d) \times 
\dot w(\{e:e \in b_2,\ e \in d$ and $\psi_2(e,d)\})\le$\\

$\sum\limits_{e \in b_2} \dot w(e) \times \frac{\bold p}{1-\bold p}\times
\dot w_p(\{d \setminus\{e\}:d \in b_1,\ e \in d$ and 
$\psi_1(d\setminus\{e\}, e)\})$\\

$\qquad + \sum\limits_{d \in b_1} \dot w_{\bold p}(d) \times
c^2_{i_2}\le$\\

$\frac{\bold p}{1-\bold p}\sum\limits_{e \in b_2} \dot w(e)\times c^1_{i_1}+
\dot w_{\bold p}(b_1) \times c^2_{i_2} =$\\

$\frac{\bold p}{1-\bold p}\times c^1_{i_1} \times
\sum\limits_{e \in b_2} \dot w(e) + c^2_{i_2} \times \dot w_{\bold p}(b_1)=
\frac{\bold p}{1-\bold p} \times c^1_{i_1} \times \dot w(b_2)+c^2_{i_2}\times
\dot w_{\bold p}(b_1)$\\

$\le c^1_{i_1} \times (\frac{\bold p}{1-\bold p} \times 2 
\times c^2_{j_2})+ c^1_{j_1}\times 2\times c^2_{i_2}$
$< \frac{1}{8}\times \bold p \times c^1_{j_1} \times 2 \times c^2_{j_2}+
\frac{1}{8}\times
c^1_{j_1}\times \bold p \times c^2_{j_2}=\frac{1}{4}\times \bold p \times
c^1_{j_1}\times c^2_{j_2}$.
\end{enumerate}
\mn
Now $(\oplus)_1+(\oplus)_2$ gives contradiction.
\end{PROOF}

\begin{remark}
\label{c74}
If in \ref{c71} we agree to have $\delta^1 = \delta^2$, we
can weaken $(iii)$ of $(*)$ demanding $j_2<j_1$.
\end{remark}

\begin{claim}
\label{c77}
Assume $a \in \gC$ is pseudo finite and infinite, if
$\gC \models ``c \in (0,1)_{\bbR}$ and $\frac{n}{|a|} < c 
< \frac{1}{n}"$ for $n<\omega$ and $\delta,\delta_1,\delta_2$ are 
limit ordinals \then \,

\noindent
1)  If $n \times \ln(|a|)/|a|<c$ for $n<\omega$ \then \, we
can find $\bold p,{\bar c} =\langle c_i: i<\delta\rangle$ and 
$\bar c' = \langle c'_i:i< \delta\rangle$ such that:
\mn
\begin{enumerate}
\item[$(*)$]  $a,\bold p,\bar c,\bar c'$ are as in \ref{c71}(1),
$\lg(\bar c)=\delta= \lg (\bar{c'})$, $\bold p = c=c_0$.
\end{enumerate}
\mn
2)  We can find $\bold p,{\bar c}^1,{\bar c}^2$ as in $(*)$
of \ref{c71}(2), and $\bar c^1$, $\bar c^2$ wide, $\lg(\bar
c^\ell) = \delta_\ell$ and $\bold p,c^1_i,c^2_j \ge c$. 

\noindent
3)  In part (2), moreover if $\delta_1 = \omega=\delta_2$, we can
choose (in ${\gC}$) any $\bold p \in ({}^n\sqrt{c}, 1/n)$, 
for every $n<\omega$, and choose (for each $n$)
$c^2_m=8\times c \times \bold p^{-m-2}$ and $c^1_0=\frac{1}{3}\times
\bold e !^{-c\times p^d\times |a|}$ when ${\gC} \models ``d>n \, \& \, d
\in \bbN",d$ small enough, and $c^1_m=\sqrt[m]{c^1_0}$.
\end{claim}

\begin{PROOF}{\ref{c77}}
By compactness \wilog \, $\delta = \delta_1=\delta_2=\omega$.

\noindent
1) In $\gC$, first choose $\bold p = c$; recall that the function
$x\ln(1/x)$ is strictly increasing for $x \in (0, 1/\bold e !)$ 
(as the derivative is
$-(\ln x)-1$ which is positive), has values in $(0,1/{\bold e})$ and
$\lim\limits_{x\rightarrow 0}(x\ln(1/x))=0$; the same is true for
$\sqrt{x}$ except having values in $(0,\sqrt{1/\bold e !})$. 
Choose by induction on $n$: $c_0=\bold p$,
$c_{n+1}=\sqrt{c_n}$; clearly $(\forall m<\omega) [c_n<1/m]$ and
$2c_n<c_{n+1}$ and even $m<\omega\Rightarrow m \times 
c_n < c_{n+1}$ (by induction on $n$), also $c_0 = \bold p > m/|a|$ for
$m<\omega$. 

Also
\mn
\begin{enumerate} 
\item[$(*)_1$]  $\ln (1/c_n)< c_{n+1}/c_n$. 
\end{enumerate}
\mn
[Why? As this means $c_n \ln(1/c_n)<
c_{n+1}$, but $x\in (0, 10^{-3})\Rightarrow
x\ln(1/x)<\sqrt{x}$ as $\sqrt{x}-x\ln(1/x)$ is increasing in this domain
and $0=\lim\limits_{x\rightarrow 0)} (\sqrt{x}-x\ln (1/x))$, so it holds.]

Let
\mn
\begin{enumerate} 
\item[$(*)_2$]  $c_0'=:\ln(c_0)/(|a|\times \ln(1-\bold p))$ 
\sn
\item[$(*)_3$] $(\forall m<\omega)\bigl[c_0'\ge \frac{m}{|a|}\bigr]$.
\end{enumerate}
\mn
[Why? As $(\forall m<\omega)[\ln (c_0)/ \ln(1- \bold p)>m]$ because $(\ln
c_0)/\ln(1-\bold p) \sim \frac{1}{\bold p} \times \ln(1/c_0) =
\frac{1}{c} {\rm ln} (1/c) > m\times 1=m$].

Lastly 
\mn
\begin{enumerate}
\item[$(*)_4$]  $(\forall n<\omega) [c'_0<1/n]$.
\end{enumerate}
\mn 
[Why?  Clearly $c'_0=\ln (c_0)/(|a|\times \ln(1- \bold p)) \sim 
\frac{1}{\bold p}\times \ln(\frac{1}{\bold p})/|a|$ so the 
requirement mean $\frac{1}{\bold p}\ln \frac{1}{\bold p} <
|a|/m$ for $m<\omega$. Now $c=\bold p$ and by
assumption $n \times \ln (|a|)/ |a| <c$ which mean $\frac{1}{c}< \frac{1}{n}
\times \frac{|a|}{\ln (|a|)}$  hence $\frac{1}{\bold q} \ln
\frac{1}{\bold q} = \frac{1}{c} \ln (\frac{1}{c}) < 
(\frac{1}{n} \times \frac{|a|}{\ln (|a|)}) \times \ln
(\frac{|a|}{\ln (|a|)} \times \frac{1}{n}) < \frac{1}{n} 
\times \frac{|a|}{\ln(|a|)} \times \ln
(|a|)=\frac{|a|}{n}$ as required.]

We have finished as there is no problem to define $c'_n$ for $n\in [1,
\omega)$ by induction on $n$ as before e.g. as $\sqrt[n]{c'_o}$. 

\noindent
2) Note that $\bold e !^{-x\times |a|}$ is decreasing with $x$, and

\[
(\forall x)\bigl(\bigwedge\limits_{n} x>\frac{n}{|a|}\Leftrightarrow
\bigwedge\limits_{n} \bold e !^{-x\times |a|}<\frac{1}{n}\bigr).
\]

\mn
First choose $\bold p \in (0,1)_{\bbR}$ such that

\[
\bigwedge_n [c < \bold p^n \, \& \, \bold p<\frac{1}{n}]
\]

\mn
(so clause $(ii)$ of $(*)$ of \ref{c71}(2) holds).
Second choose $c^2_m=8 \times c\times \bold p^{-m-2}$ (so clause $(iv)$ of 
$(*)$ of \ref{c71}(2)) hold and 
$\bar c^{2}=\langle c^2_n:n<\omega\rangle$
is wide and O.K. which give half of (i) of (*) of
\ref{c71}(2)).

Third choose $c^1_0$ such that

\[
n< \omega \Rightarrow \frac{n}{|\cP(a)|}< c^1_0 <\frac{1}{n}
\]

\mn
and

\[
\bigwedge_n \big[\bold e^{-\frac{1}{8} \bold p^2\times c^2_n\times |a|}<
\frac{1}{2}c^1_0\big]
\]

\mn
equivalently

\[
\bigwedge_n \big[\bold e^{-c\times \bold p^{-n}\times |a|} < 
\frac{1}{2}c^1_0\big].
\]

\mn
For this to be possible we need $\bigwedge\limits_n [\bold e !^{-c\times
\bold p^{-n} \times |a|}< 1/n]$ equivalently 
$\bigwedge\limits_n \big[\bold e !^{-8c\times \bold p^{n}\times |a|} <
\bold e !^{-n}\big]$ hence equivalently
$\bigwedge\limits_n [8c \times \bold p^{-n}> \frac{n}{|a|}]$ i.e.
$\bigwedge\limits_n [\frac{|a|}{n}\times c> \bold p^n.]$, 
now as $\frac{|a|}{m}\times c> 1,1 > \bold p$ this holds.] 
This guarantee clause $(iii)$ of $(*)$ of \ref{c71}(2) 
if we shall have $c^1_0 \le c^1_n$ for $n<\omega$. 

Lastly choose e.g. $c^1_n=\sqrt[n]{c^1_0}$ so clearly ${\bar c}^1
=\langle c^1_n:n<\omega\rangle$ is wide and O.K. which give the second
fall of (i) of (*) of \ref{c71}(2).

\noindent
3) We are left with proving the ``moreover". 

As as $\bigwedge\limits_n [\frac{n}{|a|}< c< \frac{1}{n}]$,
clearly $\bigwedge\limits_n \sqrt[n]{c}<\frac{1}{n}$ so there is
$\bold p$ satisfying $\bigwedge\limits_n \bold p \in [\sqrt[n]{c},
\frac{1}{n}]$, so $\bigwedge\limits_{n} [c < \bold q^n \, \& \, 
p<\frac{1}{n}]$ as required in the proof of part (2).

Now we continue as in part (2).

Now check the requirement in $(*)$ of \ref{c71}(2).
\end{PROOF}

\noindent
The claim we shall mostly use in this context is
\begin{claim}
\label{c80}
Assume that
\mn
\begin{enumerate}
\item[$(a)$]   $\gC$ is as in \ref{a2}(2)
\sn
\item[$(b)$]  $a$ is pseudo finite (in $\gC$ so $|a|>n$
for $n<\omega$ of course)
\sn
\item[$(c)$]   $\delta_1$, $\delta_2$ limit ordinals
\sn
\item[$(d)$]  $\gC \models \frac{n}{|a|} < c < 
\frac{1}{n}$ for every true natural number $n$.
\end{enumerate}
\mn
\Then \,  we can find $\bold p,\bar c^1,\bar c^2$ such that:
\mn
\begin{enumerate}
\item[$(A)$]   the triple $(a,\bold p,\bar c^1)$ and the pair 
$(a,{\bar c}^2)$ are is as in \ref{c32} (so as in $(*)(i)$ 
of \ref{c71}(2))
\sn
\item[$(B)$]  if $\{\bold p,c\},\bar c^1,\bar c^2\subseteq A =
\acl(A)$ and $\tp(d^*,A)$ is $\Gamma^{\wmg}_{a,\bold p,\bar c^1}$-big and
$\tp(e^*,A)$ is $\Gamma^{\ms}_{a,\bar c^2}$-big \then \,  we
can find $d',e'$ ($\in \gC$ of course) such that
\sn
\begin{enumerate}
\item[$(i)$]  $\tp(e',A+d')$ is $\Gamma^{\ms}_{a,\bar c^2}$-big extending
$\tp(e^*, A)$
\sn
\item[$(ii)$]  $e' \in d'$
\sn
\item[$(iii)$]  $\tp(d' \setminus \{e'\},A+e')$ is
  $\Gamma^{\mg}_{a,\bold p,\bar c^1}$-big
\sn
\item[$(iv)$]  $\tp(d',A+ e')$ nicely extend $\tp(d^*,A)$
\end{enumerate}
\mn
\item[$(C)$]  $n/ |a| < \bold p < c^1_0< 1/n$ for $n< \omega$
\sn
\item[$(D)$]  $c^2_i < c$.
\end{enumerate}
\end{claim}

\begin{PROOF}{\ref{c80}}
By compactness \wilog \, $\delta_1 = \delta_2 = \omega$. Work in
$\gC$.  First, we find a non-standard integer $\dot n$ small enough,
i.e.
\mn
\begin{enumerate}
\item[$(*)_1$]  $\gC \vDash ``n < \dot n \, \& \, 2^{\bold n}<
c\times |a|"$ for $n< \omega$.
\end{enumerate}
\mn
We let $\bold p = 1/ \dot n,c^2_i = \bold p^{-i} 2^{\dot n}/ |a|$ and
lastly $c^1_i = (\ln \dot n)^i / \dot n$. Now clauses (A), (C), (D) are
immediate. For clause (B) we have to check the demand $(*)$ in
\ref{c71}(2).
 There clause (i) holds by clause (A), clause (ii) holds by 
clause (C). Clause (iv) holds by the choice of $c^2_i$ and as 
for clause (iii) for $i,j<\omega$ we have

\[
\bold e^{-\frac{1}{8} \bold p^2 \times c^2_j \times |a|} =
\bold e^{-\frac{1}{8} \bold p^{2-j} 2^{\bold n}} \le 
\bold e^{-2^{\bold n/2}} < 1/2 \bold n < \frac{1}{2} c^1_0 
\le \frac{1}{2}c^1_i.
\]
\end{PROOF}

\begin{definition}
\label{c83}
Let $T$ be a complete first order theory, ${\gp} = {\gp}(\bar x)$ a 
type definition (see \ref{x2}(8), say with parameters in $A \subseteq
M^*$, see more in \cite{Sh:384}.  \Then \, we let $\Gamma = 
\Gamma^{\gp}$ be the following bigness notion:
if $A \subseteq M,\bar a \subseteq M$ \then \,:
$\varphi(\bar x,\bar a)$ is $\Gamma$-big iff $\varphi(\bar x,\bar a)
\in {\gp}^M$.
\end{definition}

\begin{claim}
\label{c86}
1) For $T$, ${\gp},M^*$ as in Definition \ref{c83}, $\Gamma^{\gp}$ is
a $\ell$-bigness notion.

\noindent
2)  For $T^*$ as in \ref{a2}(2), if $\Gamma$ is an instance of
$\Gamma^{\ms}$,or $\Gamma^{\wm}$, or any $\Gamma$ such that for some
$a \in A_\Gamma$, a pseudo finite set and every $\Gamma$-big $p\in
{\bold S}(A_\Gamma,M^*)$ we have 
$[x_i \in a]\in p$ (for $i<\lg (\bar x_\Gamma)$),
${\gp}$ is a type definition not increasing ``finite'' sets (for example
$\Gamma^{\gp}_{\uf}$ see \cite[2.11=L2.6]{Sh:384}) \then \, 
$\Gamma \perp \Gamma^{\gp}$.

\noindent
3) Assume $p_2$ has a unique extension in ${\bold S}(A',M^*)$ which 
is necessarily $\Gamma_2$-big so $\Gamma_1 \bot \Gamma_2$ \when \,:
\mn
\begin{enumerate}
\item[$(a)_\ell$]  $\Gamma_1$, $\Gamma_2$ are bigness notions, 
\sn
\item[$(b)_\ell$]  $p_{\ell}(\bar x^{\ell}) \in 
\bold S^{\lg(\bar x^\ell)}(A,M^*)$ is $\Gamma_{\ell}$-big (for $\ell=1,2$), 
\sn
\item[$(c)\hphantom{_\ell}$]  there is $d \in A,M^* \models ``d$
  finite", $[{\bar x}^1\subseteq d]\in p_1$, and
for every $\varphi({\bar x}^1,{\bar x}^2)$ with parameters from $A$
from some $e_\varphi\in A$ we have $(\forall {\bar z}^1)[{\bar z}^1\subseteq d
\rightarrow {\bar z}^1 \in e_\varphi \equiv\varphi ({\bar x}^2,
{\bar z}^1)] \in p_2$ and $A \subseteq A'\subseteq A\cup d(M^*)$
where $d(M^*) = \{c\in M^*:M^* \models ``c \dot e d"\}$.
\end{enumerate}
\end{claim}

\begin{PROOF}{\ref{c86}}
Straightforward.
\end{PROOF}

\begin{definition}
\label{c89}
We say that a bigness notion $\Gamma$ for models of 
$T^*$ is orthogonal to pseudo finite if: for any 
$\Gamma$-big $p(\bar x)$ and pseudo finite $d \in \gC$, there 
is an extension $q(\bar x)$ of $p(\bar x)$ which is 
$\Gamma$-big and satisfies clause (c) of \ref{c86}(2).
\end{definition}

\begin{claim}
\label{c92}
(For $T^*$ as in \ref{a2}(B)). If $\Gamma_1,\Gamma_2$ are 
bigness notion, $\Gamma_1\bot\Gamma_2$ if:
\mn
\begin{enumerate}
\item[$(a)$]  $\Gamma_1$ is orthogonal to pseudo finite
\sn
\item[$(b)$]  if $\Gamma(\bar x)$ is $\Gamma_2$-big in $M^* \supseteq
  A_{\Gamma_2}$ then for some sequence 
$\langle d_i:i < \lg(\bar x)\rangle$ of pseudo finite sets we have 
$P(\bar x_{\Gamma_2}) \cup \{x_i \in d_i:i < \lg(\bar x_\Gamma)\}$ 
is $\Gamma_2$-big (we say: pseudo finitary).
\end{enumerate}
\end{claim}

\noindent
The following generalizes $\Gamma^{\na}$.
\begin{definition}
\label{c95}
Let $T^*$ be as in \ref{a2}(B) and $M^*$
be a model of $T^*$.  Let $a,\dot{\bold d},c_i(i < \delta,\delta$ is 
limit ordinal) be members of $M^*$ such that (in $M^*$):
\mn
\begin{enumerate}
\item[$(a)$]  $a$ is a set
\sn
\item[$(b)$]  $\dot{\bold d}$ is a distance function on $a$ in the
  sense of $M^*$,  i.e. for $b_1,b_2,b_3 \dot e^{M^*} a$ we have
\sn
\begin{enumerate}
\item[$\bullet$]   $\dot{\bold d}(b_1,b_2)$ is a non-negative real 
which is positive \underline{if and only if} $b_1 \ne b_2$,
\sn
\item[$\bullet$]  $\dot{\bold d}(b_1,b_2) = \bold d(b_2,b_1)$,
\sn
\item[$\bullet$]   $\bold d(b_1,b_3) \le \bold d(b_1,b_2) + \bold
  d(b_2,b_3)$
\end{enumerate}
\sn
\item[$(c)$]  $c_i$ a positive real
\sn
\item[$(d)$]  $2c_i \le c_{i+1}$ and $i<j \Rightarrow c_i< c_j$
\sn
\item[$(e)$]  for every $i < \delta$ and $n$ there are
  $b_0,\dotsc,b_{n-1} \dot e a$ such that $\ell < k < n \Rightarrow
\dot{\bold d}(b_\ell,b_k) \ge c_i$.
\end{enumerate}
\mn
Let $\bar c=\langle c_i:i<\delta\rangle$. We define the $\ell$-bigness 
notion $\Gamma^{\mt}_{a,\dot{\bold d},\bar c}$ ($\mt$ for metric) as follows:
$\varphi(x,\bar b)$ is $\Gamma$-big \underline{if and only if} in $M^*$ 
there are $n$ members of $a$ satisfying $\varphi(-,\bar b)$ 
pairwise of distance $\ge c_i$, that is: 

\[
\psi = \psi_{\varphi(\bar x,\bar b),n} := (\exists x_1,\ldots x_n)
[\bigwedge\limits_{\ell=1}^{n} \varphi(x_\ell,\bar b) \wedge
\bigwedge\limits_{\ell<m} \dot{\bold d}(x_\ell,x_m)\ge c_i \wedge
\bigwedge\limits^{n}_{\ell=1} x_\ell \dot e a]
\]

\mn
for every $i<\delta,n<\omega$. Let $\Gamma_\delta^{\mt}$ be the
corresponding bigness scheme for $T^*$. We may omit $a$ if it is defined as
$\{x:(\exists y)[\dot{\bold d}(x,y)$ well defined$\}$; we may write
``$\dis$" instead of $\dot{\bold d}$.
\end{definition}

\begin{claim}
\label{c98}
Assume $T^*,M^*,a,\dot{\bold d},\bar c = \langle c_i:i <
\delta\rangle$ are as in \ref{c95}.

\noindent
1)  $\Gamma^{\mt}_{a,\dot{\bold d},\bar c}$ is an invariant 
$\ell$-bigness notion, co-simple, $\aleph_1$-presentable and 
orthogonal to every invariant $\Gamma$ in particular
to $\Gamma^{\ms},\Gamma^\wm$ (and to $\Gamma^{\na},\Gamma^{\ids}$ of
course).

\noindent
2) Suppose $A\subseteq B\subseteq \gC,A \subseteq A' \subseteq M^*,
M^*$ a $\kappa$-saturated model of $T^*$.  We can find elementary
mappings ${\bold f}_n$ (for $n<\omega$) such that:
\mn
\begin{enumerate}
\item[$(i)$]   $\Dom(\bold f_n)=B,{\bold f}_n \restriction A=\id_A$
\sn
\item[$(ii)$]   if $\Gamma=\Gamma^\mt_{a,\dot{\bold d},\bar c}$ for some $a$, 
$\dot{\bold d},\bar c$ from $A$ and $b \in B,\tp(b,\acl(A),M^*)$ 
is $\Gamma$-big \then \, for any $n<\omega,\tp(\bold f_n(b),A'\cup
\bigcup\limits_{\ell \ne n} A_{\ell},\gC)$ is $\Gamma$-big.
\end{enumerate}
\end{claim}

\begin{remark}
\label{c99}
1) Compare with \ref{e27}-\ref{e44} below.

\noindent
2) In \ref{c98}(2) we can replace $w$ by any $\Gamma$.
\end{remark}

\begin{PROOF}{\ref{c98}}
1) Left to the reader e.g. use \ref{b8}(2) 
(note: if $\tp(d,A)$ is $\Gamma^{\mt}_{a,\dot{\bold d},\bar{c}}$-large
then for every $\lambda$ we can find $\dot{\bold d}_\alpha$ \st \, 
each $\dot{\bold d}_\alpha$ realizes $\tp(a,A),N^*$) 
(for $\alpha<\lambda$) and $N^*$ such that $M^* \prec N^*$ and
$\alpha < \beta < \lambda \, \& \, i < \lg(\bar{c})\Rightarrow N^*\models$
``$\dot{\bold d}(d_{\alpha},d_{\beta})\ge c_i \wedge d_\alpha \dot e a"$.

\noindent
2) For one $\Gamma,b$ this should be clear by the definition of
$\Gamma^{\mt}_{a,\dot{\bold d},\bar c}$. Generally use compactness; in
more detail assume to show that it is enough to prove, for any 
$n<\omega$ that if 
\mn
\begin{enumerate}
\item[$(*)_1$]   for $\ell<m,a_\ell,\dot{\bold d}_\ell,\bar c^\ell$ are 
from $A$ and as in Definition \ref{c95}, $b_\ell \in B$ and 
$\tp(b_\ell,A)$ is $\Gamma^{\mt}_{a_\ell,\dot{\bold d}_\ell,\bar
  c^\ell}$ big  clearly
\sn
\item[$(*)_2$]  it is enough to prove
\sn
\begin{enumerate}
\item[$\bullet$]  for every $m,m(*)$ and find elementary 
mapping ${\bold f}_n$ for $n<\omega$ \st \, $\Dom(\bold f_n)=B,\bold
f_n \rest A = \id_A$ and for each $\ell<m$, and $n_1<n_2<\omega$ 
and $i < \lg(\bar c_\ell)$ we have $\dot{\bold d}_\ell(\bold f_{n_1}
(b_\ell),\bold f_{n_2}(b_2)) \ge c_{\ell,i}$.  
\newline
So assume that
\sn
\item[$\bullet$]  from $(*)_2$ fail, so for some $m,m(*)$ and $\psi$
  we have (supervising parameters from $A$)
\end{enumerate}
\sn
\item[$(*)_3$]  $(a) \quad M^* \models \psi[b_0,\dots,b_{m-1}]$
\sn
\item[${{}}$]  $(b) \quad M^*\models (\exists
  y_{\ell,k},\ldots)_{\ell<m,k<n(*)} (\forall x_0,\ldots,x_{m-1})
[\psi(x_0,\ldots,x_{m-1})\rightarrow$

\hskip25pt $\bigvee\limits_{\ell<m} \bigvee\limits_{k<n(*)} \dot{\bold d}_\ell 
(x_\ell,y_{\ell,k})<c_{\ell,i_\ell}]$.
\end{enumerate}
\mn
Hence (as we are dealing with models of $T^*$)
\mn
\begin{enumerate}
\item[$(*)_4$]   we can find $b_{k,\ell}\in \dcl(A)$ for
  $\ell<m,k<n(*)$ \st \,

\[ 
M^* \models (\forall x_0,\ldots,x_{m-1})[\psi (x_0,\ldots,x_{m-1})
\rightarrow \bigvee\limits_{\ell<m}  \bigvee\limits_{k<n(*)}
\dot{\bold d}_\ell(x_2,y b_{\ell,k}) < c_{\ell,i_\ell}].
\]
\end{enumerate}
\mn
Substituting $b_\ell$ for $x_\ell$ we get that for some 
$\ell<m,k<m(*) M^* \models ``\dot{\bold d}_\ell
(b_i,b_{\ell,k})<c_{\ell,i_\ell}"$, hence clearly $\tp(b_\ell,A,M^*)
\vdash \dot{\bold d}_\ell (x,b_{\ell,k}) < c_{i,i_\ell}$, a 
contradiction to the assumption ``$\tp(b_\ell,A,M^*)$ is 
$\Gamma_{a_\ell,\dot{\bold d}_\ell,\bar c_\ell}$-big". 
\end{PROOF}

\noindent
Trivial but useful in the proof of \ref{e24} is:
\begin{observation}
\label{c102}
Let $M$ be a model of $T^*$ from \ref{a2}(B).
Let $A=\acl(A) \subseteq M$, $M \restriction A \prec M, e\in M,
\{a,\dot{\bold d}\} \cup \Rang(\bar c)\subseteq A$. 
Now $\tp(e, A, M)$ is $\Gamma^{\mt}_{a,\dot{\bold d},\bar c}$-big
\underline{if and only if} for every $e' \in A,i < \lg(\bar c)$ we have $M
\vDash ``e' \in a \rightarrow \bold d(e',e) \ge c_i"$.
\end{observation}

\begin{PROOF}{\ref{c102}}
Easy as in the proof of \ref{c98}(2)
(see more in \ref{e27}-\ref{e51}).
\end{PROOF}

\begin{definition}
\label{c105}
Let $T^*$ be as in \ref{a2}(B) and ${\gC} \models ``a$ an
infinite set".

\noindent
1)  We define a local bigness notion scheme $\Gamma^{\ms}_a$:
\mn
\begin{enumerate}
\item[$\bullet$]  $\varphi(x,\bar b)$ is $\Gamma^{\tms}_a$-big \Iff \, 
$\gC \models ``\{x\subseteq a: \varphi(x,\bar b)\}$ is not a null
subset of $\cP(a)$".
\end{enumerate}
\mn
[Note: $\gC$ ``think" itself a model of set theory,
hence for $A\subseteq \cP(a)$ we can define its outer Lebesgue measure
identifying $\cP(a)$ with ${}^a 2$.]

\noindent
2) We define a local bigness notion scheme $\Gamma^{\inf}_a$ by:
\mn
\begin{enumerate}
\item[$\bullet$]  $\varphi(x,\bar b)$ is $\Gamma^{\inf}_a$-big 
if $\gC \models ``\{x \dot e a:\varphi(x,\bar b)\}$ is infinite".
\end{enumerate}
\end{definition}

\begin{claim}
\label{c108}
Assume $\gC \models ``a,a_1$ are infinite, $a_2$ is finite".

\noindent
1)  $\Gamma^{\tms}_a$ is a simple invariant $\ell$-bigness notion.

\noindent
2) $\Gamma^{\tms}_a$ is orthogonal to $\Gamma^{\tms}_{a_1}$ and
to $\Gamma^{\wm}_{a_2,\dot w,\bar c}$ if $(a_2,\dot w,\bar c)$ are as in
Definition \ref{c32}.

\noindent
3) $\Gamma_a^{\inf}$ is an invariant $\ell$-bigness notion, uniformly
$\aleph_1$-simple (hence co-simple) orthogonal to any invariant local 
and even global bigness notion.
\end{claim}

\begin{PROOF}{\ref{c108}}
1)  Easy.

\noindent
2)  Use Fubini theorem.

\noindent
3)  Easy.
\end{PROOF}
\newpage

\section {General Construction for $T$}

In this section we think on building ${\gB}_\alpha$ a model of $T$ of
cardinality $\lambda$ by induction on
$\alpha <\lambda^+$ representing ${\gB}_\alpha$ as the increasing union 
of $A^\alpha_j$ ($j< \lambda$) and having special 
$\bar a^\alpha_i \subseteq A^\alpha_{i+1}$, we better demand $\lambda$
is regular uncountable.  On the one hand
constructing ${\gB}_{\alpha+1}$ we do it by approximations which
are types over ${\gB}_\alpha$ of cardinality $< \lambda$, 
restricting ourselves to appropriate $\Omega^\alpha$-big types. 
On the other hand for $\beta< \alpha$ we demand that for ``many" 
$i$, $A^\beta_i\subseteq A^\alpha_i$ and $\tp(\bar a^\beta_i,
A^\alpha_i,\gB_\alpha)$ is a $\Gamma^\alpha _i$-big type. 
To be able to carry this we need the orthogonality of the
$\Omega$'s with the $\Gamma$'s. We look at $\langle A^\alpha_j: j<
\lambda\rangle$ as increasing vertically and at $\langle {\gB}_\alpha:
\alpha< \lambda^+\rangle$ as increasing horizontally.

\begin{context}
\label{d2}
1)  $T$ a complete first order theory (usually as in
\ref{a2}(B)), $\gC$ a monster for $T$.

\noindent
2)  $\lambda$ a regular cardinal $\ge |T|$ (and $\chi >\lambda$).

\noindent
3)  ${\bf \Upsilon}_\hor$ ($\hor$ -short for horizontal) is a set of
$\le \lambda^+$ global-bigness notions and schemes of $g$-bigness 
notions for $T$ such that $\Gamma^{\tr} \in {\bf \Upsilon}_\hor$ 
or $\Gamma^\na \in {\bf \Upsilon}_\hor$ and
if such scheme has $\kappa$ parameters then $\lambda^\kappa \le 
\lambda^+$ (or we do not use all instances of the scheme).

\noindent
4)  ${\bf \Upsilon}^\ver$ ($\ver$ - short for vertical) is a set of
$g$-bigness notions and schemes of $g$-bigness notions for $T$ 
such that: if $\Gamma \in {\bf \Upsilon}^{\ver}$ is a scheme with
$\kappa$ parameters then $\lambda^\kappa = \lambda$.

\noindent
5) We assume: if $N \models T$, $\Gamma_1$ is an instance of
${\bf \Upsilon}_\hor$ for $N$, $\Gamma_2$ is an instance of 
${\bf \Upsilon}^\ver$ for $N$ \then \, $\Gamma_1 \perp \Gamma_2$ 
(at least for those actually used), in fact
nicely orthogonal (used only in niceness of (D)(7) below, in the present
context is not an extra assumption by \ref{a53}(3),(4)).

\noindent
6)  For a given model $N$ of $T$, an instance of ${\bf \Upsilon}_\hor$ 
for $N$ mean $\Gamma \in {\bf \Upsilon}_\hor$ or $\Gamma$ a
case of a scheme $\Gamma \in {\bf \Upsilon}_\hor$ with parameters from
$N$, similarly for ${\bf \Upsilon}^{\ver}$.

\noindent
7) Let for ${\bf \Upsilon}$ as above $\Gamma \in c \ell_N({\bf
  \Upsilon})$ means $\Gamma= \langle \Gamma_i:i < \alpha\rangle$ 
for some $\alpha< \lambda$, each $\Gamma_i$ an instance of ${\bf \Upsilon}$ 
for $N$, see \ref{a37}- \ref{a53}.
\end{context}

\begin{discussion}
\label{d5}
If $\lambda$ is a regular uncountable cardinal
$S^{\lambda^+}_{\lambda} := \{\delta<\lambda^+:\cf(\delta)
=\lambda\} \in \check{I}[\lambda^+]$, see \cite[3.4=Lcdl.1]{Sh:E62}
 or at least some stationary $S \subseteq S^{\lambda^+}_{\lambda}$ is
 in $\check{I}[\lambda]$, things are nicer. Assuming there are 
$\lambda^+$ almost disjoint stationary subsets of $\lambda$ 
(a very weak assumption, see \cite{Sh:247} or \cite[4.1=Ld4]{Sh:E62}
and Gitik-Shelah \cite{GiSh:577}), simplifies (can use 
$u^\alpha_i=\{i\}$) but till now was not really necessary. 
Then below $\bar S = \langle S_\alpha:\alpha <
\lambda^+\rangle,S_\alpha \subseteq \lambda$ is stationary, $\beta <
\alpha \Rightarrow |S_\beta \cap S_\alpha| < \lambda$.

We shall describe a construction of a model of $T$ of
cardinality $\lambda^+$ by an increasing continuous sequence $\langle
\gB_{\alpha}:\alpha<\lambda^+\rangle$ of $\lambda^+$ approximations:
models of $T$ of cardinality $\lambda$, and for $\alpha= \beta+1$,
$\gB_\alpha$ is constructed in $\lambda$ steps; in step $i<\lambda$, we have
already constructed a type $p^\alpha_i$ over some $A^\alpha_i
\subseteq \gB_\beta$ of cardinality $<\lambda$ (stipulating 
$\gB_0$ is empty so $A^0_i$ is empty),
and a $g$-bigness notion, $\Omega^\alpha_i$, such that $p^\alpha_i$ is
$\Omega^\alpha_i$-big and $p^\alpha_i,\Omega^\alpha_i$ are increasing with $i$.
We described the construction by assigning some
persona called {\it contractor} to perform various tasks. Each
contractor may play the major role for some $\alpha<\lambda^+$, so it is
called ``the contractor at $\alpha$", but it is also  assigned some $i$'s
for every $\alpha$.  For each $\alpha<\lambda^+$ a contractor
$\zeta_\alpha$ plays the major role, in particular chooses a set of
permissible sequences $\langle\Omega^\alpha_i:i<\lambda\rangle$ (see below) and
possibly a linear ordering $<_i^\alpha$ of $i$ with 
$j<i\Rightarrow <_j^{\alpha}= <_i^{\alpha}\restriction j$ 
(in \S5 we choose such sequences, generally this choice has to be closed
under limits, has no maximal member; the default value is the usual
order). A simple case of the $<^{\lambda}_i$ (the one, which we
already use) is when the contractor $\zeta_{\alpha}$ choose a 
linear order $<_{\alpha}$ of $\lambda$ and
let $<_i^{\alpha}=<_{\alpha}\restriction i$). If not said otherwise we
allow to add instances of $\Gamma^{\na}$.  We demand that all the 
bigness notions are nice: also we can replace
$\langle\Omega^\alpha_j:j<i\rangle$ by one bigness notion $\Omega^\alpha$. We
may use games to describe the constructions as in \cite{Sh:107},
\cite{HLSh:162}, see also  \cite[AP]{Sh:326} or \cite[AP]{Sh:405}.
\end{discussion}

\begin{tcc}
\label{d8}
We have $\bar S = \langle S_i:\ i<\lambda\rangle$, a partition of
$\lambda$ to stationary sets, and $\bar W = \langle
W_\alpha:\alpha<\lambda^+\rangle$ a partition of $\lambda^+$
such that for every regular $\theta \le \lambda$ the set
$\{\delta< \lambda^+:\cf(\delta)=\theta$ and $\delta\in W_\alpha\}$ is
a stationary subset of $\lambda^+$ (actually here we use $W_\alpha$ only for
$\alpha<\lambda$). In applications for each
$\alpha<\lambda^+$ we assign a ``contractor" who can
make sure the model $\gB_{\alpha+1}$ of $T$ which we shall construct
will have some properties.  Let $W^\theta_\alpha = \{\delta \in W_\alpha:
\cf(\delta)=\theta\}$.
\end{tcc}

\noindent
Now we start with the formal description.
\begin{preliminaries}
\label{d11}
We choose by induction on $\alpha<\lambda^+,\bold c_i^\alpha$ for
$i < \lambda$ such that:
\mn
\begin{enumerate}
\item[$(a)$]   $\bar{\bold c}_\alpha = \langle \bold c_i^\alpha:
i<\lambda\rangle$ is an increasing continuous sequence of subsets of
$\alpha$
\sn
\item[$(b)$]  $\alpha=\bigcup\limits_{i<\lambda} \bold c_i^\alpha$
\sn
\item[$(c)$]  $|\bold c_i^\alpha|<\lambda$
\sn
\item[$(d)$]  $\beta\in \bold c_i^\alpha\Rightarrow \bold c_i^\beta= 
\bold c_i^\alpha\cap\beta$
\sn
\item[$(e)$]  $\bold c^0_i=\emptyset,\alpha \in \bold c_0^{\alpha+1}$,
$0 \in \bold c_0^{1+\alpha},[\aleph_0 \le \cf(\alpha)<\lambda \, \& \,
\bold c^\alpha_i \ne \{0\}\Rightarrow \alpha=\sup(\bold c_i^\alpha)]$.
\end{enumerate}
\end{preliminaries}

\begin{tcd}
\label{d14}
We define a game $\Game = \Game_{\bar{\bold c}}$ between the
portagonist and antagonist player, the antagonist choices are divided
to the work of various so called contractors and they are actually
independent sub-players.  All the other choices are of the protagonist;
the protagonist wins a play when always there is a legal move.  The
order of the choices is first by $\alpha < \lambda^+$ and then by
$\varepsilon < \lambda$.  During a play the following are chosen.

For $\alpha<\lambda^+,\gB_\alpha$ and given $\alpha$ and 
${\gB}_\alpha$ toward with choosing ${\gB}_{\alpha+1}$ 
by induction on $\epsilon<\lambda$, ordinal we choose an $i^\alpha_\epsilon =
i_{\alpha,\epsilon}$, set $u^\alpha_{i_{\alpha,\epsilon}}$, and for $j\in
u^\alpha_{i_{\alpha,\epsilon}}$ a type $p^\alpha_j$ and a set 
$A^\alpha_j$, and (for $\epsilon$ and $\alpha$) $E_\alpha\cap 
i^{\alpha}_{\epsilon+1}$ and $\langle <^\alpha_j:j\in 
u^\alpha_{i_{\alpha,\epsilon}}\rangle,
\bar u^\alpha_\epsilon=\langle u_j^\alpha:j\in E_\alpha\cap
i^{\alpha}_{\epsilon+1}\rangle,E^+_\alpha \cap 
i^\alpha_{\epsilon+1},\langle \Gamma^\alpha_i,\bar c_i^\alpha:i 
\in E_\alpha^+\cap i^\alpha_{\epsilon+1}\rangle$ and
$\bar\Omega^\alpha=\langle\Omega^\alpha_j:j\in E^+_\alpha\cap
i^\alpha_{\epsilon+1}\rangle$ (and for some $\alpha$'s also $<_0^\alpha$) 
and $\langle \bar a_{\alpha,i}:i<\lambda\rangle$ such that:
\mn
\begin{enumerate}
\item[$(A)$]  $(a) \quad \gB_\alpha$ is a model of $T$ with universe
$\lambda\times\alpha$ (so we stipulate $\gB_0$ 

\hskip25pt  is an empty model and 
for notational simplicity ignore the case

\hskip25pt  $|{\gB}_{\alpha+1} \setminus {\gB}_\alpha| < \lambda$;
alternatively you may ask that the universe of 

\hskip25pt $\gB_\alpha$ is $\gamma_\alpha<\lambda\times \alpha$ with
no serious changes)
\sn
\item[${{}}$] $(b) \quad \beta < \alpha \Rightarrow \gB_\beta
\prec \gB_\alpha$ (so $\gB_\alpha$ is $\prec$-increasing continuous)
\sn
\item[${{}}$]  $(c) \quad$ essentially ${\gB}_{\alpha+1} \setminus 
{\gB}_\alpha=\cup\{\bar a_{\alpha,i}:i<\lambda\}$, more exactly as 
we would

\hskip25pt  like to allow elements appearing $\bar a_{\alpha,i}$ to be
equal to a member of 

\hskip25pt  ${\gB}_\alpha$ and as we
may like not to use $\Gamma^{\tr}$, we demand just

\hskip25pt  ${\gB}_{\alpha+1} = \acl_{\gB_{\alpha+1}}\bigl(\gB_\alpha\cup
\bigcup\limits_{i<\lambda} \bar a_{\alpha, i}\bigr)$
\sn
\item[${{}}$]  $(d) \quad \langle A^\alpha_j:j<\lambda\rangle$ is
an increasing sequence of subsets of ${\gB}_\alpha$ each of 

\hskip25pt  cardinality
$<\lambda$  and $\bigcup\limits_{j<\lambda} A^\alpha_j = \gB_\alpha$
\sn
\item[$(B)$]  $(a) \quad E_\alpha$ is a club of $\lambda$
\sn
\item[${{}}$]  $(b) \quad$ if $\beta \in \bold c_i^\alpha$ and 
$i \in E_\alpha$ then $E_\alpha \setminus (i+1) \subseteq E_\beta$
\sn
\item[${{}}$] $(c) \quad E_\alpha=\{i^\alpha_\epsilon:
\epsilon<\lambda\},i_\epsilon^\alpha=i_{\alpha, \epsilon}$ is 
increasing continuous with $\epsilon$
\sn
\item[${{}}$]  $(d) \quad 0\in E_\alpha$ (so $i^\alpha_0=0$)
\sn
\item[$(C)$]   $(a) \quad \bar u^\alpha=\langle u_i^\alpha:i\in
E_\alpha\rangle,u_i^\alpha$ a closed set of ordinals $<\lambda$ with a last
element, 

\hskip25pt  $\min(u_i^\alpha)=i$, $\max(u_i^\alpha)<\lambda$
\sn
\item[${{}}$]  $(b) \quad E^+_\alpha = \bigcup\limits_{i\in E_\alpha}
u_i^\alpha$
\sn
\item[${{}}$] $(c) \quad i \in E_\alpha \Rightarrow
\min(E_\alpha\setminus(i+1))>\max (u_i^\alpha)$
\sn
\item[${{}}$]  $(d) \quad i \in E_\alpha,\beta \in \bold c^\alpha_i,
j\in u^\alpha_i$ then $A^\alpha_j\cap {\gB}_\beta=A^\beta_j$
\sn
\item[${{}}$]  $(e)(\alpha) \quad$  if $i \in E_\alpha,\beta \in \bold
  c_i^\alpha$, then $i \in E_\beta$ and $u_i^\beta$ is an initial 
segment of $u_i^\alpha$
\sn
\item[${{}}$]  $\,\,\,\,\,(\beta) \quad$ if $i \in E_\alpha,j \in u_i^\alpha$,
$\beta\in \bold c_j^\alpha,j \notin u_i^\beta$ then $j \in E_\beta$
\sn
\item[${{}}$]  $\,\,\,\,\,(\gamma) \quad <^\alpha_j$ is a linear, 
well ordering of $\{i:i \in j \cap E^+_\alpha\}$ 

\hskip35pt increasing with $j$ for $j \in E^+_\alpha$, i.e. if $j_1 <
j_2$ then $\dom(<^\alpha_{j_1})$

\hskip35pt is an initial segment of $(\dom(<^\alpha_{j_2}),<^\alpha_{j_2})$
\sn
\item[$(D)$]  $(a) \quad {\bar p}^{\alpha}=\langle p_j^\alpha:j\in
E_\alpha^+\rangle$ is an increasing continuous sequence of
types over 

\hskip25pt  ${\gB}_\alpha$ (for $\alpha=0$ this means consistent with $T$)
\sn
\item[${{}}$]  $(b) \quad p_ j^\alpha$ is a type in the variables $\{\bar
x_{\alpha,j_1}:j_1 \in j \cap E^+_\alpha\}$
\sn
\item[${{}}$]  $(c) \quad p_j^\alpha$ is a complete type over $A^\alpha_j$
(in ${\gB}_\alpha$)
\sn
\item[${{}}$]  $(d) \quad \langle\bar a_{\alpha,j_1}:j_1<j\rangle$ realizes in
${\gB}_{\alpha+1}$ the type $p_j^\alpha$
\sn
\item[${{}}$]  $(e) \quad \Omega_j^\alpha$ a bigness notion, an instance
 of ${\bf \Upsilon}_{\hor}$ with parameters from $A^\alpha_j$ or

\hskip25pt  at least $A^\alpha_j\cup \bigcup \{\bar x_{\alpha,j_1}:
j_1<^\alpha_{j+1} j\}$
\sn
\item[${{}}$]  $(f) \quad$ the type $p_j^\alpha$ is big for $\langle
\Omega_{j_1}^\alpha:{j_1}< j\rangle$ by the order $<^\alpha_{j}$
\sn
\item[${{}}$]  $(g) \quad$ for $j_1<j_2$ in $E^+_{\alpha}$ we have :
$p^\alpha_{j_2}\restriction\bigcup\limits_{j<j_1} {\bar x}_{\alpha,
j}$ is a nice extension of $p^{\alpha}_{j_1}$

\hskip25pt  (i.e. in $\gB_{\alpha+1}$,
$A^\alpha_{j_2} \cap \acl \bigl( A^\alpha_{j_1}\cup
\bigcup\limits_{j<j_1} \bar a_{\alpha, j}\bigr) = A^\alpha_{j_1}$); the
niceness 

\hskip25pt  require nice  orthogonality in \ref{d2}(4)$\otimes$,
but in our present proof 

\hskip25pt  this is automatic
\sn
\item[$(E)$]  $(a) \quad \bar c_j^\alpha \subseteq \gB_\alpha$ 
is a sequence of length
$\lg (\bar x_{\Gamma_j^\alpha})$ and $\tp(\bar{c}_j^\alpha, A_j^\alpha,
{\gB}_\alpha)$ is $\Gamma_j^\alpha$-big 

\hskip25pt  (for $j\in E^+_\alpha$)
\sn
\item[${{}}$]  $(b) \quad \Gamma_j^\alpha$ is an instance of
$c \ell_{\gB_\alpha} ({\bf \Upsilon}^{\ver})$ with parameters from $A_j^\alpha$
\sn
\item[${{}}$]  $(c) \quad$ if $i\in E_\alpha,j \in u_i^\alpha,\beta \in
\bold c_i^\alpha$, $j\in u_i^\beta$ \then \,:
\sn
\begin{enumerate}
\item[${{}}$]  $(\alpha) \quad \bar c^{\beta}_j$ is an 
initial segment of $\bar c^\alpha_j$
\sn
\item[${{}}$]  $(\beta) \quad$ if $\bold c^\alpha_i$ has no 
last element then $\bar c^\alpha_j$ is the limit of $\langle\bar c^\gamma_j:
\beta \le \gamma \in \bold c^\alpha_i\rangle$
\sn
\item[${{}}$]  $(\gamma) \quad \Gamma^\beta_j$ is a restriction (to initial
segment of the variables) of $\Gamma^\alpha_j$
\end{enumerate}
\mn
(note: if $j\in u^\alpha_i\setminus\bigcup\{u^\beta_i :\beta\in
\bold c^\alpha_i\}$, the only restriction on $\bar c^\alpha_j$ is: $\in
\gB_\alpha$); the \underline{easy}\ case is $\bar c^\alpha_i =
\bar c^\beta_i$, and this is the one used, in the general case we need
to ensure that the limit in clause $(\beta)$ exists.
\sn
\item[$(F)$]  $(a) \quad \cP_\alpha$ is a family of cardinality $\le \lambda$
of: types over ${\gB}_\alpha$, subsets, relations 

\hskip25pt  on $|{\gB}_{\alpha}|$ 
and partial function from $|{\gB}_\alpha|$ to
$|{\gB}_\alpha|$; and this family 

\hskip25pt  is increasing with $\alpha$, with
reasonable closure conditions.  E.g. choose

\hskip25pt  $\gA_\alpha \prec (\cH(\chi),
\in,<^*_\chi)$, increasing with $\alpha$ of cardinality
$\lambda$ with 

\hskip25pt  $\lambda+1\subseteq {\gA}$ 
such that the construction so far belong to it, 

\hskip25pt   including $\langle {\gA}_\beta: \beta< \alpha\rangle$,
$\langle E_\beta,E^+_\beta$:

\hskip25pt $\beta<\alpha\rangle$, $p^\beta_j,A^\beta_j<^\beta_j, 
u^\beta_j,\Gamma^\beta_j,\bar c^\beta_j,\bar a^\beta_\eta): \beta<\alpha,
j \in E^*_\beta)$

\hskip25pt $\langle i^\beta_\epsilon,u^\beta_{i_{\beta,\epsilon}}:
\beta<\alpha,\epsilon<\lambda)$
\sn
\item[${{}}$]  $(b) \quad$ in ${\gB}_{\alpha+1}$ all types
over $\gB_\alpha$ from $\cP_\alpha$ of cardinality $<\lambda$
are realized 

\hskip25pt in $\gB_{\alpha+1}$ 
\sn
\item[$(G)$]  In stage $\alpha$, the construction is done by induction on
$\epsilon<\lambda$.  First we decide (or are given) what is
$i_\epsilon^\alpha$ and then $\Omega_j^\alpha,p_j^\alpha$ for
$j\in\cup\{u_{i_\epsilon}^\beta:\beta\in \bold c_{i_\epsilon}^\alpha\}$
by induction on $j$, and then we continue adding more elements to what
will be $u_i^\alpha$ and the corresponding $\Gamma^\alpha_j,
\Omega^\alpha_j$ and $p_j^\alpha$ $(j\in u_i^\alpha)$ and lastly
choose $i_{\epsilon+1}^\alpha$, $p^\alpha_{i_{\alpha,\epsilon+1}}$.
The decisions are distributed among the various contractors, for
$j\in E^+_\alpha$ it will be $\zeta^\alpha_j=\zeta_{\alpha, j}$, but
$\zeta^\alpha_{i_{\alpha, \epsilon}}$ will have a say on every $j\in
u^\alpha_\epsilon$ (and is called ``the (major) contractor for 
$(\alpha,\epsilon)$ or the major $\epsilon$-substage in the stage
$\alpha$"), and $\zeta^\alpha_0$ on all $j\in E^+_\alpha$, in
particular $\zeta_{0}^\alpha$ (which is called ``the (major) 
contractor for $\alpha$)", and written also as $\zeta_\alpha$
 decide what will be the family of permissible
$\langle\Omega_j^\alpha:j<i\rangle$ for $i<\lambda$ (usually unique, in
particular whether we have $<_i^\alpha$'s).
Let $E_\alpha^0\subseteq\lambda$ be a thin enough club such that $0\in
E^0_\alpha$. For $\epsilon=0$, $i_{\alpha,\zeta}$ is zero and but 
$\zeta_{i_{\alpha, \epsilon}}^\alpha=\zeta^\alpha_0$ is the
$\zeta$ such that $\alpha\in W_\zeta$.
For $\epsilon$ a limit ordinal let $i=i_\epsilon^\alpha$ be
$\bigcup\limits_{\xi<\epsilon} i^\alpha_\xi$ (necessarily it is in
$E^0_\alpha$ and $[\xi<\epsilon\Rightarrow\max (u_\xi^\alpha)<i]$) and let
$\zeta_{i_{\alpha, \epsilon}}^\alpha$ be the unique $\zeta$
such that $i_\epsilon^\alpha\in S_\zeta$. For $\epsilon$
successor it will be a member of $E^0_\alpha \cap 
\bigcap\limits_{\beta\in \bold c_{\alpha,\epsilon-1}} E_\beta$ 
(which is $>\max (u_{\epsilon-1}^\alpha)$) as decided in stage
$\epsilon-1$.
\sn
\item[$(H)$]   Also in stage $\alpha$, in the induction on $\epsilon$ we
choose ${\cP}_{\alpha,\epsilon}$ as in (F)(1) increasing with
$\epsilon$, ${\cP}_\alpha\subseteq {\cP}_{\alpha,\epsilon}$,
$\bigcup\limits_{\epsilon<\lambda} {\cP}_{\alpha,\epsilon}\subseteq
{\cP}_{\alpha+1}$, $|{\cP}_\alpha| \le \lambda$  and the
construction up to $\epsilon$ belongs to it.
Also (essentially) all types in ${\cP}_{\alpha,\epsilon}$ of
cardinality $<\lambda$ over ${\gB}_\alpha \cup A^\alpha_j$ will 
be realized in ${\gB}_{\alpha+1}$ (as in \ref{d20}(2) below).
\sn
\item[$(I)$]  the division of the decisions:
\sn
\begin{enumerate}
\item[$(a)$]   for any $\alpha< \lambda^+$ the antagonist chooses 
an index family
called ``the horizontal contractor" $\bar \zeta^\hor_\alpha
= \langle \zeta^\hor_{\alpha, x}: x\in X^\hor_\alpha\rangle$ (or we
call $\zeta^\hor_{\alpha, x}$ a case of $\bar \zeta^\hor_\alpha$) 
where $X^\hor_\alpha$ has cardinality $\le \lambda^+$
\sn
\item[$(b)$]  for any $\alpha<\lambda^+,i^\varepsilon < \lambda$ the
  antagonist chooses a non-empty index family called 
``the vertical contractor" 
$\bar \zeta_\alpha^{\ver} = \langle \zeta^\ver_{\alpha,y}:y \in 
Y_{\alpha}^\ver\rangle$ where $Y^\ver_\alpha$ has cardinality $\le
\lambda$.
\sn
\item[$(c)$]  By bookkeeping for each stage $\alpha < \lambda^+$ of
the construction, exactly one of the cases of the vertical
contractors $\zeta^{\hor}_{\beta,x}$are active as the major contractors where
$\beta \le \alpha, x \in X^{\hor}_\alpha$; for each 
$\zeta^\hor_{\gamma, x}$ (so $\gamma< \lambda^+,x \in X^\hor_\alpha
\cap {\cP}_\gamma$) for some $\beta< \lambda^+$, in all 
$\alpha\in W_\beta$ we have $\zeta^\hor_{\gamma,x}$ is ``the
contractor": this contractor chooses the $\Omega$'s
\sn
\item[$(d)$]  By bookkeeping for each $\zeta^{\ver}_{\alpha, y}$ (so $\alpha<
\lambda^+,y \in Y^{\ver}_\alpha$), for every $y \in Y^{\ver}_{\alpha,}$ 
and $\beta\in (\alpha,\lambda^+)$  we have: for a club of 
$\epsilon< \lambda$  for some $j\in u^\beta_\epsilon \setminus
\{\max(u^\beta_\epsilon)\}$ we have: this contractor choose the
$\Gamma^\alpha_j$.
\sn
\item[$(e)$]  some contractor in $Y^{\ver}_{\alpha,\varepsilon}$ 
choose $\Gamma^{\na}$ which belong to ${\bf \Upsilon}_{\ver}$
\sn
\item[$(f)$]  the antagonist chooses also the $\cP_\alpha$ and the
  $\cP_{\alpha,\varepsilon}$.
\end{enumerate}
\end{enumerate}
\end{tcd}

\begin{theorem}
\label{d17}
The protagonist wins the game $\Game_{\bold c}$.
\end{theorem}

\begin{PROOF}{\ref{d17}}
Note:
\mn
\begin{enumerate}
\item[$(\alpha)$]  for $\alpha$'s with no $<_i^\alpha$'s we can carry
the construction by the orthogonality of ${\bf \Upsilon}_\hor,
{\bf \Upsilon}^\ver$
\sn
\item[$(\beta)$]  for $\alpha$'s with the $<_i^\alpha$ not well ordered
(case not used) we need stronger demand on the bigness notion 
$\Omega^\alpha_j$: it is simple.
\sn
\item[$(\gamma)$]  the contractor $\zeta_\epsilon^\alpha$ can choose a
$j_1$ to be in $u_\epsilon^\alpha$ though $j_1$ belongs some $S_\xi$,
$\xi\not=\zeta_\epsilon^\alpha$.
\end{enumerate}
\end{PROOF}

\begin{observation}
\label{d20}
1)  ${\gB}^* = \bigcup\limits_{\alpha<\lambda^+}{\gB}_\alpha$ is 
a model of $T$.

\noindent
2)  $\|{\gB}^*\|=\lambda^+$.
\end{observation}

\begin{PROOF}{\ref{d20}}
1)  As $\gB_\alpha$ is a model of $T$, $\prec$-increasing with
$\alpha$.

\noindent
2) As $\Gamma^\na \in {\bf \Upsilon}_\hor$ and some $\zeta_\alpha$ 
contractor allows it.
\end{PROOF}

\begin{discussion}
\label{d23}
1) We can demand only $\gB_{\alpha+2}$ to be quite
saturated model, while for limit ordinal 
$\delta,\gB_{\delta+1}$ can be any algebraically closed set.

\noindent
2)  We can weaken the demand on ${\bf \Upsilon}_\hor,{\bf \Upsilon}^\ver$ to be
sets of $\lambda$-weak, $\lambda^+$-weak $g$-bigness notion
respectively, see Definition \ref{b54}. 

If we assume that
\mn
\begin{enumerate}
\item[$(*)$]   there is $f^* \in {}^\lambda \lambda$ such that for
$\alpha<\lambda^+$ if $f_\alpha$ is the $\alpha$ the function in
${}^\lambda \lambda$ (e.g. $f_\alpha(i) = \otp(\bold c^\alpha_i)$) then
$f_\alpha<_{\cD_\lambda} f^*$ (where $\cD_\lambda$ is the club 
filter on $\lambda$):
then it is enough to demand also on ${\bf \Upsilon}^\ver$ that it is
$(<\lambda)$-weak $g$-bigness notion.
\end{enumerate}
\end{discussion}
\newpage

\section {Proving the compactness}

Here we prove our main theorem.
For suitable expansion ${\gC}_*$ of ${\cH}(\chi^*,\in)$, and
$\kappa\geq |\tau({\gC})|$, $\lambda = (2^\kappa)^+$, there is a
$\kappa^+$-saturated model ${\gB}^*$ of $\Th(\gC_*)$ of cardinality $\lambda^+$
(or $\lambda^{++}$) in which if $\bold b_1,\bold b_2$ are 
Boolean algebras (or rings) in 
$\gB^*$'s sense, any complete embedding of ${\bold b}_1^{{\gB}^*}$
into ${\bold b}_2^{{\gB}^*}$ is one from ${\gB}^*$; this gives
compactness of appropriate logics.

\begin{question}
\label{e2}
{\rm Characterize the first order theories ${\gt}$ such that for any $T$ (as
in \ref{a2}(A) or as in \ref{a2}(B)) there are
models of $T$ in which each automorphism of any instance of ${\gt}$
(interpreted there) is definable (by a first order formula with
parameters) hence is represented.}
\end{question}

We carry this here for Boolean rings (in \ref{e8}).
\bigskip

\subsection {Explanation of the proof of \ref{e8}:}\
\bigskip

We shall prove here, in particular the compactness of the logic
${\bbL}(\dot{\bold Q}^{\isb})$-where $\dot{\bold Q}^\isb$ is  
quantification over isomorphism of one atomic Boolean ring
onto another atomic Boolean ring (not atomless ones as in
\cite{Sh:384}; [we later in the section deal with any Boolean algebra]. 
In fact we deal with a more general case which says something for any case
of the independence property, but here we try to
explain the proof for this specific case.
Of course we do it in the framework of \S4 showing
that for any ``positive" set of moves, i.e. a strategy for the
antagonist, there is a strategy for them, i.e. for the protagonist
guaranteeing all such isomorphisms are definable (with parameters).
Below we shall survey the proof so we oversimplify in some points.
In particular assume $2^\kappa$ is regular so we could let $\lambda=2^\kappa$.
Let $T$ be a first order complete theory satisfying
$|T|<\lambda=(2^\kappa)^+$.

We build by induction on $\alpha<\lambda^+$ model $\gB_\alpha$ of 
$T$ of cardinality $\lambda$ such that $\gB_\alpha$ is 
$\prec$-increasing, continuous in
$\alpha$ and ${\gB}={\gB}_{\lambda^+}=\cup \{{\gB}_\alpha:
\alpha<\lambda^+\}$ should serve.  We consider
${\bold b}_1, {\bold b}_2$ which are definitions by
first order formulas with parameters of atomic Boolean rings in $\gB$.
For stationarily many $\alpha$, we think there will be an ``undesirable"
isomorphism ${\bold f}$ from one atomic Boolean ring ${\bold b}_1 [{\gB}]$ 
onto the other, ${\bold b}_2 [{\gB}]$  such that 
$(\gB_\alpha,\bold f \rest \gB_\alpha) \prec (\gB,\bold f)$ 
so ${\bold b}_\ell[\gB_\alpha] = {\bold b}^{\gB_\alpha}_\ell$ is the
Boolean ring ${\bold b}_\ell$ as interpreted in $\gB_\alpha$.
We cannot list and treat all such possibilities and we do not know to guess
then (note that G.C.H. may fail here), so we try to add few elements
 such that the restriction of ${\bold f}$ to them will 
suffice to reconstruct ${\bold f} \restriction {\bold b}_1 [\gB_\alpha]$.
A first approximation is to add $ \langle a^\alpha_i:i <
2^\kappa\rangle$ with each $a^\alpha_i \dot e {\bold b}_1$
such that  for no two disjoint equivalently distinct
atoms $b_1,b_2$ of ${\bold b}_1^{\gB_\alpha}$ is
$\{i < 2^\kappa:b_1 \le_{{\bold b}_1} \, a^\alpha_i\} = \{i <2^\kappa:b_2
\le_{{\bold b}_1} \, a^\alpha_i\}$, for this use the bigness notion from
Definition \ref{c65}.
So from ${\bold f} \restriction \{a^\alpha_i:i < 2^\kappa\}$ which has
fewer possibilities  we will be able to reconstruct 
${\bold f} \restriction \gB_\alpha$. So consider $\beta>\alpha$ 
such that $\gB_\beta$ is closed under $\bold f,\bold f^{-1}$.
Now add new $a^{\alpha,\beta}_\xi$
(for $\xi <\kappa$) such that for every distinct $b_1,b_2 \in {\bold b}^\at_1
[\gB_{\beta, j}]$  we have $\{i:b_1 \le_{{\bold b}_1}
a^{\alpha,\beta}_i\} \ne \{i:b_2 \le_{{\bold b}_1} a^{\alpha,\beta}_i\}$, 
such a sequence exists as suitable types appear in the ${\bf \Upsilon}_\hor$.
So if we list the possible ${\bold f} \restriction
\{a^{\alpha,\beta}_i:i < \kappa\}$, together with $\alpha,\beta$
we have essentially listed the possible ${\bold f} \restriction 
\bold b_1[\gB_\alpha]$!  In limit $\delta < \lambda^+$ such that 
$\lambda^{\cf(\delta)} = \lambda$, we can list $\le \lambda$ candidates to 
${\bold f} \restriction {\bold b}_1 [\gB_\delta]$
(those such that for unboundedly many $\alpha < \delta,\bold f \restriction
\gB_\alpha$ was listed before $\delta$) so ${\bold f} \rest \bold b_1 
[\gB_\delta]$ is listed in $\delta$, for a club of such $\delta$'s.

The next stage - assume for simplicity that on $\delta < \lambda^+$ we guess
${\bold b}_1$, ${\bold b}_2$, ${\bold f} \restriction {\bold b}_1 [\gB_1]$.
As in older proofs we add $x_\delta \in {\bold b}_1[{\gB}_1]$ 
and try to omit the type $\{y \in \gB_\delta \, \& \,  
[{\bold f}(a) \le_{{\bold b}_2} y]^{\iif[a \le_{{\bold b}_1}x_\delta]}: 
a \in {\bold b}_1^\at [\gB_\delta]\}$.
It is hard to omit types without $\diamondsuit_\mu$, so we use the special
types (as in \cite{Sh:72}) by preserving bigness for vertical bigness notions
as explained below. For $x_\delta$, we define by induction on 
$j < \lambda$, $p_j \in {\bold S}^1 (\gB_{\delta, j})$ increasing 
continuous in $j<\lambda$ such that ``$x \dot e {\bold b}_1" \in p_j$, 
and $p_j$ is big in the sense that for any pairwise distinct $a_0,
a_1,\dots, a_n \in {\bold b}^\at_1 [\gB_\delta] \setminus 
\gB_{\delta,j}$, the type $p_j \cup 
\{[a_\ell \le x_\delta]^{\iif[\ell\,\even]}:\ \ell \le n \}$ is consistent
 in ${\gB}_\delta$, this is a case of $\Gamma^{\ind}$.
For each $\zeta < \lambda$ we would like to choose $a^\delta_\zeta 
\in {\bold b}^\at_1 [\gB_\delta]$ and an infinite 
$\dot{\bold J}^\delta_{j} \subseteq \gB_\delta$ to which
${\bold f}(a^\delta_j)$ belongs, indiscernible over $\gB_{\delta,j} 
\cup \{a^\delta_j\}$, and to make it indiscernible over
$\gB_{\delta,j} \cup \{a^\delta_j\} \cup \{x_\delta\}$.
More generally, we would like to promise that for $\beta > \delta,\{j <
\lambda:\dot{\bold J}^\delta_j$ is an indiscernible sequence over 
$\gB_{\beta,j} \cup \{a^\delta_j\}\ ({\rm in}\ \gB_\delta) \}$ is stationary.
To preserve this in limit we promise that this occurs for $j_\epsilon\in
u_{\delta,\epsilon}$ for ``almost" all $\epsilon \in S$ where
$S \subseteq \lambda$ is stationary, almost means except a non-stationary set
(so it is clear that having $\lambda$ almost disjoint stationary subsets of
$\lambda$ is helpful though not actually used).

However we have outsmarted ourselves: if we add
$\langle a^\alpha_i:i < 2^\kappa\rangle$ as above,
this does not let us fulfill the obligation we have intended to add in order 
to omit the type - omitting types by the indiscernibility is a strong 
commitment. 
There are various directions to try to solve the dilemma, our choice is to
 weaken the demand on $\langle a^\alpha_i:\ i < 2^\kappa \rangle $ - we
demand just that:
\mn
\begin{enumerate}
\item[$(*)$]    letting ${\cE}_{\alpha}$ be the following equivalence relation
on ${\bold b}_1^{\at} [\gB_\alpha]$
\[ 
b_1 {\bold E}_{\alpha} b_2 =: \bigwedge\limits_i [b_1 \le a^\alpha_i
\equiv b_2 \le a^\alpha_i]
\]
\end{enumerate}
\mn
we demand:
\mn
\begin{enumerate}
\item[$(*)$]  for an unbounded set ${\cU}_\delta \subseteq \lambda$,
for every $j \in {\cU}_\delta$ and $b \in {\bold b}_1^{\at}$  
$[\gB_{\alpha,j +1}] \setminus {\bold b}^\at_1 [\gB_{\alpha, j}]$,
the set $b/\cE_\alpha$ is a singleton.
\end{enumerate}
\mn
So from ${\bold f} \restriction \{a^\alpha_i:\ i < 2^\kappa\}$ 
we can reconstruct ${\bold f} \restriction 
\bigcup\limits_{j \in \cU_\delta} ({\bold b}^\at_1[\gB_{\alpha,
j +1}] \setminus {\bold b}^\at_1 [\gB_{\alpha,j}])$,
more exactly ${\bold f} \restriction \{b \in {\bold b}^\at_1 (\gB_\alpha):
b/\cE_{\alpha}$ a singleton$\}$ which includes the mapping above.  Now it is
natural to demand on $\gB_\alpha$ that every definable (with parameters)
infinite set has cardinality $\lambda$, and when it is a subset of
$\bold b^{\at}_1[\gB_\alpha]$ then it
has members in ${\bold b}^\at_1[\gB_{\alpha,j +1}] \setminus
{\bold b}^\at_1[\gB_{\alpha,j}]$ for every $j < \lambda$ large enough.
So if $c_1$, $c_2 \in {\bold b}_1 [\gB_\alpha]$ and  the symmetric difference
is not (really) finite union of atoms, then we can distinguish between
${\bold f}(c_1)$, ${\bold f}(c_2)$. From the definition of 
${\bold f} \restriction \{b \in {\bold b}_1^\at [\gB_\alpha):
b /\cE_{\alpha}$\ a singleton $\}$ we can reconstruct the isomorphism
${\bold f}$ induce on ${\bold b}_1 [\gB_\alpha]/$ (truly finite union of atoms)
onto ${\bold b}_2[\gB_\alpha]/$ (truly finite union of atoms) though
if we like to assume just ``${\bold f}$ is a complete embedding" 
we have to use a larger ideal.
For our purpose it is enough to show that ${\bold f} \restriction \gB_\alpha$
can be reconstructed up to having $\le \lambda$ possibilities.
So assume that ${\bold f}_1,{\bold f}_2$ are isomorphisms 
from ${\bold b}_1 [\gB_\alpha]$ onto ${\bold b}_2[\gB_\alpha]$ 
inducing the same isomorphism above and let $A=:\{b \in {\bold b}^\at_1
[\gB_\alpha]:\ {\bold f}_1(b) \ne {\bold f}_2(b)\}$. 
In the case $A$ is infinite, and
$\gB_\alpha$ is $\aleph_1$-saturated we get contradiction, how? there is
$A' \subseteq A$ infinite such that for $b_1 \ne b_2$ in $A',\langle
{\bold f}_1(b_1),{\bold f}_2(b_1),{\bold f}_1(b_2),{\bold f}_2(b_2)
\rangle$ is with no repetition, so there is $a \in 
{\bold b}_1[{\gB}_\alpha]$ such that
$b \in A' \Rightarrow b \le a\, \& \, {\bold f}^{-1}_2 {\bold f}_1(b) 
\cap a = 0_{{\bold b}_1}\, \& \, {\bold f}_1 {\bold f}^{-1}_2 (b) \cap a = 
0_{{\bold b}_1}$, and ${\bold f}_1(a) \triangle {\bold f}_2(a)$ 
is not finite union of atoms).
We assume $\lambda^{\aleph_0} = \lambda$ and $\cf(\delta)> \aleph_0$
- the latter can be waived.
All this is not the end - we have just succeed to have  for
stationarily many $\delta < \lambda^+$ such that 
$\lambda^{\cf(\delta)} = \lambda,\cf(\delta) > \aleph_0$ and 
${\bold f} \restriction {\bold b}_1 [\gB_\delta ]$ being among our
guesses and for some $j < \lambda$ for every $a \in {\bold b}^{\at}_1
[\gB_\delta] \setminus {\bold b}_1^\at[{\gB}_{\delta,j}]$ there is no
infinite set indiscernible over $\gB_{\delta,j} \cup \{a\}$ to which $h(a)$
belongs.
We would like to deduce $h(a) \in \acl_{\gB_{\delta}} (\gB_{\delta,j}
\cup \{a\})$. If this does not occur we cannot immediately
add $x_\delta$ and promise to omit a type as above (for possible
lack of $\dot{\bold J}^\delta_\zeta$'s) but we can add such
indiscernibles and then have $x_\delta$ (can do it all in $\gB_{\delta + 1})$.

Above we were obscure on which bigness notions we use. Actually these
come from ``random enough sets" (like  $\Gamma^{\ms}_{a,\bar c}$).

So we have accomplished two things. First, for every $\alpha<
\lambda^+$, ${\bold f}\restriction {\bold b}_1 [\gB_\alpha]$ appear in
$\cP_{\lambda^+}$ hence for a club of $\delta< \lambda^+$ of
appropriate cofinality $\theta$ ($<\lambda$), ${\bold f} \restriction
{\bold b}_1[\gB_\delta]$ appear in $\cP_\delta$ (i.e. immediately). Second,
for a club of $\delta < \lambda^+$ of cofinality $\lambda$, in $\cP_\delta$ we
have ${\bf f}\restriction B^{3,\bold f}_{\alpha,\delta}$, where $B^{3,
{\bold f}}_{\alpha, \delta}$ is a suitable ``large" subset of 
${\bold b}^{\at}_1[\gB_\delta]$ (i.e. have a member below $d \in
{\bold b}_1[\gB_\delta]$ if $[x \dot e d]$ is $\Gamma^{\ms}_{\bar
  c}$-big). So then the pre-killer contractor ``acts". 
He tries to promise that for stationary many $\epsilon < \lambda$ 
for some $j \in u^\delta_{\dot e} \setminus \{\max(u^\delta_{\dot e})\}$ the
sequence $\langle (d, d_n): n<\omega\rangle$ is indiscernible 
over $A^\delta_j$ and $\Gamma^\delta_j = \Gamma^{\ids}$ 
(so this will be preserved) where $d \in {\bold
  b}_1^{\at}[\gB_\delta],d_n \in {\bold b}_2[\gB_\delta]$ and
${\bold f}(d) = d_0$ (so $d \in B^{3,\bold f}_{\alpha,\delta}$) but
$d_0 \notin \acl_{\gB_\delta}(A^\delta_j +d)$.

In later stage $\beta\in \cC^3_{\bold f},\cf(\beta)= \theta$, we will know
${\bold f} \restriction {\bold b}_1[\gB_\beta]$ (i.e. it belongs 
to $\cP_\beta$), and so we ``promise" that for stationary 
many $\epsilon < \lambda$ for some $j \in u^\beta_\epsilon$ the 
sequence $\langle (d^{\beta,j},d^{\beta,j}_n:n < \omega\rangle 
= \langle(d^{\delta,j}, d^{\delta,j}_n): n<\omega\rangle$ is 
indiscernible over $A^\beta_j$ and
$e^\delta_j \in {\bold b}^{\at}_1[\gB_\beta]$ satisfies 
${\bold f}(e^\delta_j) \cap d^{\delta,j}_1 > 0_{\bold b_2}$ and 
$a^\delta_\omega \in {\bold b}_1[\gB_{\beta+1}]$ satisfies 
$e^\beta_j\cap a^\delta_\omega = 0_{\bold b_1},d^{\beta,j} 
\le_{\bold b_1} a^\delta_\omega$ (this is the old way to kill).
So we get that for $b \in B^{3,\bold f}_{\alpha,\lambda^+},{\bold f}(b)\in
\acl_{\gB_{\lambda^+}} (A+b)$ where $A \in \cP_{\lambda^+} \cap 
[\gB_{\lambda^+}]^{<\lambda}$, so ${\bold f}(b) \in \{f_i(b):i<
i(*)\}$, where

\[
\gB_{\lambda^+} \models ``\breve{f}_i \text{ is a partial function
  from } \bold b^{\at}_1 \text{ to } \bold b_2",
\]

\mn
so let ${\bold f}(b) = f_{i(b)}(b)$.

Now as $B^{3,\bold f}_{\alpha,\lambda^+}$ is large enough, we can show
that \wilog \, $i(b)$ depends just on $\tp(b,A,\gB_{\lambda^+})$. 
Then we show that $\breve{f}_i$ (\wilog \,) satisfies

\[
(\forall x, y \in \Dom(\breve{f}_i)) [x \ne y \rightarrow
\breve{f}_i(x) \cap f_i(y)= 0_{\bold b_2}].
\]

\mn
Next, we have one such $\dot f$ which will define ${\bold f} \restriction
({\bold b}_1 \restriction(-d))$ where the formula $[x \dot e d]$ is
$\Gamma^{\ms}_{a,\bar c}$-small (note: as ${\bold b}_1$ is a
``finite", it is a Boolean algebra, so $-d$ is legal).

Now for $\delta$ in which the relevant  contractor
works, let $\bar c=\langle a^\delta_n:n < \omega\rangle$, we get $f$
as above.

So this translates into: we have a tree whose levels are non-standard
integer ${\bold n}^*\in \gB_{\lambda^+}$, with inverse order and we have
to show that also such trees have no undefinable branch.
This needs: replacing nice by strictly nice (see \ref{e24}).  
We succeed to deal with this thus at last we finish the proof.

\begin{ml}
\label{e8}
Let $T=T^*$ be as in \ref{a2}(B), $\lambda \ge |T|$ regular,
$\kappa > \aleph_0$ regular $\lambda^+ \ge \lambda^{<\kappa}$ and:
 $\lambda^{\aleph_0}=\lambda,\theta=\aleph_0$.
\Then \, in the framework of \ref{d2} we
can get that the model $\gB = \gB^* = \gB_{\lambda^+}=
\bigcup\limits_{\alpha<\lambda^+} \gB_\alpha$ satisfies:

\noindent
1)  $\gB^*$ is a model of $T$ of cardinality $\lambda^+$.

\noindent
2) ${\gB}$ is $\kappa$-compact that is every type over $\gB^*$ of
cardinality $< \kappa$ is realized, even by $\lambda^+$ elements.

\noindent
3)  There are pseudo finite sets $a_\alpha \in {\gB}$ (for
$\alpha<\lambda^+$) increasing by $\subseteq^{\gB}$ such that for 
every pseudo finite $b \in {\gB}$, for every large enough $\alpha,
\gB \models$ ``$b \subseteq a_\alpha$".
Also, if $Y \subseteq {\gB}$ and each $Y \cap a_\alpha$ is
represented in ${\gB}^*$ (i.e. for some $b \in \gB$ we have, 
for every $c \in \gB$: 

\[
\gB \models ``c \dot e b" \Leftrightarrow c \in Y \, \& \,
\gB \models ``c \dot e a_\alpha")
\]

\mn
\then \, $Y$ is definable (with parameters, by a first order
formula) in ${\gB}$.

\noindent
4) If $(*)$ below holds and ${\bold b}_1, {\bold b}_2$ are
atomic Boolean rings in ${\gB}$ (so their set of members is a ``set" 
of ${\gB}$ not just a definable subset (with parameters))
\then \, every isomorphism from ${\bold b}_1$ onto ${\bold b}_2$ 
is represented in ${\gB}$; where 
\mn
\begin{enumerate}
\item[$(*)$]  for some $n(*)<\omega,\kappa_{n(*)} < \kappa_{n(*)-1} <
\ldots < \kappa_0=\lambda$ we have $(2^{\kappa_{\ell+1}}) \ge \kappa_\ell$ and 
$\cf\bigl([\kappa_\ell]^{\le \kappa_{\ell+1}},\subseteq\bigr) \le \lambda^+$
and $\lambda^{\kappa_{n(*)}} \le \lambda^+$.  Let 
$\kappa^*_\ell$ be $\kappa^+_\ell$ if $\ell \in
\{1,\ldots,n(*)\}$ and $\lambda$ if $\ell=0$ 
\newline
for example
\sn
\item[$(**)$]   $2^\kappa$ regular, 
$\lambda=(2^\kappa)^+,\kappa_0=\lambda,\kappa_1=2^\kappa,
\kappa_2=\kappa$ \underline{or} $\lambda = (2^\kappa)^{++},\kappa_0 =
\lambda,\kappa_1 = (2^\kappa)^+,\kappa_2 = 2^\kappa,\kappa_3 = \kappa$.
\end{enumerate}
\end{ml}

\begin{remark}
\label{e11}
1) From $(*)$ of (4), the demand $\lambda^{\kappa_{n(*)}} \le
\lambda^+$ is used only in proving $\otimes_4$ during Stage E. 
We can weaken it to:
\mn
\begin{enumerate}
\item[$(*)$]   We can find $Y_i\subseteq \kappa_{n(*)}$ for $i<\theta, 
\dot{\bold I}$ an ideal on $\theta$ such that: $2^\theta \le
\lambda^+,T_{\dot{\bold I}} (\lambda)\le \lambda^+$ (see
  \cite[3.7=Lc18]{Sh:E62}) and $(\forall Z\in \cP(\theta) \setminus 
\dot{\bold I})(\exists^{<\sigma}\alpha)(\exists \beta)(\alpha \ne 
\beta<\theta\, \& \, \bigwedge\limits_{i\in Z}[\alpha \in Y_i \equiv 
\beta\in Y_i])$ and $\lambda^{<\sigma} \le \lambda^+$ 
(mostly it suffice $\theta=\cf(\theta), \lambda^{\langle\theta\rangle}
=\lambda,\lambda^\theta \le \lambda^+$,
which means every tree with $\le \lambda$ nodes has $\le \lambda$
$\theta$-branches (no much harm done if we demand
$\lambda=\lambda^{<\kappa}\, \& \, \kappa >|T|$)).
\end{enumerate}
\mn
2)  Recall from \cite[0.12=L2.8A]{Sh:384} that ${\bold f}$ is a 
complete embedding of the Boolean ring $\bold B_1$ into the Boolean 
ring $\bold B_2$ \If \, it is an embedding and maps every maximal 
antichain of $\bold B_1$ to a maximal antichain of $\bold B_2$; 
equivalently if $a \in \bold B_2 \smallsetminus \{0_{\bold B_2}\}$
then there is $b_1 \in \bold B_1,b \ne 0_{\bold B_1}$ such that 
$\bold B_1 \models 0 <_{\bold B_1} c \le b^1 \Rightarrow \bold B_2 
\models a \cap {\bold f}(c) \ne 0^1$.

\noindent
3)  Recall a Boolean ring is like an ideal of a Boolean algebra.
\end{remark}

\begin{PROOF}{\ref{e8}}
\underline{Stage A}:   We use \ref{d2}-\ref{d14} 
(almost as in \ref{d17})   for $T^*$ and
${\bf \Upsilon}_{\hor}$, ${\bf \Upsilon}^{\ver}$ as in 
\ref{d2} such that:
\mn
\begin{enumerate}
\item[$(*)_0$]  $\Gamma^{\na},\Gamma^{\mt},\Gamma^{\wm}$ (wide cases only!),
$\Gamma_{<\kappa}^{\av},\Gamma^{\gp_{\uf}}$
(see \ref{c83} or see
\cite[2.11=L2.6]{Sh:384}) are in ${\bf \Upsilon}_{\hor}$ and
\sn
\item[$(*)_1$]  $\Gamma^{\na},\Gamma^{\ids}_{\omega+1},\Gamma^{\wm}$ are
in ${\bf \Upsilon}^{\ver}$. 

Also 
\sn
\item[$(*)_2$]  ${\cP}_\alpha,{\cP}_{\alpha,i}$ (${\cP}_\alpha$ 
increasing in $\alpha$), $|{A}_{\alpha,i}| \le \lambda$ (see \ref{d14}(F)),
${\cP}_{\alpha+1} \supseteq \bigcup\limits_{i<\lambda}
{\cP}_{\alpha,i},\cP_\alpha \subseteq \cP_{\alpha,0}$ (of course 
$\cP_{\alpha,i}$ increasing in $i$); they will be, for some 
$\chi$ large enough, the set of objects definable in 
$({\cH}(\chi),\in,<^*_\chi)$ from the construction
up to this point and finitely many members of
$\alpha \cup |\gB_\alpha| \cup (\lambda+1)$, (for any $\alpha$ and $i$ this
include $p^\alpha_i$).
\end{enumerate}
\mn
For this to make sense we have to check the orthogonality condition which
holds: we check each one in ${\bf \Upsilon}^{\ver}$: for $\Gamma^\na$ by
\ref{b8}(2), for $\Gamma^{\ids}_{\omega +1}$
by \ref{b14}(2); for $\Gamma^{\wm}$,
orthogonality:  to $\Gamma^\na$ by \ref{b8}(2), to
$\Gamma^{\mt}$ by \ref{c98}, to
$\Gamma^{\wm}$ by \ref{c38}(2), to $\Gamma_{<\kappa}^{\av}$ by
\ref{c38}(3) and to $\Gamma^{\gp_{\uf}}$ by
\ref{c86}; may compare with \cite[2.17=L2.8A]{Sh:384}.

We fix a winning strategy for the protagonist and then we decide
various things for the antagonist in the form of decisions for various
contractors.  Each such commitment implies that $\gB$, (the outcome of
a play under the restrictions above) satisfies more.

This essentially fits in \ref{d17} but
in some cases things are more complicated. One contractor, $\zeta_{\sat}$ (the
saturator) acts in every stage $\alpha< \lambda^+$ but he uses cases of
$\Gamma^{\na}$ only (actually we can let it act in stage $\alpha$ 
only for $\alpha$ successor of successor ordinals).
Another contractor $\zeta_{\br}$ (the branch killer) do not
need to add $\Omega^\alpha_i$'s but acts for every 
$\alpha< \lambda^+$ guaranteeing amalgamation of certain kind exist.

The third real deviation from \ref{d17} is the coder (see Stage
D). We should be careful and show that the demands can be fulfilled.

The fourth deviation is that in the end we have to use a refinement of
``niceness", from \ref{e24}-\ref{e57}.
\medskip

\noindent
\underline{Stage B}:  We assign a contractor called the
\underline{saturator} $\zeta_{\sat}$, for $\alpha\in W_{\zeta_{\sat}}
\subseteq \lambda^+$, he chooses $\Omega_0^\alpha = \Gamma^{\na}$, 
and a type $p(\bar x)$ to be $\subseteq p_0^\alpha$ if possible 
such that every member of ${\cP}_\alpha$ which is a type over 
$\gB_\alpha$ of cardinality $<\lambda$ is
eventually chosen. However here we also need for every $\alpha$ for
every $i<\lambda$, every type over $A_i^{\alpha+1}$ from
${\cP}_{\alpha,i}$ of cardinality $<\lambda$ which is consistent with
${\gB}_{\alpha+1}$, is realized in $\gB_{\alpha+1}$
(note: there are just $\le \lambda$ many).
Moreover, if the type is $p(\bar x)\in {\bold S}^m(B),B \subseteq
A^{\alpha+1}_i$, and $p(\bar x)$ ``say" $\bar x \cap \acl(B)=\emptyset$
he can demand such type to be realized by a sequence disjoint to
${\gB}_\alpha$: using $\Gamma^{\na}$ he can; moreover, for stationarity
many $i \in E_\alpha\cap S_{\zeta_\sat}$, for some $j_1<j_2$ successive member
of $u^\alpha_i$ there is such a sequence $\subseteq 
A^{\alpha+1}_{j_1} \setminus A^{\alpha+1}_{j_2}$.
This is help for (1)+(2) of \ref{e8} in the absence of Skolem
functions.  We have Skolem functions here because we use \ref{a2}(B), but if
we like to use \ref{a2}(A) and in some continuations we seem to have to be
more careful, for this end we list the cases.
\medskip

\noindent
\underline{Stage C}:  We turn to (3) of \ref{e8} to
which we assign two contractors: the end extender contractor
$\zeta_{\ex}$ and the branch killer contractor $\zeta_{\br}$. 
In order to satisfy the first phrase of (3) of \ref{e8} we 
can for $\alpha\in W_{\zeta_{\ex}}$, add
an element $a_\alpha=\bar{a}_{\alpha,0}$ realizing
$\gp_{\uf}^{\gB_\alpha}$, i.e. the end extender contractor decides that
$\Omega_{1+i}^\alpha$ is a case of $\Gamma^{\na}$ and
$\Omega_0^\alpha=\Gamma^{\gp_{\uf}}$ (see \ref{c86}(1) or see
\cite[2.11]{Sh:384}). 
But for $\alpha\in W'_{\zeta_{\br}} =: \{\alpha \in
W_{\zeta_{\br}}:\cf(\alpha)=\lambda$, and
$\alpha = \sup(W_{\zeta_{\ex}} \cap \alpha)\}$ we demand more:
\mn
\begin{enumerate}
\item[$\otimes_{1,\alpha}^c$]  for every $\beta>\alpha$ and $b \in
  \gB_\beta$ we have $(\alpha)$ or $(\beta)$ where:
\sn
\begin{enumerate}
\item[$(\alpha)$]  $(\exists\gamma\in W_{\zeta_{\ex}} \cap \alpha)[b\cap
a_\gamma \notin \gB_\alpha]$ (so $a_{\gamma'}$ for every large enough
$\gamma' \in W_{\zeta_{\ex}} \cap \alpha$ or even pseudofinite set $a\in
\gB_\alpha$ extending $a_\gamma$ is o.k. instead of $a_\gamma$, i.e.
$b \cap a \notin \gB_\alpha,a \in \gB_\alpha$)
\sn
\item[$(\beta)$] for some $\bar b \subseteq \gB_\alpha$ and formula
$\psi = \psi(x,\bar b)$, for every $c \in \gB_\alpha$:
\[
\gB_\beta \models [c \dot e b \equiv \psi(c,\bar b)].
\]
\end{enumerate}
\end{enumerate}
\mn
However to do this, the branch killer, for $\alpha$ as
above $(\in W'_{\zeta_{\br}})$ at stage $\alpha$, guarantee:
\mn
\begin{enumerate}
\item[$\otimes_{2,\alpha}^c$] for some club $E=E_\alpha^{\br}$ of
$\lambda$, for every $i \in S_{\zeta_{\br}} \cap E$ there is 
$j=j_i=j^{\br}(\alpha,i) \in u_i^\alpha,j < \max(u_i^\alpha)$, such 
that $\Gamma_j^{\alpha+1} = \Gamma^{\ids}_{\omega+1},\bar
c_j^{\alpha+1} = \langle c_n^{\alpha+1,j}: n \le \omega \rangle,
c_\omega^{\alpha+1,j}$ is $a_{\gamma(\alpha+1,j)}$ ($\in \gB_\alpha$) 
for some $\gamma(\alpha+1,j)\in W_{\zeta_{\rm ex}} \cap \alpha,
c_0^{\alpha+1,j} \in \gB_{\alpha+1}$ realizes ${\gp}_{\uf}^{\gB_\alpha}$
(we could have used $\Gamma_\omega^{\ids}$ instead) and $\langle
\gamma(\alpha+1,j_i):i\in S_{\zeta_{\br}} \cap E \rangle$ is strictly
increasing with limit $\alpha$.
\end{enumerate}
\mn
In stage $\alpha$ itself (i.e. defining $\gB_{\alpha+1}$),
for this the branch killer contractor $\zeta_{\br}$ chooses $<^\alpha$ (so
$<^\alpha_i = <^\alpha \restriction i$) as follows:
 $j_1 <^\alpha j_2$ \underline{if and only if} both $j_1$ and $j_2$ 
are $<\lambda$ and exactly one of the following occurs:
\mn
\begin{enumerate}
\item[$(a)$]  $j_1$ even, $j_2$ odd
\sn
\item[$(b)$]  both even, $j_1<j_2$
\sn
\item[$(c)$ ] both odd, $j_1<j_2$
\end{enumerate}
\mn
Further he decrees $\Omega_j^\alpha$ is: $\Gamma^{\gp_{\uf}}$ for
$j$ even, instance of $\Gamma^{\na}$ for $j$ odd.

Now to make $\otimes^c_{1,\alpha}$ and $\otimes^c_{2,\alpha}$ true in
stage $\alpha$ is straightforward, but preserving $\otimes^c_{1,\alpha}$
needs care.  In stage $\alpha$ itself we can think we first add 
$\langle\bar{a}_{\alpha,j}:j<\lambda$ even $\rangle$,
then the others.  When adding $\bar a_{\alpha,0}$ use the properties
of $\gp_\uf$ (see \ref{c83}(2), compare with \cite[2.20(2)]{Sh:384}) adding 
$\bar a_{\alpha,2+2i}$ no subset of $\bar a_{\alpha,0}$ is added
(\cite[2.20(2)]{Sh:384}) so it preserves the old 
$\otimes^c_{1,\alpha'}$ for $\alpha'\in W_{\zeta_{\br}} \cap \alpha$, 
as for $\bar a_{\alpha, 2i+1}$ we use only $\Gamma^{\na}$ so it is 
like the successor case below.

As for the case $\alpha$ is limit the preservation is automatic; we
are left with the successor case. So now suppose we are in stage
$\alpha$ and we would like to define $\gB_{\alpha+1}$ etc. and to preserve
$\otimes_{1,\beta}^c$ for $\beta \in \alpha \cap W'_{\zeta_{\br}}$. In step
$\epsilon<\lambda$ from $S_{\zeta_{\br}}$, after defining
$i=i_\epsilon^\alpha,u^\alpha_i = u_{i_\epsilon^\alpha}^\alpha$ 
and $p_{j_1}^\alpha,j_1=j(1) := \max(u^\alpha_{i^\alpha_\epsilon})$, 
by some bookkeeping we choose $\beta\in
W'_{\zeta_{\br}} \cap \bold c_i^\alpha$, and $y \in 
\bigcup\{\bar x^\alpha_\xi:\xi<j_1\}$ or just $y=\sigma(\bar x^\alpha_{\xi_1},
\ldots,\bar x^\alpha_{\xi_n})$ for some $n<\omega$, term $\sigma$
with parameters in $A^\alpha_{i_{\alpha,\epsilon}}$ and
$\xi_1,\ldots,\xi_n < j_1$ (if $T$ has no Skolem function: just in
their algebraic closure, no real difference).
We can find, by $\otimes^c_{2,\beta}$ and $\gB_\alpha$ satisfying
clause (E) of \ref{d14}, ordinals $i(*) \in E_\alpha^{\br} \cap 
S_{\zeta_{\br}}$ such that letting $j_2=j(2) := j^{\br}(\beta,i(*))$, 
the sequence $\bar c_{j_2}^\beta=
\langle c_{nj(2)}^\beta:n \le \omega\rangle$ is indiscernible over 
$A^\alpha_{j_2}$ which include $A^\alpha_{j_1}$. We
extend $p^\alpha_{j_1}$ to a complete type $q$ over
$\acl_{\gB_\alpha}(A^\alpha_{j_2} \cup \{c_\omega^{\beta,j_2}\})$ satisfying
the required bigness conditions (concerning $\Omega^\alpha)\epsilon,
\epsilon<j_1$ recall $\Gamma^{\ids}$ is orthogonal to every instance of 
${\bf \Upsilon}^{\ver}$). Remember: $c^{\beta,j(2)}_\omega \in 
\gB_\beta$ and $c^{\beta,j(2)}_0 \in \gB_{\beta+1}$
does $\subseteq^{\gB_{\beta+1}}$ -- extend every pseudo-finite set 
of $\gB_\beta$.  Choose appropriate $i_{\varepsilon +1}^\alpha$, i.e.
$A_{i\alpha,\epsilon+1}^\alpha \supseteq \dom(q)$.
\medskip

\noindent
\underline{First case}:  $q$ says that 
$y \cap c_\omega^{\beta,j(2)}$ is not equal to any member of the 
domain of $q$ or just $\dom(q)\cap \gB_\beta$.
Using niceness of the bigness demand on $q$ (for $\langle
\Omega_\xi^\alpha:\xi<j_1\rangle$), we can extend $q$ 
to an appropriate complete type over $A_{j_2+1}^\alpha$ and 
by clause (D)(7) of \ref{d14} we get the
desired contradiction, $\bigcup\{p^\alpha_i: i<\lambda\}$ will say $y\cap
c^{\beta,j(2)}_\omega$ is $\notin \gB_\alpha$. So clearly we succeed in
guarantying $\otimes^c_{1,\beta}$.
\medskip

\noindent
\underline{Second case}: $q$ says $y\cap c_\omega^{\beta,j(2)}=d$
for some $d \in \gB_\beta\cap \dom(q)$.  As $c_\omega^{\beta,j(2)}, 
c_0^{\beta,j(2)}$ realize the same type over $A_{j_2}^\alpha$,  
there is an elementary mapping $\bold g$ from 
$\acl_{\gB_\alpha}(A_{j_2}^\beta \cup \{c_\omega^{\beta,j(2)}\})$ onto
$\acl_{\gB_\alpha}(A_{j_2}^\beta \cup \{c_0^{\beta,j(2)}\})$ satisfying 
$\bold g \restriction A_{j_2}^\beta=$ the identity, now use 
$\bold g(q)$ instead of $q$ above, so we know $y \cap c_0^{\beta,j(2)}\in
\acl_{\gB_\alpha}(A_j^\beta \cup \{c_0^{\beta,j(2)}\}) \subseteq 
\gB_\alpha$ and as $\otimes^c_{1,\beta}$ holds for 
$\alpha$ and $c^{\beta,j(2)}_0$  does
$\subseteq^{\gB_{\beta+1}}$-extends 
every pseudo-finite member of $\gB_\beta$ we are done.

This argument works in both cases (as we use nice types), and ``moving by
$\bold g$" preserved the relevant properties.
\medskip

\noindent
\underline{Stage D}:   We now start dealing with part (4),
 but meanwhile, more generally 
we deal with complete embedding of a pseudo-finite 
Boolean Algebra into a Boolean ring (both represented in the model);
 the more restricted case from part (4) of \ref{e8} will 
use this. We assign a contractor, the coder $\zeta_{\cd}$, to try 
to code a complete embedding of one ``finite" Boolean ring 
(hence algebra) $\bold b_1[\gB]$ to another not necessarily ``finite"
Boolean ring, $\bold b_2[\gB]$ so, both ${\bold b}_1,\bold b_2$ are 
Boolean rings in the sense of $\gB$, note that as $\bold b_1$ is 
pseudo-finite it is atomic and call its set of atoms $\bold b_1^{\at}$.
 
By normal bookkeeping we assign to every such pair $(\bold b_1,\bold b_2)$ a
stationary subset $W_{({\bold b}_1, {\bold b}_2)}$ of $\{\delta \in
W_{\zeta_{\cd}}:\cf(\delta)=\lambda\}$. For $\delta \in 
W_{({\bold b}_1,{\bold b}_2)}$, the coder decrees that
$\Omega_n^\delta = \Gamma^{\av}_{\dot D,\langle n:n<\omega\rangle}$
($\dot D \in {\cP}_{0,0}$ any non principal ultrafilter on $\omega$), and
$\Omega_{\omega+i}^\delta = \Gamma^{\ms}[{\bold b}_1^{\at},\langle
1/x^\delta_n : n<\omega\rangle]$ for $i<\kappa_1$ (see $(*)_1$ of
\ref{e8}(4); this is an instance of $\Gamma^{\wm}$)
(so $u^\delta_0=\omega+\kappa_1$). All $\Omega_i^\delta$
(for $i \ge \omega+\kappa_1$)  will be instances of
$\Gamma^{\na}$.  Now for $\epsilon<\lambda$ such that
$i=i_\epsilon^\delta\in S_{\zeta_{\cd}}$,
$j_0=\sup\bigcup\{u_i^\beta:\beta\in \bold c_i^\alpha\})$, choose
$j_1=j_1^{\delta,\epsilon}$, $j_2=j_2^{\delta,\epsilon}$, such that
$j_0<j_1^{\delta,\epsilon}<j_2^{\delta,\epsilon}<\lambda$ and:
\mn
\begin{enumerate}
\item[$(*)_3$]  $(\gB_\delta \restriction
A_{j_1}^\delta,A_{j_0}^\delta) \prec (\gB_\delta,A_{j_0}^\delta)$
and
\sn
\item[$(*)_4$]   $(\gB_\delta \rest A_{j_2}^\delta,A_{j_0}^\delta,
A_{j_1}^\delta) \prec (\gB_\delta,A_{j_0}^\delta,A_{j_1}^\delta)$
define $u^\delta_i$ as $\bigcup\{u^\beta_i: \beta \in \bold c^\alpha_i\}\cup
\{j_0, j_1,j_2\}$ and demand:
\sn
\item[${{}}$]  $\otimes_0 \quad$ if $d \in \gB_\delta,\gB_\delta
  \models ``d \dot e \bold b_1^{\at}",d \in A_{j_2}^\delta 
\setminus A_{j_1}^\delta$ then there is no $d' \in \gB_\delta,
\gB_\delta \models$

\hskip25pt $``d' \dot e \bold b_1^{\at}",d' \ne d"$
such that 

\hskip25pt $\{d \le x_{\omega+i}^\delta \equiv d' \le x^\delta_{\omega+i}:
i<\kappa_1\}\subseteq \bigcup\limits_{j<\lambda} p_j^\delta$ and 

\hskip25pt for simplicity there is $i<\kappa_1$ such that 
$[d \le x_{\omega+i}^\delta]\in p_{j_2}^\delta$.
\end{enumerate}
\mn
Note that \underline{this is not} a part of the general machinery of \ref{d17},
but we shall see that it is compatible with it, i.e. this is
part of contractor $\zeta_\cd$'s work i.e. he overtake more control
this ``at the expense of" ``$\zeta_{\cd}$ is the main contractor for $\delta$".
Now it is reasonable to demand that when $p_{j}^\delta$ is defined, the
condition holds for $d\in A_{j_2^{\delta,\epsilon}} \setminus
A_{j_1^{\delta,\epsilon}},d' \in A_{j}^\delta$ (for every $\epsilon$ such that
$j_1^{\delta,\epsilon}<j_2^{\delta,\epsilon} \le j)$. So how can
we preserve this condition when defining $p^\delta_j$ 
(for $j \in E_j^+$)? For $j=0$ no problem. For limit $j$ there is no
problem. So assume, $i_\epsilon^\delta$ is defined, 
$j \in u^\delta_{i^\delta_\epsilon}$,
$p_j^\delta$ is defined and we have to define $p_{j(*)}^\delta$ where
$j(*)>j$ is the successor of $j$ in $u^\delta_{i^\delta_\epsilon}$ or 
$j$ is the last member of $u^\delta_{i^\delta_\epsilon}$ and
$j(*) = i^\delta_{\epsilon+1}$. We have to consider what is
the constraint.

Note: as the $1/x_n^\delta$ in $p^\delta_j$ satisfies
\mn
\begin{enumerate}
\item[$(*)_4$]  for every $n$, $m<\omega$, $[\bigwedge\limits_n
\gB_\delta\models 0<c<1/n \Rightarrow \bigwedge\limits_{n,m}
\gB_{\delta+1} \models ``c < 1/x^\delta_m <1/n"]$
\end{enumerate}
\mn
we really have freedom.
\medskip

\noindent
\underline{The First Case}:  No constraint.

So we have to extend $p_j^\delta$ in a nice way (to preserve clause (D7)
and $\otimes_0$).  By induction on $i<j$ we choose $\bar a_{\delta,i}$
to realize over $A^\delta_{j(*)}$ the right type (say
in some saturated $M$, $\gB_\delta\prec M$); i.e. we preserve (remember
we can look at $\bold b_1$ as the family of subsets of $\bold b^{\at}_1$, in
$T^*$'s sense):
\mn
\begin{enumerate}
\item[$(i)$]   $A^\delta_{j(*)} \cap \acl_M (A^\delta_j + 
\{\bar a_{\delta,\xi}:\xi<i\}) = A^\delta_j$
\sn
\item[$(ii)$]  $\tp(\bar a_{\delta,i},A^\delta_{j(*)} \cup 
\{\bar a_{\delta,\xi}:\xi<i\})$ is $\Omega^\delta_i$-big
\sn
\item[$(iii)$]  $a_{\delta,i} = \bar a_{\delta,i} \in \bold b_1[M]$ when
$\omega \le i<\omega+\kappa_1$
\sn
\item[$(iv)$]  $\bigl[ d \dot e^M \bold b^\at_1 \, \& \, d \in
(A^\delta_{j(*)}\setminus A^\delta_j)\Rightarrow \neg d \dot e^M
a_{\delta, i}\bigr]$ when  $\omega \le i<\omega+\kappa_1$
\sn
\item[$(v)$]  if $i \ge \omega$ and $d \in A^\delta_{j(*)}\setminus 
A^\delta_j$, $e\in \acl(A^\delta_j\cup \bigcup\limits_{\epsilon<i} 
\bar a_{\delta,\epsilon})$,
and $M \models e \subseteq \bold b^\at_1, e \le \log_2(a_{\delta,n})$
(equivalently if $2^{-|e|} \ge 1/a_{\delta,n}$)" for each $n$
\then \, $M \models ``\neg d \dot e e$".
\end{enumerate}
\mn
[Why we need (v)?  After defining $a_{\delta,i}$ for $i \le \omega$,
applying the relevant claim from \S3, in order to  have freedom for 
(iv) we need (v).]

This is possible by \ref{c53}(2) (check conditions).
If $i<\omega$ by clause (ii) we have exactly one choice, clause (i) is
easy by niceness (and uniqueness), clauses (iii)+(iv)+(v) are irrelevant, and
clause (v) is easy. For $i=\omega$ clauses (i)-(iv) are immediate, 
to assumption of (v) implies $e$ is truly finite hence all 
$\dot e$-members are by (i) not in $A^\delta_{j(*)},\setminus
A^\delta_j$.  If $\omega<i \le \omega+\kappa_1,i$ limits we have 
no problem, if $i$ is a successor ordinal we use \ref{c53}(2). 

Lastly,if $i \ge \omega+\kappa_1$ we just use niceness. 
\medskip

\noindent
\underline{The Second Case}:  We have obligation from Stage C.

I.e. $i^\delta_\epsilon \in S_{\zeta_{\br}}$ and let $j_1=$ max 
$(u^\delta_\epsilon)$, choose $\beta,i(*),j_2,y$ as in Stage C, 
``for the successor case" (with $j_1,j_2$ here standing for 
$j(*),j$ there).  Choose the type $q$ over $\acl_{\gB_\alpha}
(A^\alpha_{j_1} \cup \{c^{\beta,j_2}_\omega\})$ as the first case
(in our present stage) with $A^\alpha_{j_1},\acl_{\gB_\alpha}
(A^\alpha_{j_1} \cup \{c^{\beta,j_2}_\omega\})$ here standing for 
$A^\alpha_{j(*)}, A^\alpha_j$ there.

If in Stage C, first case apply, then we can choose appropriate
$i^\delta_{\epsilon+1}$ and extend $q$ to a type as required by 
the proof of first case (in our present stage).
If the second case in Stage C apply, the
elementary mapping ${\bold g}$ preserve the right things so no
problems: just like when the first case applies.
\medskip

\noindent
\underline{The Third Case}:  $\Gamma_j^\delta$ is defined and
equal to $\Gamma_{\omega+1}^{\ids}$.
We just first lengthen the indiscernible set to $\dot{\bold I}$,
$|\dot{\bold I}| > \beth_{(2^\lambda)^+}$ say in $M^*$ where 
${\gB}_\alpha \prec  M^*$ this extend $p_j^\delta$ to $p^+$, complete type
over $\acl(A_j^\delta \cup \dot{\bold I})$ which is a nice extension of $p$
big for $\langle \Omega^\delta_\zeta:\zeta \in u^\delta_{i,\dot e}\rangle$
such that for some $u \subseteq \kappa_1$ for all $d \in
\acl(A_j^\delta \cup \dot I) \setminus A_j^\delta$ (note
$A_j^\delta=\acl (A_j^\delta$)) we have: $\{(d \le 
x_{\omega+i}^\delta)^{[\iif(i\in u)]}:i<\kappa_1\}\subseteq p^+$, 
where $u$ is chosen such that for no $d\in A_j^\delta$, $\{(d \le 
x_{\omega+i}^\delta)^{\iif(i\in u)]}:i<\kappa_1\}\subseteq p$ 
(note: $p^+$ exists by \ref{c53}(1)),
then use the proof in \ref{b17}.
(i.e. $\Gamma^{\ids} \bot \Gamma^{\ms}$, in fact $u=\emptyset$ is O.K.)
\medskip

\noindent
\underline{The fourth case}:   $\Gamma_j^\delta$ is defined and
equal to $\Gamma^{\wm}_{a,\dot w,\bar{c}}$.

Similar to first case; for $n< \omega$, we choose $a_{\delta,n}$ (corresponding
 to  $x_n^\delta$) by \ref{c29}(3), i.e. as $\Gamma_\omega^{\av}$,
$\Gamma_{a,\dot w,\bar c}^{\wm}$ are orthogonal.

For $x_{\omega+i}^\delta$ ($i<\kappa_1$) we can use \ref{c59}.
\medskip

\noindent
\underline{The fifth case}: $\Gamma_j^\delta$ is defined and equal 
to $\Gamma^{\na}$.

Easy.
\medskip

\noindent
\underline{The sixth case}: $j,j(*)$ are like $j_1,j_2$ above
 in $(*)_2+\otimes_0$ of this stage. 

Again by claim \ref{c53}(2) choosing the function $h$ carefully
enough remembering $|A^\delta_{j(*)}| < \lambda \le 2^{\kappa_1}$.
\medskip

\noindent
\underline{Stage E}:  Assume
\mn
\begin{enumerate}
\item[$\boxtimes$]  ${\bold b}_1$, ${\bold b}_2$ as in stage D, for
$\gB_{\lambda^+}$, ${\bold f}$ is a complete embedding of 
${\bold b}_1[\gB_{\lambda^+}]$ into ${\bold b}_2[\gB_{\lambda^+}]$. 
\end{enumerate}
\mn
${\bold f}$ will be fixed for stages E---J. 

So

\[
\cC^1_{\bold f} = \{\delta<\lambda^+:(\gB_\delta,\bold f \rest 
\gB_\delta) \prec (\gB,{\bold f})\}
\]

\mn
is a club of $\lambda^+$; let $W_{\bold f} = \cC^1_{\bold f} \cap
W_{({\bold b}_1,{\bold b}_2)} \subseteq \cC^1_F \cap W_{\zeta_{\cd}}$.

Clearly by \ref{e11}(2)
\mn
\begin{enumerate}
\item[$\otimes$]  for $\delta \in \cC^1_{\bold f}$, ${\bold f} \rest
  \bold b_1[\gB_\delta]$ is a complete embedding of ${\bold b}_1[\gB_\delta]$
into ${\bold b}_2[\gB_\delta]$.
\end{enumerate}
\mn
For $\alpha \in W_{\bold f}$ let $B^*_\alpha := \{b \in 
\bold b^{\at}_1[\gB^*]: b \in \gB_\alpha$ and for no $b' \ne b$ do we have

\[
b' \dot e^{\gB_\alpha} {\bold b}^{\at}_1 \, \& \, b'\in \gB_\alpha \,
\& \, \bigwedge\limits_{i<\kappa_1}[b \dot e^{\gB^*} a^\alpha_{\omega+i}\equiv
b' \dot e^{\gB^*} a^\alpha_{\omega+i}]\}.
\]

\mn
Note that $B^*_\alpha$ depends on $\alpha$, $\gB_\alpha$,
${\bold b}_1$, $\langle a^\alpha_{\omega+i}:i < \kappa_1\rangle$ but 
not on ${\bold f}$. For the rest of stage E we fix $\alpha$.

For $\delta \in \cC^1_{\bold f}$ for awhile we shall try to show that
${\bold f} \restriction B^*_\delta$ belongs to ${\cP}_{\lambda^+}$, 
this in the following substages:
\mn
\begin{enumerate}
\item[$\otimes_1$]  from ${\bold f}^0_\delta = {\bold f} \restriction
\{a_{\omega+i}^\delta:i<\kappa_1\}$ we can reconstruct 
${\bold f}^1_\delta = {\bold f} \restriction B_\delta^*$,
so if ${\bold f}^0_\delta = {\bold f} \restriction 
\{a_{\omega+i}^\delta:i<\kappa_1\} \in \bigcup\limits_{\alpha} \cP_\alpha$
then ${\bold f}^1_\delta = {\bold f} \restriction B_\delta^* \in 
\cP_{\lambda^+}$.
\end{enumerate}
\mn
[How?  For $d \in B^*_\delta$, ${\bold f}(d)$ is the maximal
member of ${\bold b}_2[\gB_\delta]$ which is $\le {\bold f}(c)$ 
whenever $c\in\{a_{\delta, \omega+i}:i<\kappa_1,\ d \le
a_{\delta, \omega+i}\} \cup \{\bold f(-a_{\delta, \omega+i}):i< \kappa_1$ and
$d \le -a_{\delta,\omega+1}\}$. Note that $-a_{\delta,\omega+i}$ is
well defined as ${\bold b}_1 [\bold B]$ is a Boolean Algebra, by the first
paragraph of stage D. Also note: ${\bold f}(d)$
satisfies this as ${\bold f} \restriction {\bold b}_1[\gB_\delta]$ is 
a complete embedding of ${\bold b}_1 [\gB_\delta]$ into 
${\bold b}_2 [\gB_\delta]$ and in ${\bold b}_1 [\gB_\delta],d$ is 
a maximal member of ${\bold b}_1 [\gB_\delta]$ which is $\le c$
whenever $c \in \{a_{\delta,\omega+i}:i < \kappa_1,d \le a_{\delta,\omega+i}\}
\cup \{-a_{\delta, \omega+i}:i< \kappa_1,d \le
-a_{\delta,\omega+i}\}$).]
\mn
\begin{enumerate}
\item[$\otimes_2$]  if $\ell<n(*)$, $A\subseteq \gB_{\lambda^+}$,
$|A|<\kappa_\ell$, \then \, we can find $A' \subseteq \gB_{\lambda^+}$ 
satisfying $|A'|<\kappa_\ell,A' \in \bigcup\limits_{\alpha<\lambda^+} 
\cP_\alpha$ and $A \subseteq A'$.
\end{enumerate}
\mn
[Why? By $(*)$ of \ref{e8}(4) we prove by induction on
$\ell$. For $\ell =0$ immediate by $\kappa_\ell$ being regular (as
$\kappa_0= \lambda= \cf(\lambda)$ and have 
$\cf([\lambda]^{<\lambda},\subseteq)=\lambda$).
For $\ell+1$, by the induction hypothesis we can find $A'' \in 
\bigcup\limits_{\alpha< \lambda^+} {\cP}_\alpha$ such that 
$A \subseteq A'',|A''| < \kappa^*_\ell$ and now recall
$\cf([\kappa_\ell]^{\le \kappa_{\ell+1}}, \subseteq) \le \lambda^+$
so a cofinal subset of $[\kappa_\ell]^{\le \kappa_{\ell+1}}$ has
cardinality $\le \lambda^+$ and belongs to $\bigcup\limits_{\alpha<
\lambda^+} {\cP}_\alpha$ hence is included in 
$\bigcup\limits_{\alpha< \lambda^+} {\cP}_\alpha$.]
\mn
\begin{enumerate}
\item[$\otimes_3$]  if $\ell< n(*),A \subseteq {\bold b}^{\at}_1
[\gB_{\lambda^+}],|A| < \kappa_\ell,A \in 
\bigcup\limits_{\alpha<\lambda^+} \cP_\alpha$ \then \, 
we can find $A'$, $A_1$, 
$A_2\in\bigcup\limits_{\alpha<\lambda^+}
{\cP}_\alpha$, such that $A \subseteq A_1$, $\rang({\bold f} 
\restriction A_1)=A_2,|A_1 |+ |A_2|<\kappa_\ell,A_i \subseteq \bold b_i
[\gB_{\lambda^+}]$ for $i=1,2$ and $A' \subseteq {\bold b}_1[\gB_\delta],
|A'|<\kappa_{\ell+1}$ and from ${\bold f} \restriction
A'$ and $A_1,A_2$ we can reconstruct ${\bold f} \restriction A_1$; 
i.e. it belongs to $\bigcup\limits_{\alpha<\lambda^+} \cP_\alpha$.
\end{enumerate}
\mn
[Why ? By $\otimes_2,\cf([\lambda^+]^{< \kappa_\ell},\subseteq) =
\lambda^+$ so by \ref{e8}(4)$(*)$ and 
\cite[3.11]{Sh:E62} there is a stationary
${\cS} \subseteq [\lambda]^{\le \kappa_\ell}$ of cardinality
$\le \lambda^+$ hence there is a model $N,|N| \in
\bigcup\limits_{\alpha<\lambda^+} \gA_\alpha$,
 such that $N \prec (\gB_{\lambda^+},\bold f),A \subseteq N,
\|N\| < \kappa_\ell$ and ${\bold b}_1,{\bold b}_2 \in N$. 
Let $A_i= N \cap {\bold b}_i[\gB_{\lambda^+}]$ for $i = 1,2$, 
so ${\bold f} \restriction N$ is a complete embedding
of ${\bold b}_1[\gB_{\lambda^+}] \restriction N$ into 
${\bold b}_2[\gB_{\lambda^+}]\restriction N$; hence to reconstruct it it
suffice to reconstruct ${\bold f} \restriction 
({\bold b}^\at_1[\gB_{\lambda^+}]\cap N)$. But the saturator guarantee the
existence of $a_i\in \gB_{\lambda^+}$ satisfying 
$a_i \dot e {\bold b}_1[\gB_{\lambda^+}]$ for $i<\kappa_{\ell+1}$ such that:
if $d' \ne d'' \in {\bold b}_1[N]$ then 
$\bigvee\limits_{i<\kappa_{\ell+1}} [d' \le a_i \,\& \, 
d''\cap a_i=0_{{\bold b}_1}]$. Let $A'=\{a_i: i<\kappa_{\ell+1}\}$.]
\mn
\begin{enumerate}
\item[$\otimes_4$]  if $A \subseteq {\bold b}_1[\gB_{\lambda^+}],
|A| \le \kappa_{n(*)}$ then ${\bold f} \restriction A' \in
\bigcup\limits_{\alpha<\lambda^+}{\cP}_\alpha$ for some 
$A'$ satisfying $A\subseteq A' \subseteq {\bold b}_1 [{\gB}_{\lambda^+}]$
and $|A'|<\lambda$.
\end{enumerate}
\mn
[Why?  As $\lambda^{\kappa_{n(*)}} \le \lambda^+$.]
\mn
\begin{enumerate}
\item[$\otimes_5$]  if $A \in {\cP}_{\lambda^+},A \subseteq
{\bold b}_1[\gB_{\lambda^+}],|A| < \lambda$ then ${\bold f} \rest A \in
\cP_{\lambda^+}$ (remember that $[{\bold b}_1[\gB_{\lambda^+}]]^{<\lambda}\cap
\cP_{\lambda^+}$ is cofinal in $[{\bold b}_1[\gB_{\lambda^+}]]^{<\lambda}$).
\end{enumerate}
\mn
[Why?  Put together $\otimes_2,\otimes_3,\otimes_4$; i.e. we can prove
by induction on $\ell \le n(*)$ that if $|A| < \kappa_{n(*)-\ell}$ 
then the conclusion holds; now for $\ell=0$ use $\otimes_4$ and 
for $\ell+1$ use $\otimes_3$.]
\mn
\begin{enumerate}
\item[$\otimes_6$]   ${\bold f} \restriction
B_\alpha^* \in \bigcup\limits_{\alpha<\lambda^+}{\cP}_\alpha$ for
$\alpha \in W_{\bold f}$.
\end{enumerate}
\mn
[Why?  Put together $\otimes_1$, $\otimes_5$.]

Let ${\cI}^1_\alpha= {\cI}^{1,{\bold f}}_\alpha$ be the ideal of
${\bold b}_1[\gB_\alpha]$ generated by $\{c:c \in 
\bold b^{\at}_1[\gB_\alpha]$ but $c \notin B^*_\alpha\}$, clearly
${\cI}^1_\alpha\in {\cP}_{\lambda^+}$ but we do not claim any
definability in $\gB_{\lambda^+}$. Also clearly $c \in \cI^1_\alpha$
\underline{if and only if} for some $n<\omega$ and $c_1,\ldots c_n \in 
{\bold b}^{\rm at}_1[{\gB}_\alpha] \setminus B^*_\alpha$
we have ${\bold b}_1[{\gB}_\alpha] \models ``c=c_1\cup \ldots \cup c_n"$.a

As $B^*_\alpha$ has members in every infinite subset of ${\bold b}_1
[{\gB}_\alpha]$ definable in ${\gB}_\alpha$ with parameters, clearly
\mn
\begin{enumerate}
\item[$\otimes_7$]  for $b \in {\bold b}_1[\gB_\alpha]$ we have $b \in
{\cI}^1_\alpha \Leftrightarrow \neg (\exists c \in 
\bold b_1^{\at}[\gB_\alpha])(c \in B^*_\alpha \, \& \, c \le b)$.
\end{enumerate}

Let

\[
{\cI}^2_\alpha = {\cI}^{2,\bold f}_\alpha = 
\{c \in {\bold b}_2[\gB_\alpha]: \text{ for every } b \in 
B^*_\alpha \text{ we have } \bold f(b) \cap c=0_{\bold b_1}\}.
\]

\mn
Clearly ${\cI}^2_\alpha \in {\cP}_{\lambda^+}$ is an ideal of
${\bold b}_2 [{\gB}_\alpha]$ . We define a function
${\bold f}^0_\alpha$ with domain ${\bold b}_1[\gB_\alpha]$:

\begin{equation*}
\begin{array}{clcr}
\bold f^0_\alpha(d) = \{c \in {\bold b}_2[\gB_\alpha]: &\text{ for
  every } b \in B^*_\alpha \text{ we have}\\
  &b \le_{\bold b_1} d \Rightarrow \bold f(b) \le_{\bold b_2} c) \text{ and}\\
  & b \cap d= 0_{\bold b_1} \Rightarrow \bold f(b) \cap c= 0_{\bold b_2})\}.
\end{array}
\end{equation*}

\mn
It is easy to see that ${\bold f}^0_\alpha \in {\cP}_{\lambda^+}$ and
${\bold f}^0_\alpha$ is a homomorphism  from the Boolean Algebra 
${\bold b}_1[\gB_\alpha]$ into the Boolean ring ${\bold
  b}_2[\gB_\alpha] \cI^2_\alpha$, and ${\bold f}(d) \in 
\bold f^0_\alpha(d)$, and ${\cI}^1_\alpha$ is the kernel of 
${\bold f}^0_\alpha$.

We now show (recall that by assumption $\lambda=
\lambda^{\aleph_0}, \theta=\aleph_0$):
\mn
\begin{enumerate}
\item[$\otimes_8$]  Assume $\cf(\alpha) \ne \theta$ and the saturator
works for unboundedly many $\beta < \alpha$. If ${\bold f}',\bold f''$
are two functions satisfying the information on
${\bold f}\restriction \gB_\alpha$ gathered so far 
(more exactly: ${\bold f}' \restriction B^*_\alpha = \bold f \rest 
B^*_\alpha = {\bold f}'' \restriction B^*_\alpha$, and ${\bold
  f}',\bold f''$ are complete embeddings of ${\bold  b}_1[\gB_\alpha]$ 
into ${\bold b}_2[\gB_\alpha]$ hence looking at the definitions of 
$\cI^{1,{\bf f}}_\alpha,\cI^{2,{\bold f}}_\alpha, {\bold f}^0_\alpha$, clearly 
$d \in {\bold b}_1[\gB_\alpha] \Rightarrow \bold f'(d) /
{\cI}^2_\alpha = {\bold f}(d) / {\cI}^2_\alpha = {\bold f}''
(d)/\cI^2_\alpha$) \then \, $\{e:e \in {\bold b}_1^{\at}[\gB_\alpha]$ and
${\bold f}'(e) \ne {\bold f}''(e)\}$ has  $<\theta$ members.
\end{enumerate}
\mn
[Why?  Assume $e_i \in {\bold b}^{\at}_1[\gB_\alpha]$ for $i<\theta$ are
pairwise distinct and ${\bold f}'(e_i) \ne {\bold f}''(e_i)$. So
\wilog \, ${\bold f}'(e_i) - {\bold f}''(e_i)> 0_{{\bold b}_2[\gB_\alpha]}$. As
${\bold f}''$ is a complete embedding of ${\bold b}_1[\gB_\alpha]$ into
${\bold b}_2[\gB_\alpha]$ clearly for $i<\theta$ there is $e'_i \in
\bold b_1^\at [\gB_\alpha]$ such that ${\bold f}''(e'_i) \cap
({\bold f}'(e_i) - \bold f''(e_i)) > 0_{\bold b_2[\gB_\alpha]}$.
This implies $\bold f''(e'_i) \ne \bold f''(e_i)$ hence 
$e'_i \ne e_i$, and as $\langle e_i:i< \theta\rangle$ is without repetitions,
\wilog \, $e'_i \notin \{e_j:j<\theta\}$ for $i<\theta$. As
$\cf(\alpha) \ne \theta$ \wilog \, for some $\beta< \alpha,
\{e_i,e'_i:i<\theta\} \subseteq \gB_\beta$, hence (as $\lambda = 
\lambda^{\langle \theta\rangle}$ and less suffice) for some
countable $u \in [\theta]^{\aleph_0},\{(e_i,e'_i):i < \omega\}
\in {\cP}_{\beta}$, \wilog \, $u=\omega$; hence by the saturater work 
there is $d \in {\bold b}_1[\gB_\alpha]$ such that $e_n \le d,e'_n \cap d =
0_{{\bold b}_1}$.  Hence by $(*)_3$ for some $c \in {\cI}^1_\alpha,
\bold f'(d-c) \le {\bold f}''(d),{\bold f}''(d-c) \le {\bold f}'(d)$ and so
in ${\bold b}_2[\gB_\alpha]$ we have

\begin{equation*}
\begin{array}{clcr}
({\bold f}'(d) - {\bold f}''(d)) \cap {\bold f}''(e'_n) \ge &\bold
f''(-d) \cap {\bold f}'(d) \cap {\bold f}''(e'_n) = \\
   &{\bold f}''((-d) \cap e'_n) \cap {\bold f}'(d) \ge \\
   & {\bold f}''(e'_n)\cap {\bold f}'(d) \ge\\
   & {\bold f}''(e'_n) \cap {\bold f}'(e_n) \ge \\
   & {\bold f}''(e'_n) \cap ({\bold f}'(e_n) - {\bold f}''(e_n)) 
 > 0_{{\bold b}_1}
\end{array}
\end{equation*}

\mn
[Why?  As ${\bold b}_1[\gB_\alpha]$ is a Boolean algebra 
(${\bold b}_1^{\at}$ being ``finite"); as ${\bold f}''$ is an embedding 
as $e'_n \cap d =0$ so $e'_n \le -d$; as $e_n \le d$; by Boolean
rules; and by the choice $e'_i$.]

Choose $n$ such that $\neg (e'_n \le_{{\bold b}_1} c)$, so we got
contradiction to ${\bold f}'(d) / {\cI}^2_\alpha = {\bold f}''(d)/
{\cI}^2_\alpha$.  So $\otimes_8$ really holds.]

\relax From now on we assume the conclusion of $\otimes_8$ holds which suffice for
\ref{e8}.

As ${\bold f} \restriction {\bold b}_1[\gB_\alpha]$ satisfies the requirements
in $\otimes_8$, (there is a least one such ${\bold f}'$, and) by
$\otimes_8$ there are $\le \lambda^{\langle\theta\rangle} 
\le \lambda^+$ such functions ${\bold f}'$ so we conclude
\mn
\begin{enumerate}
\item[$\otimes_9$]  $({\bold f} \restriction {\bold b}_1[\gB_\alpha])
 \in \cP_{\lambda^+}$ hence $\beta< \alpha \Rightarrow {\bold b}\restriction
{\bold b}_1 [\gB_\beta] \in {\cP}_{\lambda^+}$.
\end{enumerate}
\mn
So

\begin{equation*}
\begin{array}{clcr}
\cC^2_{\bold f} = \{\gamma< \lambda^+: &\gamma 
\text{ a limit ordinal such that if}\\
  & \alpha < \gamma \text{ then } {\bold f} \restriction 
{\bold b}_1[\gB_\alpha] \in {\cP}_\beta \text{ for some }
\beta \in (\alpha, \gamma)\}
\end{array}
\end{equation*}

\mn
is a club of $\lambda^{+}$ hence (as $\lambda= 
\lambda^{\langle \theta\rangle}$)
\mn
\begin{enumerate}
\item[$\otimes_{10}$]  if $\alpha \in \cC^2_{\bold f}$ and 
$\cf(\alpha)=\theta$ then $f_\alpha := {\bold f} \restriction 
\bold b_1 [\gB_\alpha] \in {\cP}_\alpha$ [remember ${\cP}_\alpha \ne 
\bigcup\limits_{\gamma<\alpha} {\cP}_\gamma$!].
We shall not ${\bold f}_\alpha \in \cup \{{\cP}_\alpha:\alpha<\lambda^+\}$
when not necessary i.e. $\bold q^0_\alpha$ suffice.
\end{enumerate}
\medskip

\noindent
\underline{Stage F}:  We have a contractor, the pre-separator,
$\zeta_{\ps}$ which acts for any (fixed for this stage):
\mn
\begin{enumerate}
\item[$\otimes_{11}$]  ${\bold b}_1,{\bold b}_2$ as
in Stage D, and $\bold p,c,\bar c^1,\bar c^2$ as in
\ref{c80} of length say $\omega$ with ${\bold b}^\at_1$ for $a$
\end{enumerate}
\mn
(for example those which the contractor
$\zeta_{\ps}$ posed, they will be fixed in this stage),  for
stationarily  many $\alpha \in W_{\zeta_{\ps}}$ such that
$\cf(\alpha)= \lambda$ (and the contractor for this $\alpha$ 
choose ${\bold b}_1,{\bold b}_2,\bar c^1,\bar c^2,\bold p$ which are 
from ${\gB}_\alpha$ and so $\bar c^1,\bar c^2 \in {\cP}_\alpha$), and 
is quite closed under the saturator work. 

Now
\mn
\begin{enumerate}
\item[$(*)_5$]  the pre-separator takes care that for every $i \in
E_\alpha\cap S_{\zeta_{\ps}}$, there is $j_i\in u^\alpha_i$, $j_i <
\max (u^\alpha_i)$,\ such that $\Gamma^\alpha_i$ has the form
\[
\Gamma^{\wm}[{\bold b}_1,\dot w_{\bold p},\bar{c}^1]=
\Gamma^{\wmg}[{\bold b}^\at_1,\bold p,\bar c^1]
\]
\end{enumerate}
\mn
(note: ${\bold b}_1$ can be considered the power set of 
${\bold b}_1^{\at}$, recall $\gB_\alpha$ ``think" that ${\bold b}^\at_1$ is
finite hence ${\bold b}_1$ is a Boolean Algebra).
The set of such $\alpha$'s will be called $W_{\zeta_{\sep},\bold
  b_1,\bold b_2,c,\bar c^1,\bar c^2,\bold p}$.

We have another contractor, the separator, $\zeta_{\sep}$, such that:
stationarily many $\alpha\in W_{\zeta_{\sep}},\cf(\alpha)=\lambda$,  are
assigned to ${\bold b}_1$, ${\bold b}_2$, $c$, $\bar c^2,\bar
c^1,\bold p$ and $p \in \cP_{\alpha,0}$ which is a 
$\Gamma^{\ms}_{{\bold b}^\at_1,\bar c^2}$-big type over
$\gB_\alpha$ of cardinality $<\lambda$, and the pre-separator has
acted in some $\alpha_1 < \alpha$, for the relevant parameters; 
the separator chooses $\Omega^\alpha_0 = \Gamma^{\ms}_{{\bold b}^{\at}_1,
\bar c^2}$  and make $p\subseteq p^\alpha_0$
(so $x^\alpha_0$ will be in ${\bold b}^\at_1 [\gB_{\alpha+1}]$). Now for
every $i^\alpha_{\epsilon} \in S_{\zeta_{\sep}} \cap E_{\alpha_1}$ with
$\alpha_1 \in \bold c^\alpha_{i^\alpha_\epsilon}$,
he took care to have $j \in u^\alpha_{i^\alpha_{\epsilon}}$, $j<
\max(u^\alpha_{i_{\epsilon}})$ such that for some $d^\alpha_j \in
{\bold b}_1 [\gB_{\alpha_1}]$, the type $\tp(d^\alpha_j,A^\alpha_j,
\gB_\alpha)$ is $\Gamma^{\wmg}_{{\bold b}^{at}_1,\bold q,\bar{c}^1}$-big and
$[x^\alpha_0 \dot e d^\alpha_j] \in p_{j'}$ for $j' > j$, but
$\Gamma^{\alpha+1}_j = \Gamma^{\wmg}_{{\bold b}^\at_1,\bold p,\bar c^1}$ and
$c^{\alpha+1}_j = d^\alpha_j - a_{\alpha,0}$ (in $\bold b_1[\gB_{\lambda^+}]$'s
sense!); this is possible by the pre-separator
work for some higher $j$ and \ref{c80} (and the assumption on
$\bold p,c,\bar c^1,\bar c^2$ in $\otimes_{11}$).

Let $\alpha < \lambda^+$ be such that ${\gB}_\alpha$ is closed under
${\bold f}$ and we fix $\alpha$ for a while. Let 
$\beta \le \lambda^+$ be such that $\beta \ge \beta_\alpha^*
= \min\{\beta:\beta>\alpha$ and ${\bold f}_\alpha \in {\cP}_\beta\}$ 
where below $\alpha$ there is $\alpha_1$ as above. 
We define an equivalence relation
${\cE}_{\alpha,\beta}$ on ${\bold b}^{\at}_1 [\gB_{\beta}]$:
\, $d_1 {\cE}_{\alpha,\beta} d_2$ \underline{if and only if} 
for every $b \in {\bold b}_1 [\gB_{\alpha}]$, we have $d_1 \le 
b \Leftrightarrow d_2 \le b$.

Let $B^1_{\alpha,\beta} = \{d \in {\bold b}^{\at}_1[\gB_{\beta}]:
d/{\cE}_{\alpha,\beta}$ is a singleton$\}$.

Note: as the separator do his for stationarily many $\alpha' <
\lambda^+$ of cofinality $\lambda$, we can use $\alpha$ such that
$\gB_\alpha$ is closed under ${\bold f}$. Let $B^2_{\alpha, \beta}=
B^1_{\alpha,\beta} \setminus \gB_\alpha$,

\begin{equation*}
\begin{array}{clcr}
B^3_{\alpha, \beta} = B^{3,{\bold f}}_{\alpha,\beta} = 
\{a_{\gamma,0}: &\gamma \in [\alpha, \beta) \text{ and in stage } \gamma
\text{ the separator acts}\\
  &(\text{for the parameters } {\bold b}_1, {\bold b}_2,c,\bold p,
\bar c^1,\bar c^2)\}.
\end{array}
\end{equation*}
\mn
Recall that there is $\alpha_1< \alpha$ 
such that the pre-separator act for those parameters.
Clearly $B^3_{\alpha,\beta} \subseteq B^2_{\alpha,\beta}\subseteq
B^1_{\alpha, \beta}$, moreover
\mn
\begin{enumerate}
\item[$(*)_6$]  for every $\beta_1\in [\beta,\lambda^+)$ we 
have $B^3_{\alpha, \beta} \subseteq B^3_{\alpha,\beta_1}$
and $B^3_{\alpha,\beta} \subseteq B^2_{\alpha,\beta_1}$.
\end{enumerate}
\mn
[Why?  $B^3_{\alpha,\beta}\subseteq B^3_{\alpha,\beta_1}$ holds trivially.
As for $B^3_{\alpha,\beta} \subseteq B^2_{\alpha,\beta}$                       
(this is the whole point of the work of the
separators, that is, assume $a_{\gamma,0} \in 
B^3_{\alpha,\beta}$ and \wilog \, $\beta=\beta_1$ 
if $a_{\gamma,0} \notin B^2_{\alpha,\beta}$ then there is
$d \in a_{\gamma,0}/\epsilon_{\alpha,\beta} \smallsetminus 
\{a_{\gamma,0}\}$,  but for unboundedly many $j \in 
E^+_{\gamma+1}$  the separator in
stage $\gamma$, for $j$ choose $c^{\gamma+1}_j= d^\alpha_j- a_{\gamma,0}$
and $A^{\gamma+1}_j \cup
\{a_{\gamma,0},d\} \subseteq A^\beta_j$, and $\tp(c^{\gamma+1}_j,
A^\beta_j,\gB_\beta)$ is $\Gamma^{\wmg}_{{\bf b}^{\at}_1, 
\bold q,\bar c^1}$-big hence $d' \in {\bold b}^{\at}_1
[{\gB}_\beta] \cap A^\beta_j \Rightarrow \neg (d' \dot e c^{\gamma+1}_j)$ hence
in particular $\neg (d \dot e d^\alpha_j)$ but $(a_{\gamma,0} \dot e
d^\alpha_j)$ and $d^\alpha_j \in {\bold b}^{\at}_1
[\gB_\alpha]$ so $d^\alpha_j$ exemplify 
$\neg(d \dot e_{\alpha,\beta} a_{\gamma,0})$ so $(*)_5$ holds].

For $\beta\in (\beta^*_\alpha, \lambda^+)$ we define a relation $R_\beta=
R_{\alpha,\beta}$:
\mn
\begin{enumerate}
\item[$\boxtimes$]  $c R_\beta b$ iff:
$c \in B^3_{\alpha, \beta},b \in {\bold b}_2[\gB_\beta]$,  and for every
$d \in {\bold b}_1[\gB_{\alpha}]$ and $d' \in {\bold
  b}_2[\gB_{\alpha}]$ such that $\bold f(d) = d' \mod
{\cI}^2_{\alpha}$ equivalently, $d' \in {\bold f}^0_\alpha(d))$
we have $c \le_{{\bold b}_1} d  \Rightarrow b \le_{{\bold b}_2} d'$
and $\bold b_1 \models c \cap d =0_{\bold b_1} \Rightarrow \bold b_2 
\models b \cap d' = 0_{\bold b_2}$.
\end{enumerate}
\mn
Now
\mn
\begin{enumerate}
\item[$(*)_7$]  if $\beta_1 <\beta_1 \le \alpha$ are as above then
$R_{\alpha,\beta_1}= R_{\alpha,\beta_2} \cap ({\gB}_\beta 
\times {\gB}_\beta)$
\sn
\item[$(*)_8$]  $R_\beta \in {\cP}_\beta$.
\end{enumerate}
\mn
[Why?  As ${\bold f}^0_\alpha$, belongs to ${\cP}_\beta$ etc.]
\mn
\begin{enumerate}
\item[$(*)_9$]  if $c R_\beta b$ then $b \notin 0_{\bold b_2}$.
\end{enumerate}
\mn
[Why?  Think].
\mn
\begin{enumerate}
\item[$(*)_{10}$]  If $({\bold f}(c) = b)\in {\gB_\beta}$, and 
$c \in B^3_{\alpha, \beta}$ then $c R_\beta b$.
\end{enumerate}
\mn
[Why?  Let $d \in {\bold b}_1[\gB_{\alpha}]$ and $d' \in 
{\bold b}_2[\gB_{\alpha}]$ be such that ${\bold f}(d) = d' \mod \cI^2_\alpha$. 
Now assume $c \le d$ (in ${\bold b}_1[\gB_\beta]$) then ${\bold f}(c)
\le \bold f(d)$ (in $\bold b_2[\gB_\beta]$) and let $\bold f(d) - d' \le
{\bold f}(d_0) \cup {\bold f}(d_1)\cup \ldots \cup {\bold f}(d_{n-1})$,
for some true natural number $n$ and $d_\ell \dot e \bold b^\at_1
[\gB_\alpha]$ hence

\[
{\bold f}(c) - d' \le ({\bold f}(c) - {\bold f}(d)) \cup 
({\bold f}(d) - d') = ({\bold f}(d) - d') \le {\bold f}(d_0) \cup 
\ldots \cup {\bold f}(d_{n-1})
\]

\mn
but $c \in {\bold b}^\at_1 [\gB_\beta] \smallsetminus \gB_\alpha$, so $c\cap
d_\ell =0$ in ${\bold b}_1 [\gB_\beta]$ hence ${\bold f}(c)
\cap {\bold f}(d_\ell) =0$ in ${\bold b}_2[{\gB}_\beta]$ 
so at last ${\bold f}(c) - d' = 0$ in ${\bold b}_2[{\gB_\beta}]$. 
So $c \le_{{\bold b}_1} d \Rightarrow b \le_{{\bold b}_2} d'$. Similarly
${\gb}_2[!!{\gB}_\beta] \models c \cap d =0 \Rightarrow 
{\bf b}_2 [{\gB}_\beta] \models b \cap d' = 0$.]
\mn
\begin{enumerate}
\item[$(*)_{11}$]  If ${\bold f}(c) \in \gB_\beta$ and 
$c R_\beta b$ then $b\le {\bold f}(c)$.
\end{enumerate}
\mn
[Why?  If not, $b- {\bold f}(c)> 0$ in ${\bold b}_2[\gB_{\lambda^+}]$, so as
 the embedding ${\bold f}$ is complete and ${\bold b}_1 [{\gB}_{\lambda^+}]$
is atomic, for some $e \in {\bold b}^{\at}_1[\gB_{\lambda^+}], 
0 <_{{\bold b}_1} e$ and we have (in this formula 
${\bold f}(e),\bold f(c)$ are members of ${\bold b}_2[\gB_{\lambda^+}]$

\[
{\bold b}_2[\gB_{\lambda^+}] \models ``\bold f(e) \cap (b-\bold
f(c))>0".
\]

\mn
Now if $e \cap c > 0$ then $c=e$ so ${\bold f}(e) {\bold f} (c)$ hence 
${\bold f}(e) \cap (b- {\bold f}(c)) = 0$, contradiction.

So clearly $e \cap c =0$ in ${\bold b}_1[\gB_{\lambda^+}]$, so as $c \in
B^3_{\alpha,\beta} \subseteq B^1_{\alpha,\beta}$ there is
$d \in {\bold b}_1[\gB_{\alpha}]$ such that $c \le d,e \cap d= 0$
(in ${\bold b}_1[\gB_{\lambda^+}]$), hence ${\bold f}(c) \le \bold
f(d),\bold f(e) \cap \bold f(d) =0$ in ${\bold b}_2[\gB_{\lambda^+}]$,
and by the definition of $R_\beta$, as $c R_\beta b$ holds. 

Also $b \le {\bold f}(d)$, hence $b \cap {\bold f}(e)=0$
in ${\bold b}_2$ so $b-{\bold f} (e) =0$ in ${\bold b}_2[!{\gB}_{\lambda^+}]$ 
so we have gotten a contradiction to the choice of $e$.
Hence $b \le {\bold f}(c)$.]
\mn
\begin{enumerate}
\item[$(*)_{12}$]  $b= {\bold f}(c)$ if $c R_\beta b$ and 
$(\gB_\beta,{\bold f}) \prec (\gB_{\lambda^+},{\bold f})$.
\end{enumerate}
\mn
[Why?  If not, by $(*)_{11}$ above we have $b< {\bold f}(c)$ so by 
the ``${\bold f}$ is a \underline{complete} embedding" 
for some $e$ we have $\gB_{\lambda^+}\models$ ``$e \dot e 
{\bold b}^\at_1 \, \& \, {\bold f}(e) \cap ({\bold f}(c)-b)> 0_{{\bf b}_2}$". 
So as $e$ is an atom of ${\bold b}_1 [{\gB}_{\lambda^+}]$ clearly
$e \le_{{\bold b}_1} c e \cap c = 0_{{\bold b}_1}$. First assume 
$e \cap c = 0_{{\bf b}_1}$ hence $0_{{\bold b}_2}=
{\bold f}(e \cap c) = {\bold f}(e) \cap {\bold f}(2) \ge \bold f(e) 
\cap ({\bold f}(c)-b) > 0_{{\bf b}_2}$ contradiction. 
Second assume if $e \le_{{\bold b}_1} c$ recall that $c$
$\in B^3_{\alpha,\beta} \subseteq {\bold b}^{\at}_1$ [$\gB_\beta] 
\smallsetminus \gB_\alpha$ so $e=c$ contradiction to $c R_\beta b$ as the later
implies $b \ne 0_{{\bold b}_2}$].

We can conclude
\mn
\begin{enumerate}
\item[$(*)_{13}$]  if $\beta\in (\beta^*_\alpha, \lambda^+)$ and 
${\bold f} \rest {\bold b}_1 [{\gB}_\beta]$ is a complete embedding of 
${\bold b}_1[\gB_\beta]$ into ${\bold b}_2[\gB_\beta]$ (which occurs for a club
of $\beta$'s) then ${\bold f} \restriction B^3_{\alpha, \beta} \in \cP_\beta$.
\end{enumerate}
\mn
What have we gained compared to $\otimes_10$ in the end of stage E?
Here this works for all confinalities (for a club of $\beta$-s). 
Let $\bar c = \langle c_i:i < \omega \rangle$ be an increasing sequence for
${\bold b}^{\at}_1$ and let $\cC^3_{\bold f} = \bar{\cC}_{\bold f,\bar c}
= \{\delta: \delta< \lambda^+$ is a limit ordinal $> \beta^*_\alpha$ and 
${\bold f}$ is a complete embedding of $\bold b_1[\gB_\delta]$ into
$\bold b_2[\gB_\beta]$ and if $\gB_\delta \models ``d \dot e {\bold
  b}_1$ and $|d|/|a| \ge c_i"$  for every  $i < {\lg}(\bar c)$ 
then there is $e \in {\bold b}^\at_1[\gB_\delta],e \le d$ such that $e
\in B^3_{\alpha,\beta}\}$.

Clearly it is a club of $\lambda^+$. Now we can note more (but shall not
use it below).
\mn
\begin{enumerate}
\item[$(*)_{14}$]  for every $c$, $\bar c^1$, $\bar c^2$ as in \ref{c11}
for some club $\cC^4_{\bold f} \subseteq \cC^3_{\bold f} \cap 
\cC^1_{\bold f}$ of $\lambda^+$ we have
\sn
\item[$(*)_{15}$]  for $\delta \in \cC_{\bold f}$ there is 
${\bold f}^2_\delta \in {\cP}_\delta$ which is a function from 
${\bold b}_1[\gB_\delta]$ into ${\bold b}_2[\gB_\delta]$ and $e \in
\bold b_1[\gB_\delta] \Rightarrow \bold f(e) \vartriangle 
\bold b^2_\delta(e) \in \cJ^\delta_{\bold f,\bar c}$ where
\[
{\cJ}^\delta_{{\bold f},\bar c} = \{\bold f(b):\bold b_1[\gB_\delta] \models
``b \subseteq a \text{ has } \le c_i|a| \text{ elements" for some }
i< \lg(\bar c)\}.
\]
\end{enumerate}
\mn
Note that now $\bold p,\bar c^2$ disappear, as by
\ref{c80} such $\bold p,\bar c^2$ exists.
\medskip

\noindent
\underline{Stage G:}  Assume
\mn
\begin{enumerate}
\item[$\otimes$]   ${\bold b}_1$, ${\bold b}_2$, $\alpha$, $\bar c$,
${\bold f} \restriction {\bold b}^\at_1[\gB_\alpha]$, and ${\bold f}
\rest B^{3,{\bold f}}_{\alpha, \beta} \in {\cP}_\beta$
for $\beta \in \cC^3_{\bf f}$ as above.
\end{enumerate}
\mn
We have a pre-killer contractor $\zeta_{\pk}$ such that: for stationarily
many (even for every) $\delta \in \cC^3_{\bold f} \cap 
W_{\zeta_{\pk}}$ of cofinality $\lambda$,
(for each candidate ${\bold f}$ for ${\bold f} \rest B^3_{\alpha,\delta})$ for
stationarily many $\epsilon \in S_{\zeta_{\pk}}$, if possible we ensure
that for some $j \in u^\delta_{i_{\delta,\dot e}} \setminus
\{\max(u^\delta_{i_{\delta,\dot e}})\}$, we have
$\Gamma^{\delta+1}_j = \Gamma^{\ids}_{\langle
  (d_{\delta,j},d_{\delta,j,n}):
n<\omega\rangle}$ so $\bar c^\delta_j = \langle (d_{\delta,j},
d_{\delta,j,n}):n<\omega\rangle$
where $d_{\delta,j} \in B^{3}_{\alpha,
\delta}\subseteq {\bold b}^\at_1[\gB_\delta],d_{\delta,j,n} \in
{\bold b}_2[\gB_\delta],\bold f(d_{\delta,j}) = d_{\delta,j,0}$ and $\langle
d_{\delta,j,n}:n < \omega\rangle$ is indiscernible over
$A^\delta_j \cup \{d_{\delta, j}\}$ with $d_{\delta,j,0} \ne
d_{\delta,j,1}$, of course, and \wilog \, $\neg (d_{\delta,j,0}
\le_{{\bold b}_2} d_{\delta,j,1})$ as we can find $d_{\delta,j,-1},
d_{\delta,j,-2},\ldots$ \st \, $\langle d_{\delta,j,n}:n \in \bbI\rangle$
is an indiscernible over $A^\delta_j \cup \{d_{\delta,j}\}$.
Later the automorphism killer contractor $\zeta_{\ak}$ is active for
our case for stationarily many $\beta\in \cC^3_{\bold f} \cap
W_{\zeta_{\ak}}$ with $\cf(\beta)=\theta$ (so ${\bold f} \rest 
\bold b_1[\gB_\beta]\in {\cP}_\beta$ by $\otimes_{10}$, 
of course,  the automorphism killer contractor deal there 
with all such candidates as he does not know which
is really necessary). For stationarily many $\epsilon \in S_{\zeta_{\ak}}
\subseteq \lambda$ he ensure for  $\delta\in \beta \cap 
\cC^3_{\bold f} \cap W_{\zeta_{\pk}}$ as above,
if possible, that for some $j \in u^\beta_{i_{\beta,\dot e}} \setminus
\{\max(u^\beta_{i_{\beta,\dot e}})\}$
and $j_1< \lambda$,  $d_{\delta, j_1} \notin
\acl_{\gB_\beta}(A^\beta_j)$, $\langle d_{\delta, j_1, n}:
n<\omega\rangle$ is indiscrenible over $A^\beta_j+ d_{\delta, j_1}$ and
$\Gamma^\beta_j = \Gamma^{\ids}_{\bar c^\delta_j}$.
He choose $e^\beta_j\in {\bold b}_1^\at [\gB_\beta]$ such that
${\bold b}^\at_1 [{\gB}_\beta] \models ``e^\beta_j \cap d_{\delta,j}
=0"$ and $({\bold f}\restriction {\bold b}^\at_1[\gB_\beta])(e^\beta_j)
\cap d_{\delta,j,1} > 0_{{\bold b}_2}$ (if there is no one, his candidate for
${\bold f} \restriction {\bold b}_1[\gB_\beta]$ is faked, failing coming from 
a complete embedding, so can be forgotten) and he
make $e^\beta_j \in A^\beta_{\min(u^\beta_j \setminus (j+1))}$. Then
he let $\Omega^\beta_n= \Gamma^{\av}_{\dot D,\langle m: m<\omega\rangle}$
(for $n<\omega$), $\Omega^\beta_\omega = \Gamma^{\ms}_{{\bold b}^{\at}_1,
\langle 1/x_{\beta,n}:n<\omega\rangle}$ be as in stage D  (alternatively
$a_{\beta, n} > n/|{\bold b}^\at_1|$ for $n<\omega,a_{\beta, n}< b$ if $b\in
\acl_{\gB_{\beta+1}}(\gB_\beta\cup\{a_{\beta, m}: m<\omega\})$ and
$n<\omega \Rightarrow n/|{\bold b}^\at_1| <b$ (check in stage D)), and he
demand:
\mn
\begin{enumerate}
\item[$(*)_{16}$]  for a club of $\epsilon<\lambda$, for $j$ as above:
\sn
\begin{enumerate}
\item[$(\alpha)$]  $d_{\delta,j} \le_{{\bold b}_1} a_{\beta, \omega}\,
  \& \, e^\beta_j\cap a_{\beta, \omega} =0_{{\bold b}_1}$ in ${\bold b}_1 
[{\gB}_{\beta_1}]$
\sn
\item[$(\beta)$]  $\langle (d_{\delta,j},d_{\delta,j,n}):n<\omega\rangle$ is
indiscernible over $A^{\beta+1}_j$.
\end{enumerate}
\end{enumerate}
\mn
No problem and ${\bold f}(a_{\beta, \omega})$ will give a contradiction.

So
\mn
\begin{enumerate}
\item[$(*)_{17}$]  for every $\delta \in \cC^3_{\bold f} 
\cap W_{\zeta_{\pk}}$ of cofinality $\lambda$, for some
$j<\lambda$, the pre-killer contractor cannot choose $d_{\delta, j},
d_{\delta,j,n}$ ($n<\omega$) as above. Which means: if $d \in B^{3}_{\alpha,
\lambda^+} \setminus \acl_{\gB_\delta}(A^\delta_j)$ then ${\bold f}(d)
\in \acl_{\gB_\delta}(A^\delta_j +d)$ which is equal to $\dcl_{\gB_\delta}
(A^\delta_j+d)$ as $T^*$ has Skolem functions.
\end{enumerate}
\mn
Now \wilog \, $\lambda^2$ divides $\delta$ hence there is $\delta' <
\lambda^+$ such that $\cf(\delta')= \lambda$ and $A^\delta_j \subseteq
\gB_\delta'$ hence for some $j'< \lambda$, $A^\delta_j \subseteq
A^{\delta'}_{j'}$ and of course $A^{\delta'}_{j'} \in
\bigcup\limits_{\beta< \delta} {\cP}_\beta$. 

So using Fodor lemma 
\mn
\begin{enumerate}
\item[$(*)_{18}$]  for some $A^*$ for stationarily many $\delta$ as above 
satisfying $A^*\in \cup \{\cP_\gamma:\gamma < \delta\}$, 
in $(*)_{17}$ we can replace $A^\delta_j$ by $A^*$.
\end{enumerate}
\medskip

\noindent
\underline{Stage H}:

We let the automorphism killer contractor act also for $\breve{f} \in A^*$
$\beta\in \cC^3_{\bold f} \cap W_{\zeta_{\ak}}$ of cofinality $\lambda$ for
stationarily many $\epsilon  \in S_{\zeta_{\pk}}$, 
where $\gB_\beta \models ``\breve{f}$ a partial function from 
${\bold b}^\at_1$ to ${\bold b}_2"$ to ensure that for some
$j \in u^\beta_{i_{\beta,\epsilon}},j<
\max(u^p_{i_{\beta,\epsilon}})$, he choose in ${\gB}_\delta$, if possible 
$\langle d_{\beta, j, n}: n<\omega\rangle$, indiscernible over
$A^\delta_j + \breve{f}$ such that $d_{\beta,j,0} \ne d_{\beta, j, 1} \in
B^{3}_{\alpha, \beta}$, ${\bold f}(d_{\beta,j,0}) =
\breve{f}(d_{\beta,j \ge})$ and ${\bold f}(d_{\beta,j,1}) \ne 
\breve{f}(d_{\beta,j,1})$ and are, of course, well defined. 
If so \wilog \, add $\breve{f}(d_{\beta,j,1}) \notin 
\{\dot f(d_{\beta,j,n}):n<\omega\}$, let
$\Gamma^\beta_j = \Gamma^{\ids}_{\bar c^\beta_j}$, $\bar c^\beta_j
=\langle (\breve{f}, d_{\beta, j, n}):
n<\omega\rangle$ and get contradiction as above.

Now easily (possibly shrinking $\cC^3_{\bf f}$, using the freedom in choosing
a type for $\zeta_{\ak}$; using that $T$ has Skolem functions)
\mn
\begin{enumerate}
\item[$(*)_{19}$]   $\delta \in \cC^3_{\bold f}$, and $A\in
[\gB_\delta]^{<\lambda}\cap \cup\{{\cP}_\alpha:\alpha<\delta\}$ 
contains the relevant parameters and $d^1,d^2 \in
B^{3}_{\alpha,\delta}$ realizes the same non-algebraic type
over ${\acl}_{{\gB}_\delta} (A^*)$, \then \, for some $d^* \in B^{3}_{\alpha,
\delta}$ we have that for $\ell=1,2$ there is an infinite 
indiscernible sequence to which $d^*,d^\ell$ belong.
\end{enumerate}
\mn
Hence (together with stage G)
\mn
\begin{enumerate}
\item[$(*)_{20}$]  for some $\alpha < \lambda^+$ and $A^*\in 
\cP_{\lambda^+} \cap [\gB_{\lambda^+}]^{<\lambda}$,
(\wilog \, $A^* = \acl_{\gB_{\lambda^+}}(A^*) \prec \gB_{\lambda^+}$ and of
course $A^* \in {\cP}_{\lambda^+} = 
\bigcup\{{\cP}_\gamma:\gamma<\lambda^+\}$),
for every $\Gamma^{\ms}_{{\bold b}^{\at}_1,\bar c}$--big type 
$p\in {\bold S}(A^*,\gB_{\lambda^+})$ from ${\cP}_{\lambda^+}$
such that $[x \dot e {\bold b}^{\at}_1]\in p$ for some 
$\breve{f}_p \in A^*$ we have
\sn
\begin{enumerate}
\item[$(\alpha)$]  $\gB_{\lambda^+} \models ``\breve{f}_p$ is a (partial)
function from ${\bold b}^\at_1$ to ${\bold b}_2"$
\sn
\item[$(\beta)$]  for every $d \in B^{3}_{\alpha,\lambda^+}$ realizing
$p$, we have $\breve{f}_p(d) = {\bold f}(d)$
\sn
\item[$(\gamma)$]  $x \ne y \wedge x \dot e \Dom(\breve{f}) \subseteq 
{\bold b}_1^{\at} \wedge y \dot e \Dom(\breve{f}) \subseteq {\bold b}^\at_1
\Rightarrow \breve{f}_p(x) \cap \breve{f}_p(y) = 0_{{\bold b}_2}$ or
at least
\sn
\item[$(\gamma)'$]  for some $i$ and $\psi(y,x)$, \st \, 
$\forall x [|\{y \dot e {\bold b}^{\at}_1:\psi(y,x)\}| \le c_i \times 
|{\bold b}^{\at}_1|]$ we have $x \ne y \wedge x \dot e \Dom(f)
\wedge y \dot e \Dom(\dot{f}) \subseteq {\bold b}^{\at}_1$.
\end{enumerate}
\end{enumerate}
\mn
[Why? First note that if $b \in \acl_{{\gB}_{\lambda^+}} (A^*+c), A^*$
as above then for some $\breve{f} \in A^*,\gB_{\lambda^+} \models ``\breve{f}$
is a partial function and $\breve{f}(c)=b"$; also if $c \in {\bold b}^\at_1
[{\gB}_{\lambda^+}],b \in {\bold b}_2 [{\gB}_{\lambda^*}]$ \wilog \,
${\gB}_{\lambda^+} \models ``\breve{f}$ is a function from ${\bold b}_1^\at$
into ${\bold b}_2$''; this takes care of clause ($\alpha$). Second by the first
paragraph of this stage we get clause $(\beta)$. 

Third,  concerning clause $(\gamma)$ we get $p(x) \cup p(y) 
\cup\{\neg \psi(y,x,\bar a):\bar a \subseteq A$, and for some $i <
\lg(\bar c)$ we have ${\gB}_{\lambda^+} \vdash |\{v \dot e 
\bold b^{\at}_1:\psi(b,x,\bar a)\}| \le c_1 x |\bold b^{\at}_1|\} 
\cup \{\breve{f}_p (x) \cap \breve{f}_p (y) \ne 0_{\bold b_2}$ in 
${\bold b}_2,x \ne y\}$ belong to ${\cP}_{\lambda^+}$ and is not 
realized in $\gB_{\lambda^+}$ hence by the saturator 
we can change $\breve{f}_p$ to make clause $(\gamma)$ true. Now we can extend
$\breve{f}_p$ by giving the value $0_{{\bold b}_2}$ for $x \dot e 
\bold b^{\at}_1$ on which it is not defined,
so we can make $\Dom(\breve{f}_p) = {\bold b}^{\at}_1$].

\Wilog \, $c_a \times |{\bold b}^{\at}_1|$ is an integer.  So for some 
$\bar a,\gB_{\lambda^+} \models ``\langle b_{\dot n}:\dot n < c_1 \times
|{\bold b}^{\at}_1|\rangle$ is a partition of ${\bold b}^{\at}_1$, and 
for each $\dot n,a' \ne a'' \wedge a'' \dot e b_{\dot n} \Rightarrow
\breve{f}_p(a') \cap \breve{f}_p(a'') = 0_{\bold b_2}"$.

Now we can find for each $\breve{f}_p$ (for $p \in \bold S(A^*,
\gB_{\lambda^+}) \cap {\cP}_{\lambda^+}$ which is 
$\Gamma^{\ms}_{{\bold b}^{\at}_1,\bar c}$-big), in $\gB_{\lambda^+}$ a
``finite" sequence $\langle c^{\ell,p}_k:k < k^*_p,\ell<2\rangle$ such that (in
$\gB_{\lambda^+}$), $2k^*_p \le |{\bold b}^\at_1|$ and $c^{\ell,p}_k
\dot e {\bold b}^\at_1$ are pairwise distinct and $\breve{f}_p(c^{0,p}_k) \cap 
\breve{f}_p(c^{1,p}_k)>0_{{\bold b}_2}$ and

\[
(\forall x,y \dot e {\bold b}_1^{\at})(x \ne y \, \& \, x, y \notin
\{c^{\ell,p}_k: k< k^*_p, \ell<2\} \rightarrow \breve{f}_p(x) \cap 
\breve{f}_p(y) =0_{{\bf b}_2}).
\]

\mn
Let $d^\ell_p = \bigcup \{c^{\ell, p}_k: k< k^*_p\} \dot e {\bold b}_1$ for
$\ell=1,2,d^2_p = 1_{{\bold b}_1} - d^0_p - d^1_p$,
so $d^0_p$, $d^1_p$, $d^2_p$ are pairwise disjoint and let $e^\ell_p =
{\bold f}(d^\ell_p) \dot e {\bold b}_2$ and $e^{\ell,p}_k = 
\breve{f}_p (c^{\ell,p}_k)$ (all in $\gB_{\lambda^+}$'s sense). 

Lastly, let

\[
d^{\ell,*}_p = \bigcup \{c:c \dot e {\bold b}^{\at}_1 \text{ and } c 
\le_{{\bold b}_1} d^\ell_p \text{ and } \breve{f}_p(c) \le e^\ell_p\}
\]

\mn
so $d^{\ell,*}_p \dot e {\bold b}_1$.

Now for $\ell=0, 1, 2$
\mn
\begin{enumerate}
\item[$(*)_{21}$]  $p(x) \cup [x \le d^\ell_p - d^{\ell,*}_p]$ is
$\Gamma^{\ms}_{\bar c}$-small.
\end{enumerate}
\mn
[Why?  Clearly this type belong to ${\cP}_{\lambda^+}$ hence if $(*)$
fail this there is $b \in B^{3}_{\alpha,\lambda^+}$
realizing it so ${\bold f}(b) = \breve{f}_p(b)$ and as ${\bold f}$ 
is a homomorphism we get contradiction.]
\mn
\begin{enumerate}
\item[$(*)_{22}$]  for every $k< k^*_p$, we have 
$(e^0_k \le_{{\bold b}_1} (d^0_p - d^{0, *}_p)) \, \vee \, 
(e^1_k \le_{{\bold b}_1} (d^1_p - d^{1, *}_p))$.
\end{enumerate}
\mn
[Why?  As $f_p(e^0_k) \cap f_p(e^1_k)> 0_{{\bold b}_1}$ whereas 
${\bold f}(d^0_p) \cap F(d^1_p) =0_{{\bf b}_2}$.]

Hence [just proving that (in $\gB_{\lambda^+} \vee \kappa^*_p$ is
``small"]
\mn
\begin{enumerate}
\item[$(*)_{23,\ell}$]  For some $\vartheta_p(x) \in p(x)$ we have
$\vartheta_p(x) \, \& \, [x \le_{{\bold b}_1} (d^0_p \cup d^1_p) = 
\bigcup \{e^\ell_k:\ell< 2, k< k^*_p\}]$ is $\Gamma^{\ms}_{\bar c}$-small.
\end{enumerate}
\mn
Let $a_p =\{x \dot e {\bold b}_1^\at: \vartheta(x)\}$ and $a^-_p = \{
x \dot e {\bold b}^\at_1: \dot e(x) \, \& \, x \le_{{\bold b}_1} d^0_p \cup
d^1_p\}$ and let $a^+_p = \{x \dot e {\bold b}^\at_1: \vartheta(x)
\, \& \, \neg x \le d^0_p \cup d^1_p\}$ so $a_p$ is the disjoint union
of $a^-_p$ and $a^+_p$. Note that $\langle c^{\ell,p}_\kappa:\ell < 2,
\kappa< \kappa^*_p\rangle,\langle e^\ell_k: \ell, k\rangle,
d^\ell_p,d^{\ell,*}_p,a^-_p,a^+_p$ depend on $\breve{f}_p$ 
(and not on $p$) and $\breve{f}_p \in A^*$ so \wilog \, $\in A^*$ for every 
$\Gamma^{\ms}_{\bar c}$-big $p \in \bold S(A,\gB_{\lambda^+})
\cap {\cP}_{\lambda^+}$ we have $[x \le_{{\bold b}_1} a^-_p] \notin p$ 
and for some $j(p) < \lg(\bar c)$ we have ``$|a^-_p| \le c_{j(p)}
\times |{\bold b}^\at_1|$ and $\breve{f}_p \restriction a_p$ maps
${\bold b}^{\at}_1$ into ${\bold b}_2,x \ne y \in a^+_p 
\Rightarrow \breve{f}_p(x) \cap \breve{f}_p(y) = 0_{{\bold b}_2}$ 
and so $\breve{f}_p$ induce an embedding of ${\bold b}_1 \restriction
a^+_p$ into ${\bold b}_2$ called $\breve{f}_p"$.  (We are identifying 
${\bold b}_1$ with ${\cP}({\bold b}^\at_1)$ where $a^+_p := \cup \{x
\dot e {\bold b}^{\at}_1:\breve{f}_p(x)>0_{{\bold b}_2}$ in ${\bold b}_2$-s
sense$\}$ satisfy $a^+_p \dot e \bold b_1$ so for some 
$\breve{f}_p$ we have $\gB_{\lambda^+} \models$
``$\breve{f}_p$ is an embedding of ${\bold b}_1 \rest a^+_p$ into 
${\bold b}_2, {\bold b}_1 \rest a^+_p$ is the sub-boolean ring of 
${\bold b}_1$ with set of elements $\{d \in {\bold b}_1:d\le
a^+_p\}"$.
\medskip

\noindent
\underline{Stage I}:  Let $c^*\in \gB_{\lambda^+}$ be such that 
$c_i < c^* \dot e (0,1)_{\bbR},c^*<1/n$ and $c^* |{\bold b}^{\at}_1|$ 
integer $>n$ for every base $n$.
Let $\dot k^* \in \bbN[{\gB}_{\lambda^+}]$ be small enough.

In $\gB_{\lambda^+}$ let $\langle \breve{f}_{\dot n}:\dot n \le \dot
k^* \rangle$ be a list of functions from ${\bold b}^{\at}_1$ 
into ${\bold b}_2$ satisfying $x \in {\bold b}_1^{\at}\, \& \, y \in 
{\bold b}^\at_1\, \& \, x \ne y \Rightarrow \breve{f}_{\dot n}(x) \cap 
\breve{f}_{\dot n}(y) = 0_{{\bold b}_2}$ including all such 
members of $A^*$ (exist by the saturator work). Let $e \in
\gB_{\lambda^+}$ be an equivalence relation on ${\bold b}^{\at}_1$ with
$< {\dot k}^*$ equivalence classes such that
$\gB_{\lambda^+} \vDash ``x e y"$ implies $\tp(x,A^*,\gB_{\lambda^+}) 
= \tp(y,A^*,\gB_{\lambda^+})$.

Clearly 
\mn
\begin{enumerate}
\item[$(*)_{24}$]  if $\dot n_1,\dot n_2 < \dot k^*,\dot m_1 e \dot m_2$ then 
$\breve{f}_{\dot n_1} (\dot m_1) \le \breve{f}_{\dot n_2} (\dot m_1)
\Leftrightarrow \breve{f}_{\dot n_2} (\dot m_2) \le \breve{f}_{\dot n_2} 
(\dot m_2)$
\sn
\item[$(*)_{25}$]  if $a' \ne a'', a' e a''$ and $\dot n < \dot k^*$ then 
$\breve{f}_{\dot n}(a') \cap \breve{f}_{\dot n}(a'') = 0_{{\bold b}_2}$.
\end{enumerate}
\mn
In $\gB_{\lambda^+}$ let $b^* = \{x \dot e {\bold b}^\at_1: |x/e|
\ge {\bold k}^*_1\}$, where e.g. ${\gB}_\lambda = ``\dot m k^*_1=
2 \dot m k^*"$ and we can choose $\langle a_{\dot m}:\dot m < 
\dot m^*\rangle$ ``randomly in ${\gB}_{\lambda^+}$ sense such that:
\mn
\begin{enumerate}
\item[$(\alpha)$]  $\dot m_1 < \dot m_2 < \dot k^*_1 \Rightarrow  
\neg (a_{{\dot m}_1} e a_{{\dot m}_2})$
\sn
\item[$(\beta)$]  $a_{\dot m} \dot e b^*,b^* = 
\bigcup\limits_{\dot m < \dot k^*_1} (a_{\dot m}/e)$.
\end{enumerate}
\mn
So almost surely 
\mn
\begin{enumerate}
\item[$(\gamma)$]  if $\dot m_1 < m^*,\dot n_1,\dot n_2 < \dot k^*$
 and for every $\dot n_2 < \dot k^*$, if $\breve{f}$ 
satisfies $\Dom(\breve{f}) = \{0,\ldots,\dot m^*_1\} \setminus \{\dot
m_1\},(\forall \dot m) \breve{f}(\dot m) \subseteq \{0,\ldots \dot k^*\}$ 
and $b =\cup\{\breve{f}_{\dot n}(a_{\dot m})\dot m \in
\Dom(\breve{f}),\dot n \in \breve{f}(\dot m)\} \dot e {\bold b}_2$ and 
for every a $\dot e(a_{\dot m}/e) \setminus \{a_{\dot m}\}$ 
we have $\breve{f}_{\dot n_1} (a) - \breve{f}_{\dot n_2}(a)$ 
is disjoint to $b$ (in ${\bold b}_2$), \then \, $a_{\dot m_1}$ satisfies this.
\end{enumerate}
\mn
So 
\mn
\begin{enumerate}
\item[$(*)_{27}$]   for each ${\dot m} < {\dot m}^*$, for some $\dot n <
  \dot k^*$ we have $\breve{f}(a_{\dot m}) = \breve{f}_{\dot
    n}(a_{\dot m})$.
\end{enumerate}
\mn
Hence by $(*)_{27}$ we have $a \in
(a_{\dot m}/e) \Rightarrow \breve{f}(a) = \breve{f}_{\dot n}(a)$. 
Our next aim is to show that the choice of ${\dot n}$ can be done uniformly: 
in a way represented in ${\gB}_{\lambda^+}$.
For each ${\dot m}_1 \ne \dot m_2 < \dot m^*$ and $\dot n_1,\dot n_2 < 
\dot k^*$ choose if possible a member $b = b_{\dot n_1,\dot n_2,
\dot m_1,\dot m_2} \dot e(a_{{\dot m}_2}/ e)$ such that
$\breve{f}_{{\dot n}_1} (a_{{\dot m}_1}) \cap \breve{f}_{{\dot n}_2} (b) 
\ne 0_{{\bold b}_2}$,  otherwise we let 
$b_{{\dot n}_1,{\dot n}_2,{\dot m}_1,{\dot m}_2}$ be $0_{{\bold b}_1}$; 
so of course, \wilog \, the function
$({\dot n}_1,{\dot n}_2,{\dot m}_1,{\dot m}_2) \mapsto 
b_{{\dot n}_1,{\dot n}_2,{\dot m}_1,{\dot m}_2}$ is represented in 
$\gB_{\lambda^+}$. 

Let $a^1 =: \{a_{\dot m}:{\dot m} < {\dot m}^*\}$ and $a^2 =: 
\breve{f}(a^1) \in {\bold b}_2 [\gB_{\lambda^+}].$
Now by stage H for each ${\dot m}_1 < {\dot m}^*$ for some 
${\dot n}_0 < {\dot k}^*$ we have $\breve{f}(a_{{\dot m}_1}) =
\breve{f}_{\dot n_0} (a_{{\dot m}_1})$ hence $\breve{f}_{\dot n_0}
(a_{{\dot m}_1}) = \breve{f}(a_{{\dot m}_1}) \le \breve{f}a^1) = a^2$, 
also we have
\mn
\begin{enumerate}
\item[$(*)^{28}_{\dot n_1,\dot m_1}$]  $a \in a_{{\dot m}_1} / e 
\Rightarrow \tp(a,A^*,\gB_{\lambda^+})= \tp(a_{{\dot m}_1},A^*,
\gB_{\lambda^+}) \Rightarrow \breve{f}(a) = \breve{f}_{\dot n_0}(a)$.
\end{enumerate}
\mn
We like to define $\dot n_0$ from $a$, or just from $a_{{\dot m}_1}$ (inside
$\gB_{\lambda^+})$. 

Clearly 
\mn
\begin{enumerate}
\item[$(*)_{29}$]  $\breve{f}_{{\dot n}_0} (a_{{\dot m}_1}) \le_{{\bold b}_2}
a^2,\dot n_0 < \dot k^*$ and $a \in a_{\dot m_1}/e \setminus 
\{a_{\dot m_1}\} \Rightarrow \breve{f}_{\dot n_0}(a) \cap a^2=
0_{{\bold b}_2}$ 
\sn
\item[$(*)_{30}$]  for no ${\dot n}_1 < \dot k^*,\dot n \ne \dot n_0$
  do we have 
\[
\dot f_{{\dot n}_1} (a_{{\dot m}_1}) \le_{{\bold b}_2}
 a^2 \, \& \, \neg \big(\breve{f}_{{\dot n}_1} (a_{{\dot m}_1}) 
\le_{{\bold b}_2} \breve{f}_{{\dot n}_0} (a_{{\dot m}_1})\big).
\]
\end{enumerate}
\mn
[Why?  Assume ${\dot n}_1$ is a counterexample,
necessarily by ${\bold f}$ being a complete embedding using the maximal 
antichain ${\bold b}_1^{\at} - \bigcup\limits_{a \in b^*} a\} \cup
\{a:a \in b^* \}$ of ${\bold b}_1 [{\gB}_{\lambda^+}]$,
for some $\dot m_2 < \dot m^*$ and $a''_{\dot m_2} \dot e 
a_{\dot m_2} / e$ we have $\breve{f}_{{\dot n}_1} (a_{{\dot m}_1}) 
- \breve{f}_{{\dot n}_0} (a_{{\dot m}_1})$ computed in ${\bold b}_2
[\gB_{\lambda^+}]$ is not disjoint to
${\bold f}(a''_{{\dot m}_2})$, but for some ${\dot n}_2 < \dot m^*$ we have
$\breve{f}(a_{{\dot m}_2}) = \breve{f}_{{\dot n}_2} (a_{{\dot m}_2})$, 
so necessarily $b_{{\dot n}_1,{\dot n}_2,{\dot m}_1,{\dot m}_2}
\ne 0_{{\bold b}_1}$ and is a member of ${\bold b}_1$ disjoint to $a^1$ so
${\bold f}(b_{{\dot n}_1,{\dot n}_2,{\dot m}_1,{\dot m}_2})$ is
disjoint (in ${\bold b}_2)$ to $a^2= {\bold f}(a^1)$.

Also as $b_{{\dot n}_1,{\dot n}_2,{\dot m}_1{\dot m}_2} 
\dot e a_{{\dot m}_2} / e$ and ${\bold f} (a_{{\dot m}_2}) =
\breve{f}_{{\dot n}_2} (a_{{\dot m}_2})$ clearly 
${\bold f}(b_{{\dot n}_1,{\dot n}_2,{\dot m}_1,{\dot m}_2}) = 
\breve{f}_{{\dot n}_2}(b_{{\dot n}_1,{\dot n}_2,{\dot m}_1,
{\dot m}_2})$ so by the previous sentence $\breve{f}_{{\dot n}_2} 
(b_{{\dot n}_1,{\dot n}_2,{\dot m}_1,{\dot m}_2}) \cap a^2
=0_{{\bold b}_2}$, but by the choice of $b_{{\dot n}_1,{\dot n}_2,
{\dot m}_1,{\dot m}_2} (> 0_{{\bold b}_1})$, i.e. the previous 
sentence we have $f_{{\dot n}_2} (b_{{\dot n}_1,{\dot n}_2,{\dot m}_1,
{\dot m}_2}) \cap \breve{f}_{{\dot n}_1} (a_{{\dot m}_1})
> 0_{{\bold b}_2}$ hence $\neg(\breve{f}_{{\dot n}_1} (a_{{\dot m}_1})
\le_{{\bold b}_2} a_2)$, contradicting the choice of ${\dot n}_1$].

So
\mn
\begin{enumerate}
\item[$(*)_{31}$]  in $\gB_{\lambda^+}$, if $x \dot e {\bold b}^\at_1
  \, \& \, x \le_{{\bold b}_1} {\bold b}^*$ then ${\bold f}(x)$ is 
$\breve{f}_{{\dot n}(x)}(x)$ where ${\dot n}(x)$ is the unique 
${\dot n} < {\dot k}^*$ such that
\begin{equation*}
\begin{array}{clcr}
(\exists {\dot m} < {\dot m}^*)(x \dot e a_{\dot m}/ e \, \& \, 
\breve{f}_{\dot n} (a_{\dot m}) &\le_{{\bold b_2}} a^2 \, \& \, 
 (\forall \dot n' < \dot n)[\breve{f}_{\dot n'},(a_{\dot m}) \\
  &\le a^2 \rightarrow \breve{f}_{\dot n} (a_{\dot m}) \le 
\breve{f}_{{\dot n},1} (a_{\dot m})]),
\end{array}
\end{equation*}
\sn
so $x \mapsto  {\dot n}(x)$ is a function in $\gB_{\lambda^+}$.
\end{enumerate}
\medskip

\noindent
\underline{Stage J}:

Let

\[
a_\lev := \{b: b \text{ a natural number } <\log \log |{\bold b}^\at_1|\}\ \ \
\text{(in }\gB_{\lambda^+})
\]

\mn
though really we are interested just in

\[
a^-_\lev := \{b \dot e a_\lev: b > n \text{ for every (true) natural
  number \,(not in }\gB_{\lambda^+}!)
\]

\mn
for $b \dot e a_\lev$,  in $\gB_{\lambda^+}$ let

\begin{equation*}
\begin{array}{clcr}
I_b := \{\breve{f}: &\breve{f} \text{ a partial function from } 
\bold b^\at_1 \text{ into } {\bold b}_2 \setminus \{0_{{\bold b}_2}\} 
\text{ such that} \\
  &e_1 \ne e_2 \dot e {\bold b}^\at_1 \rightarrow \breve{f}(e_1) 
\cap \breve{f}(e_2) = 0_{{\bold b}_2} \text{ such that} \\
  &|{\bold b}^\at_1 - \Dom(\breve{f})| \le 2^b\}
\end{array}
\end{equation*}

\mn
($b$-th level). We define a distance function on $I= 
\bigcup\{I_b:b \dot e a_\lev\}$:

\[
\dis(\breve{f}_1,\breve{f}_2) = |\{x \dot e {\bold b}^{\at}_1:x \notin
\Dom(\breve{f}_1) \text{ or } x \notin \Dom(\breve{f}_2) 
\text{ or } \breve{f}_1(x) \ne \breve{f}_2(x)\}|.
\]

\mn
Next define a branch,  it is an (outside) function $\hat H$, 
satisfying $\Dom(\hat H) = a^-_\lev$, $\hat H(b) \in I_b,
b_1 < b_2$ (in $a^-_\lev$) $\Rightarrow \dis(\hat H(b_1),\hat H(b_2))
< 10^{b_2}$.

Now ${\bold f}$ induce a branch $\hat H$ as $\boxtimes$ in the 
end of stage I holds for $b \in a^-_{\lev}$.

By \ref{e24} below we have ``no undefinable branch" so there is an
equivalent branch $H'$ (see \ref{e24} below) which is definable 
in ${\gB}_{\lambda^+}$ hence is represented say by $\dot f$.

Let $\dot f_n = \dot f(n)$ ($n$ a true natural number in 
$\gB_{\lambda}$'s sense).

For each $b \dot e a^-_\lev$, for some $n_b<\omega$ we have
$\dis(\hat H(b),\hat H(n)) \le |10^b|$ for every $n \ge n_b$ 
(otherwise $\{n: n$ true natural number$\}$ is
definable in $\gB_{\lambda^+}$) so some $n^*$ is $n_b$ for
arbitrarily small $b \dot e a^-_\lev$ hence (changing $10^b$
slightly) \wilog \, this holds for every small enough $b \dot e
a^-_\lev$.

[Why?  As the cofinality of $(\{\dot n \in {\bbN}:n
<_{{\gB}_{\lambda^+}} {\dot n}$ for $n<\omega\},\ge)$ is uncountable
(in fact is $\lambda^+$ as $\Omega_{\alpha,0} = 
\Gamma^{\ver}_{\langle n:n<\omega\rangle}$ for unboundedly many $\alpha
<\lambda^+$].  In particular $|{\bold b}^\at_1 \setminus 
\Dom(\dot f_{n^*})|$ is really $<10^b$ for each $b \in a_{\ell ev}$ 
hence is finite. Assume $\{d \in {\bold b}^{\at}_1[\gB_{\lambda^+}]: 
d \dot e \Dom (\dot f_{n^*}),\ \dot f_{n^*}(d) \ne {\bold f}(d)\}$
is infinite.  So we can find $d^1_n,d^2_n \dot e {\bold b}_1^{at}$, 
pairwise distinct hence disjoint, satisfying ${\bold f}
(d^2_n) \cap \dot f_{n^*}(d^1_n) \ne 0_{{\bold b}_2}$. 
We can find $b \dot e {\bold b}_1[\gB_{\lambda^+}]$,
such that $d^1_n \le_{{\bold b}_1} b,d^2_n \cap b=0_{{\bold b}_1}$, let
$b' = {\bold f}(b)$, and in $\gB_{\lambda^+}$:

\[
b'' := \{d':d' \dot e {\bold b}^\at_1,d' \cap_{{\bold b}_1} b
= 0_{{\bold b}_1} \text{ but } \dot f_{n^*}(d') \cap b' \ne 0\},
\]

\mn
so necessarily $\gB_{\lambda^+} \models ``|b''|$ is $>n"$ for each
$n$, hence for some $c \in  a^-_{\lev}$ we have $\gB_{\lambda^+}\models$
``$|b''|> c"$ and we get easy contradiction.

Thus modulo \ref{e24} we have finished proving:
\mn
\begin{enumerate}
\item[$(*)$]  every complete embedding of ${\bold b}_1[\gB_{\lambda^+}]$ into
${\bold b}_2[\gB_{\lambda^+}]$ appear in $\gB_{\lambda^+}$ where
$\gB_{\lambda^+} \models ``{\bold b}_1$ is a finite Boolean ring hence algebra,
${\bold b}_2$ is a Boolean ring (e.q. algebra)".
\end{enumerate}
\medskip

\noindent
\underline{Stage K}:  I assume ${\gB}_{\lambda^+} \models 
``{\bold b}_1,{\bold b}_2$ are atomic Boolean rings" and 
${\bold f}$ is a isomorphism from ${\bold b}_1 
[{\gB}_{\lambda^+}]$ onto ${\bold b}_2 [{\gB}_{\lambda^+}]$, and let
${\cY} =\{a:{\gB}_{\lambda^+} \models ``a \in {\bold b}_1$
is $\ne 0_{{\bold b}_1}$ and is a finite union of atoms"$\}$. 
So for every $a \in {\cY}, {\bold f}$ is an isomorphism from Boolean Algebra
${\bold b}^1_a [{\gB}_{\lambda^+}]$ where ${\bold b}^a_1 := {\bold b}_1
\rest \{x:x \le_{{\bold b}_1} a\}$ onto ${\bold b}_2 [{\gB}_{\lambda^+}]$
where ${\bold b}^{{\bold f}(a)}_2 = {\bold b}_2 \rest \{y:y 
\le_{{\bold b}_2} {\bold f}(a)\}$, i.e. ${\bold f}_a = {\bold f}
\rest {\bold b}^a_1[{\gB}_{\lambda^+}]$ is in ${\gB}_{\lambda^+}$ 
hence by Stage C we are done. 
\end{PROOF}

\begin{remark}
We may try to replace the proof from $(*)_{27}$ till here by:
\mn
\begin{enumerate}
\item[$(*)$]  for $\dot m_1 < \dot m^*,\dot n_1 < \dot k^*$, we have $\bold f
(a_{\dot m_1})= \breve{f}_{\dot n_1} (a_{\dot m_1})$ \underline{if and
  only if} $(\dot m_1,\dot n_1)$ satisfies $\breve{f}_{\dot
  n_1}(a_{\dot m_1}) \le_{{\bold b}_2} a^2$ and if 
$\breve{f}_{\dot n_2} (a_{\dot m_1} \le_{{\bold b}_1} a^2$ and 
$\breve{f}) \dot m) \{\dot n < \dot k^*:\breve{f}_{\dot n} (a_{\dot m})
\le a^2\}$ when $\dot m < \dot m^*\, \& \, \dot m \ne \dot m^*$, 
then $\breve{f}_{\dot n_1} (a)- \breve{f}_{\dot n_2}(a)$ is 
disjoint to $a^2$ for every $a e a_{\dot m_1}$ and $a^2 =
\cup\{\breve{f}_{\dot n} (a_{\dot m}):\dot n \in \breve{f}'(\dot m)\}
\cup \breve{f}_{\dot n_1}(a_{\dot m_1})$.
\end{enumerate}
\mn
Hence playing with $\dot k^*$ we get: 
\mn
\begin{enumerate}
\item[$\boxtimes$]  for every ${\dot k}^*$ such that 
$\gB_{\lambda^+} \models ``n < \dot k^* <|{\bold b}^{\at}_1|$ 
for $n<\omega$, there is $b_{\dot k^*}^* \in {\bold b}_1$ such that
 ${\bold f} \rest \{a \dot e {\bold b}^\at_1:
a \cap b_{\dot k^*} = 0_{{\bold b}_1}\}$ is represented in 
${\gB}_{\lambda^+}$ hence ${\bold f} \rest \{a \dot e {\bold b}_1: 
a \cap b_{{\dot k}^*} = 0_{{\bold b}_1}\}$ is represented in 
${\gB}_{\lambda^+}$ and $|\{a \dot e {\bold b}^\at_1:a 
\le_{{\bold b}_1} b_{{\dot k}^*}\}| < {\dot k}^*$.
\end{enumerate}
\end{remark}

\begin{discussion} 
\label{e21}
1)  In \ref{e8} we can add 

\noindent
5) Assume $N = \gB^{[\bar\varphi]}_*$ and 
\mn
\begin{enumerate}
\item[$(i)$]  $N$ is a model of ${\gt}^{\ind}_n$ and
\sn
\item[$(ii)$]   $(\forall y_1 \ne y_2 
\in \theta^N) \Rightarrow (\exists^{\ge n} x \in P^N)
[x R y_1 \equiv \neg x R y_2]$.
\end{enumerate}
\mn
\Then \, any auto of $N$ is 
represented in ${\gB}^*$. The proof is just easier. 
Without (ii) we get a weaken result. 
If above we have dealt with complete embedding rather then just isomorphism
onto, we can get more. May like to allow such ${\bf f}$'s (and get
the same result).

\noindent
2)  We may like in \ref{e8} to allow ${\bold b}$ to be non-atomic.
One way is combining our proof with \cite{Sh:384}.  We shall give a
complete proof elsewhere.

Suppose $\gB_{\lambda^+} \models ``{\bold b}_1$, ${\bold b}_2$ are 
Boolean rings", ${\bold f}$ a complete embedding of 
${\bold b}_1[\gB_{\lambda^+}]$ into ${\bold b}_2[\gB]$, we would like to
show that it is representable.

Let us define, inside $\gB_{\lambda^+}$, $Y = \{a: a$ a ``finite"
subset of ${\bold b}_1$ consisting of pairwise disjoint elements such that
${\gB}_{\lambda^+} \models$ ``for every $c \in a$, either $c$ is 
an atom of ${\bold b}_1$ or below $c,{\bold b}_1$ is atomless"$\}$ and
for $a \dot e Y$ let $b_a= \bigcup \{x:x \dot e a\} \in {\bold b}_1$.

For every (internally) finite $a \dot e Y$ let
${\bold b}_1[a] := $ ``the sub-ring of ${\bold b}_1$ generated by $a$'', it
is, in $\gB_{\lambda^+}$, a finite sub-Boolean ring of ${\bold b}_1$ and 
itself is a Boolean algebra and inside $\gB_\lambda$
\mn
\begin{enumerate}
\item[$\otimes$]  if $a' \subseteq {\bold b}_1[a]$ is a maximal antichain
of ${\bold b}_1[a]$, then in ${\bold b}_1,a'$ is a maximal family of
pairwise disjoint elements which are $\le b_a$.
\end{enumerate}
\mn
Hence ${\bold f} \restriction {\bold b}_1[a]$ is a complete
embedding of ${\bold b}_1[a]$ into ${\bold b}_2 \restriction \{x \in
{\bold b}_2:x \le_{{\bold b}_2} {\bold f}(b_a)\}$ which also is 
a Boolean algebra (and sub-Boolean ring of ${\bold b}_2$), hence 
is represented in $\gB_{\lambda^+}$.  So by Stage C, i.e. as we 
have proved part (3) of \ref{e8}, we have
finished proving part (4) too.  For ${\bold b}$ a not necessarily atomic 
Boolean ring. (Note: if we have earlier the result only 
for isomorphism (onto ${\bold b}_2$), here we would have to work harder.)
\end{discussion}

\begin{claim}
\label{e24}
For $\gB_* = \gB_{\lambda^+}$ as in \ref{e8} we can add:
\mn
\begin{enumerate}
\item[$\bullet$]  every branch is equivalent to a definable 
branch $\hat H$ \when \, we assume (compare with \ref{c98}(2)):
\sn
\begin{enumerate}
\item[$(*)$]  $(a) \quad d^*,a \in {\bbN}^{\ns}[\gB_{\lambda^+}] 
= \{x:\gB_{\lambda^+} \vDash$ ``$x \dot e {\bbN}\, \& \, n < x"$ 
for every 

\hskip25pt  truly finite $n\}$
\sn
\item[${{}}$]  $(b) \quad \gB_{\lambda^+} \models ``\bar I = \langle I_d:d \in 
\bbN^{d^*}\rangle$ is a sequence of sets, $\dis$ is a 

\hskip25pt  symmetric two
place function from $I = \bigcup\{I_d:d \in \bbN^{d^*}\}$ 
into ${\bbN}$, 

\hskip25pt  satisfying $\dis(x, z)\le \dis(x,y) + \dis(y,z)"$ where

\hskip25pt $\bbN^d := \{x \in \bbN:x < d\}$
\sn
\item[${{}}$]  $(c) \quad \gB_{\lambda^+}$, $c$ is a function (monotonic for
simplicity) from ${\bbN}^{d^*}$ to ${\bbN}^a$ 

\hskip25pt  such that the value ``converge to finite" when 
the argument does, 

\hskip25pt i.e.:

\begin{equation*}
\begin{array}{clcr}
b \dot e \bbN^a \, \& \, \bigwedge\limits_{n} n<b \Rightarrow &(\exists
d)[\bigwedge\limits_n n < d \dot e {\bbN}^{d^*} \, \& \\
 &\forall x \forall d'\ (x \dot e I_{d'} \wedge 
\bigwedge\limits_{n} n<d' \le d \Rightarrow c(d')<b)].
\end{array}
\end{equation*}
\sn
\item[${{}}$]  $(d) \quad$ We call a function $\hat H$ (generally, not
  necessary in $\gB_*$) a $c$-branch 

\hskip25pt (but may omit $c$)
$\Dom(\hat H) = {\bbN}^{d^*} \setminus \{n:n<\omega\}$
\sn
\item[${{}}$]  $(e) \quad \hat H(d) \in I_d$  if
\sn
\item[${{}}$]  $\quad (\alpha) \quad
\bigwedge\limits_{n} n<d_1<d_2 \le d^* \Rightarrow \dis(\hat H(d_1),
\hat H(d_2)) \le 2 \times \max\{c(d_1),c(d_2)\}$
\sn
\item[${{}}$]  $(f) \quad$ branches 
$\hat H_1,\hat H_2$ are $i$-equivalent if 

\hskip35pt $\bigwedge\limits_{n} n<d \le d^* \Rightarrow \dis(\hat H_1(d),
\hat H_2(d)) \le 4 \times c(d)$
\sn
\item[${{}}$]  $(g) \quad$ a branch 
$\hat H$ is $c$-definable if for some 

\hskip35pt  $f \in \gB_{\lambda^+}$,
$d \in \Dom(\hat H) \Rightarrow \hat H(d)=f(d)$.
\end{enumerate}
\end{enumerate}
\end{claim}

\noindent
The proof of \ref{e24} is broken to some definition and claims.


\begin{fact}
\label{e25}
Let $\Gamma$ be a $g$-bigness notion, $\gk = (\bold K_{\gk},\le_{\gk})$ 
elementary class. \Then \,  $p = \tp(\bar a, A, M)$ is 
$\Gamma$-big iff $p^* = \tp(\bar a,\acl (A),M)$ is $\Gamma$-big.
\end{fact}

\begin{PROOF}{\ref{e25}}
The ``if" is by monotonicity. For the ``only if" assume the left,
let $M \le_{\gk} N \in {\bold K},N$ is strongly $|A|^+$-saturated, $\bar a'
\in {}^\alpha N,p' = \tp(a',\acl A,N)$ extend $p$ and is
$\Gamma$-big. But there is an automorphism $f$ of $N$ mapping $\bar a'$
to $\bar a,f \restriction A = \id_A$ so $f^{-1}(\acl_M(A)) = \acl_M(A)$.
\end{PROOF}

\begin{definition}
\label{e27}
1)  We say a bigness notion $\Gamma$ is strict (or strictly nice) \If
\, (where $\alpha = \lg(x_\Gamma)$)
\mn
\begin{enumerate}
\item[$(*)$]  if $p \in {\bold S}^\alpha(A,M)$ is $\Gamma$-big, $A=\acl_M(A)$
and $\beta$ is an ordinal \then \, for some $N,\bar a_i$ ($i< \beta$)
we have
\sn
\begin{enumerate}
\item[$(\alpha)$]  $M \le_{\gk} N$
\sn
\item[$(\beta)$]  $\bar a_i \in  {}^\alpha N$
\sn
\item[$(\gamma)$]  $\tp(\bar a_i,A \cup \bigcup\limits_{j< i} \bar
  a_j,N)$ is $\Gamma$-big
\sn
\item[$(\delta)$]  if $\bar b \in {}^m(\acl(A\cup \bar a_i)),\Omega =
\Gamma^\mt_{a,\dot{\bold d},\bar c}$ is an $a$-bigness notion 
(see Definition \ref{c95}), $A_\Omega
\subseteq A,\tp(\bar b,A,N)$ is $\Omega$-big then $\tp(\bar b,A \cup
\bigcup\limits_{j<i} \bar a_j, N)$ is $\Omega$-big; hence 
$\tp(\bar b, A\cup \bigcup\limits_{j \ne i} a_j,N)$ is $\Omega$-big.
\end{enumerate}
\end{enumerate}
\mn
2)  We say $p=\tp(\bar a, A_2, M)$ is a strict (or strictly nice) extension of
$p\restriction A_1$, where $A_1\subseteq A_2\subseteq M$
\when \,: if $m<\omega$, $\bar b \in {}^m(\acl_M(A_2+\bar a)),
\Omega = \Gamma^\mt_{a,\dis,\bar c}$ is a bigness notion 
(see Definition \ref{c95}), $A_\Omega \subseteq A_1,A_2 \subseteq A_1,
\tp(\bar b,\acl(A_1),M)$ is $\Omega$-big \then \, 
$\tp(\bar b,\acl (A_2),M)$ is $\Omega$-big.

\noindent
3)  A global bigness notion $\Gamma$ is strict (or strictly nice) 
when: if $p = \tp(\bar a,A_2,M)$ is $\Gamma$-big and $A_\Gamma \subseteq A_1$
and $A_1 \subseteq A_2 \subseteq M$ \then\  $p$ has a strictly nice extension
$\tp(\bar a',A_2,N),M \le_{\gk} N$, which is 
$\Gamma$-big. We may omit the ``nice" one 
and leave the ``strictly".
\end{definition}

\begin{definition}
\label{e31}
For $g$-bigness notions $\Gamma_1,\Gamma_2$, we say $\Gamma_1$ is
strictly orthogonal to $\Gamma_2$ or $\Gamma_1 \bot_s \Gamma_2$ when:
\mn
\begin{enumerate}
\item[$(*)$]  if $A_{\Gamma_1} \cup A_{\Gamma_2} \subseteq A\subseteq
M \in {\bold K},\tp(\bar a_\ell,A,M)$ is $\Gamma_\ell$-big for $\ell=1,2$
\then \,  we can find $N$, $\bar a'_1$, $\bar a'_2$ such that:
\sn
\begin{enumerate}
\item[$(\alpha)$]  $M \le_{\gk} N \in {\bold K}$
\sn
\item[$(\beta)$]  $\bar a'_1,\bar a'_2 \subseteq N$
\sn
\item[$(\gamma)$]  $\tp(\bar a'_\ell,A,N)$ extend $\tp(\bar a_\ell,
A,N)$ for $\ell = 1,2$
\sn
\item[$(\delta)$]  $\tp(\bar a'_\ell,A+ \bar a'_{3-\ell},N)$ is
$\Gamma_\ell$-big for $\ell=1,2$
\sn
\item[$(\varp)_1$]  $\tp(\bar a'_1,A+ \bar a'_2,N)$ is a
strictly nice extension of $\tp(\bar a'_1,A)$ 
that is: if $\Omega = \Gamma^{\mt}_{a,\dis,\bar c}$, a bigness notion,
$A_\Omega \subseteq \acl(A)$ and $\bar b \subseteq \acl (A+\bar a'_1)$
and $\tp(\bar b,A,N)$ is $\Omega$-big then $\tp(\bar b, A+\bar
a'_2,N)$ is $\Omega$-big.
\end{enumerate}
\end{enumerate}
\end{definition}

\begin{claim}
\label{e41}
Assume $T$ has Skolem function  (as in \ref{a2}(B))
If $\Gamma_1 \bot_s \Gamma_2$ \then \, $\Gamma_2 \bot_s \Gamma_1$, 
moreover we have
\mn
\begin{enumerate}
\item[$\otimes$]  in $(*)$ inside Definition \ref{e31} above it follows that
\mn
\begin{enumerate}
\item[$(\varp)_2$]  $\tp(\bar a'_2, A+ \bar a'_1, N)$ is strictly nice
extension of $\tp(\bar a'_2, A, N)$, i.e. if $\Omega =
\Gamma^{\mt}_{a,\dis,\bar c}$ a 
bigness notion, $A_\Omega\subseteq \acl(A)$ and $\bar b
\subseteq  \acl(A+ \bar a'_2),\tp(\bar b, A, N)$ is $\Omega$-big
\then \,  $\tp(\bar b, A+\bar a'_1, N)$ is $\Omega$-big.
\end{enumerate}
\end{enumerate}
\mn
So we can say ``$\Gamma_1$, $\Gamma_2$ are strictly nicely orthogonal",
i.e. this is symmetric relation.
\end{claim}

\begin{PROOF}{\ref{e41}}
Use \ref{c95}(2). Let us prove $(\epsilon)_2$ from $\otimes$, so
assume that $\Omega = \Gamma^{\mt}_{\dis,\bar c}$, $\bar b \subseteq
\acl(A+\bar a'_2)$ is a counterexample, in particular $A_\Omega
\subseteq A$ and
\mn
\begin{enumerate}
\item[$(*)$]  $\bar b \subseteq \acl(A + \bar a'_2)$ is a counter
  example, so
\sn
\begin{enumerate}
\item[$(a)$]  $\tp(\bar b,\acl_{N_*}(A),N)$ is $\Omega$-big
\sn
\item[$(b)$]  $\tp(\bar b,\acl_N(A + \bar a'_1),N)$ is $\Omega$-small.
\end{enumerate}
\end{enumerate}
\mn
Now $\tp(\bar b,A_{\bar{a}_1})$ is $\Omega$-small, but
$\Gamma^{\mt}_{a,\dis,\bar c}$ is a local bigness notion hence by
\ref{e25} there is a formula $\vartheta (\bar{y},\bar{a}'_1) \in
\tp(\bar{b},A+ \bar a'_1)$, possibly with parameter from $A$, which is
$\Omega$-small, so there are $n<\omega$ and $i < \lg(\bar{c})$ such 
that $\{\dis(\bar{y}_\ell,\bar{y}_m) \ge c_i:\ell < m<n\} \cup \{\vartheta
(y_\ell,\bar a'_1):\ell \le n\}$ is inconsistent and \wilog \, 
$n$ is minimal.

Clearly, it follows that (we can add $\bigwedge\limits_{k < \ell < n}
\dis(\bar y_k,\bar y_\ell) \ge c_i$ but no need)

\[
N \models (\exists \bar y_0\ldots \bar y_{n-1} [\bigwedge\limits_{\ell<n}
\vartheta (\bar y_\ell,\bar a'_1) \wedge (\forall \bar y) 
(\vartheta_\ell(\bar y,\bar a_1) \rightarrow \bigwedge\limits_{\ell<n} 
\dis(\bar y,\bar y_\ell)<c_i)].
\]

\mn
As $T$ has Skolem functions clearly there are 
$\bar b_0,\ldots,\bar b_{n-1} \subseteq \acl(A+\bar a'_1)$ such that 
$N \models \bigwedge\limits_{\ell<n} \vartheta (\bar b_\ell,\bar a'_1)
\wedge (\forall \bar y)(\vartheta_0,(\bar y,\bar a'_1) \rightarrow 
\bigwedge\limits_{\ell<n} \dis(\bar y,\bar b_\ell)<c_i)$. 
We can substitute $\bar b$ for $\bar y$ so for some 
$\ell<n,\bar b^*\in \{\bar b_0,\ldots \bar b_{n-1}\}$ 
we have $N \models \dis(\bar b,\bar b^*)<c_i$.

Now first,
\medskip

\noindent
\underline{Case 1}: if $\tp(\bar{b}^*,\acl_N(A),N)$ is 
not $\Omega$-big 

\Then \, as in the previous sentences, for some 
$\bar{b}^{**} \subseteq \acl_N(A)$ and $j < \lg (\bar{c})$ we have

\[
N \models \dis(\bar b^*,\bar b^{**}) \le c_j,
\]

\mn
so together with the previous sentence $N \models \dis(\bar{b},
\bar{b}^{**}) \le \max\{i,j\}+1$, contradiction to $\tp(\bar{b},
\acl_N(A),N)$  is $\Omega$-big.  

Second, 
\medskip

\noindent
\underline{Case 2}: $\tp(\bar{b}^*,\acl_N(A),N)$ 
is $\Omega$-big 

We can still note that $\tp(\bar{b}^*,
\acl_N(A+\bar{a}'_1),N)$ is not $\Omega$-big contradicting 
$(\epsilon)_1$ of (*) of \ref{e31} by \ref{e25}.

Together we have gotten a contradiction.
\end{PROOF}

\begin{claim}
\label{e44}
If $\Gamma_1,\Gamma_2$ are orthogonal $g$-bigness notions and $\Gamma_1$
 is strict \then \, $\Gamma_1$, $\Gamma_2$ are 
strictly orthogonal.
\end{claim}

\begin{PROOF}{\ref{e44}}
Let $p_\ell \in \bold S^{\alpha(\Gamma_\ell)}(A,\gC)$ for $\ell=1,2$
and $A_\Gamma \subseteq A$.

Now let $\lambda = |T| + |A| + |\alpha(\Gamma_1)| + |\alpha(\Gamma_2)|
+ \aleph_0 + |\{\Omega:\Omega$ is as in \ref{e31}$(*)(\varp)_2$ for $A$ and
we choose $A_\alpha,\bar a_{1,\alpha}$ by induction on $\alpha \le
\lambda^+$ such that:
\mn
\begin{enumerate}
\item[$(*)_{1,\alpha}$]  $(a) \quad A_\alpha = \cup\{\bar
  a_{1,\beta}:\beta < \alpha\} \cup A$
\sn
\item[${{}}$]  $(b) \quad \tp(\bar a_{1,\alpha},A_\alpha,\gC)$ is a
  $\Gamma_1$-big strictly nice extension of $p_1$.
\end{enumerate}
\mn
This is possible by the assumption.  We can find $\bar a_2$ realizing
$p_2$ such that $\tp(\bar a_2,A_{\lambda^+},\gC)$ is $\Gamma_2$-big.
It is enough to prove that for some $\alpha < \lambda^+$
\mn
\begin{enumerate}
\item[$(*)_{2,\alpha}$]  for every finite $\bar b \subseteq \acl(A +
  \bar a_{1,\alpha},\gC)$ and $\Omega = \Gamma^{\mt}_{a,\dis,\bar c}$
  with parameters from $A$, if $\tp(\bar b,\acl(A),\gC)$ is
  $\Omega$-big then $\tp(\bar b,\acl(A + \bar a_2),\gC)$ is
  $\Gamma$-big.
\end{enumerate}
\mn
Toward contradiction assume this fails for every $\alpha$ and let
$(\Omega_\alpha,\bar b_\alpha,i_\alpha)$ witness this so
$\Omega_\alpha = \Gamma^{\mt}_{a,\dis,\bar c_\alpha},i_\alpha < \ell
g(\bar c_\alpha)$.  By the choice of $\lambda$ for some $\alpha <
\beta < \lambda^+$ we have $(\Omega_\alpha,\bar b_\alpha,i_\alpha) =
(\Omega_\beta,\bar b_\beta,i_\beta)$.  

By transitivity of distance we get contradiction to $(*)_{1,\beta}(b)$.
\end{PROOF}

\begin{claim}
\label{e47}
1)  If $T$ is as in \ref{a2}(2) or at least $T$ has 
Skolem function \then \, every bigness notion is a strictly nice 
bigness notion. 

\noindent
2) If the bigness notion $\Gamma_1,\Gamma_2$ are orthogonal \then \, 
they are strictly orthogonal.
\end{claim}

\begin{PROOF}{\ref{e47}}
1)  By \ref{c98}(2).

\noindent
2)  As in the proof of \ref{a53}(3).
\end{PROOF}

\noindent
Now we should check our notion, and revise the
construction in \S4 and \ref{e8}.  For local bigness notions we get
better results.
\begin{claim}
\label{e51} 
Assume $T$ is as in \ref{a2}(2) ,(or just has Skolem function 
in a strong enough sense) and $\gk = (\bold K,\le) = (\text{mod}_T,\prec)$.
Assume $\Gamma$ is a local bigness notion, $A_\Gamma \subseteq
A \subseteq B \subseteq N$ and $A = \acl_N(A),B = \acl_N(B),p=\tp(d,B,N)
\in {\bold S}^{\alpha(\Gamma)}(B,N)$ is $\Gamma$-big and 
$\acl_N (A+\bar d) \cap B = A$, i.e. niceness, 
and $\Omega= \Gamma^{\mt}_{a,\dis,\bar c}$ is a bigness notion with
$A_\Omega \subseteq A$.

\Then \, we get strict niceness, and even 
\mn
\begin{enumerate}
\item[$(\alpha)$]  if $\bar{b}_1 \in B$ and $\tp(\bar{b},A,N)$
is $\Omega$-big \then \, $\tp(\bar{b}_1,\acl(A+\bar d),N)$ is $\Omega$-big
\sn
\item[$(\beta)$]  if $\bar{b}_2 \in \acl_N (A+d)$ and $\tp(\bar{b}_2,A,N)$
is $\Omega$-big then $\tp(\bar{b}_2,B,N)$ is $\Omega$-big.
\end{enumerate}
\end{claim}

\begin{PROOF}{\ref{e51}}
For transparency we use singleton.
We can find a two-place function $\breve{f} = 
\breve{f}_\omega$ definable in $N$ such that:
\mn
\begin{enumerate}
\item[$(*)_1$]  $(\alpha) \quad \Dom(\breve{f}) = \{(d,c):d \dot e a,c 
\dot e \bbR^+\}$
\sn
\item[${{}}$]  $(\beta) \quad N \models (\forall y \dot e a) 
[\dis(y,\breve{f}(y,c))<c \, \& \, \breve{f}(y,c) \dot e \bbR^+]$ for every
$c \in (\bbR^+)^N$ 
\sn
\item[${{}}$] $(\gamma) \quad N \models (\forall y,y') 
[\breve{f}(y,c) = z \wedge {\dis}(y',z) < c/2 \rightarrow 
\breve{f}(y',\bold c,\zeta)]$ 

\hskip25pt  for every $c \in (\bbR^+)^N$.
\end{enumerate}
\mn
[Why?  Here we use ``$T$ as in \ref{a2}(B) (and $a$ is ``a set, not a
class").  That is, there is $b$ for some $F_1 \in \tau_N$ such that $N
\models$ ``for every appropriate $a,\dis,c,F_1(a,c,\dis)$ is a maximal
subset of $a$ such that $(\forall xy)[x \ne y \dot e a \rightarrow
\dis(x,y) \ge c]"$.  Clearly exist and let $F_2 \in \tau_M$ such that
$N \models$ ``for every appropriate $a,\dis,c,b \dot e a,x =
F_2(b,a,c,\dis)$ is a member of $F^N_1(a,c,\dis)$ such that $\dis(x,b)
< c$", exists by the maximality of $F_1(a,c,\dis)$.  So $(b,c) \mapsto
F_2(b,a,c,\dis)$ is a function as required.]

Now if $A_\Omega\subseteq C \subseteq N$ then
\mn
\begin{enumerate}
\item[$(*)$]  $\tp(d,b,C,N)$ is $\Omega$-big \underline{if and only if} for
every large enough $i < \lg(\bar c)$ we have $f_\Omega(b,c_i) \notin 
\acl_N (A)$.
\end{enumerate}
\end{PROOF}

\begin{claim}
\label{e54}
1)  In claim \ref{c53}(1), (2) we can strengthen clause
$(\gamma)$ to
\mn
\begin{enumerate}
\item[$(\gamma)^+$]  $\tp(b',\acl_N(A'),N)$ is a strict extension of
$\tp(a,A,M^*)$.
\end{enumerate}
\mn
2)  In claim \ref{c59} we can strengthen clause $(\beta)$ to
\mn
\begin{enumerate}
\item[$(\beta)^+$]  $\tp(b',\acl(A+b_2),\gC)$ is a
strict extension of $\tp(b_1,\acl(A),\gC)$.
\end{enumerate}
\mn
3)  In \ref{c71}(2) we can replace ``nicely" (in clause
$(\delta)$ of the conclusion) by strictly.
\end{claim}

\begin{PROOF}{\ref{e54}}
We can just use \ref{e51}, (its assumption is O.K. for our
application.

\noindent
1) In the beginning of the proof of \ref{c53} we reduce it to
the proof of $q^*(x)$ being $\Gamma^{\ms}_{a,\bar c}$-big, and there we
get (2) by applying proof to $\bigcup\limits_{i< \lambda} A'_i$; to
get $(\gamma)^+$ we just need to require there that $\tp(f_i(A'),
\acl(\bigcup\limits_{j<i} f_j(A')), M^*)$ is strict.

\noindent
2) Similarly.

\noindent
3) Similarly: when we waive ``nicely" we just replace ``$p_\zeta$
nicely extend $\tp(d^*,A)$" by ``$p_\zeta$ strictly extend $\tp(d^*,A)$".
\end{PROOF}

\begin{claim}
\label{e57}
We can repeat \S4, replace ``$g$-bigness notions'' by ``strict bigness
notions" (automatic if we use $T^*$ and in (D)(7)), ``nice extension" 
by ``strict extension".
\end{claim}

\begin{PROOF}{\ref{e57}}
Straightforward (or use \ref{e51}(4)).

Proof of \ref{e24}:

In the proof of \ref{e8} we replace nicely by strictly nice
and as we work assuming \ref{a2}(2), by \ref{e51}(4), this holds
automatically. The additional point is similar to Stage C.
We add two contractors: the pre-pseudo branch killer
$\zeta_{\pr}$, and the pseudo branch killer $\zeta_{\pk}$.

For $\beta \in W_{\zeta_{\pr}}$ of cofinality $\lambda$, assigned to
our parameters (so the relevant parameters are in $\gB_\alpha$ 
for some $\alpha<\beta$), we demand
\mn
\begin{enumerate}
\item[$\otimes_{1,\beta}$]   for stationarily many $\epsilon\in
S_{\zeta_\pr}$, for some $j \in u^{\beta+1}_{i_{\beta+1, \epsilon}}$,
non-maximal in $u^{\beta+1}_{i_{\beta+1,\epsilon}}$, satisfying
$j \ge \max(u^\beta_{i_{\beta,\epsilon}})$, we have:
\sn
\begin{enumerate}
\item[$(a)$]  $\bar c^{\beta+1} = \langle \bar c^{\beta+1, j}_n: n<
\omega\rangle$ is indiscernible over $A^{\beta+1}_j$
\sn
\item[$(b)$]  $\bar c^{\beta+1,j}_n = \langle c^{\beta+1,j}_{n,i}:i<
i(*)\rangle,i(*)< \lambda$ 
\sn
\item[$(c)$]  $\{c^{\beta+1,j}_{0, i}: i< i(*)\}$ is the universe of an
(elementary) submodel of $\gB_\beta$, including $A^{\beta+1}_j \cap
\gB_\beta$
\sn
\item[$(d)$]  $\gB_\beta \vDash$ ``$c^{\beta+1, j}_{0,\ell} \dot e
\bbN$ and $m< 2 c^{\beta+1,j}_{0, \ell} < c^{\beta+1,j}_{0,\ell+1}
< {\dot n}",\dot n \in A^{\beta+1}_j \cap \gB_\beta,\gB_\beta 
\vDash ``m < \dot n \in \bbN$" for every true natural number 
$m$, all this for any $\ell< \omega$
\sn
\item[$(e)$]  $\gB_{\beta+1} \vDash ``c^{\beta+1,j}_{0, \ell_1} >
c^{\beta+1, j+1}_{0, \ell_2}$ for $\ell_1$, $\ell_2< \omega$
\sn
\item[$(f)$] $\tp(\bar c^{\beta+1,j}_n, \bigcup\limits_{\ell<n} 
\bar c^{\beta+1,j}_\ell \cup A^{\beta+1}_j)$ is a strict extension of
$\tp(\bar c^{\beta+1}_n,A^{\beta+1}_j)$
\sn
\item[$(g)$]  $A^{\beta+1}_{j+1} \cap \gB_\beta = 
\Rang(\bar c^{\beta+1, j}_0)$.
\end{enumerate}
\end{enumerate}
\mn
For every $\alpha<\lambda^+$
\mn
\begin{enumerate}
\item[$\otimes_{2,\alpha}$]  if $\beta< \alpha$ and $\zeta_{\pr}$
acted in $\beta$ for a pseudo tree as above, \then \, for every
$\breve{f} \in I[\gB_{\alpha +1}] \setminus I[\gB_\alpha]$, 
for stationary many $\epsilon \in S_{\zeta_\pr}$, for some 
$j< \lambda$ as in $\otimes_{1, \beta}$ we have:
\sn
\begin{enumerate}
\item[$(*)$]  $(a) \quad \epsilon \in E_\beta \cap E_\alpha,\beta
\in \bold c^\alpha_{i^\alpha_\epsilon},j \in
u^\alpha_{i^\alpha_{\dot e}}$
\sn
\item[${{}}$]  $(b) \quad$  one of the following occurs:
\sn
\item[${{}}$]  $\quad (\alpha) \quad \tp(\breve{f},A^\alpha_{j+1},
\gB_{\alpha+1})$ is $\Gamma^{\mt}_{\dis,(\bar c^{\beta,j}_0 \rest
  \omega)}$-big
\sn
\item[${{}}$]  $\quad (\beta) \quad$ there are $\ell< \omega$ 
and $\breve{f}' \in I[\gB_\alpha]$ such that, in $\gB_{\alpha+1}$:
\sn
\item[${{}}$]  $\qquad (i) \quad \dis(\breve{f}',\breve{f}) \le 
c^{\beta, j}_{0, \ell}$
\sn
\item[${{}}$]  $\qquad (ii) \quad$ for every $\breve{f}'' \in 
I[\gB_\beta]$ we have $dis(\breve{f}',\breve{f}'') > c^{\beta,j}_{0,\ell+1}$
\sn
\item[${{}}$]  $\quad (\gamma) \quad$ there is $\breve{f}' \in
  A^\beta_j$ such that $\dis(\breve{f},\breve{f}') \in \bbN$ is 
$< {\dot n}$ for every 

\hskip25pt  non-standard ${\dot n} \in \bbN[\gB_\beta]$.
\end{enumerate}
\end{enumerate}
\mn
Why can we do this?   Having $p^\beta_j$ first define $p^\beta_{j+1}$ ignoring
clause (b) of $(*)$, if it holds, fine; assume not. 
As $(\alpha)$ fail by observation \ref{c102} for some 
$\breve{f}' \in A^\alpha_{j+1}$ and $\ell< \omega$ we have: 
$p^\alpha_{j+1}$ says that ``$\dis(\breve{f},\breve{f}') \le 
c^{\beta,j}_{0,\ell}"$.

If $\tp(\breve{f}', A^\beta_{j+1},\gB_\alpha)$ is 
$\Gamma^{\mt}_{\dis,\bar c^{\beta,j}_0\restriction \omega}$-big, 
possibly increasing $\ell$ we get possibility $(\beta)$ as we 
are using the strict version in $(D)(7)$ of \ref{d14}. 
So assume that $\tp(\breve{f}',A^\beta_{j+1},\gB_\alpha)$ is
$\Gamma^{\mt}_{\dis,\bar c_0^{\beta,j}} \restriction w$-small. So using
\ref{c102} again and transitivity for some $\breve{f}'' \in
A^\beta_{j+1}$ (possibly increasing $\ell$) we have
$\dis(\breve{f},\breve{f}'') < c^{\beta,j}_{0,\ell}$. We can find $i$ such that
$\breve{f}'' = c^{\beta,j}_{0,i}$ where $i<i^*$. So we can find $q$, another
candidate for $p^\alpha_{j+1}$, which say $\dis(\breve{f},
c^{\beta,j}_{1,i}) < c^{\beta,j}_{1,\ell}$ but recall $c^{\beta,j}_{1,\ell}<
c^{\beta,j}_{0,\ell}$.  So to get clause $(\beta)$ it suffices to have in
$\gB_\alpha$ that $\tp(c^{\beta,j}_{1,i},A^\beta_{j+1})$ is
$\Gamma^{\mt}_{\dis,\bar c^{\beta,j}_0 \restriction \omega}$-big. 
By the choice of $\langle\bar{c}^{\beta,j}_n:n < \omega\rangle$ 
this holds except when $\breve{f}'' = c^{\beta,j}_{0,i} \in A^\beta_j$.

But then $\dis(\breve{f},\breve{f}'')$ belongs to $A^{\beta+1}_j$, so it is a
``$\gB_{\alpha+1}$-natural number", i.e. a members of $A^{\beta+1}_j \cap
\gB_\beta \cap {\bbN}[\gB_{\beta+1}]$ smaller than all ``non standard"
member of ${\bbN}[\gB_\beta]$; as this holds for unboundedly many 
$\beta<\lambda$ clearly clause $(\gamma)$ of $(*)(b)$ holds.

So having guaranteed the relevant $\otimes_{1, \beta}$, $\otimes_{1,
\alpha}$, we apply it to a ``branch'' $\langle g({\dot n}):\dot n \in
{\bbN}^a[\gB_{\lambda^+}] \setminus \{m:m<\omega\} \rangle$, for each
$\beta< \lambda^+$ satisfying, $\cf(\beta) = \lambda$ and
$\alpha \in (\beta,\lambda^+)$ such that $\tp(a^\alpha_0,
\gB_\alpha,\gB_{\alpha+1})$ is $\Gamma^{\av}_{\dot D,\langle n:
n< \omega\rangle}$ and $\dot n \in \gB_\beta \Rightarrow (\dot n) 
\in \gB_\beta$.

So we can find $\breve{f}''_\beta \in I[\gB_\beta]$ such that
$\dis(\breve{f}''_\beta,g(a^\alpha_0)) \in {\bbN}[\gB_{\lambda^+}]$
is smaller than all non standard ${\bold n} \in {\bbN}[\gB_\beta]$. By
Fodor lemma for some $\breve{f}^*$ for stationary many $\beta$,
$\breve{f}''_\beta = \breve{f}^*$. This finish the proof.
\end{PROOF}

\begin{remark}
\label{e60}
1) It seems we can use guessing of clubs as in \cite{Sh:413} (more
\cite{Sh:572}) and ``$\check{\cD}_\lambda$ is $\lambda^+$-saturated" 
(when $\lambda> \aleph_1$) to deal with \ref{e24}, but the present 
look simpler and did not check.

\noindent
2) We can also in stage C in the proof of \ref{e8} deal with
weaker notions of trees (with distance instead equality).

\noindent
3) By \ref{e8} and the previous chapter, we can conclude the compactness
of the quantifier on complete embeddings of one boolean ring to another.

\noindent
4) It is more natural for general $T$ (in order to save
claim \ref{e41}) to replace $\Gamma^{\mt}$ in Definition 
\ref{e27},\ref{e31} by $\Gamma^{\mt}$ as in the following definition.
\end{remark}

\begin{definition}
\label{e61}
Assume
\mn
\begin{enumerate}
\item[$(a)$]  $T$ is as in \ref{a2}(1), i.e. a complete first order
  theory
\sn
\item[$(b)$]  $M^*$ is a model of $T$
\sn
\item[$(c)$]  $\bar{\varphi}=(\varphi_1(x),\varphi_1(x),
\varphi_2(x,y),\varphi_3(xy,z),\varphi_4 (x,y,z)) = 
(\varphi_{\dom}(x),\varphi_{\rang}(x),x \le y,\varphi^{\dis}_{\le z}(x,y),
\varphi_{\subadd}(x,y,z))$ are first order formulas 
(possibly with parameters),
or $\varphi$ (the intended meaning of $\varphi_{\dom} (M^*)$
is the ``space, and of $\varphi_{\rang}(M)$ the possible distances
\sn
\item[$(d)$]  the formula $\varphi^{\dis}_{\le_z} (x,y)$ is such that:
\[
\varphi^{\dis}_{\le_z} (x,y) \Rightarrow \varphi_{\dom}(x) \, \& \,
\varphi_{\dom}(y) \, \& \, \varphi_{\rang}(z) \text{ and }
 \dis_{\le z_1}(x,y) \& z_1 \le z_2 \Rightarrow \dis_{\le z_2}(x,y)
\]
(the intended meaning of $\varphi^\dis_{\le z}(x,y)$ is 
the distance from $x$ to $y$ is $\le z$)
\sn
\item[$(e)$]  $\varphi_{\subadd} (x,y,z)$ is a first order formula
such that 
\begin{equation*}
\begin{array}{clcr}
\varphi_{\subadd}(x,y,z) \Rightarrow \varphi_{\rang} (x) \, &\& \,
\varphi_{\rang}(y) \, \& \,\varphi_{\rang}(z) \,\& \, \\
  &x \le z\wedge y \le z \forall x, y(\varphi_{\rang}(x) \wedge 
\varphi_{\rang}(y) \rightarrow (\exists z) \varphi_{\subadd}{(x,y,z))}
\end{array}
\end{equation*}
the (intended meaning of $\varphi_{\subadd} (x,y,z)$ is $x+y \le z$)
\sn
\item[$(f)$]  $\varphi^{\dis}_{\le z_1}(x_1 x_2) \wedge
  \varphi^{\dis}_{\le z_2}(x_2,x_3) \wedge \varphi_{\subadd}
(z_1,z_2,z_3) \Rightarrow \varphi^{\dis}_{\le z_3}(x_1,x_3)$
\sn
\item[$(g)$]  $x_1 \le x_2$, is a partial directed order of 
$\varphi_{\Mod}(M^*)$ 
\sn
\item[$(h)$]  $\bar{c}= \langle c_i:i<\delta\rangle, c_1 \in
\varphi_{\rang}(M^*), i<j\Rightarrow \varphi_{\subadd} (x_i,x_i,x_{i+}1)$.
\end{enumerate}
\mn
We define $\Gamma=\Gamma^{\mt}_{\bar{\varphi},\bar{c}}$ as the
following local bigness notion: $\vartheta (x,a)$ is $\Gamma$-big (in $N, M^*
\prec N)$ \underline{if and only if}
 
\[
\{\vartheta(x_n):n<\omega\} \cup \{\neg \varphi^{\dis}_{\le c_i} 
(x_n,x_m): n \ne m < \omega\}
\]
\mn
is consistent.
\end{definition} 

\begin{remark}
\label{e66} 
More generally we can use any dependency relation.
\end{remark}
\newpage

\section {Constructing models in $\aleph_1$ under $\CH$}

This section has little dependence on the earlier parts. 
In Rubin Shelah \cite{RuSh:84} models in $\aleph_1$ were constructed using:
in two cases $\diamondsuit_{\aleph_1}$ and in one $\CH$. Here we prove all of
them under $\CH$ and get further results (for example the results on ordered
fields). The construction here was promised in \cite{Sh:107}.
The omission of types in \ref{f5} continue \cite{RuSh:84},
hence Keisler \cite{Ke70} or \cite{Ke71} which deal with the quantifier 
$\exists^{\ge \aleph_1}$.

\begin{context}
\label{f2}
$T^*$ as in the context \ref{a2}(B) $T^*$ countable. 
$M^*$ a countable model of $T^*$ (usually well
founded), with universe $\omega$ for simplicity.  
\end{context}

\begin{definition}
\label{f5}
1)  For a formula $\varphi(x,\bar y)$ (here first order) and term $\sigma$ of
$T^*$ we shall define the formula $(\dot{\bold Q}^{\id}_\sigma x)
\varphi(x,\bar y)$, assuming $\sigma$ is a
term whose set of free variables does not include $x$ and for notational
simplicity is $\subseteq \bar y$, so we can write $\sigma(\bar y)$; now
$(\dot{\bold Q}^{\id}_{\sigma(\bar y)}x)\varphi(x,\bar y)$ 
means ``$\sigma=\sigma(\bar y)$ is an $\aleph_1$-complete 
ideal\protect\footnote{ We can restrict ourselves to
a class $R$ (of $T^*$) or allow non-first order $\varphi$, and get other
variants. If we would like to have $|\omega^M|=\aleph_1$, we start with a
$\aleph_1$-saturated model (by section 5) and apply the theorems below to it.} 
and $\{x \dot e \Dom(\sigma):\varphi(x;\bar y)\} \dot e \sigma(\bar y)$
\underline{or} $\sigma (\bar y)$ is not an $\aleph_1$-complete ideal".

\noindent
2) If $M$ is a model of $T^*$, $p$ a type over $M$ (i.e., a set of
formulas $\varphi(\bar x,\bar a), \; \; \bar a \subseteq M$, $\bar x$ a
fixed finite sequence) we say: $p$ is suitably omitted by $M$ \If \,:
\mn 
\begin{enumerate}
\item[$(*)$]  if $\bar b\subseteq M^*$, $n<\omega$, and
$M \models (\dot{\bold Q}_{\sigma_n} y_n) \ldots 
(\dot{\bold Q}_{\sigma_1}y_1)(\exists\bar x)\psi(\bar x,y_1\ldots
y_n,\bar b)$, \then \, for some $\varphi(\bar x,\bar a)\in p,
M \models (\dot{\bold Q}_{\sigma_n} y_n) \ldots 
(\dot{\bold Q}_{\sigma_1}y_1)(\exists\bar x)[\psi(\bar x,y_1,
\ldots,y_n,\bar b)\, \& \, \neg\varphi(\bar x,\bar a)]$.
\end{enumerate}
\mn
3)  Let $M$ be a model of $T^*$, and $\Delta$ be a finite set of formulas
of the form $\varphi(\bar x;{\bar y}_{\varphi})$, $p$ a set formulas of the
form $\varphi(\bar x,\bar a)$, where $\bar a\subseteq M$, $\varphi(\bar x;
\bar y) \in \Delta$. We say that $p$ is strongly undefined over $M$
\If \, there are no sequences $\langle\psi_\varphi(\bar y_\varphi,
\bar c_\varphi):\varphi(\bar x,\bar y_\varphi) \in
\Delta\rangle$ where ${\bar c}_\varphi\subseteq M$ such that: 

\[
\varphi(\bar x,\bar a) \in  p \Rightarrow M \models \psi_\varphi[\bar a,
\bar c_\varphi],
\]

\[
\neg\varphi(\bar x,\bar a) \in p \Rightarrow M \models \neg 
\psi_\varphi[\bar a,\bar c_\varphi].
\]
\end{definition}

\begin{observation}
\label{f8}
1)  In \ref{f5}(3) if $p$ is strongly undefined type over $M$,
\then \, it is suitably omitted over $M$; if $p$ is suitable 
omitted by $M$ \then \, $p$ is omitted by $M$. 

\noindent
2)  If $a \in M$ and $A \subseteq \{b:M \models b \dot e a\}$ 
is not represented in $M$, \then \
   
\[
p := \{x \subseteq a\, \& \, [b \dot e x]^{\iif(b \in A)}:M \models b
\dot e a\} 
\]

\mn
is strongly undefined in $M$. 

\noindent
3) If $M_n \prec M_{n+1}$ for $n<\omega,p$ a type over $M_0$ suitably
omitted by $M_n$ for each $n$, \then \,  $p$ is suitably omitted by
$\bigcup\limits_{n<\omega} M_n$. 
\end{observation}

\begin{construction}
\label{f11}
We describe a construction of an elementary extension $M^{**}$ of $M^*$ of
cardinality $\aleph_1$; we leave some points for latter fulfillment.
\medskip

\noindent 
\underline{\bf Step A}:   Let $\chi>\aleph_1$, and $x\in
{\cH}(\chi)$ be given.  Let $\tau_* \in {\cH}(\chi)$ be a countable vocabulary
extending $\tau_{T^*}$, having infinitely many $n$-place predicates and
$n$-place function symbols in $\tau_{T^*}$ for each $n<\omega$.
\medskip  

\noindent
\underline{\bf Step B}:  We choose a list $\bar\eta=\langle
\eta_\alpha:\alpha<\omega_1\rangle$ of ${}^{\omega_1>}2$, such that

\[
[\eta_\alpha \vartriangleleft \eta_\beta \Rightarrow \alpha<\beta].
\] 

\mn
So $\langle\rangle=\eta_0$, and for simplicity: $\eta_\alpha \char 94 
\langle\ell \rangle\in\{\eta_{\alpha+k}:0<k<\omega\}$ for $\ell=1,2$, and 

\[
\lg(\eta_\alpha) \text{ is a limit ordinal } \equiv \alpha 
\text{ is limit.}
\]

\mn
Let $\langle S_\alpha:\alpha<\omega_1\rangle$ be a partition of $\omega_1$
to pairwise disjoint stationary sets such that
$\min(S_\alpha)>1+\alpha$ and each $S_\alpha$ is non-small, see 
\cite[3.1(2)]{Sh:E62} and history there.

By induction on $\alpha<\omega_1$ we choose $N_\alpha$ and 
$M_{\eta_\alpha}$ and $\beta(\alpha), {\bf G}_{\eta_\alpha}$
such that: 
\mn
\begin{enumerate} 
\item[$(a)$]  $N_\alpha \prec (\cH(\chi),\in,<^*_\chi)$, $\alpha
\subseteq N_\alpha$, $\{x,\bar\eta\}\in N_\alpha$, $\langle N_\gamma:\gamma
\le \beta\rangle\in N_\alpha$ for $\beta<\alpha$ and $\langle
M_{\eta_\gamma},G_{\eta_\gamma},b_{\eta_\gamma}:\gamma \le \beta\rangle\in
N_\alpha$ for $\beta <\alpha$,   
\sn
\item[$(b)$]   the sequence $\langle N_\beta:\beta \le 
\alpha\rangle$ is increasing and continuous, each $N_\alpha$ is countable,
\sn
\item[$(c)$]  $M_{\eta_\beta}$ is a model of $T^*$, with universe 
$\omega(1+ \lg(\eta_\beta)),M_{\langle\rangle}=M^*$ (but for convenience
$\beta$ as a member of those models is called $a_\beta$), 
\sn
\item[$(d)$]  for $\eta_\alpha\vartriangleleft\eta_\beta$, then
$M_{\eta_\alpha} \prec M_{\eta_\beta}$, and if $p$ belongs to 
$N_\beta$ and is a suitably omitted type over 
$M_{\eta_\alpha}$, then $M_{\eta_\beta}$
suitably omits it too,   
\sn
\item[$(e)$]  if $\lg(\eta_\alpha)$ is a limit ordinal, then $M_{\eta_\alpha}=
\bigcup\{M_{\eta_\alpha \restriction i}:i<\lg(\eta_\alpha)\}$,   
\sn
\item[$(f)$]  if $\lg(\eta_\alpha)=\gamma+1$, $\gamma\in S_{\beta(\gamma)},
M_{\eta_\alpha\restriction\gamma} \models ``a_{\beta(\gamma)}$ 
is an $\aleph_1$--complete ideal on $\Dom(a_{\beta(\gamma)})=\bigcup\{y:y
\dot e a_{\beta(\gamma)}\}"$, \then \, let $y_{\eta_\alpha} = a_{\beta(\gamma)}$; otherwise let
$y_{\eta_\alpha}$ be the ideal of non-stationary subsets of $\omega_1$ in
the sense of $M_{\eta_\alpha\restriction\gamma}$,  
\sn
\item[$(g)$]  if $\lg(\eta_\alpha)=\gamma+1$, then
\sn
\begin{enumerate}
\item[$(i)$]  ${\bold G}_{\eta_\alpha}$ is an $N_\alpha$-generic subset of
${\bbP}_{\eta_\alpha}$, where: 
\sn
\item[${{}}$]  $\bullet \quad \bbP_{\eta_\alpha} 
:= \{\varphi:\varphi=\varphi(x_{\omega(1+\alpha)},
x_{\omega(1+\alpha)+1},\ldots,x_{\omega(1+\alpha)+n};\bar b)$
\newline
\quad for some $n<\omega,\bar b \subseteq 
M_{\eta_\alpha \restriction \gamma}$ and
\newline
$\quad M_{\eta_\alpha\restriction\gamma} \models 
(\dot{\bold Q}_{y_{\eta\alpha}}x_{\omega
(1+\alpha)})(\exists x_{\omega(1+\alpha)+1}) \ldots 
(\exists x_{\omega(1+\alpha)+n})\varphi\}$
\sn
\item[${{}}$]  $\bullet \quad$ the order on $\bbP_{\eta_\alpha}$ is
  naturally defined
\sn
\item[$(ii)$]  $\bold G_{\eta_\alpha} = \{\varphi(x_{\omega(1+\alpha)},\ldots,
x_{\omega(1+\alpha)+ n},\bar b): n<\omega,\ \bar b \subseteq
M_{\eta_\alpha \restriction\gamma}$ and  $M_{\eta_\alpha}\models \varphi
[a_{\omega (1+\alpha)},\ldots,a_{\omega(1+\alpha)+n};\bar b]\}$. 
\end{enumerate}
\end{enumerate}
\mn
There is no problem to carry out the construction.

Note: in order to have a good definition of ``suitably omitted", we
restrict the family $\bbP_{\eta_\alpha}$ to be quite well defined, loosing
some cases. 

Lastly,for $\nu\in {}^{\omega_1}2$, for $M_\nu =
\bigcup\limits_{i<\omega_1} M_{\nu \restriction i}$.  
\medskip

\noindent
\underline{\bf Step C}:   We choose $N_{\omega_1} \prec
(\cH(\chi),\in,<^*_\chi),\|N_{\omega_1}\|=\aleph_1$, $\langle
N_\alpha,M_{\eta_\alpha}:\alpha<\omega_1\rangle\in N_{\omega_1}$, and
$\omega_1 \subseteq N_{\omega_1}$. Let $\langle\bar f_i:i<\omega_1\rangle$
list the sequences from $N_{\omega_1}$ of length $\le \omega$ of 
functions $f \in N_{\omega_1}$ such that $\Rang(f) \subseteq \{0,1\}$ and
(where $\tau_*$ is from Step A)

\begin{equation*}
\begin{array}{clcr}
\Dom(f) = \{(\eta_\alpha,M^+_{\eta_\alpha}):&\alpha<\omega_1,
M^+_{\eta_\alpha} \text{ is an expanion of } M_{\eta_\alpha}
\text{ by at most} \\
  &\text{ countably many relations and functions from } \tau_*\}.
\end{array}
\end{equation*}

\mn
Let $\bar{f}_i=\langle f_{i,n}: n<\alpha_i\rangle,\alpha_i \le \omega$. We
shall choose $n(i)<\alpha_i$ and let $f_i=f_{i,n(i)}$. 

Let $\langle S'_i:i<\omega_1\rangle$ be a partition of $\omega_1$ to
stationary subsets, and if $2^{\aleph_0}<2^{\aleph_1}$ non-small sets and,
if $S_\alpha$ is not small, then even $S'_i\cap S_\alpha$ is not small
(as above, \cite[3.1(2)]{Sh:E62}.

For each $i$ let $\langle l^i_\alpha:\alpha\in S'_i\rangle$ be such that for
$j<\omega_1$, $\langle\ell^i_\alpha:\alpha\in S'_i\cap S_j\rangle$ is a weak
diamond for $f_i$ (see \cite[3.1(2)]{Sh:E62}). 
So for every $\nu\in {}^{\omega_1}2$ and
expansion $M^+_\nu$ of $M_\nu$ by at most countably many relations and
functions from $\tau_*$,
  
\[
\{\beta\in S'_i:\nu\restriction\beta=\eta_\alpha\mbox{ and }f_i
(\eta_\alpha,M^+_\nu\restriction (\omega (1+\alpha)))=\ell^i_\alpha\}
\]

\mn 
is stationary. 

Choose $\nu^*\in {}^{\omega_1}2$ such that $\nu^*(\alpha)=1-\ell_\alpha^i$
if $\alpha\in S_i'$ is not a successor ordinal. Let $M^{**} = 
M_{\eta_{\nu^*}}$. 
\end{construction}

\begin{claim}
\label{f14}
1) For $\eta\in {}^{\omega_1 \ge} 2$, $M_\eta$ has the same natural
numbers as $M^*$, but when $\lg(\eta)$ is $\omega_1$ or just is large enough,
${\omega_1}^{M_\eta}$ is not well ordered. 

\noindent
2)  For $\nu \in {}^{\omega_1} 2$ and $a \in M_\nu$ we have
\mn
\begin{enumerate}
\item[$\bullet$]   $M_\nu \models$
`` $a$ is countable" \Iff \, $\{b:M_\nu \models ``b \dot e a"\}$ 
is countable (in particular $M^{**}$ satisfies this). 
\end{enumerate}
\mn
3) For $\nu \in {}^{\omega_1} 2$, if $M_\nu \models ``a$ is an uncountable
set", \then \, for stationary many $\alpha,M_{\nu\restriction(\alpha +1)}
\models ``a_\alpha$ is a countable subset of $a$ and $b \dot e a_\alpha"$,
whenever $M_{\nu\restriction\alpha} \models ``b \dot e a"$.

\noindent
4)  If $M_\nu \models ``a$ is a stationary subset of $\omega_1"$,
\then \, the set $\{\alpha<\omega_1:M_\nu \models `a_\alpha$  is an ordinal and
$a_\alpha \dot e a"$ and $\bigwedge\limits_{b \in M_\nu}
[b \dot e a_\alpha \Rightarrow b \in M_{\nu \restriction\alpha}]\}$
is a stationary subset of $\omega_1$.

\noindent
5)  Moreover, for stationary many $\alpha$, $a_\alpha$ satisfies:
$M_\nu \models``a_\alpha$ is countable", and $M_\nu \models ``b \dot e 
a_\alpha" \Leftrightarrow b \in M_{\nu \restriction \alpha}$.
\end{claim}

\begin{PROOF}{\ref{f14}}
Straightforward.
\end{PROOF}

\begin{claim}
[CH]
\label{f17}
Assume $\omega^{M_{\langle\rangle}}$ is well ordered.

\noindent
1)  If $\nu \in {}^{\omega_1}2,M_\nu \models ``|a|=\aleph_0$ and $y$ is
a family of subsets of $a$", \then \, $M_\nu \models ``y$ is
non-meagre" iff $\{\{x:M_\nu \models x \dot e b\}:b \dot e y\}$ is a non-meagre
subset of the power set of $\{x:M_\nu \models x \dot e a\}$.

\noindent
2)  Assume $\nu \in {}^{\omega_1}2$ and $M_\nu \models ``{\bold b}_1$ 
is a Boolean ring of subsets of $a$ including the singletons, 
$|a|=\aleph_0$, ${\bold b}_2$ is a Boolean ring, and ${\bold b}_1$ is 
not meagre".  \Then \, every complete embedding ${\bold f}$ of 
$\bold b_1^{M^\nu}$ into ${\bold b}_2^{M^\nu}$ is represented in
$M_\nu$.

\noindent
3)  Assume $\nu \in {}^{\omega_1}2$ and $M_\nu \models ``{\bold b}_1$
is a Boolean ring of subsets of $a_1$ including all the finite ones, 
$a \subseteq a_1,|a| = \aleph_0,\{b\cap a:b \in {\bold b}_1\}$ is a 
non-meagre family of subsets of $a$ and ${\bold b}_2$ is a Boolean
ring".  \Then \, for every embedding ${\bold f}$ of 
${\bold b}_1^{M_\nu}$ into ${\bold b}_2^{M_\nu}$ the 
following condition is satisfied:
\mn
\begin{enumerate}
\item[$(*)$]  for some $\dot g \in M_\nu$, we have: 
\sn
\begin{enumerate}
\item[$(i)$]  $\dot g$ is a function with domain $a$ (in $M_\nu$),
\sn
\item[$(ii)$] for $b \dot e^{M_\nu} a,\dot g(b)$ is an ideal 
of ${\bold b}_2$,
\sn
\item[$(iii)$]  $b_1 \dot e^{M_\nu} a,b_2 \dot e^{M_\nu}
  a, b_1 \ne b_2 \Rightarrow \dot g(b_1) \cap \dot g(b_2) =
\{0_{{\bold b}_2}\}$,
\sn
\item[$(iv)$]  for $b \dot e^{M_\nu} a,{\bold f}(\{b\}) 
\dot e^{M_\nu} \dot g(b)$. 
\end{enumerate}
\end{enumerate}
\end{claim}

\begin{PROOF}{\ref{f17}}
1) Check.

\noindent 
2) Follows by 3).

\noindent 
3) By $\CH$ (and the choice of the $N_\alpha$'s), for some
$\alpha_0,\{a,{\bold f} \restriction\{\{b\}: b \dot e a\}\} \in N_{\alpha_0}$,
$\eta_{\alpha_0} \vartriangleleft \nu$, and all parameters are in
$M_{\eta_{\alpha_0}}$. Let ${\bold b}'_1 = \{b \cap a: b \in \bold
b_1\}$, and let $\dot g \in M_\nu$ be the function from ${\bold b}_1$ 
to ${\bold b}'_1$ such that $\dot g(b) =a \cap b$. For some $\alpha
\in (\alpha_0,\omega_1)$, $y_{\nu\restriction\alpha}$ is (in $M_\nu^*$) the
ideal of meagre subsets of ${\bold b}'_1 \subseteq \cP(a)$ 
included in ${\bold b}'_1$ (${\bold b}'_1 \notin 
y_{\nu \restriction \alpha}$ as ${\bold b}_1$ is not meagre). 

We now note that if $\eta_{\beta_\ell+1} = \eta_\alpha \char 94 \langle\ell
\rangle$, then $M_{\eta_{\beta_\ell+1}}$ cannot suitably omit 

\begin{equation*}
\begin{array}{clcr}
p := \{[{\bold b}_2 &\models ``{\bold f}(b) \le y\ \&\ f(c) \cap y
=0"]\ \&\ y \dot e {\bold b}_2:{\bold b}_1[M_{\eta_{\beta_{\ell}+1}}] \\
  &\models ``b \le a_{\beta_\ell}\ \&\ a_{\beta_\ell} \cap c
=0_{{\bold b}_1}", \text{ and } b,c \in B_1[M_{\eta_{\beta_{\ell}+1}}]
\text{ are atoms } \subseteq a)\}. 
\end{array}
\end{equation*}

\mn
(as in $M_\nu,{\bold f}(a_{\beta_\ell})$ realizes it). Hence, there is a
suitable support  

\[
\dot{\bold Q}_{\sigma_0}x_0\,\ldots\,\dot{\bold Q}_{\sigma_{n-1}}
x_{n-1}\,\exists x\, \varphi(x,x_0,\ldots,x_{n-1},a_{\beta_\ell},\bar b)
\]

\mn
of $p,\bar b \subseteq M_{\eta_\beta}$. So some $t \in 
\bold G_{\eta_{\beta_\ell +1}}$ forces this 
(for $\nleq_{\eta_{\beta_\ell+1}}$ see clause (g) of 
\ref{f11}). Using this $t$ we can define $r$ as required. 
\end{PROOF}

\begin{discussion}
\label{f20}
What occurs in \ref{f17} if we omit the assumptions
``$\omega^{M_{\langle\rangle}}$ is well founded"?  We should replace
``finite" by ``finite in the sense of $M_\nu$", in particular
(see \ref{f17}(2)) for every complete embedding 
${\bold f}^*$ of ${\bold b}_1[M_\nu]$ into ${\bold b}_2[M_\nu]$ 
for some $a',\breve{f}' \in M_\nu$, $M_\nu \models ``a'$ is a finite
subset of $a$, $b \dot e {\bold b}_2,{\bold f}$ a complete 
embedding of ${\bold b}_1 \restriction \{x \dot e {\bold b}_1: x 
\cap a=0_{{\bold b}_1}\}$ into ${\bold b}_2 \restriction \{y \dot e 
\bold b_2:y \cap b = 0_{{\bold b}_2}\}$ and $\breve{f}^{M_\nu}
= \bold f \restriction \Dom({\bold f}^{M_\nu})"$.
\end{discussion}

\noindent
The next claim says that for $\nu_0 \ne \nu_1 \in {}^{\omega_1}2$, the models
$M_{\nu_0},M_{\nu_1}$ has ``very little in common over $M_{\nu_0\cap \nu_1}$".

\begin{claim}
\label{f23}
1)  Assume $\eta\in {}^{\omega_1>}2,\eta \char 94 \langle l\rangle
\vartriangleleft\nu_l\in {}^{\omega_1} 2$, and for each
$n<\omega:M_\eta \models ``\dot T$ is a tree with set of levels 
$(W,\le_W)$ which is an $\aleph_1$-directed partial order, $a_n$ a 
level of the tree, $a_n <_W a_{n+1}"$, and 
$\{a_n: n <\omega\}$ is cofinal in $(W,\le_W)^{M_\eta}$.  If, for
$n<\omega$,   

\[
M_\eta \models `b_n \dot e \dot T \text{ is in level } a_n"
\]

\mn
and, for $\ell < \omega$,  

\[
M_{\nu_\ell} \models ``b^\ell \dot e \dot T\, \& \, b_n \le_{\dot k}\
b^\ell \le_k b^{\ell+1}"
\]

\mn
\then \, for some $c$, 

\[
M_\eta \models `c \text{ is a branch of } \dot T \text{ with uncountable
cofinality and } b_n \in c"
\]

\mn
for $n<\omega$.

\noindent
2)  If $A_1$, $A_2\subseteq M_\eta$ are disjoint and no (first order)
formula (with parameters in $M_\eta$) separates them, $\eta \in {}^{\omega_1
>}2$, $\eta \char 94 \langle\ell\rangle \vartriangleleft \nu_\ell \in 
{}^{\omega_1 \ge} 2$, and $A_1 \cup A_2 \in N_{\lg(\eta)+1}$ 
(for example $A_1\cup A_2$ is represented in $M_\eta$), \then \, 
for at least one $\ell$, in $M_{\nu_l}$ no formula separates them.
\end{claim}

\begin{remark}
\label{f24}
Note: if $A_1\cup A_2 = \{b: M_\eta \models b \dot e a\}$, then: 
[$A_1$, $A_2$ not separated in $M_\nu$] means 
[$A_1$ not represented in $M_\nu$].
\end{remark}

\begin{PROOF}{\ref{f23}}
1) Let $\eta=\eta_{\alpha(0)}$. Assume that there is no $c$ as
required. We prove by induction on $\alpha\in [\alpha(0),\omega_1]$ the
statement when we replace $M_{\nu_l}$ by
$\bigcup\{M_{\eta_\beta}:\beta \le
\alpha$ and $\eta_\beta \vartriangleleft\nu_\ell\}$. This is enough -- for
$\alpha=\omega_1$ we get the result.

\noindent For $\alpha =\alpha(0)$ this is trivial.  

\noindent For $\alpha$ limit - nothing new arises.  

\noindent The only case we have to prove something is $\eta_\alpha
\vartriangleleft \nu_\ell$,  $\alpha$ a successor. We can consider all the
countably many possible $b^{1-\ell}\cup \{M_{\eta_\beta}:\beta<\alpha$ and 
$\eta_\beta \triangleleft \nu_{1-\ell}\}$, so $\langle b_n: n<
\omega\rangle$ is determined up to $\aleph_0$ possibilities, as really the
identity of $\langle b_n:n<\omega\rangle$ is not important just the branch
which $\langle b_n:n<\omega\rangle$ determines and all those branches belong
to $N_{\alpha-1}$. So the type $p_0=\{x \dot e \dot T 
\wedge b_n <_{\dot T} x: n<\omega\}\in N_\alpha$, and 
we just have to prove that
it is omitted. Let $\eta_\beta$ be the ($\vartriangleleft$--) predecessor of   
$M_{\eta_\alpha}$. By the induction hypothesis, $M_{\eta_\beta}$ omits
$p_0$; if we fail, by the construction it is not omitted by
$M_{\eta_\beta}$. But omitting $p_0$ is equivalent to omitting 

\[
p=\{x \dot e \dot T\, \& \, [b<_{\dot T} x]^{\iif[\bigvee\limits_{n}
  b<b_n]}: M_{\eta_\beta} \models ``b \dot e \dot T"\},
\]

\mn
so by \ref{f8}(1) the type $p$ is not strongly undefined. But by $T^*$'s
choice this means it is represented in $M_{\eta_\beta}$, a contradiction. 

\noindent 
2) Same proof. 
\end{PROOF}

\begin{conclusion}
\label{f26}
1) If $M_{\nu^*} \models ``\dot T \text{ is a tree with } \delta$ levels,
$\cf(\delta)$ is regular uncountable", \then \, every full 
branch of $\dot T^{M_{\nu^*}}$ (i.e., a linear ordered subset which has
members in an unbounded set of levels) is represented in $M_{\nu^*}$.

\noindent
2) The set of levels of $\dot T$ can be partially ordered as long as it is
$\aleph_1$--directed (in the sense of $T^*$), and we get the same result. 
\end{conclusion}

\begin{PROOF}{\ref{f26}}
1)  By \ref{f23} + \ref{f11} Step C, i.e., consider
expansions of $N_\nu$ by a branch $B$ of $\dot T^{N_\nu}$ (i.e., a unary
relation). Pick $i$ such that $f_i(\eta_1,(M_{\eta_\alpha},B))=0$ iff for
some $\nu\in{}^{\omega_1} 2$, $\eta_\alpha \char 94 \langle 0\rangle
\vartriangleleft\nu$ and $B'$, a full branch of 
$\dot T^{M_\nu}$, $(M_{\eta_a},B) \prec (M_\nu,B')$. 

\noindent
2) Similar.  
\end{PROOF}

\begin{definition}
\label{f29}
For an atomic Boolean ring $\bold B$:  

\noindent
1) $\bold B$ is non-meagre, if, identifying $b \in \bold B$ 
with $\{x: x \in \bold B^{\at}; x \le_{\bold B} b\},\bold B$ 
is a non-meagre family of subsets of $\bold B^{\at}$ ($\bold B^{\at}$ is
the set of ``atoms" of $\bold B$), i.e., $\bold B$ can be 
represented as a countable union $\bigcup\limits_{n<\omega} Y_n$, 
each $Y_n$ nowheredense (i.e., for every finite $a_1 \subseteq a_2
\subseteq \bold B^{\at}$ there are finite $b_1 \subseteq b_2
\subseteq \bold B^{\at}$ such that $a_1 \subseteq b_1,a_2 \setminus
a_1 \subseteq b_2 \setminus b_1$, and for no $c \in \bold B$ do we have $\cap
b_2=b_1$). 
\end{definition}

\begin{observation}
\label{f32}
1)  If $Y \subseteq {\cP}(X)$ (i.e., is a family of subsets of $X$),
$X'\subseteq X$, and $\{y \cap X':y\in Y\}$ is a meagre subset of 
$\cP(X')$, then $Y$ is a meagre subset of ${\cP}(X)$.

\noindent
2) If $Y$ is a meagre (or nowheredense) subset of ${\cP}(X)$, \then \, the
set $\{a \subseteq X: a$ is countable, $\{y \cap a:y \in Y\}$ is meagre
(or nowheredense) subset of a ${\cP}(a)\}$ is a club of 
$[X]^{\aleph_1} = \{a \subseteq X: |a|=\aleph_0\}$. 
\end{observation}

\begin{question}
\label{f35}
Phrase the statement which suffices for the proof instead $\CH$ 
(it seems the existence of a non-meagre set of
cardinality $\aleph_1$ suffices). 
\end{question}

\begin{conclusion}
\label{f38}
Assume $\omega^{M_{\langle\rangle}}$ is well ordered, and suppose $M_{\nu^*}
\models ``\bold b_1$ is a Boolean ring of subsets of $a$ including
the singletons, $|a|>\aleph_0$, ${\bold b}_1$ is non-meagre and 
${\bold b}_2$ is a Boolean algebra".  \Then \, every complete embedding of
${\bold b}_1$ into ${\bold b}_2$ is represented in $M_{\nu^*}$.  
\end{conclusion}

\begin{PROOF}{\ref{f38}}
In $M_{\nu^*}$ let $(W,\le)$ be the set of countable subsets 
of ${\bold b}^{\at}_1 \cup {\bold b}_2$ ordered by inclusion.  
We define a tree $\dot t$ with $W$ as a set of levels by: 

\[
\dot t = \{(c,f):c \dot e W \text{ and } f \text{ is a function from } c
\cap {\bold b}^{\at}_1 \text{ into } {\bold b}_2\}. 
\]

\mn
We define the order of $\dot t$ by
$(c_1,f_1) \le (c_2,f_2) \Leftrightarrow c_1 \subseteq c_2\, \& \, f_1
\subseteq f_2$ (defined in $M_{\nu^*}$).

Now, if ${\bold f}$ is a complete embedding of 
${\bold b}_1^{M_{\nu^*}}$ into ${\bold b}_2^{M_{\nu^*}}$, then for 
a club $E$ of $\alpha<\omega_1$, the restriction ${\bold f}
\restriction {\bold b}_1[M_{\nu^* \rest \alpha}]$ is a complete embedding of
${\bold b}_1[M_{\nu^*\restriction\alpha}]$ 
into ${\bold b}_2 [M_{\nu^* \rest \alpha}]$
(see \cite{Sh:384}). Let the ordinal $\gamma$ be such that 

\[
M_{\nu^*} \models ``a_\gamma = \text{ the non-stationary ideal on } 
[{\bold b}_1[M_{\nu^*}] \cup {\bold b}_2 [M_{\nu^*}]^{\aleph_0}".
\]

\mn
So for $\alpha\in E \cap S_\gamma$, $\{b:M_{\nu^*}\models b \dot e a_{\omega(1
+\alpha)}\}$ is exactly $({\bold b}_1[M_{\nu^*}] \cup {\bold b}_2 
[M_{\nu^*}]) \cap M_{\nu^*\restriction \alpha}$.
  
Now we apply \ref{f23}(2) with ${\bold f},
a_{\omega(1+\alpha)}\cap {\bold b}^\at_1[M_{\nu^*\restriction\alpha}],
\bold b^{\at}_1$, ${\bold b}_1$, ${\bold b}_2$ here standing for
${\bold f}$, $a$, $a_1$, ${\bold b}_1$, ${\bold b}_2$ there. 
We get $\dot g = \dot g_\alpha\in M_{\nu^*}$ as there. But 
$\alpha\in E$, so $\dot g = {\bold f} \restriction {\bold
  b}^{\at}_1[M_{\nu^* \restriction \alpha}]$. 
Clearly $\dot g_\alpha \dot e {\dot t}_{a_{\omega(1+\alpha)}}$
and for $\alpha<\beta$ from $E\cap S_\gamma$ we have: $(W,\le) \models
``a_{\omega(1+\alpha)} \le a_{\omega(1+\beta)}",\dot t \models 
``\dot g_\alpha \le \dot g_\beta"$; and 
$\{a_{\omega(1+\alpha)}:\alpha\in E\cap S_\gamma\}$ is
cofinal in $(W,\le)$, and $\{\dot g_\alpha:\alpha <\omega_1\}$ induce a branch
of $\dot t$. 

By \ref{f26}(2) we finish. 
\end{PROOF}

\begin{conclusion}
\label{f41}
Assume $\omega^{M_{\langle\rangle}}$ is well ordered and 

\begin{equation*}
\begin{array}{clcr}
M_{\nu^*} \models ``&(a) \quad {\bold b}_1 \text{ is a Boolean ring},\\
 &\qquad {\bf \Xi} 
\text{ is a family of maximal antichains of }{\bold b}_1,\\  
 &(b) \quad \text{for } \Xi \in {\bf \Xi}, \text{ the sub-algebra 
sub-Boolean ring}\\
 &\qquad {\bold b}_1^{[\Xi]} = \{x \in {\bold b}:\text{ for every } y \in
 \Xi, x \cap y = 0 \vee x \cap y=y\}\\ 
 & \qquad {\bold b}_1 \text{ is non meagre, i.e. essentially as a
   family of subsets of } b \equiv\\
 &(c) \quad {\bf \Xi} \text{ is } \aleph_1 \text{-directed (order: 
bigger means finer)}, \\
 &(d) \quad \text{ for every $x \in {\bold b}_1 \setminus \{0_B\}$, 
there are } \Xi \in {\bf \Xi}, \text{ and } y \in \Xi,\\
 & \qquad \text{ such that } 0 < y \le x \text{ or at least}\\
 & \qquad x = \sup_{\bold b_1}\{z:\;\text{ for every } 
\Xi \in {\bf \Xi}$ and $y \in {\bold b}_1^{[\Xi]} \text{ we have:} \\
 & \qquad x \le y \Rightarrow  z \le y\},\\ 
 &(e) \quad \bold b_2 \text{ is a Boolean ring}".
\end{array}
\end{equation*}

\mn
\Then \, every complete embedding of ${\bold b}_1[M_{\nu^*}]$ into
${\bold b}_2[M_{\nu^*}]$ is represented in $M_{\nu^*}$.
\end{conclusion}

\begin{PROOF}{\ref{f41}}
Combine \ref{f38} and \ref{f26}(2).
\end{PROOF}

\begin{definition}
\label{f44}
For a model $(A,\le,R)$ of $\gt^{\poe}$ 
(see Definition \ref{c2}) and $X\subseteq A$ we say that: 

\noindent
1) A set $X\subseteq A$ is nwd (nowhere dense) \when \, every cone has a
subcone disjoint to it (a cone is $\{x:x_0 \le x\}$). A set $X\subseteq A$
is meagre if it is a countable union of nowhere dense sets.

\noindent
2)  A set $X \subseteq A$ is non-medium meagre \when \,
if the family of countable
$a \subseteq A$ satisfying $(*)_a=(*)_a[X]$ is an unbounded subset of 
$[A]^{\le \aleph_0}$, where 
\newline
$(*)_a[X]:(A,\le,R)\restriction a$ is a model of $\gt^{\poe}$, and
we cannot find $X_n$, a nwd subset of $a$ (for $n<\omega$) such that:
$a=\bigcup\limits_{n<\omega} X_n$, and for every $c\in A$ there is
$n<\omega$ satisfying: $\{b\in a:b \le c\}\subseteq X_n$. 

If $X=A$ we may omit it. We say in this case that $X$ is non-meagre in $(A,
\le,R)$ for $a$. 

\noindent
3) $X \subseteq A$ is non-weakly meagre \when \, for a stationary set of $a\in
[A]^{\le \alpha_0}$ we have $(*)_a[X]$.
\end{definition}

\begin{remark}
\label{f47}
1)  For an ordered field or just a dense linear order ($A,<$) we use
$A$= the set of open intervals of $A_1$, with $\leq$ being subintervals,
$R$ being disjoint.

\noindent
2)  If we can get the parallel of \ref{f17}, \ref{f29},
\ref{f38} to models of $\gt^{\poe}$ (hence to ordered fields), we
later get stronger results, the missing point is \ref{f38}(1)
-- downward monotonicity of non-meagre.  
\end{remark}

\begin{claim}
\label{f50}
1)  If $(A_1,\le,R)$ is a model of ${\gt}^{\poe},X\subseteq A$
is meagre in $(A,\le,R)$, \then \, for some club $S 
\subseteq [A]^{<\aleph_0}$,  for every $a\in S$, we have: $X$ is meagre in
$(A,\le,R)$ for $a$. This in turn means: $X$ is weakly meagre.
If $X$ is medium meagre in $(A,\le,R)$, \then\ $X$ is weakly meagre. 
\end{claim}

\begin{PROOF}{\ref{f50}}
Should be clear.
\end{PROOF}

\begin{claim}
\label{f53}
[$\CH$] 

Assume $\nu\in {}^{\omega_1}2$, for $\ell=1,2$ 

\[
M_\nu \models ``(A^\ell,\le^\ell,R^\ell) \text{ is a model of }
\gt^{\poe}"
\]

\mn
and

\[
M_\nu \models ``a^\ell \subseteq A^\ell \text{ is countable}".
\]

\mn
Also $(A^1,\le^1,R^1)\restriction a^1 \models \gt^{\poe}$, and $M_\nu
\models$ ``for $(A^1,\le^1,R^1)$, $(*)_a$ from \ref{f44}(2) holds".
\Then \,  for every embedding $f$ of $(A^1,\le^1,R^1)^{M_\nu}
\restriction a$ into $(A^2,\le^2,R^2)^{M_\nu}$ mapping $a^1$ into
$a^2$ we have: 
\mn
\begin{enumerate}
\item[$\otimes$]  for every cone $C$ of $(A^1,\le^1,R^1)^{M_\nu}
\restriction a$, on some subcone $C'$ of 
\newline
$(A^1,\le^1, R^1)^{M_\nu} \rest a$, we have: 
\sn
\begin{enumerate}
\item[$(*)$]  there is $r\in M_\nu$ such that:
\sn
\item[{{}}] $(i) \quad M_\nu \models ``\dot g$ is a function 
with domain $\{x \dot e a^1:x \dot e C'\}$,  
\sn
\item[${{}}$]  $(ii) \quad$ for $x \dot e^{M_\nu} \Dom(\dot g)$ 
we have: $\dot g(x)$ is a subset of $A^2$, 
\sn
\item[${{}}$]  $(iii) \quad$ if $b_1,b_2 \dot e^{M_\nu} \Dom(\dot g)$,
  and $b_1 R^1 b_2$, then 

\hskip25pt $(\exists x_1 \dot e r(b_1))
(\exists x_2 \dot e \dot g(b_2))[x_1;R^2;x_2$ hence $x_1,x_2$ are 
$\le^1$-incomparable],
\sn
\item[${{}}$]  $(iv) \quad$ for $b \dot e^{M_\nu} \Dom(\dot g)$, 
we have: $f(b) \in \dot g(b)$.
\end{enumerate}
\end{enumerate}
\end{claim}

\begin{PROOF}{\ref{f53}}
Straightforward (like the proof of \ref{f17}(2)).
\end{PROOF}

\begin{claim}
\label{f56}
[$\CH$] 

Assume that $M_{\nu^*} \models ``(A_\ell,\le^\ell,R^\ell)$ is a model of 
$\gt^{\poe}$  and for $\ell = 1$ non-medium meagre" and ${\bold f}$ is 
a dense embedding (see \cite{Sh:384}, i.e., on branches) of
$(A^1,\le^1,R^1)^{M_{\nu^*}}$ into $(A^2,\le^2,R^2)^{M_{\nu^*}}$. 
\Then \, for a dense set of cones $C$ (of $(A^1,\le^1,R^1)$,
${\bold f} \restriction C$ is represented in $M_{\nu^*}$. 
\end{claim}

\begin{PROOF}{\ref{f56}}
Like the proof of \ref{f38} (using Fodor lemma).
\end{PROOF}

\begin{definition}
\label{f59}
1) We say that $\bold B$ is a partial Boolean algebra if the functions $(x\cap
y,x\cup y,x-y,0,1)$ are partial (but $0^{\bold B}$ 
well defined), so the identities
are interpreted as ``if at least one side is well defined then so is the
other and they are equal".  (So a Boolean ring is a partial Boolean
algebra.) Let $a \le b$ mean $a \cap b=a$, so $\neg[b\cap a=0]$ means $b\cap
a$ is an element $\neq 0$ or undefined. 

\noindent
2) Let $\bold B$ be a partial Boolean algebra. 
A set $\Xi \subseteq \bold B$ is called a maximal antichain of $\bold
B$ \when \,:
\mn
\begin{enumerate}
\item[$(*)$]  $(a) \quad a \in \Xi \Rightarrow a \ne 0$
\sn
\item[${{}}$]  $(b) \quad a \ne b \in \Xi \Rightarrow a \cap b=0$
\sn
\item[${{}}$]  $(c) \quad b \in \bold B \setminus \{0\} 
\Rightarrow \bigvee\limits_{a\in\Xi}
(\exists c)[a\cap b = c \ne 0_{\bold B}]$. 
\end{enumerate}
\mn
3) For a partial Boolean algebra $\bold B$ and a maximal antichain $\Xi$, let 
$\bold B^{\Xi}$ be a partial Boolean algebra with universe $\bold B$ and 
$\bold B^{\Xi} \models ``b \le c"$ iff for every $a \in \Xi,
[b\cap a \ne 0 \Rightarrow c \cap a \ne 0]$.

\noindent
4) For $\bold B,\Xi$ as above, and $Y \subseteq \bold B$ 
we call $Y$ nowhere dense for $(\bold B,\Xi)$, if for every 
partial finite function $h$ from $\Xi$ to $\{0,1\}$ there is a 
finite function $h^+$ from $\Xi$ to $\{0,1\}$ extending $h$
and such that for no $c \in Y$ do we have
$h^+(a)=0 \Rightarrow a \cap c=0,h^+(a)=1 \Rightarrow a \le c$. 

We say $Y$ is $\mu$-meagre for $(\bold B,\Xi)$ \If \, it is the 
union of $<\mu$ nowhere dense for $(\bold B,\Xi)$ sets; 
if $\mu=\aleph_1$ we omit it.

We say $\bold B$ is $\mu$-meagre over $\Xi$ when $\bold B$ 
is $\mu$-meagre for $(\bold B,\Xi)$ as a subset of $\bold B$. 
\end{definition}

\begin{claim}
\label{f62}
Assume $\omega^{M_{\langle\rangle}}$ is well ordered. If
${\bold b}^{M^{**}}_1$, ${\bold b}^{M^{**}}_2$ are Boolean rings in
$M^{**} = M_{\nu_*}$, $\Xi$ as in \ref{f32}(a),(c),(d) and
\mn 
\begin{enumerate}
\item[$(b)^-$]  for every $\Xi \in {\bf \Xi}$, ${\bold b}^{\Xi}_1$ is not
meagre (in the sense of $N^{**}$)
\end{enumerate}
\mn
\then \, every complete embedding of ${\bold b}^{M^{**}}_1$ into
${\bold b}^{M^{**}}_2$ is represented in $M^{**}$.
\end{claim}

\begin{PROOF}{\ref{f62}}
Straightforward.
\end{PROOF}

\begin{theorem}
\label{f65}
[$\CH$]
1)  The logic ${\bbL}$ extended by the following
quantifiers is still $\aleph_0$-compact (getting models of 
cardinality $\aleph_1$:
\mn
\begin{enumerate}
\item[$(A)$]  complete embedding of one Boolean ring to another,
\sn
\item[$(B)$]   embedding of one ordered field into another with dense range.
\end{enumerate}
\mn
2)  In the logic $\bbL_{\aleph_1,\aleph_0}$ extended by 
$\exists^{\ge \aleph_1}$ and the following 
quantifiers we still cannot characterize well ordering of order type 
$\le \omega_1$:
\mn
\begin{enumerate}
\item[$(A)$]  non-meagreness of a family of subsets of a countable
  set,
\sn
\item[$(B)$]  complete embedding of a non meagre B.A. into a Boolean
  algebra,  
\sn
\item[$(C)$]  dense embedding of a non meagre ordered field considering the
interval under inclusion as a model of $\gt^{\poe}$. 
\end{enumerate}
\mn
3) The logic $\bbL$ extended by $\exists^{\ge \aleph_1}$ and 
the quantifier from (A),(B),(C) of part (2), is $\aleph_0$-compact, getting
models of cardinality $\aleph_1$.
\end{theorem}

\begin{PROOF}{\ref{f65}}
1) Assume we are given such a theory $T$ in this
logic. First use \S5 to get an $\aleph_1$-compact model $N^*$ of $T$ 
(e.g. in ${\bold L}[A], A \subseteq 2^{\aleph_2}$
code $T$ and $\cP(\aleph_2)$, 
which satisfies (A)+(B), then create a model $\gC^*$ of $T^*$ 
in which $N^*$ is a member. Let $M^* \prec{\gC}^*,\|M^*\| =
\aleph_0$, $N^*\in M^*$, and apply this section.  

\noindent
2) Should be clear. 
\end{PROOF}

\begin{remark}
\label{f68} 
You can add in \ref{f65}(2),(3) 
also the quantifier (aa $X$), i.e., make ``for
stationary many countable $x \subseteq y$" be standard. For this in
\ref{f14}(5) we should replace ``stationary many $\alpha$" by ``club
many $\alpha$", and so restrict somewhat the $\Xi$-s which we may use.
\end{remark}

\begin{claim}
\label{f71}
Assume $\omega^{M_{\langle\rangle}}$ is well ordered. 
For ${\bold b}_1$, ${\bold b}_2\in M_{\nu^*}$ such that 
${\bold b}_\ell[M_{\nu^*}]$ is a triple $(P^\ell,Q^\ell,R,P^\ell)$ 
with the strong independence property (this means satisfying
$\aleph_0$ sentences).

\noindent
1) Claim \ref{f17}, Def.\ref{f29}, conclusion \ref{f38}
generalize naturally to dense embedding (see \cite{Sh:384}).

\noindent
2) In \ref{f65}(2) we can add:
\mn 
\begin{enumerate}
\item[$(D)$]  dense embedding of one interpretation of a model of the
strong independence property into another.
\end{enumerate}
\end{claim}

\begin{PROOF}{\ref{f71}}
No new point.
\end{PROOF}

\begin{remark}
\label{f74}
1)  We can in \ref{f5} and say that $T^*$ suitably omit $\Gamma$ \when
\, $\Gamma \subseteq \{\varphi(\bar x)\in {\bbL}(\tau_{T^*})\}$ and
\mn
\begin{enumerate}
\item[$(*)$]  if $T^*\cup\{(\dot{\bold Q}_{\sigma_n} y_n)
\ldots (\dot{\bold Q}_{\sigma_1} y_1) 
(\exists \bar x) \psi (\bar x,y_i,\ldots,y_n)\}$ is consistent then for semi 
$\varphi(\bar x)\in p, T^*\cup \{(\dot{\bold Q}_{\sigma_n} y_n)\ldots 
(\dot{\bold Q}_{\sigma_1} y_1)(\exists \bar x)(\psi(\bar x,y_1\ldots y_n)\wedge
\neg \varphi (\bar x))$ and has the ``omitting type theorem".
\end{enumerate}
\mn
2)  We can replace here $(\dot{\bold Q}_\sigma y)$ by 
$(\dot{\bold Q}_\Gamma y)$ where 
$\Gamma$ is a bigness notion and $\lg \bar y = \lg \bar x_\Gamma$. 
\end{remark}

\begin{theorem}
\label{f77}
[$\CH$]
 
Let $T$ be countable and complete first order theory. \Then \,$T$ has a
model $M^*$ of cardinality $\aleph_1$ such that: 
\mn
\begin{enumerate}
\item[$(A)$]   If ${\bold b}_1,{\bold b}_2$ are 
interpretations of Boolean rings in $M^*$, every complete embedding 
of ${\bold b}_1$ to ${\bold b}_2$ is definable (from parameters) in
$M^*$
\sn 
\item[$(B)$]  if $\bbF_1,\bbF_2$ are interpretation of dense 
linear ordered in $M^*$, every dense embedding of $\bbF_1$ into
$\bbF_2$ is on a dense set of interval definable (from parameters) in
$M^*$
\sn. 
\item[$(C)$ ] The parallel to interpretation of the theorem $\dot k^{\ind}$.
\end{enumerate}
\end{theorem}

\begin{PROOF}{\ref{f77}}
Use \cite{Sh:107}+the theorem above. 
\end{PROOF}
\newpage

\bibliographystyle{alphacolon}
\bibliography{lista,listb,listx,listf,liste,listz}

\def\germ{\frak} \def\scr{\cal} \ifx\documentclass\undefinedcs
  \def\bf{\fam\bffam\tenbf}\def\rm{\fam0\tenrm}\fi 
  \def\defaultdefine#1#2{\expandafter\ifx\csname#1\endcsname\relax
  \expandafter\def\csname#1\endcsname{#2}\fi} \defaultdefine{Bbb}{\bf}
  \defaultdefine{frak}{\bf} \defaultdefine{=}{\B} 
  \defaultdefine{mathfrak}{\frak} \defaultdefine{mathbb}{\bf}
  \defaultdefine{mathcal}{\cal}
  \defaultdefine{beth}{BETH}\defaultdefine{cal}{\bf} \def\bbfI{{\Bbb I}}
  \def\mbox{\hbox} \def\text{\hbox} \def\om{\omega} \def\Cal#1{{\bf #1}}
  \def\pcf{pcf} \defaultdefine{cf}{cf} \defaultdefine{reals}{{\Bbb R}}
  \defaultdefine{real}{{\Bbb R}} \def\restriction{{|}} \def\club{CLUB}
  \def\w{\omega} \def\exist{\exists} \def\se{{\germ se}} \def\bb{{\bf b}}
  \def\equivalence{\equiv} \let\lt< \let\gt>
  \def\implies{\Rightarrow}\def\mathfrak{\bf}\def\germ{\frak} \def\scr{\cal}
  \ifx\documentclass\undefinedcs
  \def\bf{\fam\bffam\tenbf}\def\rm{\fam0\tenrm}\fi 
  \def\defaultdefine#1#2{\expandafter\ifx\csname#1\endcsname\relax
  \expandafter\def\csname#1\endcsname{#2}\fi} \defaultdefine{Bbb}{\bf}
  \defaultdefine{frak}{\bf} \defaultdefine{=}{\B} 
  \defaultdefine{mathfrak}{\frak} \defaultdefine{mathbb}{\bf}
  \defaultdefine{mathcal}{\cal}
  \defaultdefine{beth}{BETH}\defaultdefine{cal}{\bf} \def\bbfI{{\Bbb I}}
  \def\mbox{\hbox} \def\text{\hbox} \def\om{\omega} \def\Cal#1{{\bf #1}}
  \def\pcf{pcf} \defaultdefine{cf}{cf} \defaultdefine{reals}{{\Bbb R}}
  \defaultdefine{real}{{\Bbb R}} \def\restriction{{|}} \def\club{CLUB}
  \def\w{\omega} \def\exist{\exists} \def\se{{\germ se}} \def\bb{{\bf b}}
  \def\equivalence{\equiv} \let\lt< \let\gt>
\providecommand{\bysame}{\leavevmode\hbox to3em{\hrulefill}\thinspace}
\providecommand{\MR}{\relax\ifhmode\unskip\space\fi MR }
\providecommand{\MRhref}[2]{%
  \href{http://www.ams.org/mathscinet-getitem?mr=#1}{#2}
}
\providecommand{\href}[2]{#2}
\begin{thebibliography}{}

\bibitem[Kei70]{Ke70}
H.~Jerome Keisler, \emph{{Logic with the quantifier "there exist uncountably
  many"}}, Annals of Mathematical Logic \textbf{1} (1970), 1--93.

\bibitem[Kei71]{Ke71}
\bysame, \emph{{Model theory for infinitary logic. Logic with countable
  conjunctions and finite quantifiers}}, {Studies in Logic and the Foundations
  of Mathematics}, vol.~62, North--Holland Publishing Co., Amsterdam--London,
  1971.

\bibitem[Spe87]{Sp87}
Joel Spencer, \emph{{Ten lectures on the probabilistic method}}, CBMS-NSF
  Regional Conference Series in Applied Mathematics, vol.~52, Society for
  Industrial and Applied Mathematics (SIAM), Philadelphia, PA, 1987.

\bibitem[Sh:a]{Sh:a}
Saharon Shelah, \emph{{Classification theory and the number of nonisomorphic
  models}}, {Studies in Logic and the Foundations of Mathematics}, vol.~92,
  {North-Holland Publishing Co., Amsterdam-New York, xvi+544 pp, \$62.25},
  1978.

\bibitem[Sh:c]{Sh:c}
\bysame, \emph{{Classification theory and the number of nonisomorphic models}},
  {Studies in Logic and the Foundations of Mathematics}, vol.~92,
  {North-Holland Publishing Co., Amsterdam, xxxiv+705 pp}, 1990.

\bibitem[Sh:3]{Sh:3}
\bysame, \emph{{Finite diagrams stable in power}}, {Annals of Mathematical
  Logic} \textbf{2} (1970), 69--118.

\bibitem[Sh:E58]{Sh:E58}
\bysame, \emph{{Existence of endo-rigid Boolean Algebras}}.

\bibitem[Sh:E59]{Sh:E59}
\bysame, \emph{{General non-structure theory and constructing from linear
  orders}}, arxiv:1011.3576.

\bibitem[Sh:E60]{Sh:E60}
\bysame, \emph{{Constructions with instances of GCH: applying}}.

\bibitem[Sh:E62]{Sh:E62}
\bysame, \emph{{Combinatorial background for Non-structure}},
  arxiv:math.LO/1512.04767.

\bibitem[Sh:72]{Sh:72}
\bysame, \emph{{Models with second-order properties. I. Boolean algebras with
  no definable automorphisms}}, {Annals of Mathematical Logic} \textbf{14}
  (1978), 57--72.

\bibitem[Sh:73]{Sh:73}
\bysame, \emph{{Models with second-order properties. II. Trees with no
  undefined branches}}, {Annals of Mathematical Logic} \textbf{14} (1978),
  73--87.

\bibitem[RuSh:84]{RuSh:84}
Matatyahu Rubin and Saharon Shelah, \emph{{On the elementary equivalence of
  automorphism groups of Boolean algebras; downward Skolem-Lowenheim theorems
  and compactness of related quantifiers}}, {The Journal of Symbolic Logic}
  \textbf{45} (1980), 265--283.

\bibitem[Sh:107]{Sh:107}
Saharon Shelah, \emph{{Models with second order properties. IV. A general
  method and eliminating diamonds}}, {Annals of Pure and Applied Logic}
  \textbf{25} (1983), 183--212.

\bibitem[Sh:128]{Sh:128}
\bysame, \emph{{Uncountable constructions for B.A., e.c. groups and Banach
  spaces}}, {Israel Journal of Mathematics} \textbf{51} (1985), 273--297.

\bibitem[Sh:136]{Sh:136}
\bysame, \emph{{Constructions of many complicated uncountable structures and
  Boolean algebras}}, {Israel Journal of Mathematics} \textbf{45} (1983),
  100--146.

\bibitem[HLSh:162]{HLSh:162}
Bradd Hart, Claude Laflamme, and Saharon Shelah, \emph{{Models with second
  order properties, V: A General principle}}, {Annals of Pure and Applied
  Logic} \textbf{64} (1993), 169--194, arxiv:math.LO/9311211.

\bibitem[Sh:229]{Sh:229}
Saharon Shelah, \emph{{Existence of endo-rigid Boolean Algebras}}, {Around
  classification theory of models}, {Lecture Notes in Mathematics}, vol. 1182,
  {Springer, Berlin}, 1986, arxiv:math.LO/9201238, pp.~91--119.

\bibitem[Sh:247]{Sh:247}
\bysame, \emph{{More on stationary coding}}, {Around classification theory of
  models}, {Lecture Notes in Mathematics}, vol. 1182, {Springer, Berlin}, 1986,
  pp.~224--246.

\bibitem[Sh:309]{Sh:309}
\bysame, \emph{{Black Boxes}}, {}, 0812.0656. 0812.0656. arxiv:0812.0656.

\bibitem[Sh:312]{Sh:312}
\bysame, \emph{{Existentially closed locally finite groups}}, {preprint},
  arxiv:math.LO/1102.5578v2.

\bibitem[Sh:326]{Sh:326}
\bysame, \emph{{Vive la diff\'erence I: Nonisomorphism of ultrapowers of
  countable models}}, {Set Theory of the Continuum}, {Mathematical Sciences
  Research Institute Publications}, vol.~26, {Springer Verlag}, 1992,
  arxiv:math.LO/9201245, pp.~357--405.

\bibitem[Sh:331]{Sh:331}
\bysame, \emph{{A complicated family of members of tress with $ \omega +1 $
  levels}}, {}, arxiv:math.LO/1404.2414.

\bibitem[Sh:363]{Sh:363}
\bysame, \emph{{On spectrum of $\kappa $-resplendent models}}, {preprint}.

\bibitem[Sh:384]{Sh:384}
\bysame, \emph{{Compact logics in ZFC : Complete embeddings of atomless Boolean
  rings}}, {Non structure theory, Ch X}.

\bibitem[Sh:405]{Sh:405}
\bysame, \emph{{Vive la diff\'erence II. The Ax-Kochen isomorphism theorem}},
  {Israel Journal of Mathematics} \textbf{85} (1994), 351--390,
  arxiv:math.LO/9304207.

\bibitem[Sh:413]{Sh:413}
\bysame, \emph{{More Jonsson Algebras}}, {Archive for Mathematical Logic}
  \textbf{42} (2003), 1--44, arxiv:math.LO/9809199.

\bibitem[Sh:421]{Sh:421}
\bysame, \emph{{Kaplansky test problem for $R$-modules in ZFC}}, {}.

\bibitem[Sh:482]{Sh:482}
\bysame, \emph{{Compactness in ZFC of the Quantifier on ``Complete embedding of
  BA's''}}, {}.

\bibitem[Sh:511]{Sh:511}
\bysame, \emph{{Building complicated index models and Boolean algebras}}, {}.

\bibitem[Sh:572]{Sh:572}
\bysame, \emph{{Colouring and non-productivity of $\aleph_2$-cc}}, {Annals of
  Pure and Applied Logic} \textbf{84} (1997), 153--174, arxiv:math.LO/9609218.

\bibitem[GiSh:577]{GiSh:577}
Moti Gitik and Saharon Shelah, \emph{{Less saturated ideals}}, {Proceedings of
  the American Mathematical Society} \textbf{125} (1997), 1523--1530,
  arxiv:math.LO/9503203.

\bibitem[Sh:757]{Sh:757}
Saharon Shelah, \emph{{Quite Complete Real Closed Fields}}, Israel Journal of
  Mathematics \textbf{142} (2004), 261--272, arxiv:math.LO/0112212.

\bibitem[Sh:800]{Sh:800}
\bysame, \emph{{On complicated models}}, Preprint.

\end{thebibliography}

\end{document}